\documentclass[10pt,twoside, a4paper, english, reqno]{amsart}

\usepackage{multirow}

\usepackage{tikz}
\usepackage{pgfplots}
\pgfplotsset{width=7cm,compat=newest}

\usepackage{blindtext}

\usepackage{booktabs}
\usepackage{lscape}
\usepackage{caption}
\usepackage{subcaption}
\usepackage{algorithm}
\usepackage{algorithmic}
\usepackage{amscd}
\usepackage{amssymb}
\usepackage{amsthm}
\usepackage{graphicx}
\graphicspath{{Imagestdmvareps/}}
\usepackage{epstopdf}
\usepackage{lineno,hyperref}

\usepackage{latexsym}

\usepackage{upref}
\usepackage{hyperref}
\theoremstyle{plain}
\newtheorem{thm}{Theorem}[section]
\theoremstyle{plain}
\newtheorem{lem}[thm]{Lemma}

\newtheorem{cor}[thm]{Corollary}

\theoremstyle{definition}
\newtheorem{defi}{Definition}[section]

\newtheorem*{maintheorem*}{Main Theorem}

\newcommand{\goto}{\ensuremath{\rightarrow}}
\newcommand{\grad}{\ensuremath{\nabla}}

\usepackage[toc,page]{appendix}

\newcommand{\kom}[1]{}
\renewcommand{\kom}[1]{{\bf [#1]}}
\definecolor{darkjunglegreen}{rgb}{0.1, 0.14, 0.13}
\definecolor{blau}{rgb}{0.1,0.0,0.9}

\newcounter{komcounter}
\numberwithin{komcounter}{section}

\numberwithin{equation}{section} \allowdisplaybreaks

\title[]
{Well-posedness of a Variable-Exponent Telegraph Equation Applied to Image Despeckling}

\date{}



\author[]{Sudeb Majee$^*$,  {\AA}ke Br\"annstr\"om, Niklas L.P. Lundstr\"om} 
\address[Sudeb Majee]{\newline
Department of Mathematics and Mathematical Statistics\\
Ume{\aa} University\\
SE-90187 Ume{\aa}, Sweden}
\email[]{sudebmajee@gmail.com}
\address[{\AA}ke Br\"annstr\"om]{\newline
Department of Mathematics and Mathematical Statistics\\
Ume{\aa} University\\
SE-90187 Ume{\aa}, Sweden \\
\\
Advancing Systems Analysis Program\\
International Institute for Applied Systems Analysis (IIASA)\\
Schlossplatz 1\\
2361 Laxenburg, Austria}
\email[]{ake.brannstrom@umu.se}
\address[Niklas L.P. Lundstr\"om]{\newline
Department of Mathematics and Mathematical Statistics\\
Ume{\aa} University\\
SE-90187 Ume{\aa}, Sweden}
\email[]{niklas.lundstrom@umu.se}
\thanks{$^*$ Corresponding author: Sudeb Majee}
\begin{document}

\begin{abstract}
In this paper, we present a telegraph diffusion model with variable exponents for image despeckling. Moving beyond the traditional assumption of a constant exponent in the telegraph diffusion framework, we explore three distinct variable exponents for edge detection. All of these depend on the gray level of the image or its gradient. We rigorously prove the existence and uniqueness of weak solutions of our model in a functional setting and perform numerical experiments to assess how well it can despeckle noisy gray-level images. We consider both a range of natural images contaminated by varying degrees of artificial speckle noise and synthetic aperture radar (SAR) images. We finally compare our method with the nonlocal speckle removal technique and find that our model outperforms the latter at speckle elimination and edge preservation.
\end{abstract}

\maketitle



\section{Introduction}
\label{Introduction}
In practice, due to various factors, images are often degraded by different types of noise, resulting in loss of pixel information in images. Therefore, image restoration is an essential step before beginning high-level image analysis. Hence a crucial challenge in digital image processing is eliminating noise from the acquired images by finding the best possible approximation of the unknown true image from a noisy image. One commonly employed strategy is to smooth out noise from the image while preserving essential attributes such as edges and textures. 

In this study, we focus on the removal of speckle noise. This kind of noise has a granular appearance and emerges through the interference of wavefronts in coherent imaging systems such as active radar, synthetic aperture radar, medical ultrasounds, laser, and optical coherence tomography images. Mathematically a signal-dependent image degradation model can be described as
\begin{equation*}
I_0=I\eta\,
\end{equation*}
where $I_0$, $I$, and $\eta$ indicate the speckled image, clean image, and the speckle-noise, respectively. It is commonly assumed that the speckle noise $\eta$ is Gamma$(L, L)$ distributed, where $L \in {\rm I\!N} $ is the number of ``looks" concerned to the number of spatial observations \cite{argenti2013tutorial}. Numerous studies report the fundamentals and the statistical characteristics of the speckle noise  \cite{achim2001novel,frost1982model,kuan1985adaptive,lee1980digital}.
Well-known despeckling approaches possess Bayesian approaches in the spatial domain \cite{frost1982model,kuan1985adaptive,lee1980digital},
Bayesian approaches in the transformed domain \cite{aiazzi1998multiresolution,meer1994multiresolution}, order-statistics and
morphological filters \cite{alparone1995two,crimmins1985geometric},
simulated annealing despeckling \cite{white1994simulated}, nonlocal
filtering \cite{parrilli2011nonlocal,zhong2011sar}, probabilistic
patch-based algorithm \cite{deledalle2009iterative}, homomorphic approach \cite{arsenault1984combined,huang2017multiplicative},
wavelet-based approaches \cite{achim2003sar,bhuiyan2007spatially,gagnon1997speckle,guo1994wavelet,solbo2004homomorphic}, 
nonlinear diffusion in the Laplacian pyramid domain \cite{zhang2007nonlinear}, nonlinear diffusion-based methods
\cite{jain2019non,jain2018nonlinear,majee2020gray,shan2019multiplicative,weickert1998anisotropic,yu2002speckle,zhou2015doubly},
 variational methods \cite{aubert2008variational,denis2009sar,feng2014speckle,jin2010analysis,majee2021fuzzy,rudin2003multiplicative,shi2008nonlinear}, telegraph equation \cite{majee2021fuzzy,majee2020gray,majee2022new} and, very recently, deep
learning based approaches \cite{chierchia2017sar,molini2021speckle2void,mullissa2020despecknet,wang2017sar,zhao2014adaptive}. A thorough study of these despeckling techniques is beyond the scope of this work. For a detailed explanation of these methods, we refer the readers to \cite{argenti2013tutorial,denis2021review} and the references therein. This study focuses on developing a variable-exponent telegraph diffusion model for image despeckling and provides rigorous proof of existence and uniqueness. From the beginning of the Perona-Malik (PM) \cite{perona1990scale} model, nonlinear partial differential equations (PDE) have been extensively used to develop noise reduction models. Due to the availability of well-established numerical schemes and theoretical properties, PDE-based image processing is an exciting research area for real-life applications. The PDE-based approach is well-known in the image processing community and aims to remove image noise without destroying meaningful details of the image content, typically edges, lines, or other information that is crucial for understanding an image \cite{aubert2006mathematical,weickert1998anisotropic}. The total variational (TV) based algorithms achieved remarkable results among various PDE-based models.
The first variational-based approach to deal with the multiplicative noise removal problem was proposed by Rudin et al. \cite{rudin2003multiplicative}, which is known as the RLO model.
For an image variable $I$, their model takes the form

\begin{equation}\label{fuzzy_rlo_energy}
I^{*}:=\text{arg}\min_{I\in BV \left(\Omega\right)}  \left\{ \int_{\Omega} |D I|+  \lambda_1 \int_{\Omega}  \frac{I_0}{I} dx +  \lambda_2 \int_{\Omega}\bigg(\frac{I_0}{I}-1 \bigg)^2dx   \right\}\,,
\end{equation} 
where $\lambda_1$ and $\lambda_2$ are two Lagrange multipliers, dynamically updated as explained in \cite{rudin2003multiplicative} and $I_0>0$ is the initial data. Here $BV(\Omega)$ denotes the space of functions with bounded variation \cite{evans2015measure} defined as a space of $L^1$ valued functions on an open, bounded Lipschitz domain $ \Omega \subset \mathbb{R}^n$ such that the following quantity
\begin{align*}
\int_{\Omega}|D I| := \sup_{\varphi}\left\{ \int_{\Omega} I\text{div}(\varphi) dx \big\vert \varphi \in C^1_0 (\Omega; \mathbb{R}^n), \vert \varphi \vert \leq 1 \right\}\,,
\end{align*}
is finite. The problem associated with the evolution equation of \eqref{fuzzy_rlo_energy}, along with the initial and boundary conditions, can be written as
\begin{equation*}
\left.\begin{aligned}
&{I_t} = \text{div}\bigg(\frac{\nabla I}{|\nabla I|}\bigg)+\lambda_1 \frac{I_0}{I^2} +\lambda_2 \frac{I_0^2}{I^3}\, \hspace{0.5cm} \text{in} \hspace{0.2cm} \Omega_T:= \Omega \times (0,T)\,,\\
&I(x,0)=I_0(x)	\hspace{3.35cm} \text{in} \hspace{0.2cm}\Omega\,,\\
&\partial_n I=0  \hspace{4.4cm} \text{in} \hspace{0.2cm} \partial \Omega_T:= \partial \Omega \times (0,T)\,.
\end{aligned}\right\}
\end{equation*}
Here and in the following, $\Omega$ is the domain of the image variable $I$, $T>0$ is a specified time, and $I_t$ denotes the first-time derivative of $I$. $\text{div}$ and $\nabla$ represent the divergence and gradient operator, respectively, and $\partial_n$ denotes the derivative at the boundary surface $ \partial \Omega $ in the outward normal direction $n$. In 2008, Aubert and Aujol \cite{aubert2008variational} employed the concept of speckle noise in the total variation framework. By using the maximum a-posteriori (MAP) estimator, they developed a TV model (the AA model)  for image despeckling. Their model takes the form
\begin{equation}\label{fuzzy_aa_energy}
I^{*}:=\text{arg}\min_{I\in S \left(\Omega\right)}\left\{ \int_{\Omega} |D I|+\lambda \int_{\Omega}\bigg(\log I+ \frac{I_0}{I} \bigg)dx \right\}\,,
\end{equation} 
where $S(\Omega)=\{I \in BV(\Omega), I>0 \}$ and $\lambda$ is a regularization parameter.
Later, Qiang et al. \cite{liu2013nondivergence} proposed a modified version of the AA model for image despeckling, which uses the $p$-Laplace operator with a lower order term and a maximum a posteriori estimator. The model can be expressed as
\begin{equation}\label{nondiv_energy}
I^{*}:=\text{arg}\min_{I\in S \left(\Omega\right)}\left\{ \int_{\Omega} |D I|^p+\lambda \int_{\Omega}\bigg(\log I+ \frac{I_0}{I} \bigg)dx \right\}\,.
\end{equation} 
where $1<p<2$. 
In 2013, Dong et al. \cite{dong2013convex} suggest a convex total variation model for multiplicative speckle-noise reduction with the following form
\begin{equation}\label{eq:Dong_energy}
I^{*}:=\text{arg}\min_{I\in BV \left(\Omega\right)}\left\lbrace \int_{\Omega} \alpha(x)|D I|+\lambda \int_{\Omega}\left( I+I_0 \log\frac{1}{I}\right) dx\right\rbrace.
\end{equation} 
They choose the gray level indicator function $\alpha$, as
\begin{equation}\label{eq:gray_indicator}
\left( 1-\frac{1}{1+k|G_\xi \ast I_0|^2}\right) \frac{1+kM^2}{kM^2},\,\,\,\, \text{or}\,\,\,\,\, \dfrac{G_{\xi}\ast I_0}{M}, \nonumber
\end{equation}
with $M=\sup_{x\in \Omega} (G_\xi \ast I_0)(x)$, where $\xi>0$, $k>0$, ``$\ast$" is the convolution operator, $G_\xi$ is the two dimensional Gaussian kernel and $\lambda$ is a given parameter. 
In addition to TV approaches, diffusion-based filters can also remove multiplicative noise from degraded images. One of the earliest diffusion-based models for multiplicative speckle noise removal was proposed by Yu and Acton \cite{yu2002speckle}, integrating a spatially adaptive filter with the PM model \cite{perona1990scale}. The model takes the form
\begin{equation*}
\left.\begin{aligned}
&{I_t} =\text{div}(g(q_0,q)\nabla I)  \hspace{0.85cm} \text{in} \hspace{0.2cm} \Omega_T\,,\\
&I(x,0)=I_0(x)	\hspace{1.8cm} \text{in} \hspace{0.2cm}\Omega\,,\\
&\partial_n I=0              \hspace{2.8cm} \text{in} \hspace{0.2cm} \partial \Omega_T\,.
\end{aligned}\right\}
\end{equation*}
Here $g(\cdot)$ is the diffusion coefficient, which can be defined as
\begin{equation*}
g(q,q_0)=\left( 1+\frac{q^2-{q_0}^2}{{q_0}^2(1+{q_0}^2)}\right)^{-1},
\end{equation*}
where $q$ is the instantaneous coefficient of variation (ICOV), serves as the edge detector function, and determined by the formula
\begin{equation*}
q(I, \nabla I, \nabla^2 I)=\sqrt{\dfrac{(1/2)(\nabla I/I)^2-(1/16)(\nabla^2 I/I)^2}{[1+(1/4)(\nabla^2 I/I)]^2}}\,,
\end{equation*}
and $q_0$ is the speckle scale function, serves as the diffusion threshold value determined by the ratio of the local standard deviation to mean
\begin{equation*}
q_0(I)=\frac{std(I)}{mean(I)}\,.
\end{equation*}
This filter provides significant enhancement in edge preservation and speckle suppression when compared with conventional filters. Later, based on a gray level indicator function, Zhou et al. proposed a diffusion model called ``doubly degenerate diffusion (DDD) \cite{zhou2015doubly}" for the multiplicative noise removal problem. Their model takes the form
\begin{equation*}
\left.\begin{aligned}
&I_t  = \text{div}(g(I,|\nabla I|)\nabla I),  \hspace{0.5cm} \text{in} \hspace{0.2cm} \Omega_T,  \\
&I(x,0)=I_0(x),	\hspace{1.8cm} \text{in} \hspace{0.2cm}\Omega\,,\\
&\partial_n I=0, \hspace{2.7cm} \text{in} \hspace{0.2cm} \partial \Omega_T\,.
\end{aligned}\right\}
\end{equation*}
They choose the diffusion coefficient as
\begin{align*}
g\left(I,\vert \nabla I \vert \right)=\dfrac{2\vert I \vert^\nu}{M^\nu+\vert I \vert^\nu}\cdot \dfrac{1}{\left(1+ |\nabla I|^2\right)^{(1-\beta)/2} }, 
\end{align*}
where $\nu>0,$ $0<\beta<1,$ and $M=\sup_{x\in \Omega} I$. In this case, the gray level indicator and edge detector functions are $a(I):=\dfrac{2\vert I \vert^\nu}{M^\nu+\vert I \vert^\nu}$  and  $b(I):= \dfrac{1}{\left(1+ |\nabla I|^2\right)^{(1-\beta)/2} }$ respectively. In 2017, Zhou et al. \cite{zhou2018nonlinear} proposed a variable exponent diffusion model for image despeckling. The model takes the form
\begin{equation}\label{eq:ZZDB}
\left.\begin{aligned}
&I_t  = \text{div}\left( \dfrac{\nabla I}{1+\left( \vert \nabla I_\rho \vert/K \right)^{\beta(I)}} \right),  \hspace{0.5cm} \text{in} \hspace{0.2cm} \Omega_T,  \\
&I(x,0)=I_0(x),	\hspace{3.0cm} \text{in} \hspace{0.2cm}\Omega\,,\\
&\partial_n I=0, \hspace{4.0cm} \text{in} \hspace{0.2cm} \partial \Omega_T\,.
\end{aligned}\right\}
\end{equation}
Here $\beta(I)$ is a region indicator, and the authors chose it as follows
\begin{align*}
\beta\left(I \right)=2-\dfrac{2\vert I \vert^\alpha}{M^\alpha+\vert I \vert^\alpha}, 
\end{align*}
where  $K,\alpha>0,$ $I_\rho=G_\rho \ast I,$ and $M=\sup_{x\in \Omega} I$. Also, the authors in \cite{zhou2018nonlinear} establish the existence and uniqueness results of their model. The present paper aims to propose a variable exponent telegraph diffusion model for image despeckling, and we establish the existence and uniqueness of our model.


Over the last few years, many researchers investigated the variable exponent-based diffusion models \cite{baravdish2015backward,blomgren1997total,chen2006variable,guo2011adaptive,li2010variable,theljani2019high,zhou2018nonlinear} for image denoising. Since the variable exponent models utilize the benefits of isotropic diffusion, TV diffusion, and anisotropic diffusion depending on the values of the exponent, they have some advantages over fixed exponent-based models for image restoration processes. Furthermore, telegraph diffusion-based methods \cite{cao2010class,baravdish2019damped,majee2021fuzzy,majee2020gray,majee2022new,ratner2007image,zhang2014class,zhang2015spatial} have been widely explored for removing both additive and multiplicative noise. 

Despite the promising applications of variable-exponent diffusion models in noise reduction and theoretical studies, the use of variable-exponent telegraph diffusion models remains unexplored. To the best of our knowledge, no previous research has proposed variable exponents in a telegraph equation framework for speckle noise elimination. Furthermore, the theoretical results are unexplored in this setting. Based on existing results in the literature, it seems that the variable exponent telegraph equation could be a potentially effective technique for speckle noise removal. In this study, we first analyze existing variable-exponent diffusion models for image despeckling and then extend these models to a telegraph diffusion framework. Additionally, we establish the existence and uniqueness of weak solutions for both diffusion and telegraph diffusion models with variable exponents by considering a general form for the diffusion coefficient. These results can be directly applied to similar PDE models to check their wellposedness.

The remainder of the paper is organized as follows. In Section \ref{diffusion_discussion}, we review existing diffusion-based models with variable exponents and introduce some modifications to the exponents. We then discuss existing telegraph diffusion models and propose new hybrid telegraph diffusion models for image despeckling. Section \ref{sec:analysis} presents the existence and uniqueness results for both the diffusion and telegraph diffusion models. In Section \ref{sec:numerical}, we describe the numerical implementation of our model, and Section \ref{computation_discussion} evaluates the despeckling performance of the proposed approach. Finally, Section \ref{sec:Conclusion} provides the conclusion.

\section{Models With Different Variable Exponents}
\label{diffusion_discussion}
\subsection{Diffusion Models}
\begin{figure}
       \centering
       \begin{subfigure}[b]{0.19\textwidth}           
                \includegraphics[scale=0.29]{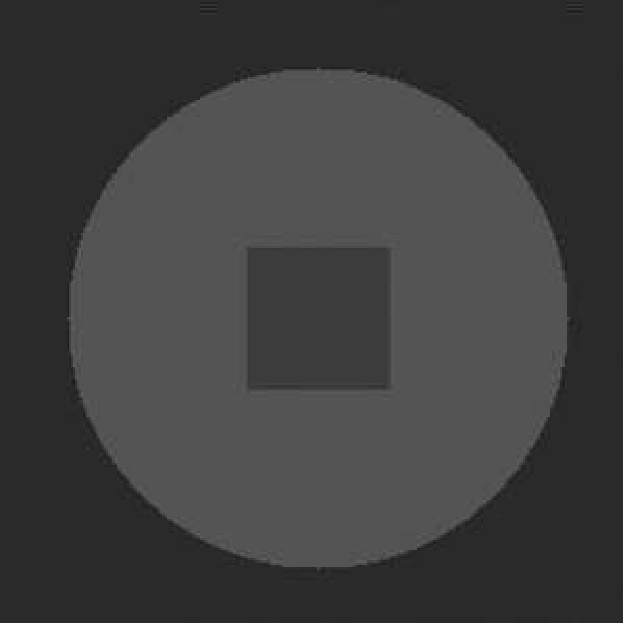}               
                \caption{Clear}
                \label{circle_clean1}
       \end{subfigure}
      \begin{subfigure}[b]{0.19\textwidth}           
                \includegraphics[scale=0.29]{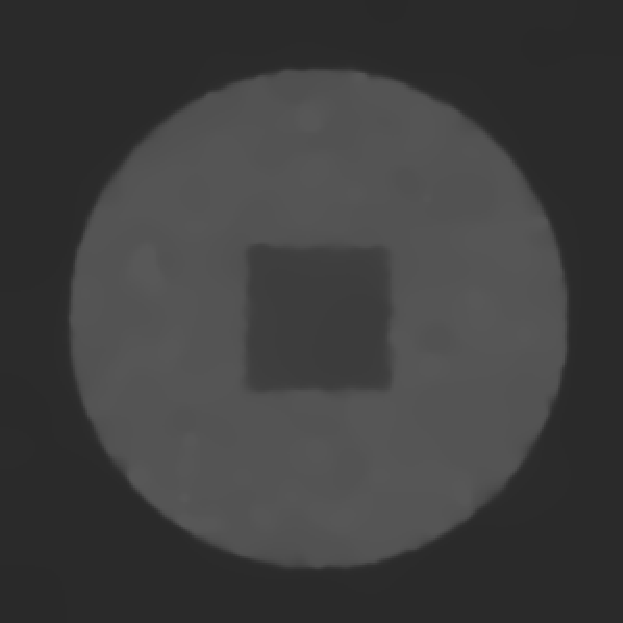}               
                \caption{$p=1.0$}
                \label{circle10_p_1.0}
       \end{subfigure}
            \begin{subfigure}[b]{0.19\textwidth}           
                \includegraphics[scale=0.29]{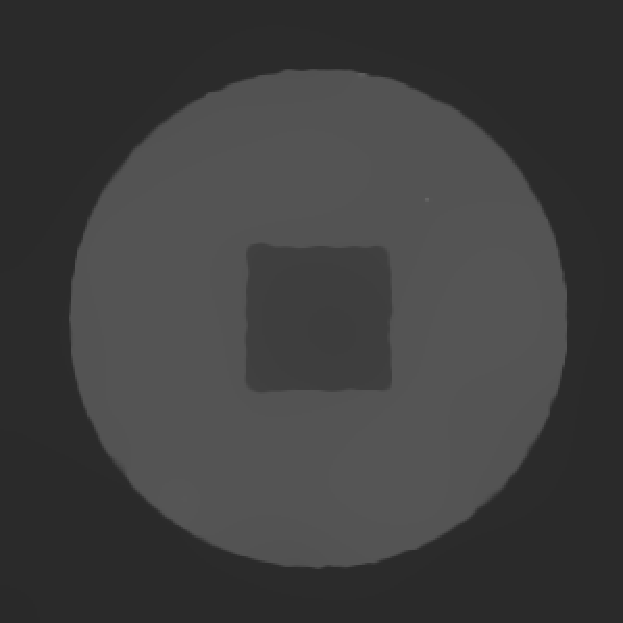}               
                \caption{$p=1.5$}
                \label{circle10_p_1.5}
       \end{subfigure}
       \begin{subfigure}[b]{0.19\textwidth}           
                \includegraphics[scale=0.29]{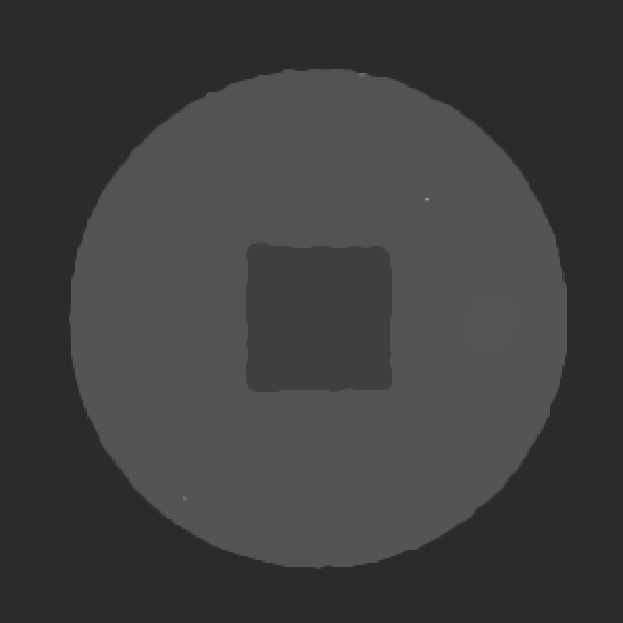}               
                \caption{$p=2.0$}
                \label{circle10_p_2.0}
       \end{subfigure}
        \begin{subfigure}[b]{0.19\textwidth}           
                \includegraphics[scale=0.29]{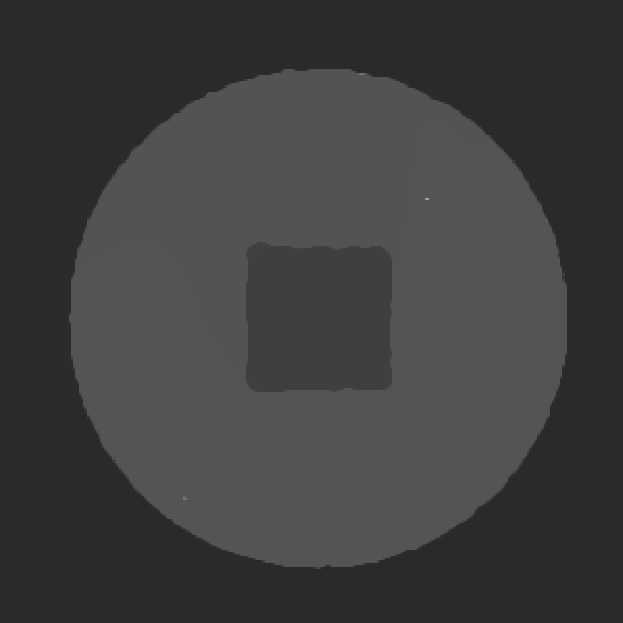}               
                \caption{$p=2.5$}
                \label{circle10_p_2.5}
       \end{subfigure}

       \begin{subfigure}[b]{0.19\textwidth}           
                \includegraphics[scale=0.29]{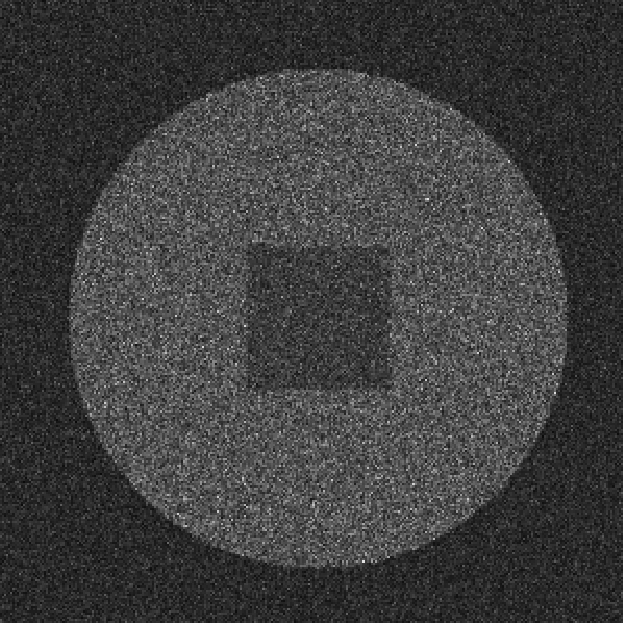}               
                \caption{Noisy}
                \label{circle_noisy101}
       \end{subfigure}
       \begin{subfigure}[b]{0.19\textwidth}           
                \includegraphics[scale=0.29]{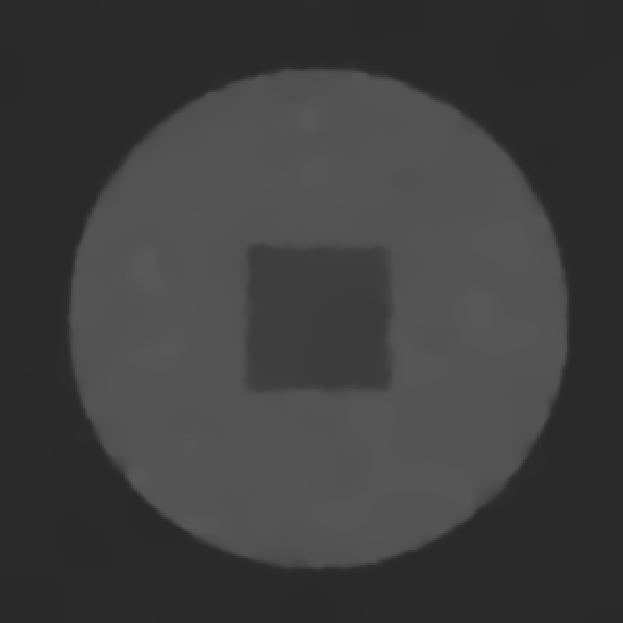}               
                \caption{$p=1.0$}
                \label{circle10_diffgrayp_1.0}
       \end{subfigure}
            \begin{subfigure}[b]{0.19\textwidth}           
                \includegraphics[scale=0.29]{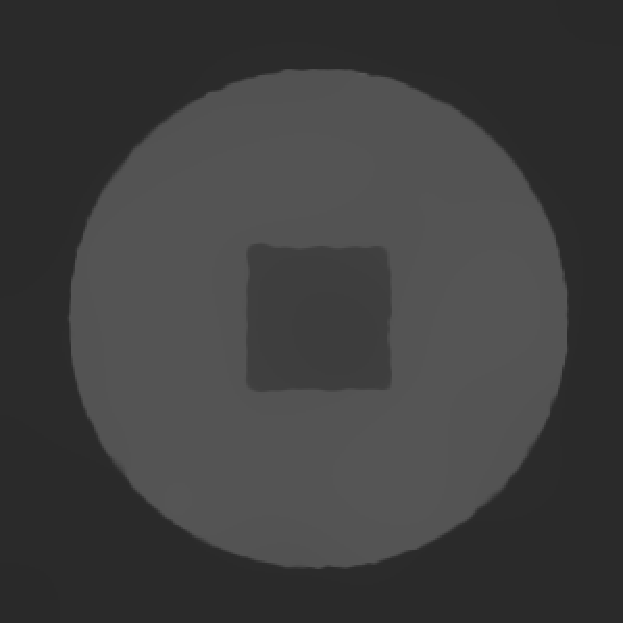}               
                \caption{$p=1.5$}
                \label{circle10_diffgrayp_1.5}
       \end{subfigure}
       \begin{subfigure}[b]{0.19\textwidth}           
                \includegraphics[scale=0.29]{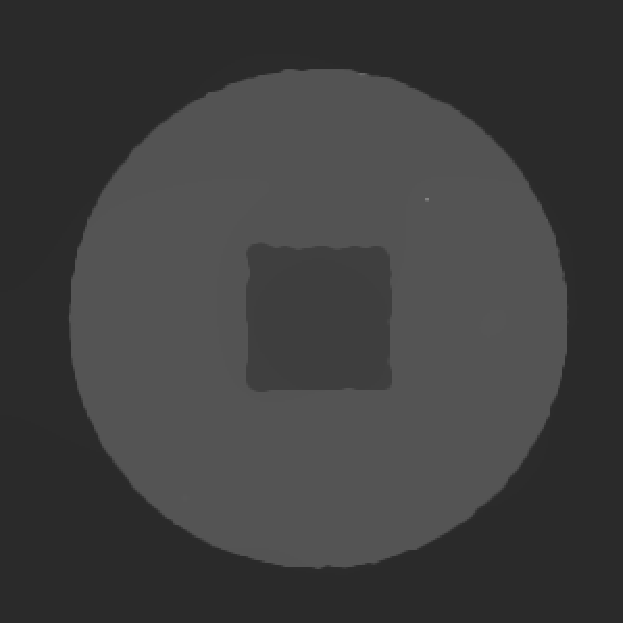}               
                \caption{$p=2.0$}
                \label{circle10_diffgrayp_2.0}
       \end{subfigure}
        \begin{subfigure}[b]{0.19\textwidth}           
                \includegraphics[scale=0.29]{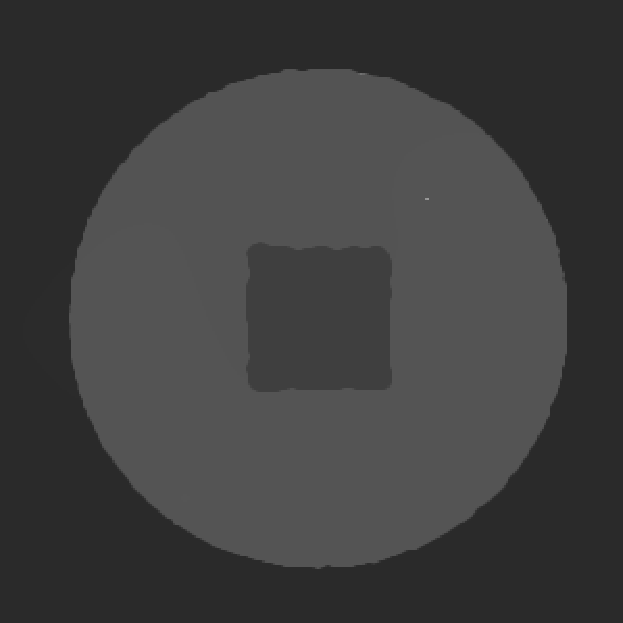}               
                \caption{$p=2.5$}
                \label{circle10_diffgrayp_2.5}
       \end{subfigure}
       
\caption{\footnotesize (a) Clean Image (f) Noisy Image ($L=10$). \textbf{First row:} (b)-(e) Restored Images using \eqref{diff_power_constant} with $\nu=0$ and different values of $p$. (b) PSNR=40.53, (c) PSNR=42.98, (d) PSNR=42.72, (e) PSNR=42.08. \textbf{Second row:}  (g)-(j) Restored Images using \eqref{diff_power_constant} with $\nu=1$ and different values of $p$. (g) PSNR=40.76, (h) \textbf{PSNR=43.47}, (i) PSNR=43.22, (j) PSNR=42.58.}\label{circle_10_restored_diffusion_constp}
\end{figure}

\begin{figure}
       \centering
      \begin{subfigure}[b]{0.19\textwidth}           
                \includegraphics[scale=0.29]{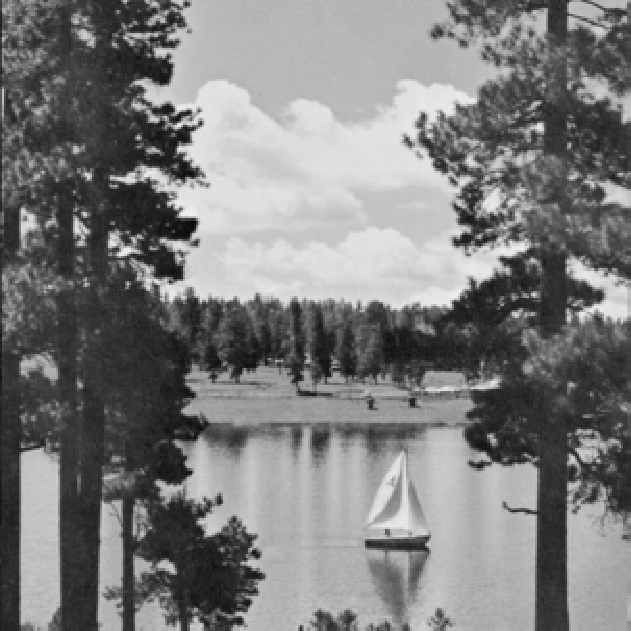}               
                \caption{Clear}
                \label{lake11_clean}
       \end{subfigure}
      \begin{subfigure}[b]{0.19\textwidth}           
                \includegraphics[scale=0.29]{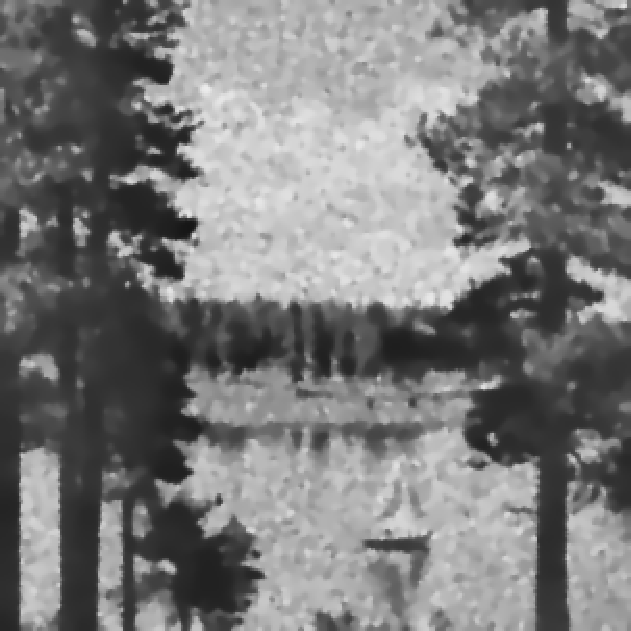}               
                \caption{$p=1.0$}
                \label{lake110_p_1.0}
       \end{subfigure}
       \begin{subfigure}[b]{0.19\textwidth}           
                \includegraphics[scale=0.29]{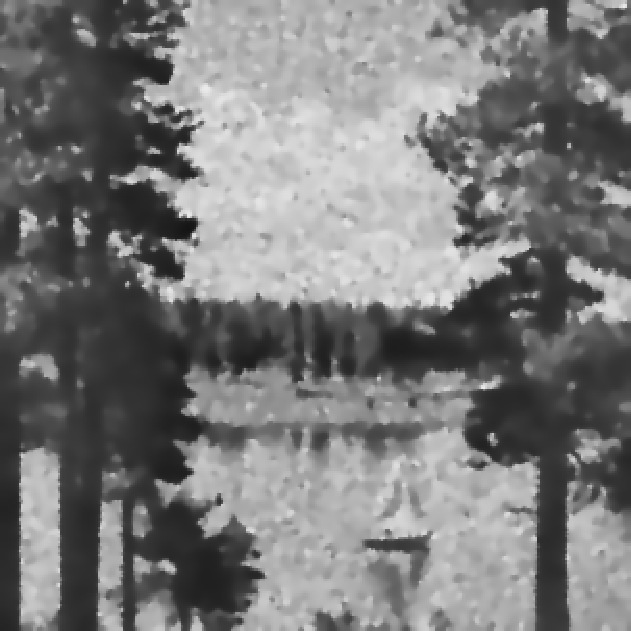}               
                \caption{$p=1.5$}
                \label{lake110_p_1.5}
       \end{subfigure}
       \begin{subfigure}[b]{0.19\textwidth}           
                \includegraphics[scale=0.29]{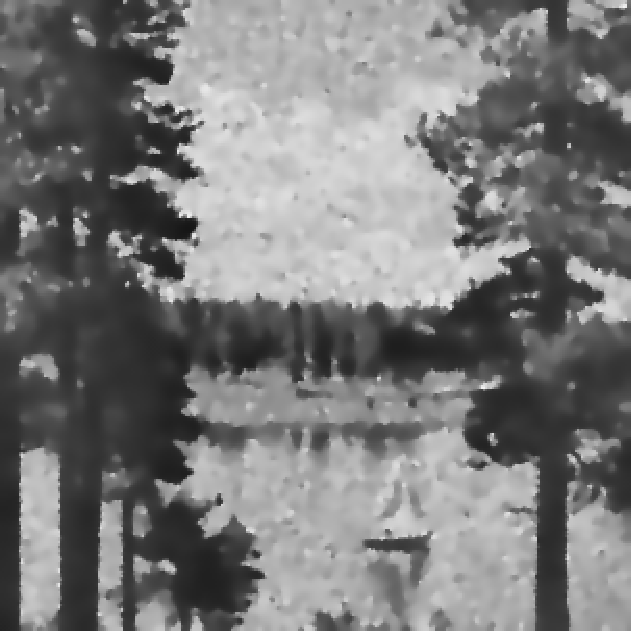}               
                \caption{$p=2.0$}
                \label{lake110_p_2.0}
       \end{subfigure}
        \begin{subfigure}[b]{0.19\textwidth}           
                \includegraphics[scale=0.29]{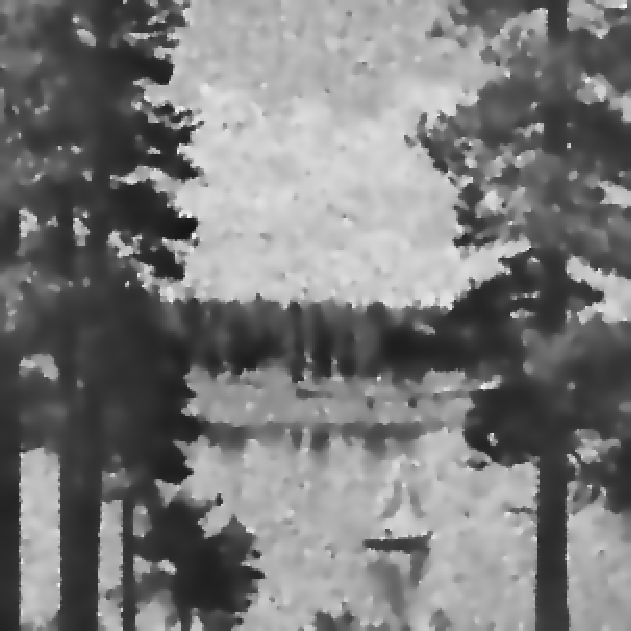}               
                \caption{$p=2.5$}
                \label{lake110_p_2.5}
       \end{subfigure}

      \begin{subfigure}[b]{0.19\textwidth}           
                \includegraphics[scale=0.29]{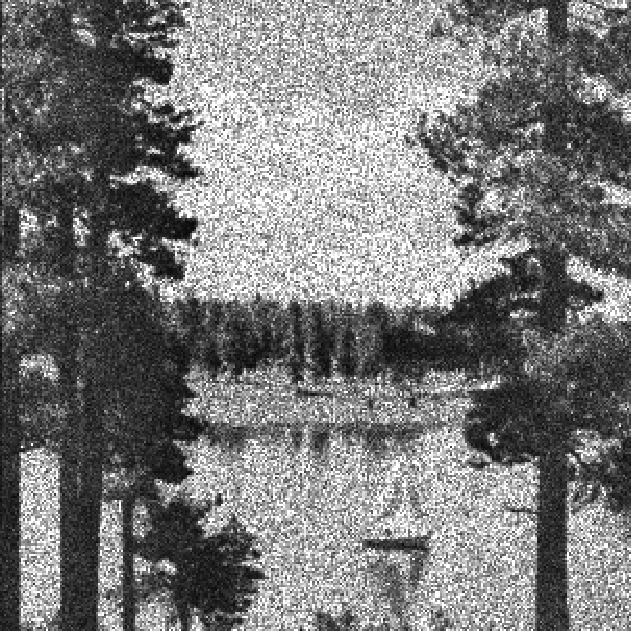}               
                \caption{Noisy}
                \label{lake11_noisy10}
       \end{subfigure}
      \begin{subfigure}[b]{0.19\textwidth}           
                \includegraphics[scale=0.29]{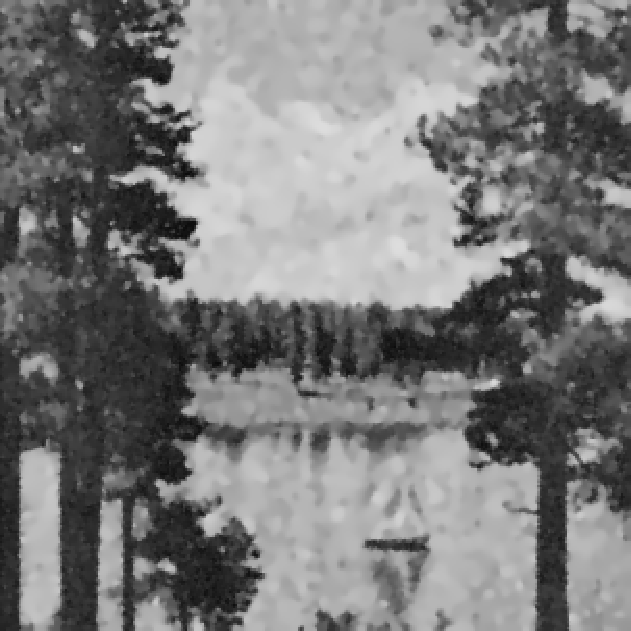}               
          \caption{$p=1.0$}
                \label{lake110_diffgray_p_1.0}
       \end{subfigure}
       \begin{subfigure}[b]{0.19\textwidth}           
                \includegraphics[scale=0.29]{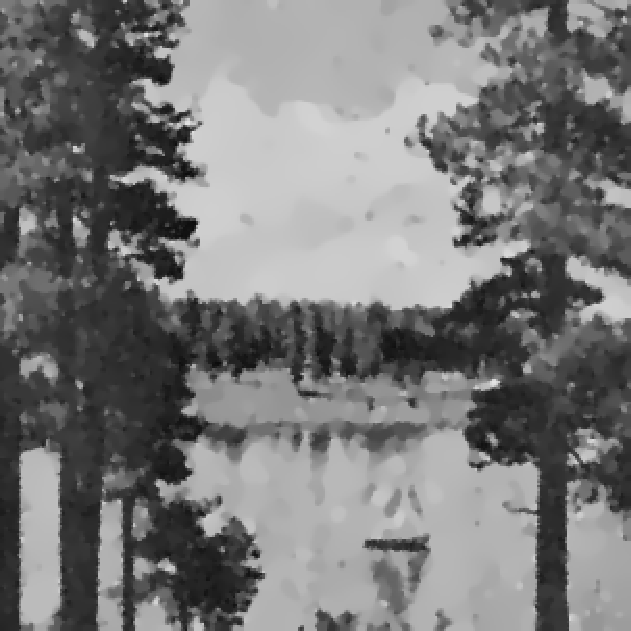}               
                \caption{$p=1.5$}
                \label{lake110_diffgray_p_1.5}
       \end{subfigure}
       \begin{subfigure}[b]{0.19\textwidth}           
                \includegraphics[scale=0.29]{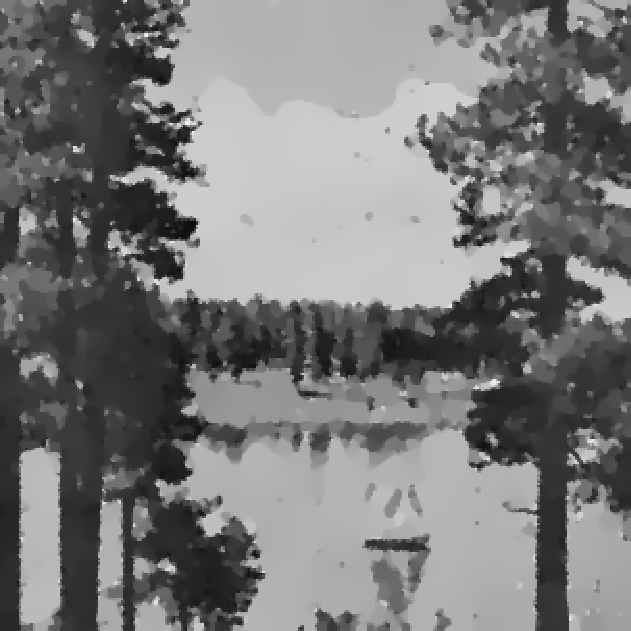}               
                \caption{$p=2.0$}
                \label{lake110_diffgray_p_2.0}
       \end{subfigure}
        \begin{subfigure}[b]{0.19\textwidth}           
                \includegraphics[scale=0.29]{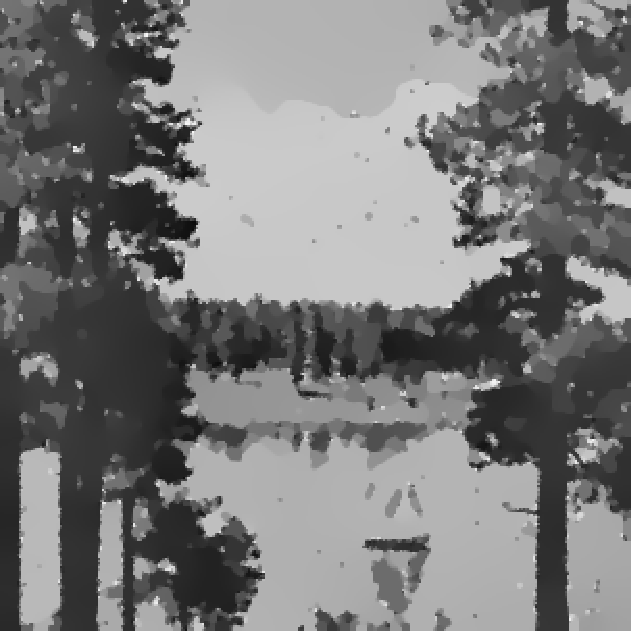}               
                \caption{$p=2.5$}
                \label{lake110_diffgray_p_2.5}
       \end{subfigure}
       
\caption{\footnotesize (a) Clean Image. (f) Noisy Image ($L=10$). \textbf{First row:} (b)-(e) Restored Images using \eqref{diff_power_constant} with $\nu=0$ and different values of $p$. (b) PSNR=24.51, (c) PSNR=24.52, (d) PSNR=24.43, (e) PSNR=24.30. \textbf{Second row:} (g)-(j) Restored Images using \eqref{diff_power_constant} with $\nu=2$ and different values of $p$. (g) PSNR=25.20, (h) \textbf{PSNR=25.26}, (i) PSNR=24.90, (j) PSNR=24.35.}\label{lake_10_restored_diffusion_constp}
\end{figure}

\begin{figure}
       \centering
       \begin{subfigure}[b]{0.19\textwidth}           
                \includegraphics[scale=0.17]{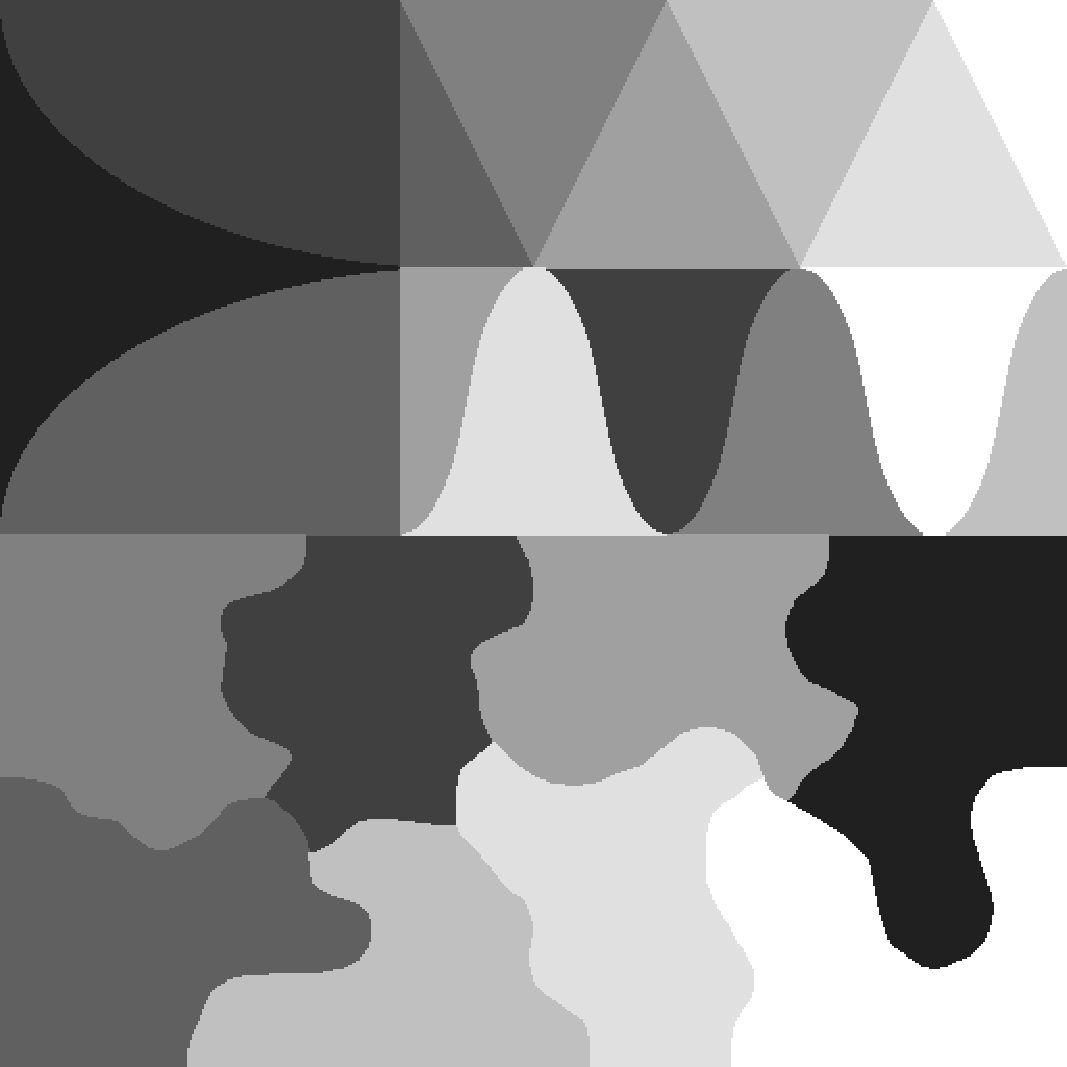}               
                \caption{Clear}
                \label{texture_clean}
       \end{subfigure}
       \begin{subfigure}[b]{0.19\textwidth}           
                \includegraphics[scale=0.17]{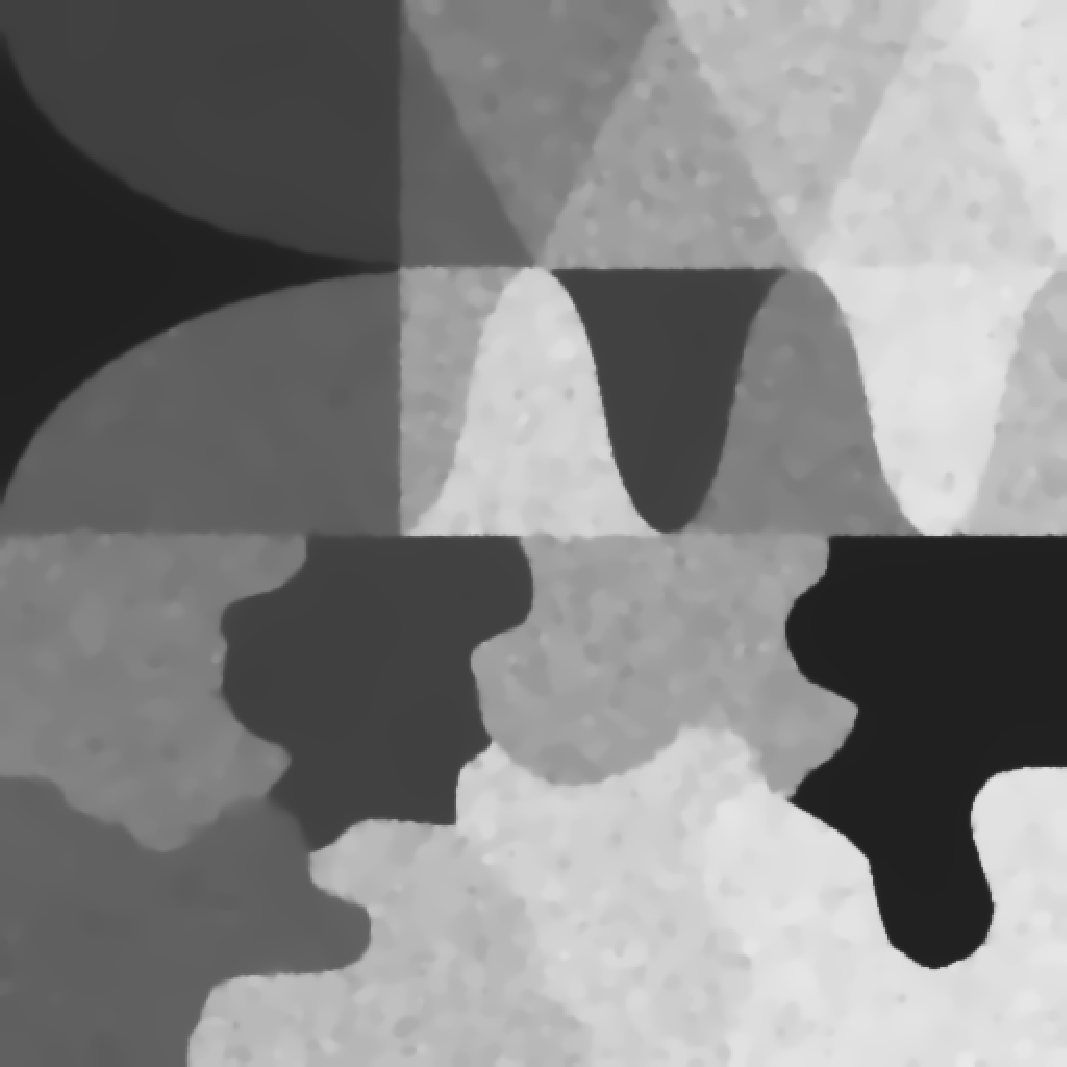}               
                \caption{$p=1.0$}
                \label{texture10_p_1.0}
       \end{subfigure}
            \begin{subfigure}[b]{0.19\textwidth}           
                \includegraphics[scale=0.17]{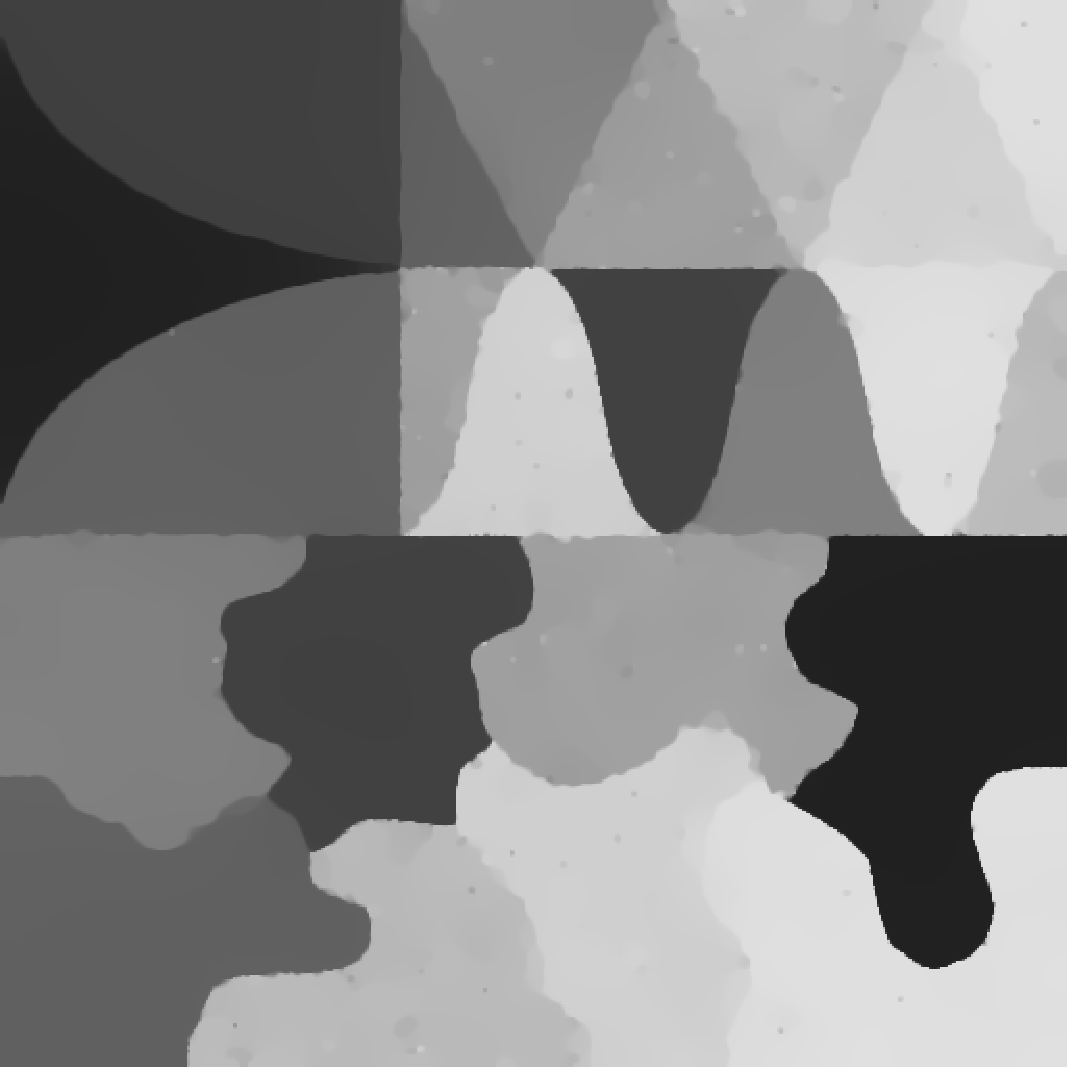}               
                \caption{$p=1.5$}
                \label{texture10_p_1.5}
       \end{subfigure}
       \begin{subfigure}[b]{0.19\textwidth}           
                \includegraphics[scale=0.17]{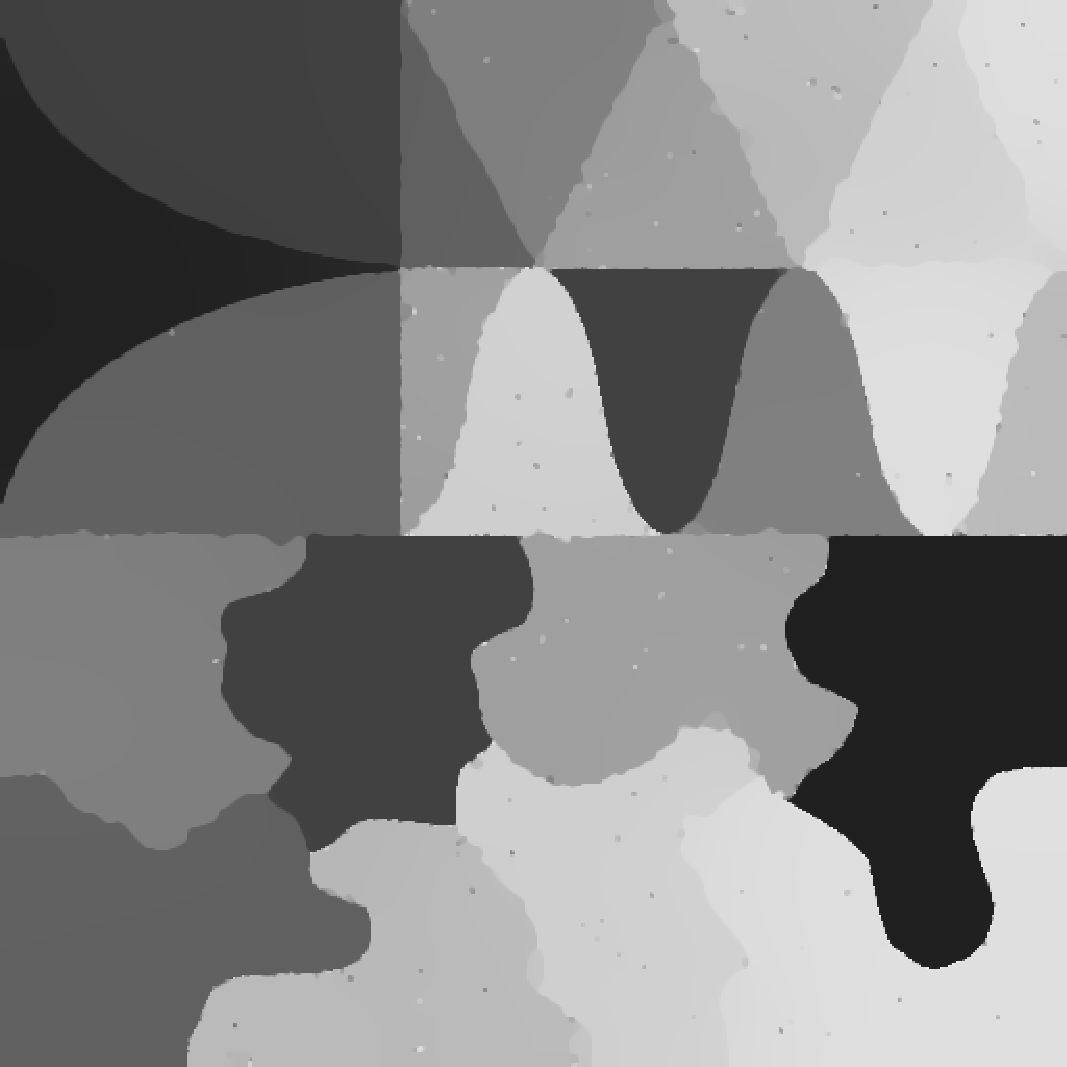}               
                \caption{$p=2.0$}
                \label{texture10_p_2.0}
       \end{subfigure}
        \begin{subfigure}[b]{0.19\textwidth}           
                \includegraphics[scale=0.17]{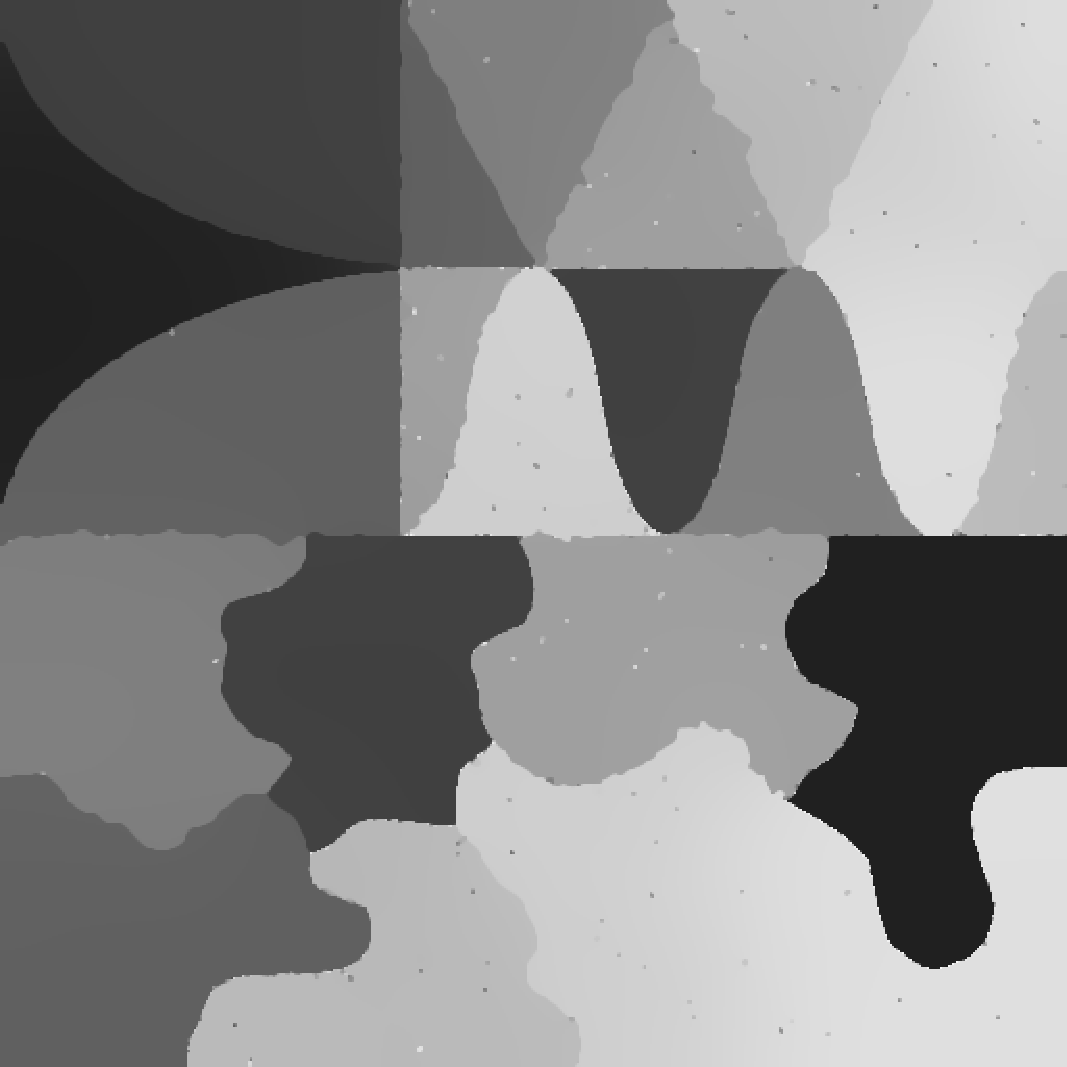}               
                \caption{$p=2.5$}
                \label{texture10_p_2.5}
       \end{subfigure}

       \begin{subfigure}[b]{0.19\textwidth}           
                \includegraphics[scale=0.17]{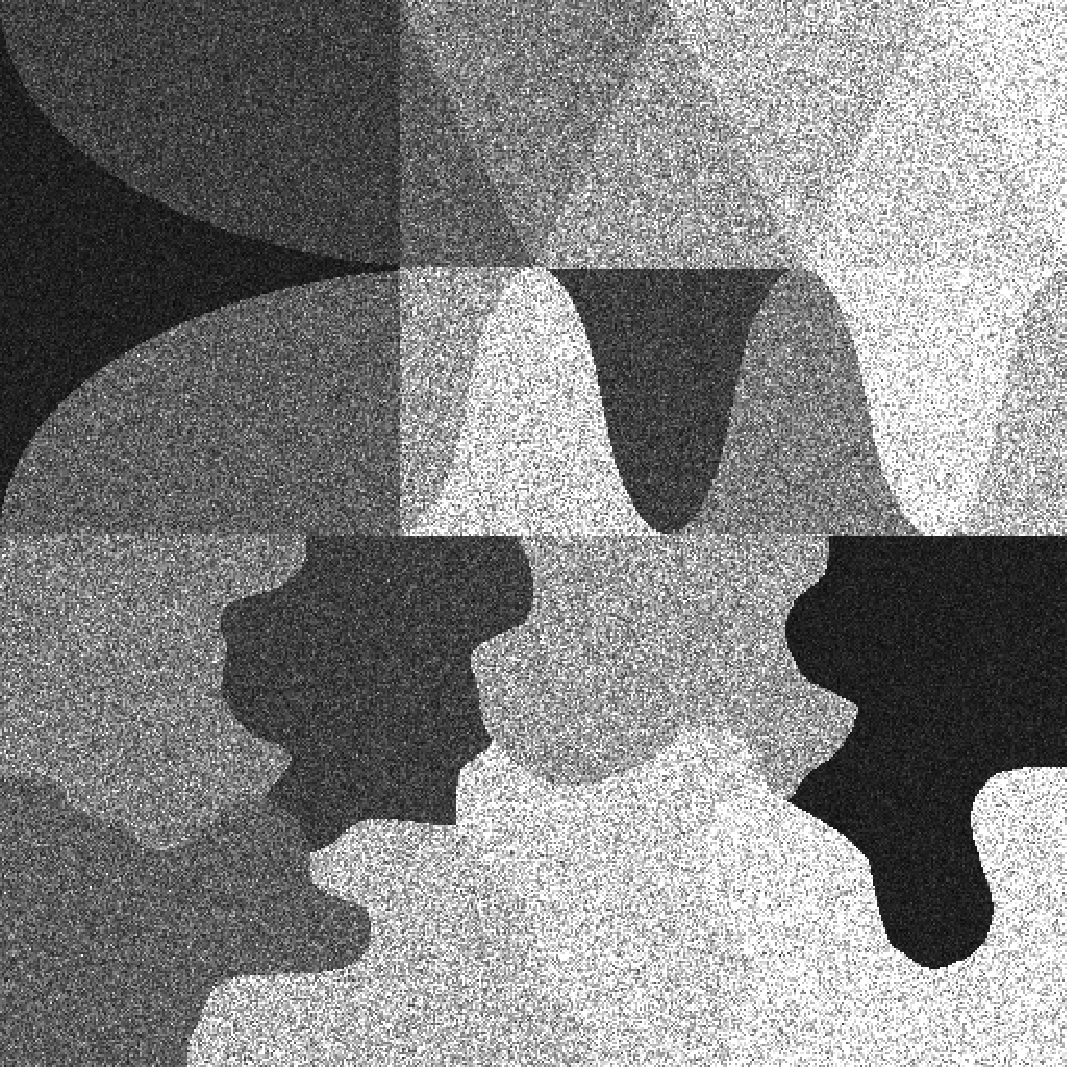}               
                \caption{Noisy}
                \label{texture_noisy10}
       \end{subfigure} 
      \begin{subfigure}[b]{0.19\textwidth}           
                \includegraphics[scale=0.17]{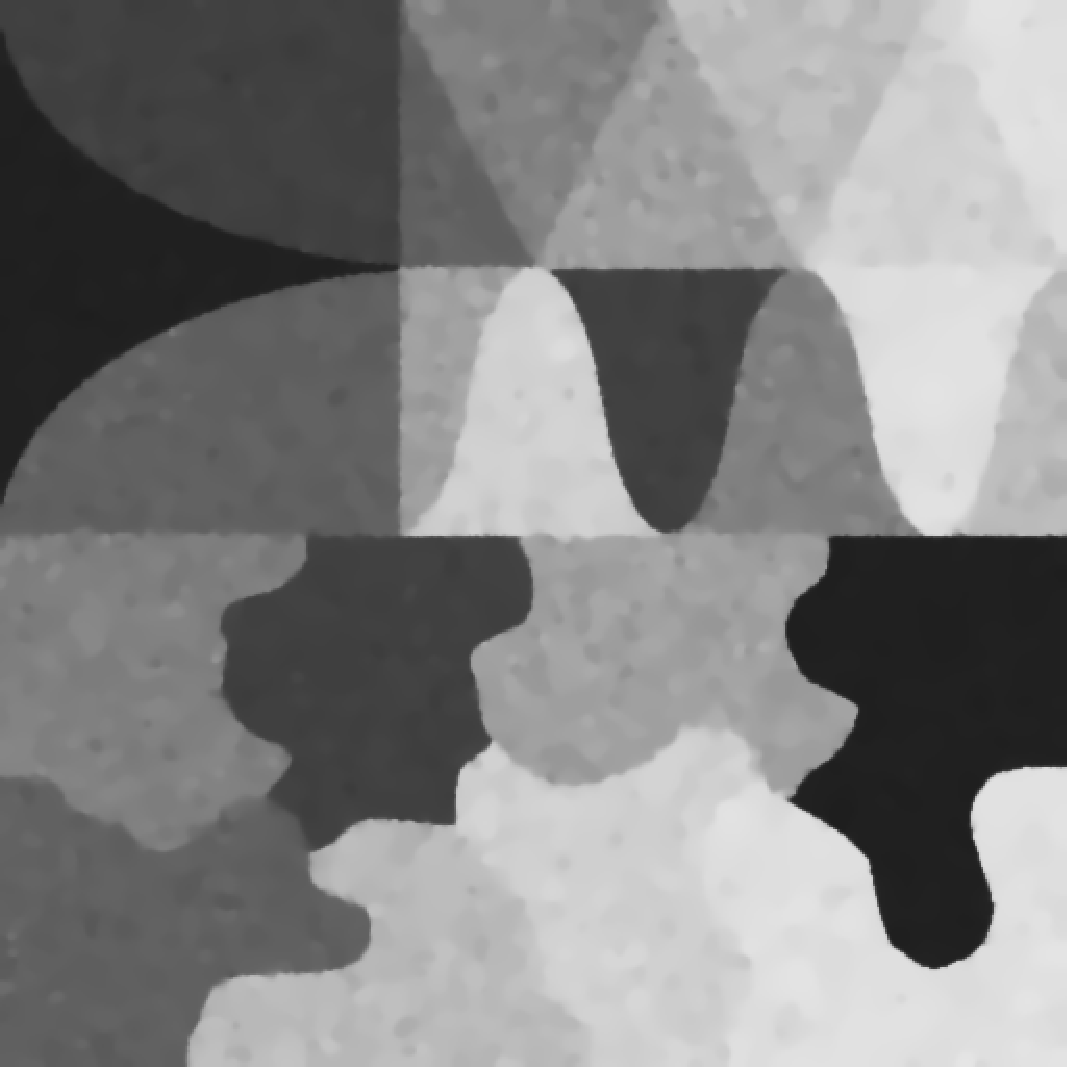}               
                \caption{$p=1.0$}
                \label{texture10_diffgrayp_1.0}
       \end{subfigure}
            \begin{subfigure}[b]{0.19\textwidth}           
                \includegraphics[scale=0.17]{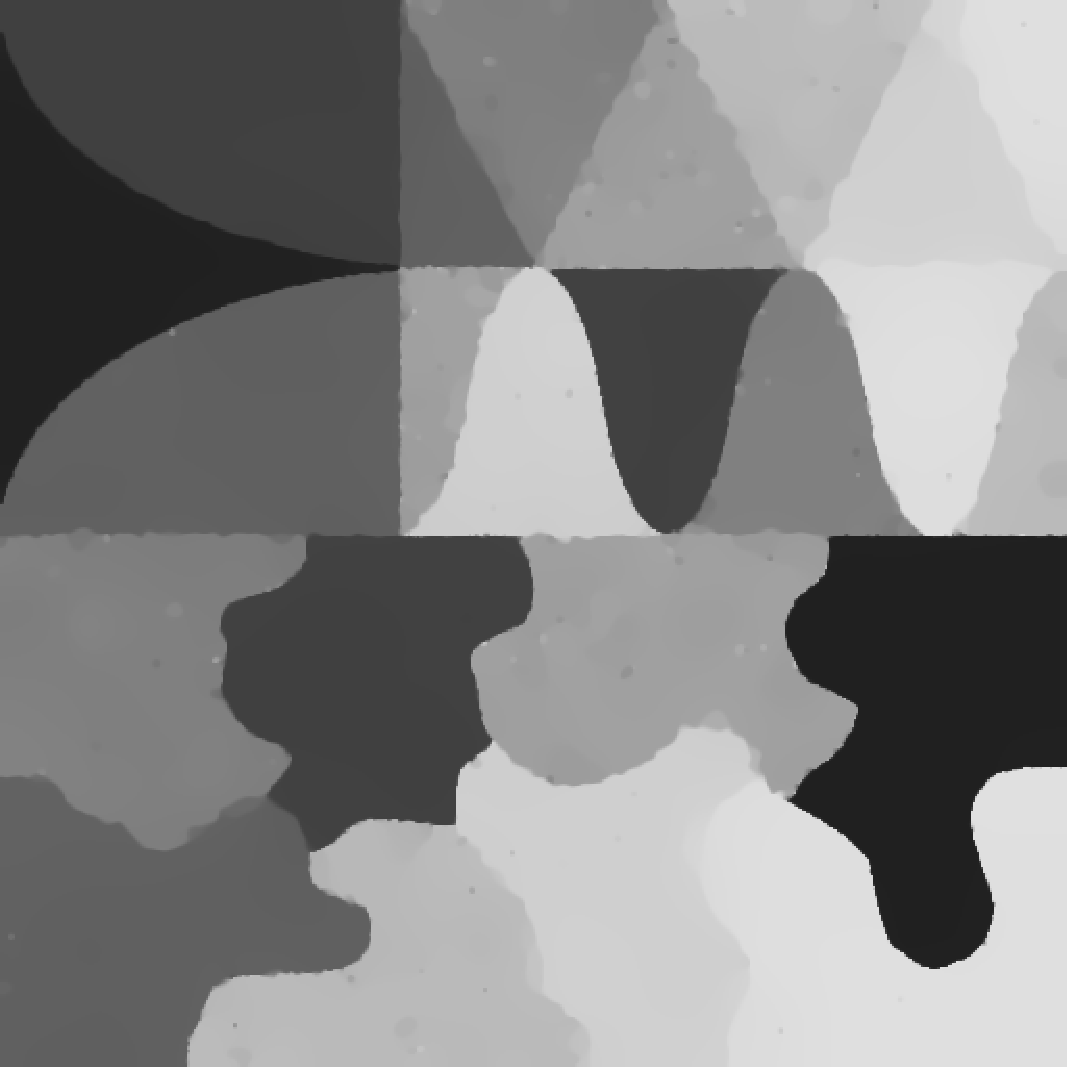}               
                \caption{$p=1.5$}
                \label{texture10_diffgrayp_1.5}
       \end{subfigure}
       \begin{subfigure}[b]{0.19\textwidth}           
                \includegraphics[scale=0.17]{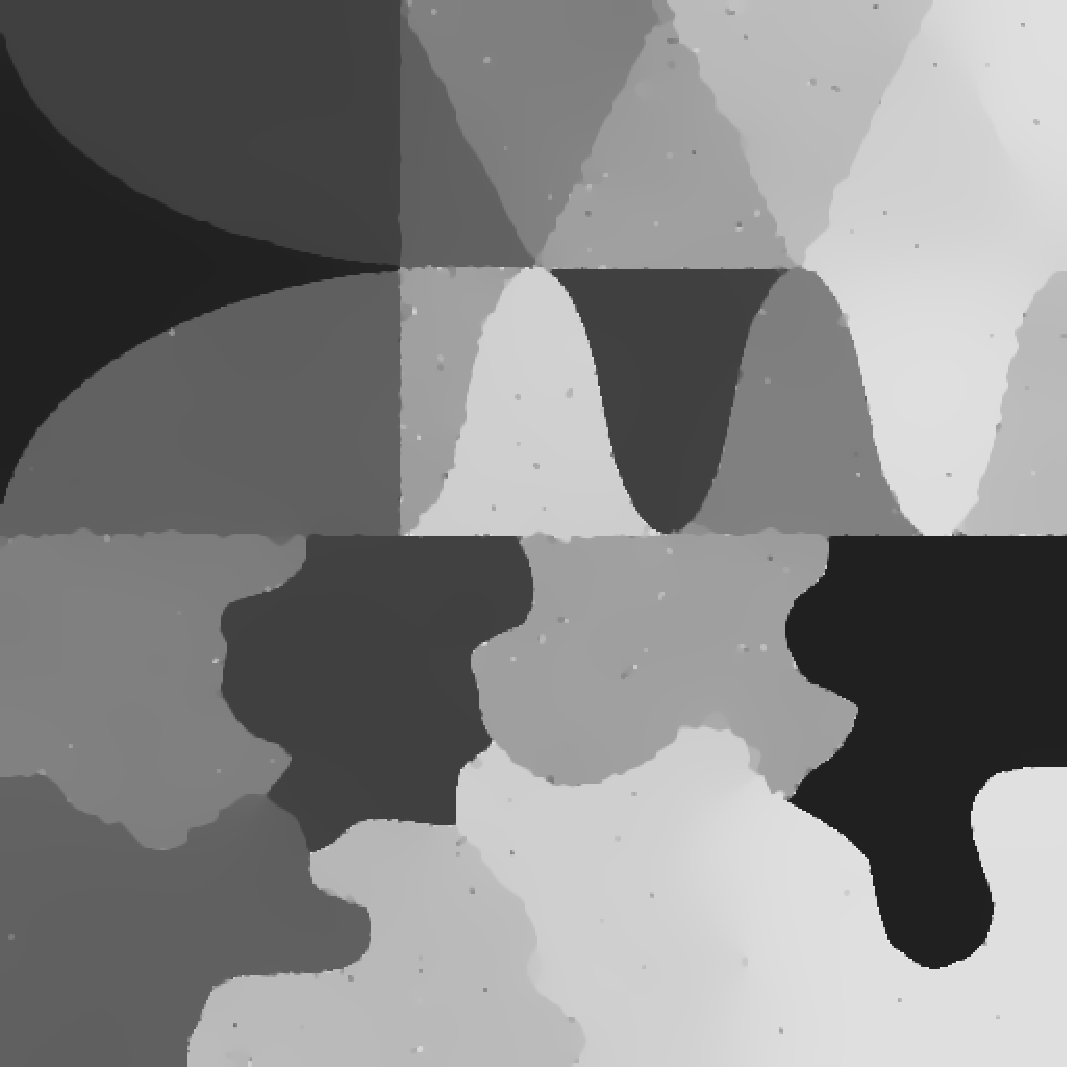}               
                \caption{$p=2.0$}
                \label{texture10_diffgrayp_2.0}
       \end{subfigure}
        \begin{subfigure}[b]{0.19\textwidth}           
                \includegraphics[scale=0.17]{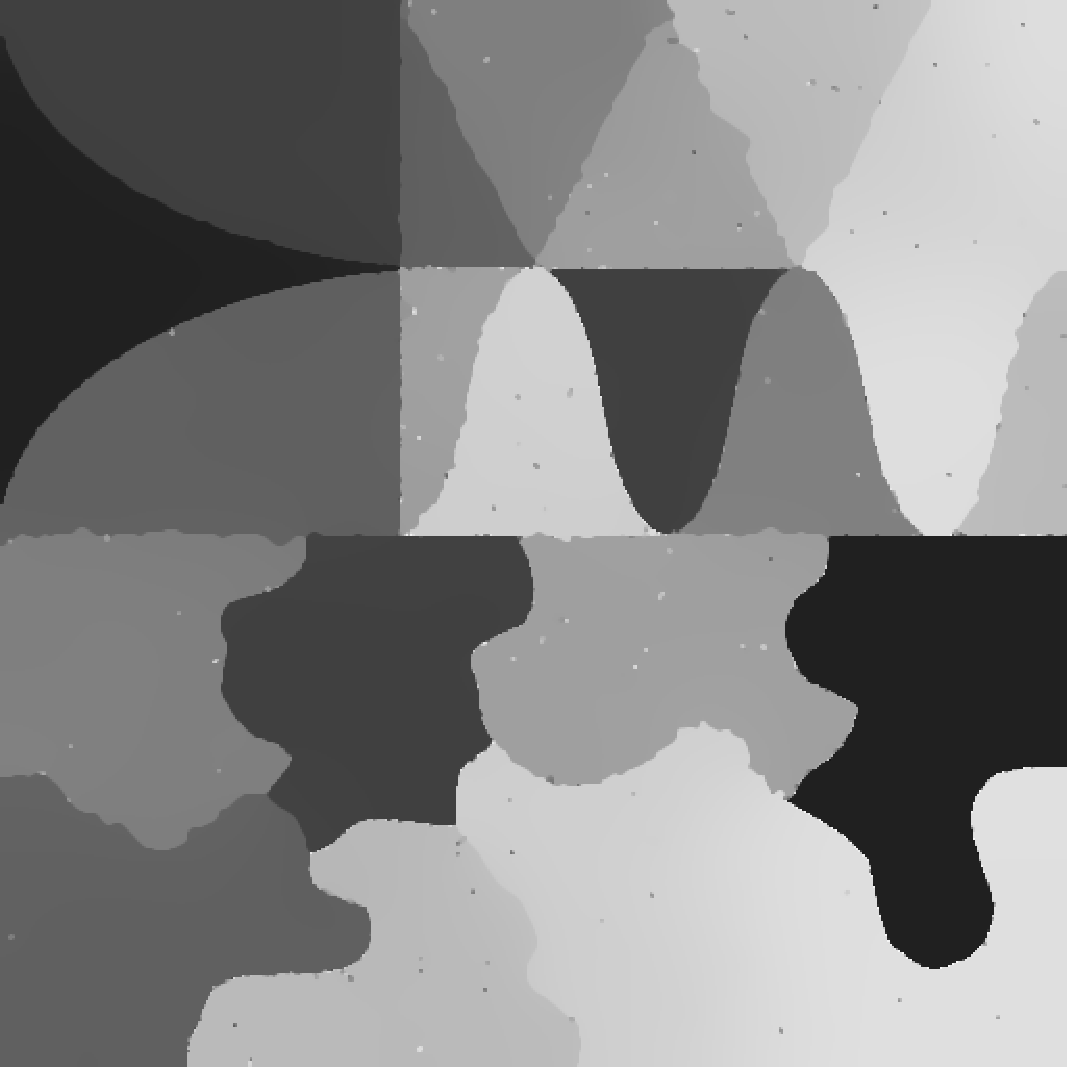}               
                \caption{$p=2.5$}
                \label{texture10_diffgrayp_2.5}
       \end{subfigure}
       
\caption{\footnotesize (a) Clean Image. (f) Noisy Image ($L=10$). \textbf{First row:} (b)-(e) Restored Images using \eqref{diff_power_constant} with $\nu=0$ and different values of $p$. (b) PSNR=25.17 (c) PSNR=25.59, (d) PSNR=25.58 (e) PSNR=25.40. \textbf{Second row:} (g)-(j) Restored Images using \eqref{diff_power_constant} with $\nu=2$ and different values of $p$. (g) PSNR=25.23, (h) \textbf{PSNR=25.61}, (i) PSNR=25.56, (j) PSNR=25.42.}\label{texture_10_restored_diffusion_constp}
\end{figure}

\noindent
In this section, we discuss the image despeckling results using the diffusion model with varying exponents. First, we discuss the image despeckling results using the following diffusion model \eqref{diff_power_constant} with the diffusion coefficient described in \eqref{diff_constatp}. We calculate the results for $\nu \geq 0$ with different values of $p$. For $\nu=0$ and $p=2$, \eqref{diff_power_constant} reduces to the Catt{\'e} model \cite{catte1992image}, which was developed for the additive Gaussian noise removal process. For a non zero $\nu$, \eqref{diff_power_constant} behaves like the DDD model \cite{zhou2015doubly}, which was developed for the speckle noise removal process, introduced in Section \ref{Introduction}. Speckle noise degraded the high gray level regions more compared to the low gray level regions \cite{zhou2015doubly}. So for speckle noise elimination, it is reasonable to smooth the high gray level regions more than the low gray level regions, which can be controlled by the gray level function $a(I_\xi)=2\vert I_\xi \vert^\nu/\left\{\left(M^{I}_{\xi}\right)^\nu + | I_\xi|^\nu\right\}$. It is easy to see that $a(I_\xi)$ becomes very small at low gray levels, which leads the diffusion coefficient close to zero and preserves low gray level image features. In high gray level regions, $a(I_\xi)$ approaches one and leads to the fact that 
$1/\left\{1+ \left(\frac{\vert \nabla I_\xi \vert }{K} \right)^{p}\right\}$ mainly controls the diffusion process.
\begin{equation}\label{diff_power_constant}
\left.\begin{aligned}
&I_t  = \text{div} \left( g\left(I_\xi,\vert \nabla I_\xi \vert \right) \nabla I\right)  \hspace{0.8cm} \text{in} \hspace{0.2cm} \Omega_T,  \\
&I(x,0)=I_0(x)	\hspace{2.6cm} \text{in} \hspace{0.2cm}\Omega\,,\\
&\partial_n I=0  \hspace{3.5cm} \text{in} \hspace{0.2cm} \partial \Omega_T\,.
\end{aligned}\right\}
\end{equation}
The diffusion coefficient $g$ is given by \begin{align}\label{diff_constatp}
g\left(I_\xi, \vert \nabla I_\xi \vert \right)=\dfrac{ 2\vert I_\xi \vert^\nu}{\big(M^{I}_{\xi}\big)^\nu+\vert I_\xi \vert^\nu}\cdot \dfrac{1}{1+ \left(\frac{\vert \nabla I_\xi \vert }{K} \right)^{p}}\,.
\end{align}
Here, $\nu \geq 0,$ $K>0$, $p>0$  are constants, $I_\xi=G_\xi\ast I$, $M_\xi^I= \max_{x\in \Omega} \vert I_\xi (x,t) \vert $.
Figures \ref{circle_10_restored_diffusion_constp}--\ref{texture_10_restored_diffusion_constp} show the restored results of three different images using \eqref{diff_power_constant}. In Figure \ref{circle_10_restored_diffusion_constp}, we observe that for both cases (\(\nu = 0\) and \(\nu > 0\)), the model \eqref{diff_power_constant} restores the image similarly as the gray level information is piecewise uniform. However, the quantitative results indicate that \eqref{diff_power_constant} performs better when \(\nu > 0\). Figure \ref{lake_10_restored_diffusion_constp} presents the results for a lake image, where the gray level information is nonuniform (with color variations across regions). Here, the performance of the two cases is more distinguishable, and the quantitative results align with the qualitative observations. Furthermore, for both images, we conclude that \eqref{diff_power_constant} yields the best results near \(p = 1.5\). Similar observations can be made for the results in Figure \ref{texture_10_restored_diffusion_constp}, as the results in Figure \ref{circle_10_restored_diffusion_constp}.

\begin{figure}
    \centering
       \begin{subfigure}[b]{0.3\textwidth}           
                \includegraphics[scale=0.45]{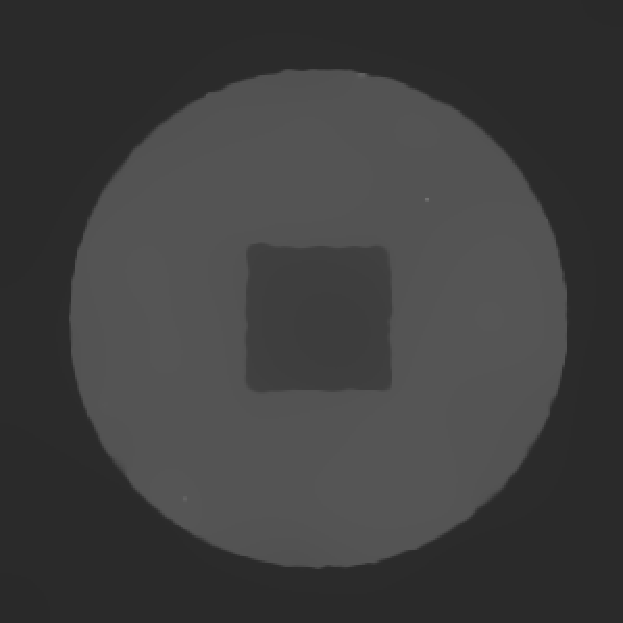}               
                \caption{\textbf{PSNR=42.83}}
                \label{circle10avggrayexp}
       \end{subfigure}
       \begin{subfigure}[b]{0.3\textwidth}           
                \includegraphics[scale=0.45]{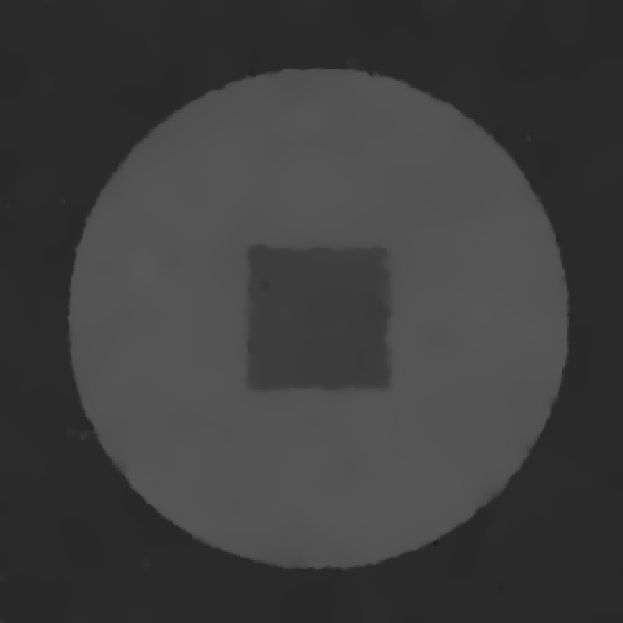}               
                 \caption{PSNR=41.79}
                \label{circle10grayexp}
       \end{subfigure}
       \begin{subfigure}[b]{0.3\textwidth}           
                \includegraphics[scale=0.45]{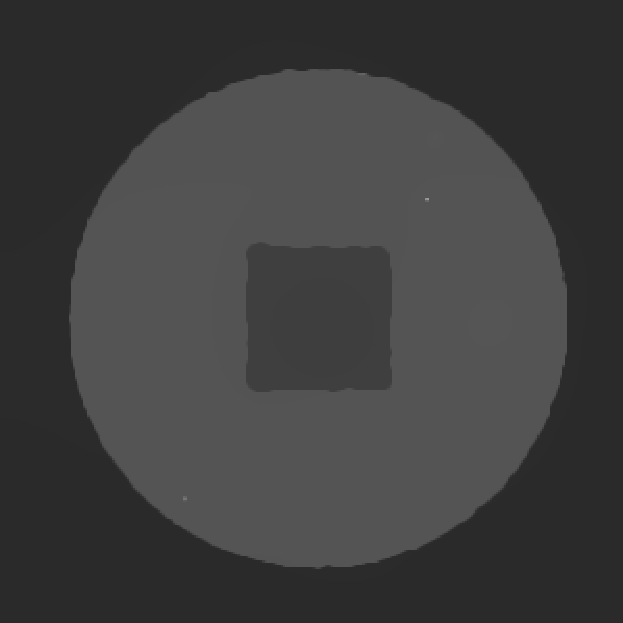}               
                \caption{PSNR=42.71}
                \label{circle10gradexp}
       \end{subfigure}
       
        \begin{subfigure}[b]{0.3\textwidth}           
                \includegraphics[scale=0.45]{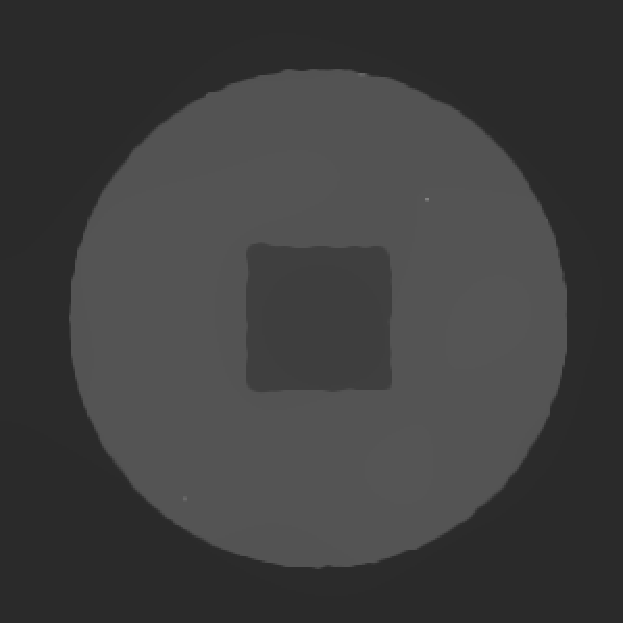}               
                \caption{PSNR=43.03}
                \label{circle10avggrayexp0}
       \end{subfigure}
       \begin{subfigure}[b]{0.3\textwidth}           
                \includegraphics[scale=0.45]{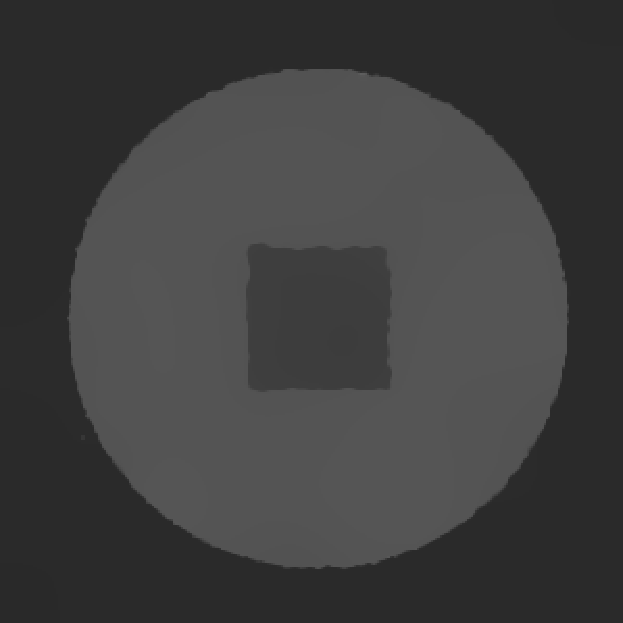}               
                 \caption{\textbf{PSNR=43.10}}
                \label{circle10grayexp0}
       \end{subfigure}
       \begin{subfigure}[b]{0.3\textwidth}           
                \includegraphics[scale=0.45]{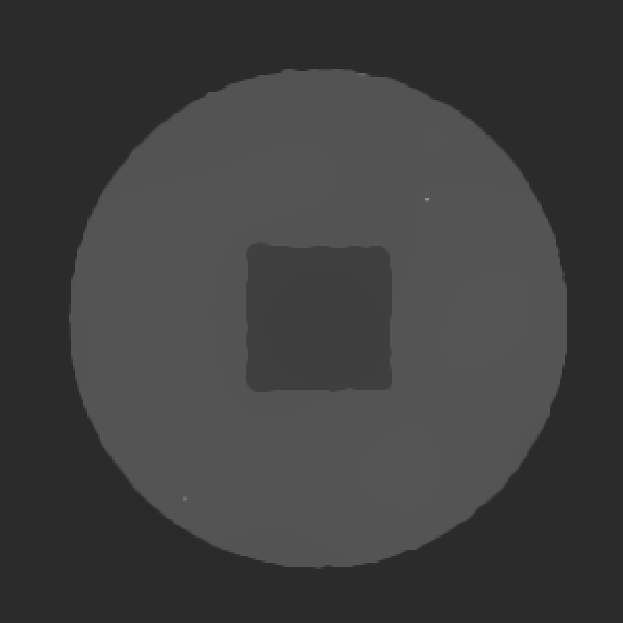}               
                \caption{PSNR=42.79}
                \label{circle10gradexp0}
       \end{subfigure}
       
        \begin{subfigure}[b]{0.3\textwidth}           
            \includegraphics[scale=0.45]{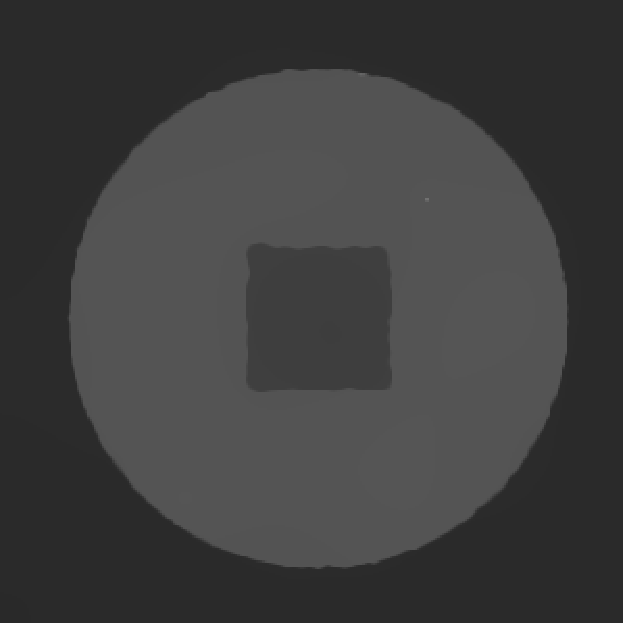}        
                \caption{\textbf{PSNR=43.53}}
                \label{circle10avggrayexp0nu1}
       \end{subfigure}
       \begin{subfigure}[b]{0.3\textwidth}           
                \includegraphics[scale=0.45]{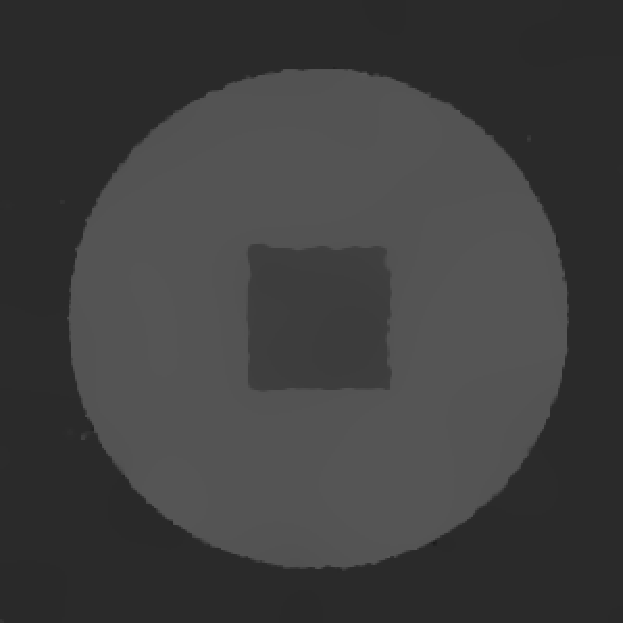}               
                 \caption{PSNR=42.74}
                \label{circle10grayexp0nu1}
       \end{subfigure}
       \begin{subfigure}[b]{0.3\textwidth}           
                \includegraphics[scale=0.45]{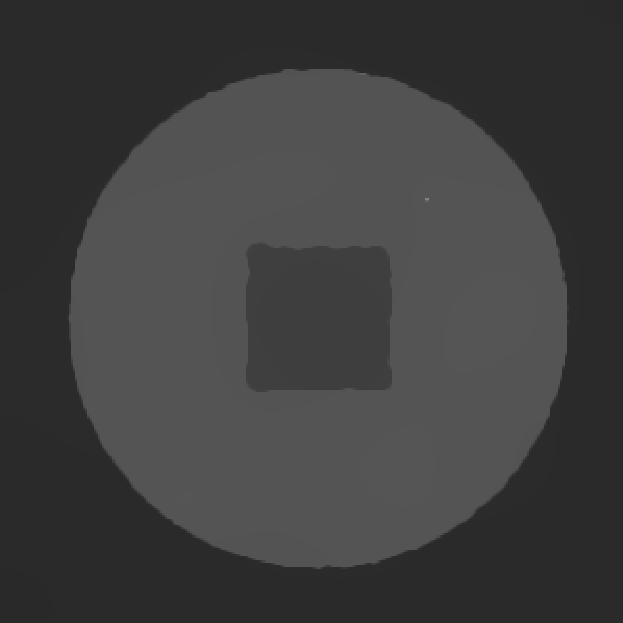}               
                \caption{PSNR=43.27}
                \label{circle10gradexp0nu1}
       \end{subfigure}

\caption{\footnotesize Restored images using \eqref{diff_power_constant} with variable exponents (Clean Image as in Figure \ref{circle_clean1}, Noisy Image as in Figure \ref{circle_noisy101}). 
\textbf{\textbf{First row}} ($p_0=2, \nu=0$): (a) $p_0-p_1$ ($K=0.1$) (b) $p_0-p_2$ ($K=0.1, \alpha=2$) (c) $p_0-p_3$ ($K=0.1, k=2$).
\textbf{Second row} ($\nu=0$): (d) $p_0-p_1$ ($p_0=2.2, K=0.1$) (e) $p_0-p_2$ ($p_0=2.6, K=0.2, \alpha=2$) (f) $p_0-p_3$ ($p_0=1.9, K=0.1, k=2$).  
\textbf{Third row} ($\nu=1$): (g) $p_0-p_1$ ($p_0=2.2, K=0.1$) (h) $p_0-p_2$ ($p_0=2.6, K=0.2, \alpha=2$) (i) $p_0-p_3$ ($p_0=1.9, K=0.1, k=2$).
 }\label{circle_10_restored_diffusion}
\end{figure}

\begin{figure}
       \centering
      \begin{subfigure}[b]{0.3\textwidth}           
                \includegraphics[scale=0.45]{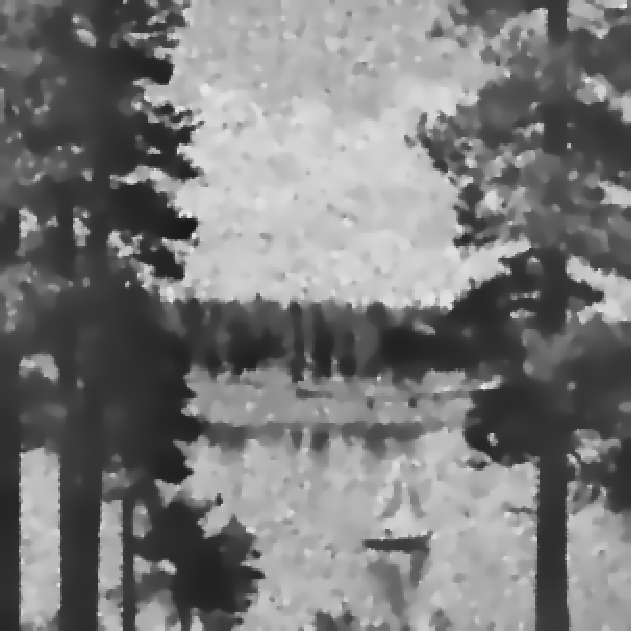}               
                \caption{PSNR=24.46}
                \label{lake110avggrayexp}
       \end{subfigure}
       \begin{subfigure}[b]{0.3\textwidth}           
                \includegraphics[scale=0.45]{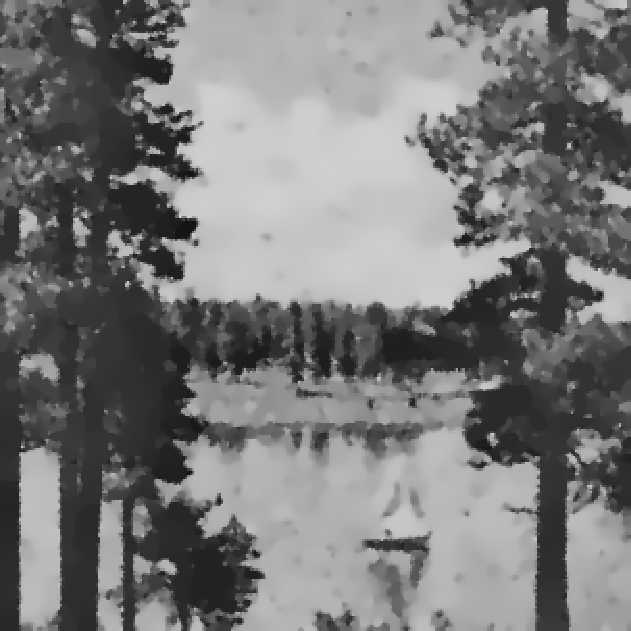}               
                 \caption{\textbf{PSNR=25.12}}
                \label{lake110grayexp}
       \end{subfigure}
       \begin{subfigure}[b]{0.3\textwidth}           
                \includegraphics[scale=0.45]{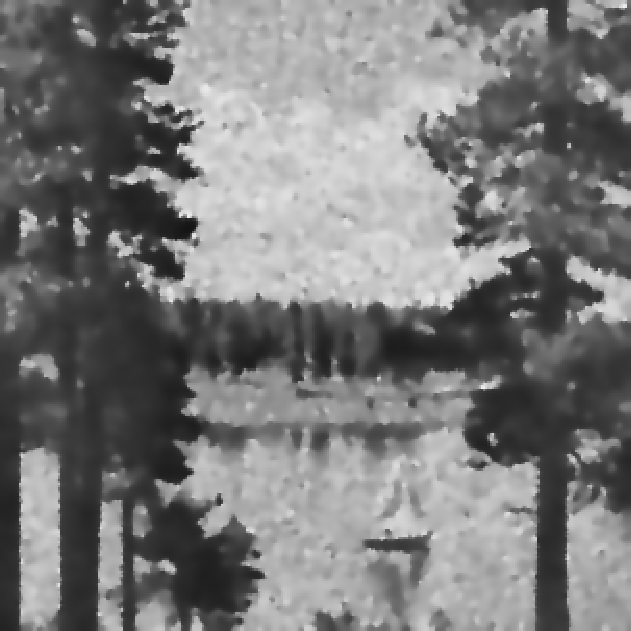}               
                \caption{PSNR=24.46}
                \label{lake110gradexp}
       \end{subfigure}
       
       \begin{subfigure}[b]{0.3\textwidth}           
                \includegraphics[scale=0.45]{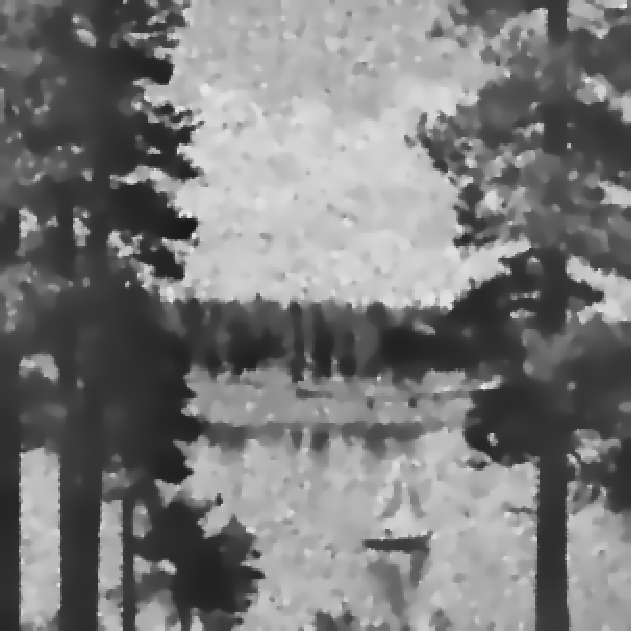}       
                \caption{PSNR=24.46}
                \label{lake110avggrayexp0}
       \end{subfigure}
       \begin{subfigure}[b]{0.3\textwidth}           
                \includegraphics[scale=0.45]{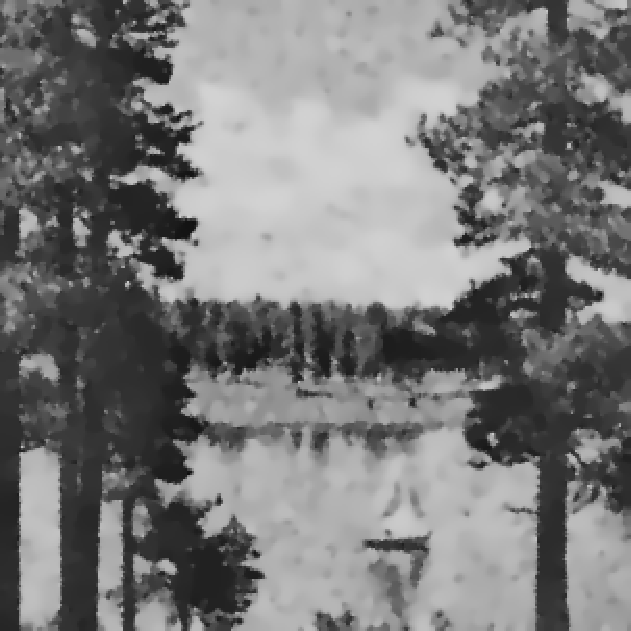}               
                 \caption{\textbf{PSNR=25.19}}
                \label{lake110grayexp0}
       \end{subfigure}
       \begin{subfigure}[b]{0.3\textwidth}           
                \includegraphics[scale=0.45]{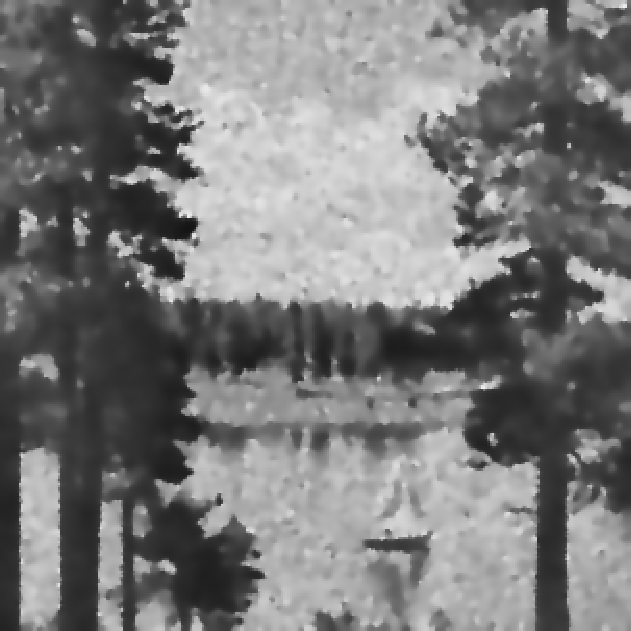}               
                \caption{PSNR=24.46}
                \label{lake110gradexp0}
       \end{subfigure}
       
       \begin{subfigure}[b]{0.3\textwidth}           
                \includegraphics[scale=0.45]{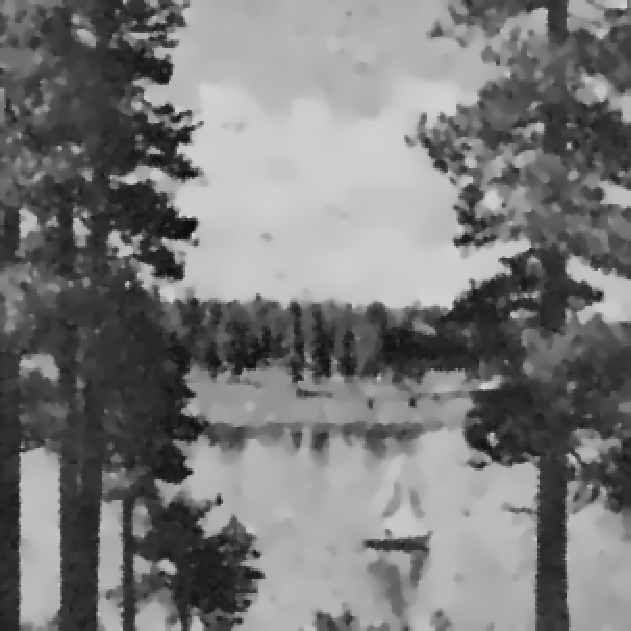}    
                \caption{\textbf{PSNR=25.29}}
                \label{lake110avggrayexp0nu1}
       \end{subfigure}
       \begin{subfigure}[b]{0.3\textwidth}           
                \includegraphics[scale=0.45]{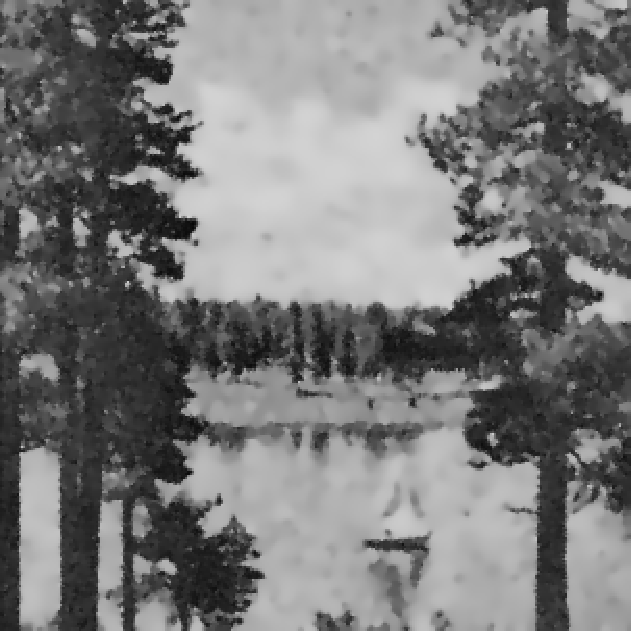}             
                \caption{PSNR=25.01}
                \label{lake110grayexp0nu1}
       \end{subfigure}
       \begin{subfigure}[b]{0.3\textwidth}           
                \includegraphics[scale=0.45]{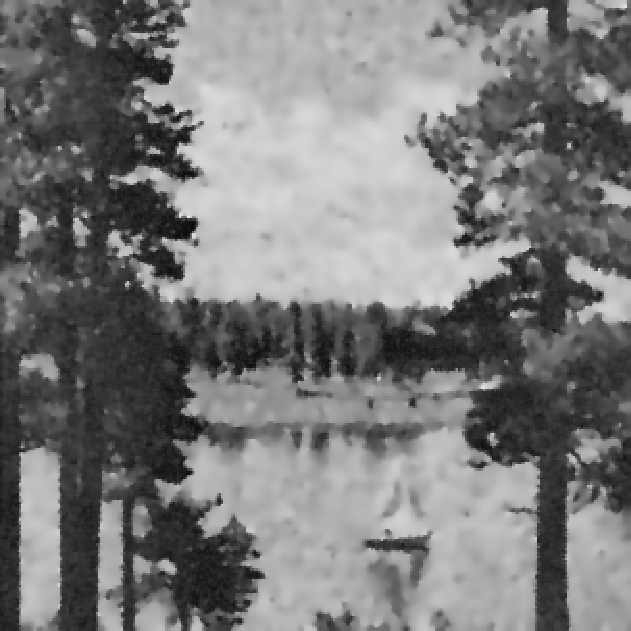}               
                \caption{PSNR=25.18}
                \label{lake110gradexp0nu1}
       \end{subfigure}
       
\caption{\footnotesize Restored images using \eqref{diff_power_constant} with variable exponents (Clean Image as in Figure \ref{lake11_clean}, Noisy Image as in Figure \ref{lake11_noisy10}). 
\textbf{First row} ($p_0=2, \nu=0$): (a) $p_0-p_1$ ($K=1$) (b) $p_0-p_2$ ($K=0.2, \alpha=2$) (c) $p_0-p_3$ ($K=4, k=2$).
\textbf{Second row} ($\nu=0$): (d) $p_0-p_1$ ($p_0=2.0, K=1$) (e) $p_0-p_2$ ($p_0=1.8, K=0.2, \alpha=2$) (f) $p_0-p_3$ ($p_0=2, K=4, k=2$).  
\textbf{Third row} ($\nu=2$): (g) $p_0-p_1$ ($p_0=1.85, K=0.40$) (h) $p_0-p_2$ ($p_0=2.0, K=1, \alpha=2$) (i) $p_0-p_3$ ($p_0=2.0, K=4, k=2$).
}\label{lake1_10_restored_diffusion}
\end{figure}

\begin{figure}
       \centering
        \begin{subfigure}[b]{0.3\textwidth}           
                \includegraphics[scale=0.26]{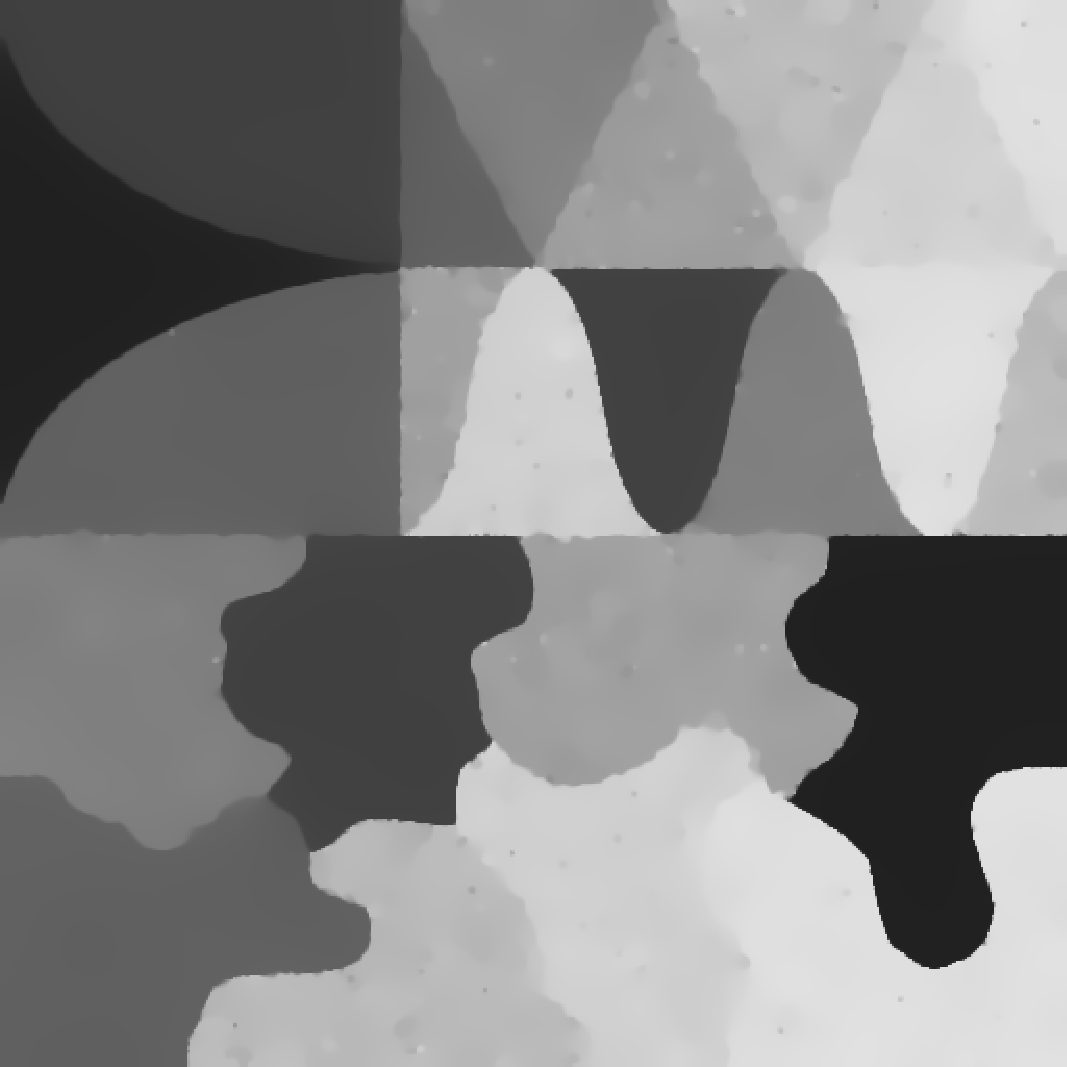}               
                \caption{\textbf{PSNR=25.57}}
                \label{texture10avggrayexp}
       \end{subfigure}
       \begin{subfigure}[b]{0.3\textwidth}           
                \includegraphics[scale=0.26]{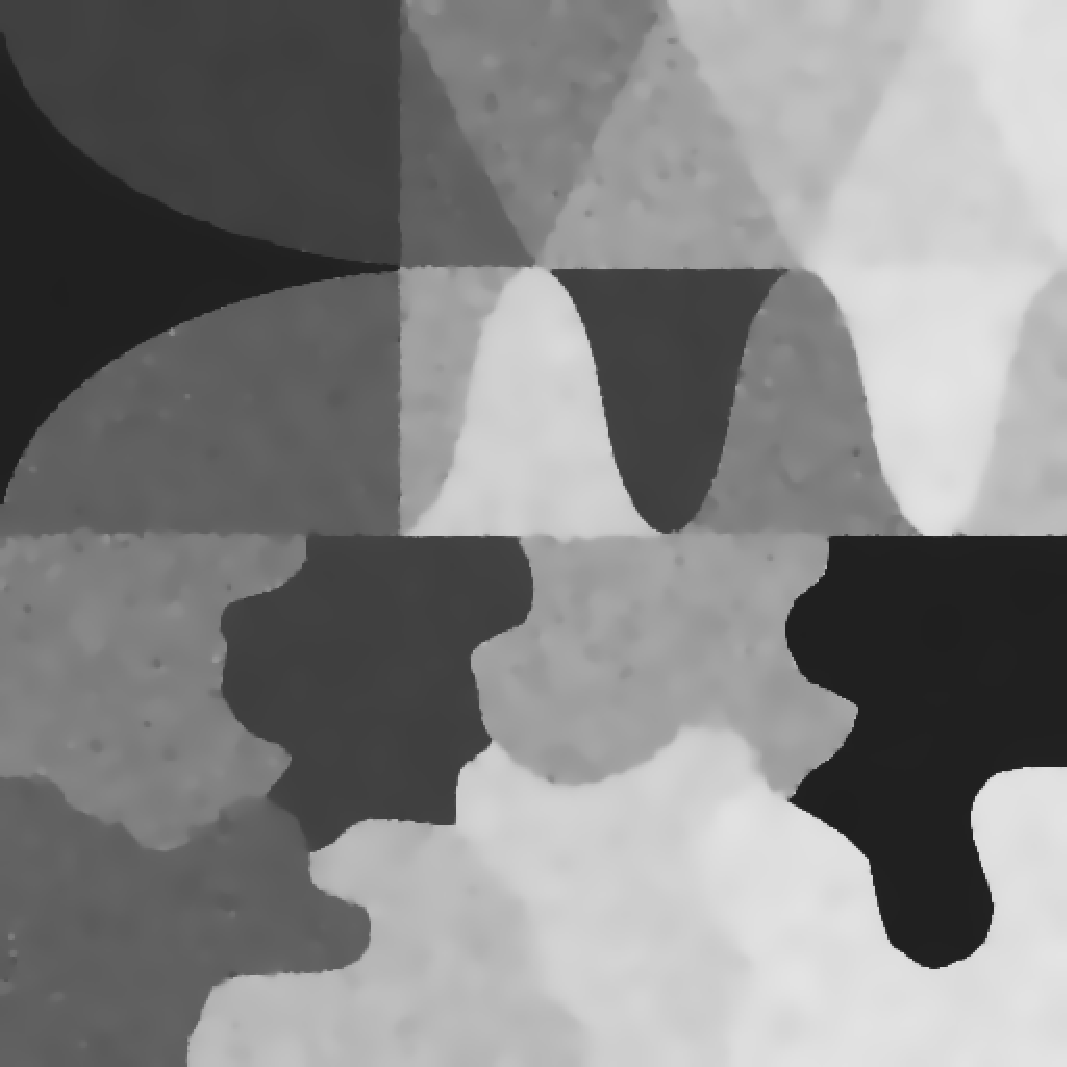}               
                 \caption{PSNR=25.47}
                \label{texture10grayexp}
       \end{subfigure}
       \begin{subfigure}[b]{0.3\textwidth}           
                \includegraphics[scale=0.26]{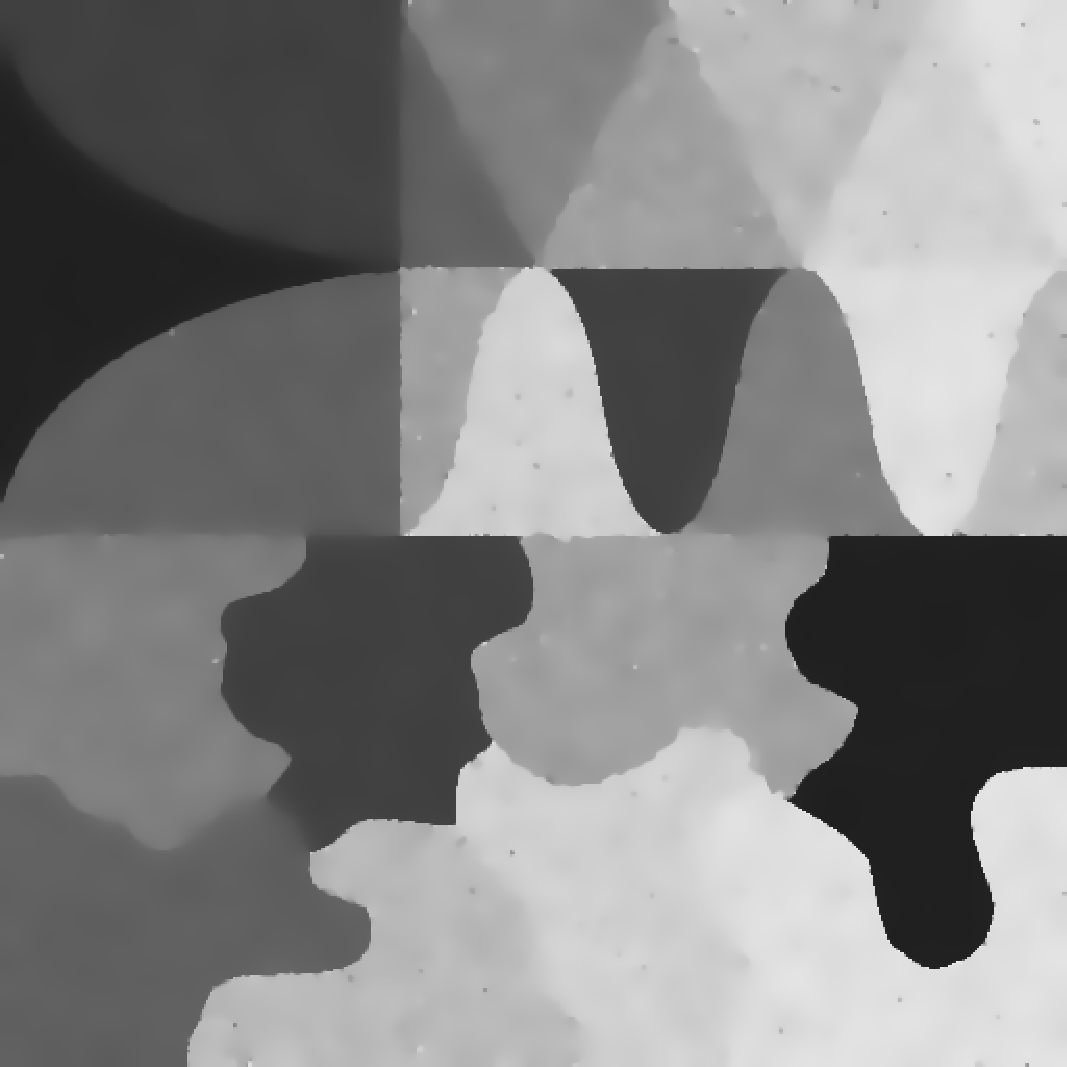}               
                \caption{PSNR=25.54}
                \label{texture10gradexp}
       \end{subfigure}
       
       \begin{subfigure}[b]{0.3\textwidth}           
                \includegraphics[scale=0.26]{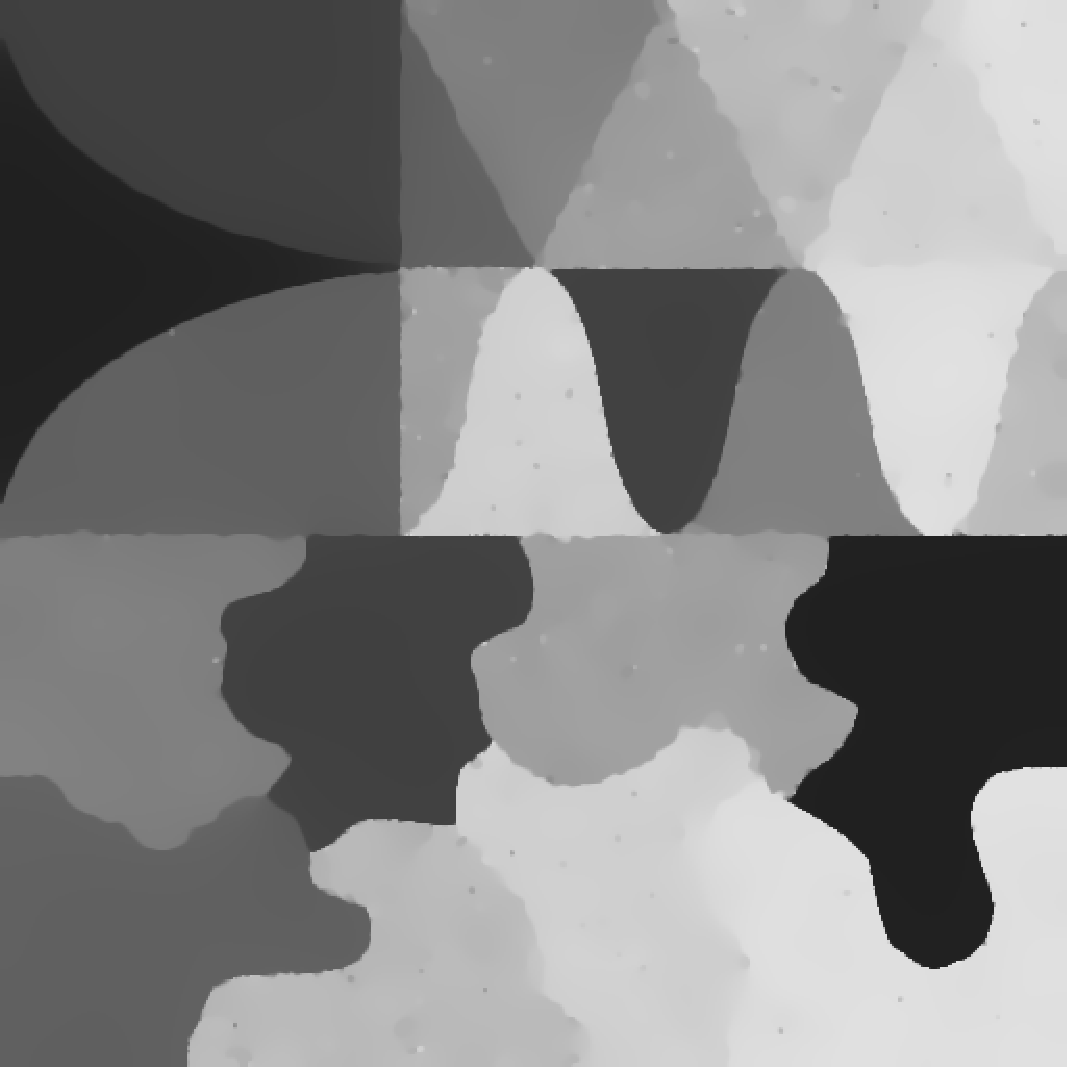}       
                \caption{\textbf{PSNR=25.61}}
                \label{texture10avggrayexp0}
       \end{subfigure}
       \begin{subfigure}[b]{0.3\textwidth}           
                \includegraphics[scale=0.26]{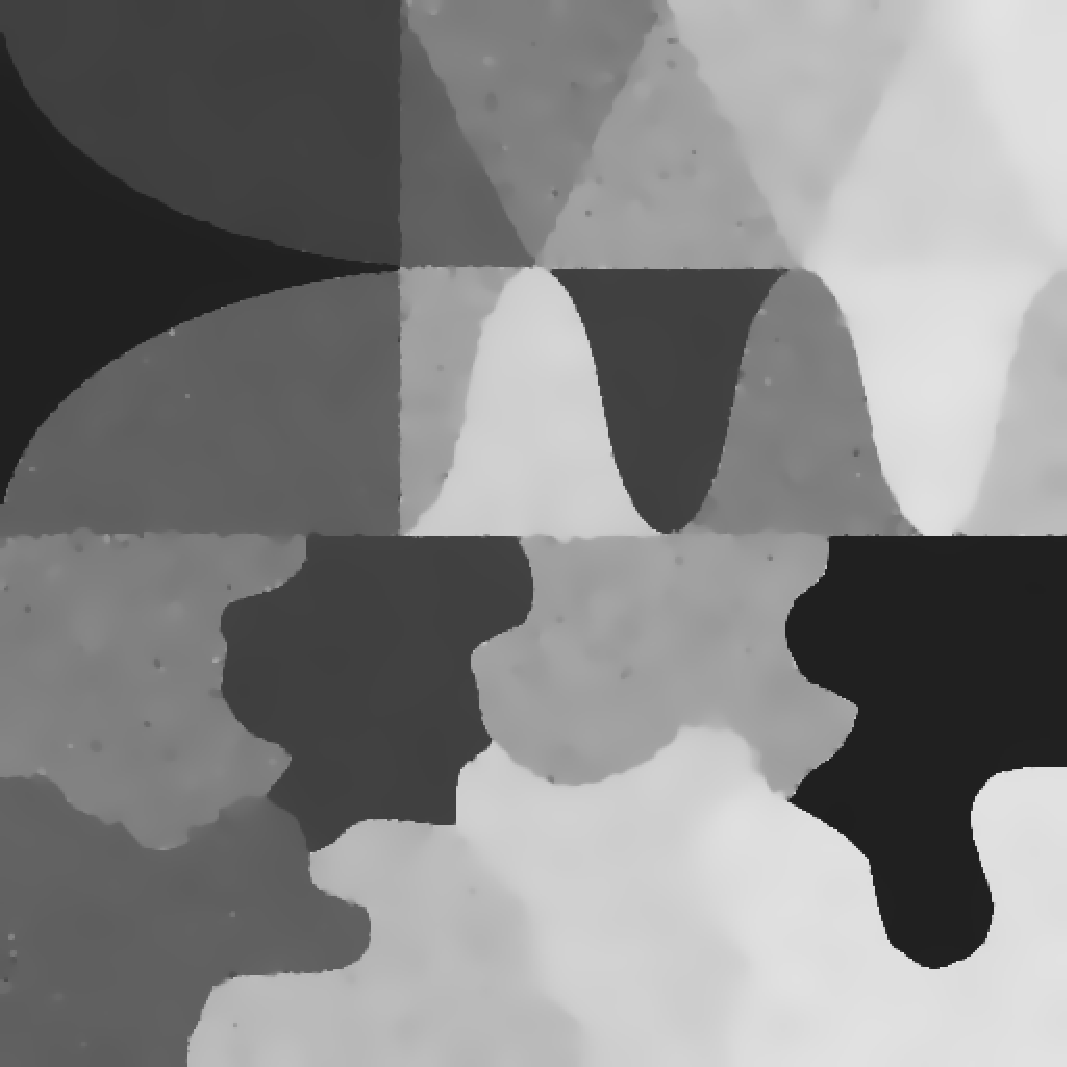}               
                 \caption{PSNR=25.57}
                \label{texture10grayexp0}
       \end{subfigure}
       \begin{subfigure}[b]{0.3\textwidth}           
                \includegraphics[scale=0.26]{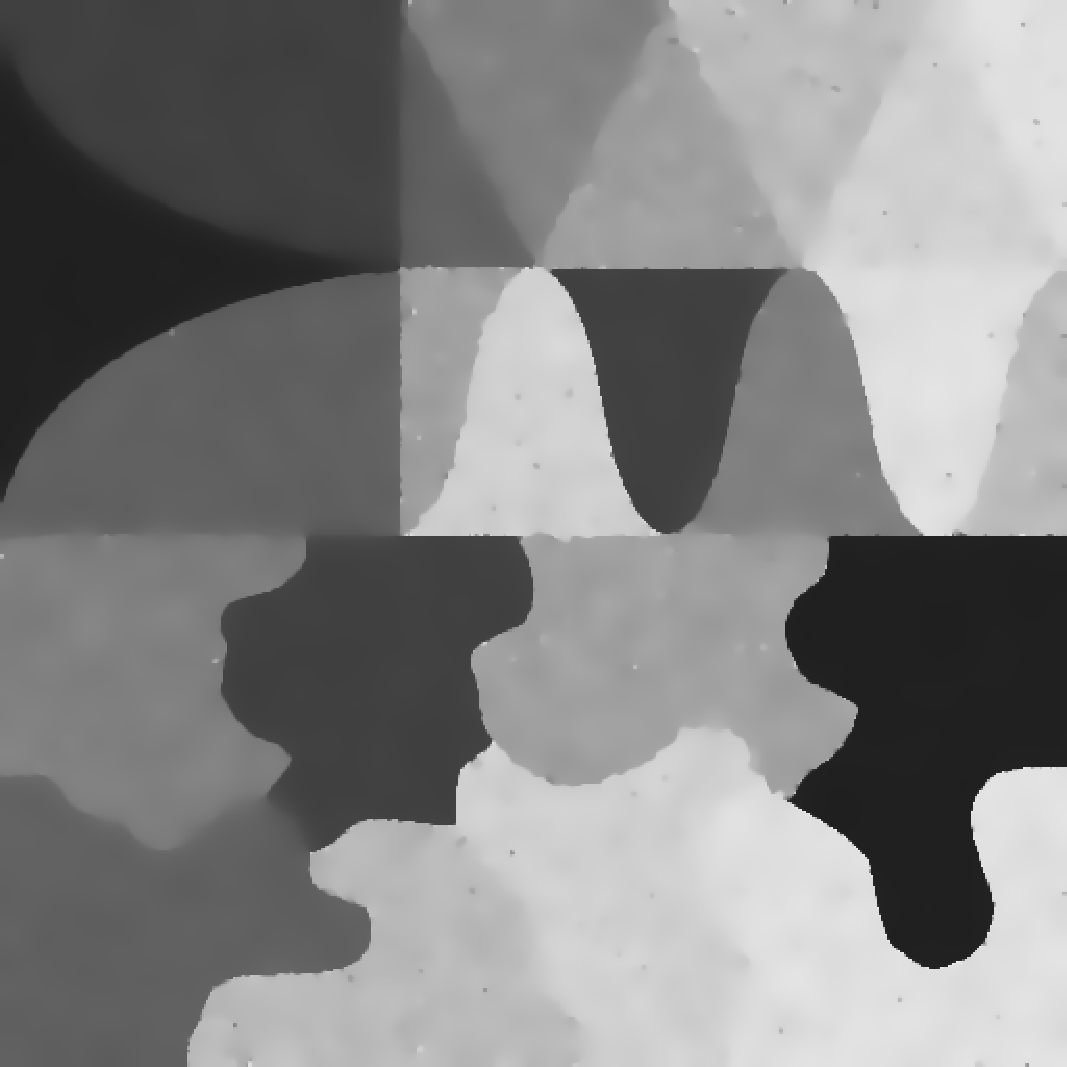}               
                \caption{PSNR=25.54}
                \label{texture10gradexp0}
       \end{subfigure}
       
       \begin{subfigure}[b]{0.3\textwidth}           
                \includegraphics[scale=0.26]{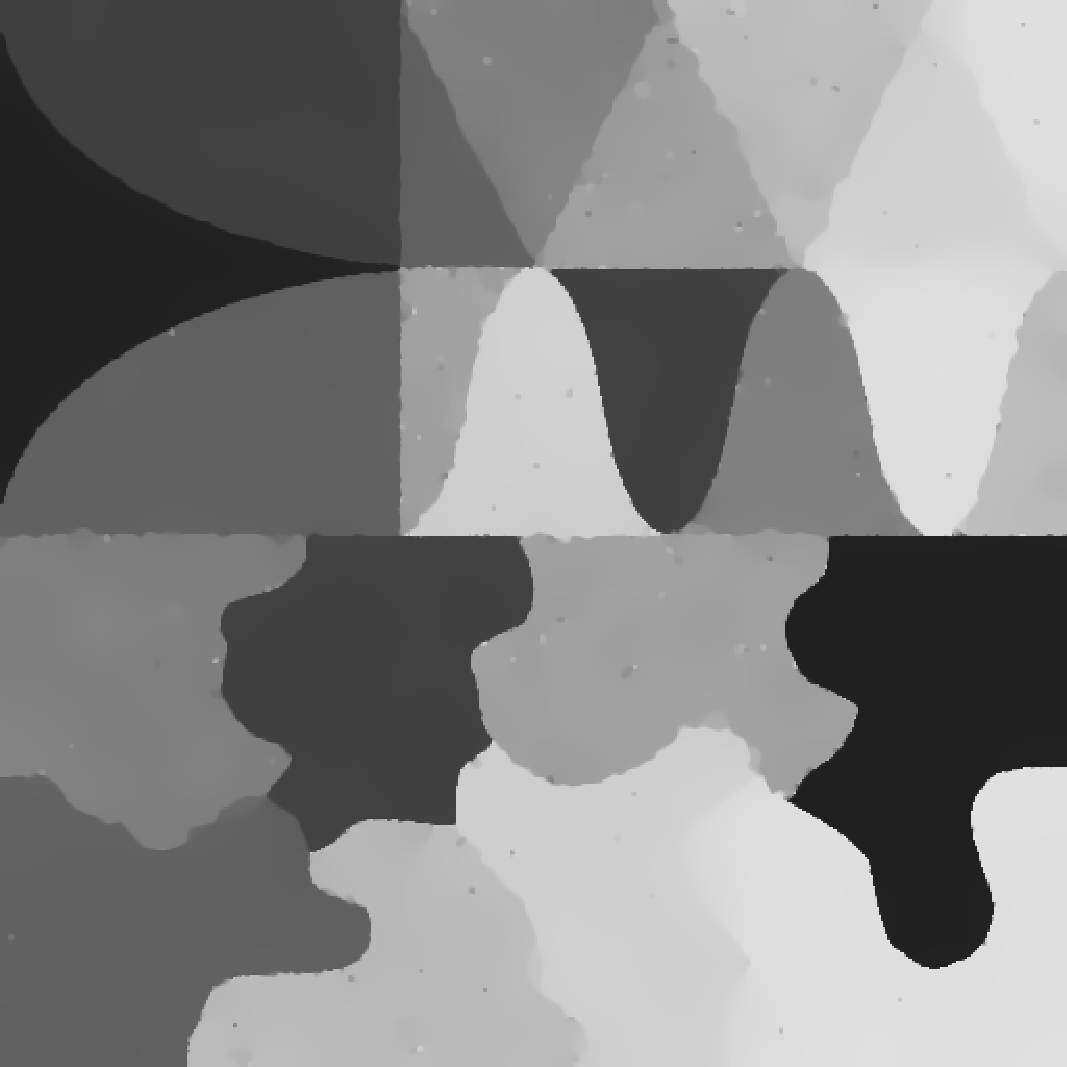}           
                \caption{\textbf{PSNR=25.63}}
                \label{texture10avggrayexp0nu1}
       \end{subfigure}
       \begin{subfigure}[b]{0.3\textwidth}           
                \includegraphics[scale=0.26]{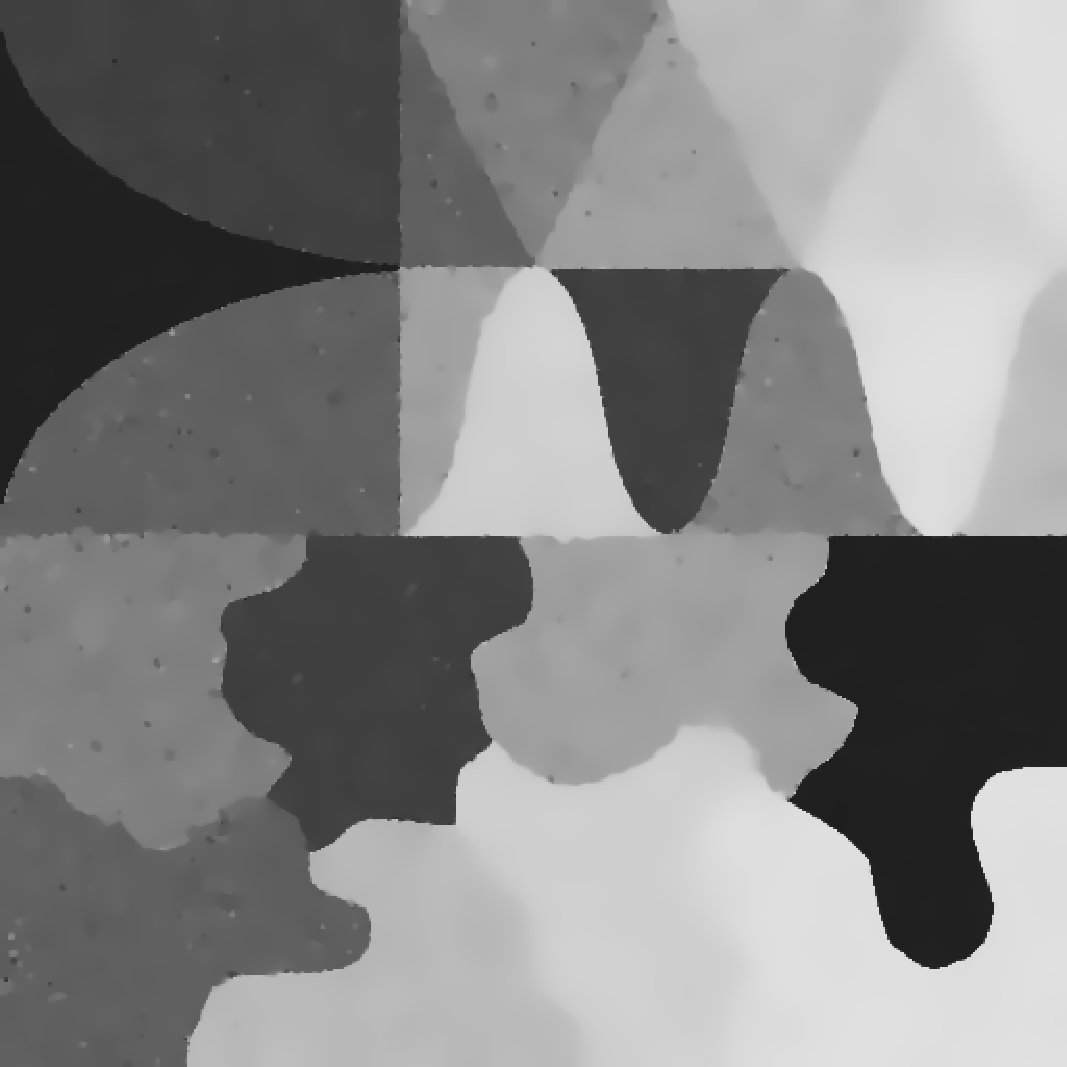}             
                \caption{PSNR=25.48}
                \label{texture10grayexp0nu1}
       \end{subfigure}
       \begin{subfigure}[b]{0.3\textwidth}           
                \includegraphics[scale=0.26]{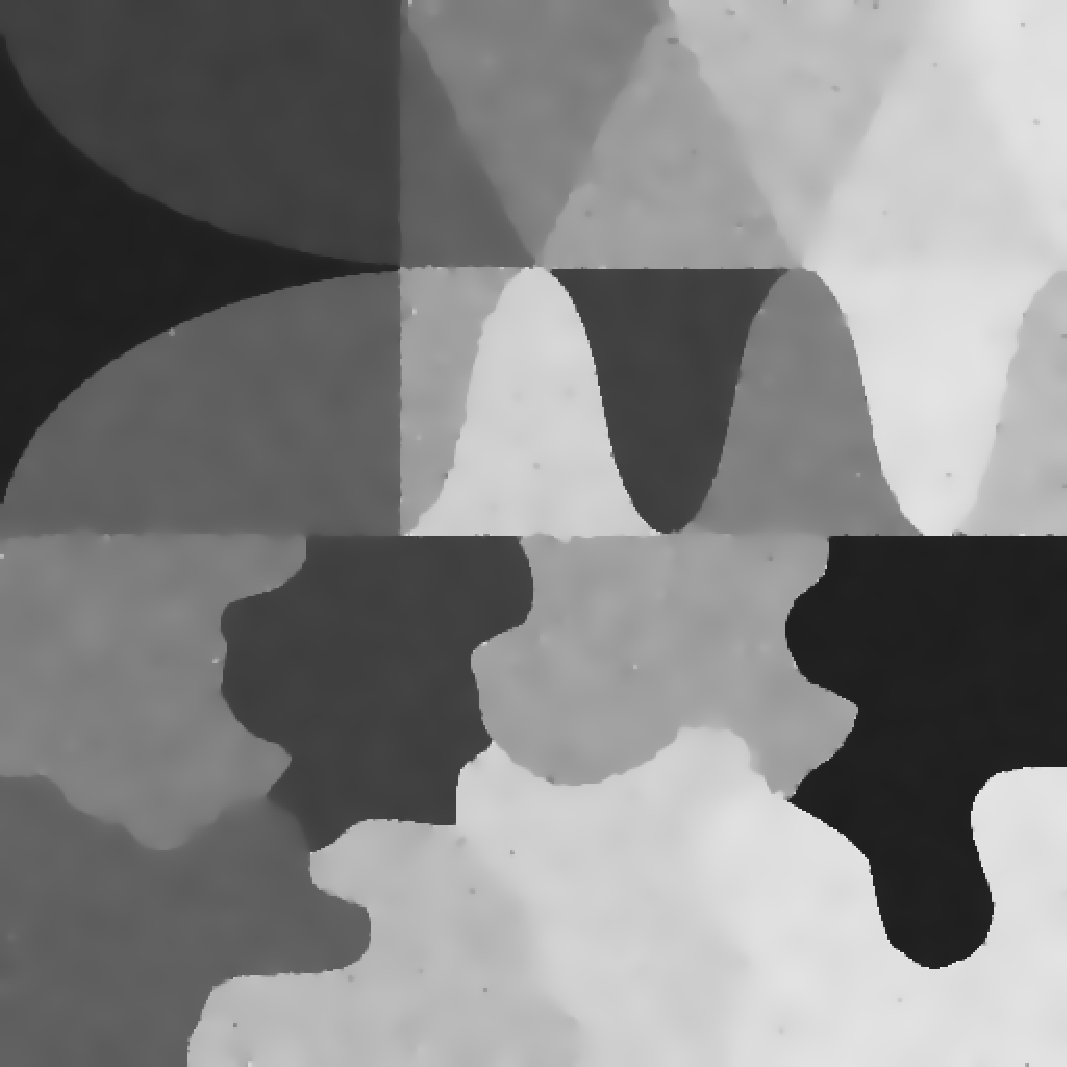}               
                \caption{PSNR=25.59}
                \label{texture10gradexp0nu1}
       \end{subfigure}
       
\caption{\footnotesize Restored images using \eqref{diff_power_constant} with variable exponents (Clean Image as in Figure \ref{texture_clean}, Noisy Image as in Figure \ref{texture_noisy10}). \textbf{First row} ($p_0=2, \nu=0$): (d) $p_0-p_1$ ($K=0.2$)  (e) $p_0-p_2$ ($K=0.5, \alpha=2$) (f) $p_0-p_3$ ($K=1, k=2$).
 \textbf{Second row} ($\nu=0$): (g) $p_0-p_1$ ($p_0=2.2, K=0.2$) (h) $p_0-p_2$ ($p_0=2.3, K=0.5, \alpha=2$) (i) $p_0-p_3$ ($p_0=2.0, K=1, k=2$). \textbf{Third row} ($\nu=1$): (j) $p_0-p_1$ ($p_0=2.2, K=0.2$) (k) $p_0-p_2$ ($p_0=2.3, K=0.5, \alpha=2$) (l) $p_0-p_3$ ($p_0=2.0, K=1.0, k=2$).}\label{texture_10_restored_diffusion}
\end{figure}

Next, due to the interest in variable exponents, we discuss the image despeckling results using diffusion model \eqref{diff_power_constant} with the diffusion coefficient as follows
\begin{align}\label{diff_edgeg}
g\left(I_\xi, \vert \nabla I_\xi \vert \right) =
\epsilon
+ \frac{ 2|I_\xi|^\nu}{\left(M^{I}_{\xi}\right)^\nu + | I_\xi|^\nu}  \cdot \dfrac{1}{1+ \left(\frac{\vert \nabla I_\xi \vert }{K} \right)^{p\left(I_\xi, \vert \nabla I_\xi \vert \right)} }\,,
\end{align}
where $\epsilon>0$ is very small parameter. 
We add $\epsilon$ to make the problem well posed, and it will not destroy the effect of edge detecting nor grey level indication whenever $\epsilon$ is very small.
In \eqref{diff_edgeg}, $p = p\left(I_\xi, \vert \nabla I_\xi \vert \right)$ is not a constant but a variable that may depend on both the gray level of the image and the image gradient. 
In order to prove well-posedness, we will assume that
\begin{align}\label{eq:ass-exponent_p}
p: \mathbf{R}\times [0,\infty) \rightarrow [1,p^+) \quad \text{is a once differentiable in both of its arguments, with $|\nabla p| \leq L$.}
\end{align}
We will use three explicit choices of exponents $p$ in our simulations as follows
\begin{align}
p\left(I_\xi, \vert \nabla I_\xi \vert \right)
&=p_0-p_1
\quad \text{where} \quad p_1 =\underset{x \in \Omega}{\mathrm{mean}} \dfrac{2\vert I_\xi \vert^\alpha}{\big(M^{I}_{\xi}\big)^\alpha+\vert I_\xi \vert^\alpha}, \label{p_1_relation}\\
p\left(I_\xi, \vert \nabla I_\xi \vert \right)
&=p_0-p_2 \quad \text{where}\quad p_2= \dfrac{2\vert I \vert^\alpha}{M^\alpha+\vert I \vert^\alpha}, \label{p_2_relation}\\
p\left(I_\xi, \vert \nabla I_\xi \vert \right)
&=p_0-p_3\quad \text{where}\quad p_3= \dfrac{2}{1+k\vert \grad G_\sigma \ast I \vert^2}, \label{p_3_relation}
\end{align}
in which $p_0$ is a parameter. We discuss the image despeckling results using the model \eqref{diff_power_constant} with the diffusion coefficient \eqref{diff_edgeg} and analyze all the distinct exponents $p$ as mentioned in \eqref{p_1_relation}, \eqref{p_2_relation}, and \eqref{p_3_relation}. Moreover, for each exponent $p$, we make three different cases depending on the choices of $p_0$ and $\nu$ in \eqref{diff_edgeg}:
\begin{itemize}
    \item[Case 1:] $p_0=2$ and $\nu=0$,
    \item[Case 2:] $p_0=2\pm\delta$ ($\delta>0$) and $\nu=0$,
    \item[Case 3:] $p_0=2\pm\delta$ ($\delta>0$) and $\nu>0$.
\end{itemize}
Figures \ref{circle_10_restored_diffusion}--\ref{texture_10_restored_diffusion} represent the results using different exponents with the above three cases. We calculate the restored results for the noise images in Figure \ref{circle_noisy101}, Figure \ref{lake11_noisy10}, and Figure \ref{texture_noisy10} and describe the results in Figure \ref{circle_10_restored_diffusion}, Figure \ref{lake1_10_restored_diffusion}, and \ref{texture_10_restored_diffusion}, respectively. In each figure, different rows represent the results for different cases. The first row shows the results for three different exponents $p$ ($p_0-p_1$, $p_0-p_2$, and $p_0-p_3$) under Case 1 ($p_0=2$ and $\nu=0$). In this case, the first term $a(I_\xi)$ in the diffusion coefficient $g$ simplifies to one. Consequently, for the second exponent (\(p_0 - p_2\)), the system reduces to the model proposed in \cite{zhou2018nonlinear}; for the third exponent (\(p_0 - p_3\)), the model behaves similarly to the approach discussed in \cite{guo2011adaptive}. The second row presents results for each of the three different exponents \( p \) under Case 2 (\( p_0 = 2 \pm \delta \) and \( \nu = 0 \)). In this case, we apply a small adjustment to \( p_0 \) (by adding or subtracting a positive value \(\delta\) from \( p_0 = 2 \) in Case 1) to explore whether this yields improved results compared to those in Case 1. Finally, the images in the third row are computed for each of the three different exponents \( p \) under Case 3 (\( p_0 = 2 \pm \delta \) and \( \nu > 0 \)). This results in a total of nine cases, with observations as follows.

Figure \ref{circle_10_restored_diffusion} shows the results for a circle image. In Case 1, the first exponent produces better results than the other two exponents. For Case 2, all exponents yield more satisfactory results than their corresponding exponents in Case 1, likely due to the adjustment in \( p_0 \). Notably, the second exponent outperforms the others in this case. In Case 3, we set \( \nu = 1 \) (which produces the best results for this image) and evaluate each of the three exponents. Here, the first and third exponents give better results compared to their counterparts in Cases 1 and 2. Interestingly, the quality of the restored image decreases for the second exponent in Case 3 compared to Case 2, possibly due to the influence of $a(I_\xi)$ in the diffusion process. Overall, for this image, the first exponent with Case 3 produces the best results among the nine cases. From the other images in Figures \ref{lake1_10_restored_diffusion} and \ref{texture_10_restored_diffusion}, we observe similar patterns as in the circle image. Additionally, Figures \ref{circle_10_restored_diffusion_constp}--\ref{texture_10_restored_diffusion_constp} show that edge information is best preserved near \( p = 1.5 \). For regions with high gray levels, the diffusion function in \eqref{diff_edgeg} (for \( p = p_0 - p_2 \)) simplifies to:
\begin{align}
\epsilon + g\left(I_\xi, |\nabla I_\xi| \right) = \dfrac{1}{1 + \left(\frac{|\nabla I_\xi|}{K}\right)}
\end{align}
Consequently, the variable exponent \( p = p_0 - p_2 \) is less effective at preserving edges in high-gray-level regions. This limitation of \( p = p_0 - p_2 \) can be easily concluded by observing the results illustrated in Figure \ref{texture_10_restored_diffusion}. We summarize the key findings for the diffusion model as follows:
\begin{itemize}
    \item If \( \nu = 0 \), the variable exponent with gray level dependence performs best in most cases.
    \item Exponents that depend on gray levels are generally more effective than those based on edge detectors.
    \item For \( \nu > 0 \), it is preferable to fix the exponent.
\end{itemize}


\subsection{Telegraph Diffusion Models}
\label{telegraph_diffusion_model}
\noindent In 2020, Majee et al. \cite{majee2020gray} introduced a telegraph-diffusion model for image despeckling, which takes the form:
\begin{align}
&I_{tt} + \gamma I_{t} - \text{div} \left( g\left(I_\xi, \vert \nabla I_\xi \vert \right) \nabla I \right) = 0 \quad \text{in} \,\, \Omega_T, \label{maina} \\
&I(x,0) = I_0(x), \quad I_t(x,0) = 0 \quad \text{in} \,\, \Omega, \label{mainc} \\
&\partial_n I = 0 \quad \text{on} \,\, \partial \Omega_T, \label{mainb} 
\end{align}
where the diffusion function \( g \) is defined by
\begin{align*}
g\left(I_\xi, \vert \nabla I_\xi \vert \right) = \frac{2 \vert I_\xi \vert^\nu}{(M^{I}_{\xi})^\nu + \vert I_\xi \vert^\nu} \cdot \frac{1}{1 + \left(\frac{\vert \nabla I_\xi \vert }{K} \right)^{2}},
\end{align*}
with \(\nu \geq 1\), and constants \( \gamma, K > 0 \). Here, \( I_\xi = G_\xi \ast I \), and \( M_\xi^I = \max_{x \in \Omega} \vert I_\xi (x,t) \vert \). Majee et al. employ a constant exponent in this model. Inspired by the effectiveness of variable exponents discussed in Section \ref{diffusion_discussion} and building on ideas from \cite{majee2020gray}, we propose a variable exponent-based telegraph diffusion model for this study. Specifically, we extend the system \eqref{maina}--\eqref{mainb} using the diffusion coefficient defined in \eqref{diff_edgeg}. It is important to note that for theoretical study, we consider the diffusion coefficient in a more general form and deal with more difficulties than the model discussed in \cite{majee2020gray}.

In Section \ref{computation_discussion}, we examine the despeckling performance of the model \eqref{maina}--\eqref{mainb} with the diffusion function \eqref{diff_edgeg}. We explore each of the three exponents \( p(I_\xi, \vert \nabla I_\xi \vert) \) as presented in equations \eqref{p_1_relation}, \eqref{p_2_relation}, and \eqref{p_3_relation}. For the first and third exponents, we use a non-zero \(\nu\) value, while for the second exponent, we set \(\nu = 0\). Results for these cases are then discussed.

\section{Existence and Uniqueness of Weak Solutions}
\label{sec:analysis}
In this section, we establish the existence and uniqueness of weak solutions of a class of parabolic and hyperbolic PDEs, including the models suggested in Section \ref{diffusion_discussion}. Let $\Omega$ denote a bounded domain in $\mathbf{R}^n$, $n \geq 1$ and let $T > 0$. For $1\le p\le \infty$ we denote by $(L^p, \|\cdot\|_{L^p})$ the standard spaces of $p$-th order integrable functions on $\Omega$.
For $r\in \mathbb{N}$, we write $(H^r, \|\cdot\|_{H^r})$ for usual Hilbert spaces on $\Omega$, and $(H^{1})^\prime$ for the dual space of $H^1$.
We will sometimes write $L^{p},$ $H^1,$ $(H^{1})^\prime$ instead of $L^{p}(\Omega),$ $H^1(\Omega),$ $(H^{1}(\Omega))^\prime$, respectively. 
When $A$ is a matrix, we adopt the matrix norm $\| A \| = \sup_{|v| = 1} |A v|$.
\subsection{Diffusion Model}
In this subsection we discuss a class of nonlinear diffusion models including the model \eqref{diff_power_constant} with the diffusion coefficients \eqref{diff_constatp} and \eqref{diff_edgeg}. 
In particular, we consider parabolic PDEs of the form
\begin{align}
 &I_t - \nabla \cdot \left( \mathcal A\, \nabla I\right) = 0 \quad (x,t) \in \Omega \times (0,T),   \tag{D.1}\label{para-general1}\\
 &I(x,0) = I_0(x)  \quad x \in \Omega,  \tag{D.2} \label{para-general2}\\
 & \left\langle\mathcal A \nabla I, {\bf n}\right\rangle = 0 \quad \text{on} \quad \partial \Omega \times (0,T),  \tag{D.3} \label{para-general3}
\end{align}
where $\mathcal A = \mathcal A(x,t,I,\nabla I, I(\cdot),\nabla I(\cdot)): \mathbf{R}^n \times \mathbf{R} \times \mathbf{R} \times \mathbf{R}^n \times \mathbf{R} \times \mathbf{R}^n \to \mathbf{R}^{n\times n}$
is a symmetric matrix satisfying some assumptions given below. Here, $I(\cdot)$ and $\nabla I(\cdot)$ indicate dependence of $I(x,t)$ and $\nabla I(x,t)$ in a non-local manner. We introduce the solution space $W(0,T)$ for the problem \eqref{para-general1}--\eqref{para-general3} as
\begin{align*}
 W(0,T)&=\Big\{w\in L^2(0,T; H^1)\,,\,w_{t} \in L^2(0,T; (H^1)') \Big\}\,,
\end{align*}
which is a Hilbert space for the graph norm, see \cite{jllions1968}.
\begin{defi}[Weak solution]\label{defi:weak_diffusion}
A function $I = I(x,t)$ is called a weak solution of \eqref{para-general1}--\eqref{para-general3} if
\begin{itemize}
 \item[a)] $I \in W(0, T) $ and \eqref{para-general2} and \eqref{para-general3} hold a.e. with derivatives in the sense of distributions.
 \item[b)] For all $\phi \in H^1$ and a.e. $t\in (0,T)$ it holds that
 \begin{align*}
\left \langle I_{t}, \phi   \right \rangle + {\displaystyle \int_{\Omega}}
 \langle\mathcal A \nabla I, \nabla \phi\rangle \,dx =0.
  \end{align*}
\end{itemize}
\end{defi}
To state our assumptions on $\mathcal A$ in problem \eqref{para-general1}--\eqref{para-general3} we define  the function space
\begin{align}\label{eq:M_1_diffusion}
 W_{1}= \big\{  v\in W(0,T):~  \|v\|_{L^\infty(0,T;L^2)}  \le  \|I_0\|_{L^2(\Omega)} \big\}.
\end{align}
We also define, for any element $w \in W_{1}$, the function $\mathcal{A}_w : \mathbf{R}^n \times \mathbf{R} \to \mathbf{R}^{n\times n}$ by
\begin{align}\label{eq:def_A_diffusion}
\mathcal{A}_w (x,t) = \mathcal A\left(x,t,w,\nabla w, w(\cdot),\nabla w(\cdot)\right).
\end{align}
To prove existence uniqueness of weak solution, we assume that there exist $\theta, C > 0$ so that  
\begin{align}
\theta |\eta|^2 \leq \langle \mathcal A_w(x,t) \eta, \eta\rangle , \quad \left(\|\mathcal{A}_w(x,t)\|\right)_{L^{\infty}(\Omega)} \leq C \quad \text{whenever}\quad (x,t) \in \Omega\times (0,T), \; \eta \in \mathbf{R}^n,\; w \in W_1.
 \tag{A}\label{eq:ass-A1-diffusion}
\end{align}
%
We also assume that there exists a constant $C$ such that whenever $I_1, I_2 \in W$ it holds that
\begin{align}
& \left( \| \mathcal A(x,t,I_1,\nabla I_1,I_1(\cdot),\nabla I_1(\cdot) ) - \mathcal A(x,t,I_2,\nabla I_2,I_2(\cdot),\nabla I_2(\cdot) ) \|\right)_{L^\infty(\Omega)} \notag \\
&\hspace{6cm} \leq C\|I_1(x,t) - I_2(x,t)\|_{L^2(\Omega)} \quad \text{a.e.}\quad  t \in (0,T).
  \tag{B}\label{eq:ass-A2-diffusion}
\end{align}
Our main results for parabolic diffusion models are summarized in the following theorem.
\begin{thm}\label{thm:existence-diffusion}
Suppose that assumptions \eqref{eq:ass-A1-diffusion} and \eqref{eq:ass-A2-diffusion} hold and that $I_0 \in H^1(\Omega)$.
Then there exists a unique weak solution of problem \eqref{para-general1}--\eqref{para-general3}in $W(0,T)$.
\end{thm}
\noindent
{\bf Proof}.
We begin with the existence part. 
Let $w \in W_1$ and $\mathcal{A}_w$ be as in \eqref{eq:def_A_diffusion}.
Define $\mathcal P$ as the problem
\begin{align}\label{eq:P_w_diffusion}
    \left\langle I_t, \varphi \right\rangle_{H^1(\Omega)' \times H^1(\Omega)}
    + \int_{\Omega} \langle \mathcal A_w\left(x,t \right)\nabla I, \nabla \varphi\rangle dx = 0,
\end{align}
whenever $\varphi \in H^1(\Omega)$ and a.e. $t \in (0,T)$, together with $I(x,0) = I_0$ and the boundary conditions \eqref{para-general3} in the sense of distribution.

Since $w \in W_{1}$ it follows by the assumption \eqref{eq:ass-A1-diffusion} that $\mathcal P$ is uniformly parabolic and thus has a unique weak solution $I_w \in W$, denoted $I=I_w$ for simplicity, which we will prove satisfies the following estimates:
\begin{align}
&  \|I\|_{L^\infty(0,T; L^2)} \le  \|I_0\|_{L^2}, \label{estimate2_diffusion}\\
&  \|I\|_{L^2(0,T; H^1)} \le C_1,\label{estimate1_diffusion} \\
&  \|I_t\|_{L^2(0,T; (H^{1})^\prime))} \le C_2, \label{estimate3_diffusion} 
\end{align}
where $C_1$ and $C_2$ are constants independent of $w$. To prove the estimates \eqref{estimate2_diffusion}--\eqref{estimate3_diffusion}, first we choose $\varphi=I$ in equation \eqref{eq:P_w_diffusion}. Then we obtain 
\begin{align}\label{diffusion-eq:stet_1}
&\frac{1}{2}\dfrac{d}{dt}\|I \|_{L^2}^2 + \large\int_{\Omega} \langle \mathcal A_w(x,t)\nabla I, \nabla I \rangle \, dx=0, \quad t \in (0,T).
\end{align}
Thanks to assumption \eqref{eq:ass-A1-diffusion} we have
\begin{align*}
&\frac{1}{2}\dfrac{d}{dt}\|I \|_{L^2}^2 + \theta \|\nabla I\|_{L^2}^2 \le 0,  \quad t \in (0,T). 
\end{align*}
Integrating the above relation w.r.t. time from $0$ to $t$, we have
\begin{align}\label{diffusion-eq:stet_2}
&\|I(t)\|^2_{L^2}+2 \theta \int_{0}^{t} \|\nabla I(\tau)\|_{L^2}^2\, d\tau \leq \|I_0\|^2_{L^2},  \quad t \in (0,T). 
\end{align}
By excluding the second term in the left-hand side and taking the supremum over $(0, T)$, we complete the proof of \eqref{estimate2_diffusion}. 

From relation \eqref{diffusion-eq:stet_2} we also have
\begin{align*}
\int_{0}^{t} \|\nabla I(\tau)\|_{L^2}^2\, d\tau \leq  \dfrac{1}{2} \theta^{-1} \|I_0\|^2_{L^2},  \quad t \in (0,T), 
\end{align*}
giving 
\begin{equation}\label{diffusion-eq:stet_3}
\|\nabla I\|^2_{L^{2}(0,T;L^2)} \leq \dfrac{1}{2} \theta^{-1} \|I_0\|^2_{L^2} \,.
\end{equation}
Now, \eqref{estimate2_diffusion} implies
\begin{align*}
\|I\|^2_{L^2(0,T,L^2)} \leq \int_0^T \|I_0\|^2_{L^2} dt = T \|I_0\|^2_{L^2} 
\end{align*}
which, together with 
\eqref{diffusion-eq:stet_3}, give
\begin{align*}
 \|I\|^2_{L^2(0,T,H^1)} =  \|I\|^2_{L^2(0,T,L^2)} +  \|\nabla I\|^2_{L^2(0,T,L^2)} \leq T \|I_0\|^2_{L^2}  + \dfrac{1}{2} \theta^{-1} \|I_0\|^2_{L^2},
\end{align*}
which completes the estimate \eqref{estimate1_diffusion}.

Next we choose $\phi \in H^1$ in \eqref{eq:P_w_diffusion} such that $\|\phi\|_{H^1}\leq 1$ and use boundedness of $\mathcal{A}_w$ together with Cauchy-Schwarz inequality to obtain
\begin{align*}
\left\vert \big \langle I_t,\varphi \big \rangle \right\vert \leq C  \|\nabla I\|_{L^2}  \|\phi\|_{H^1}\,.
\end{align*}
Hence, by the definition of norm in $(H^{1})'$, we get 
\begin{align}
\left\Vert  I_t \right\Vert_{(H^{1})'} \leq C \,.
\end{align}
Moreover, squaring both sides and integrating over $(0,T)$, we have
\begin{equation}\label{acpde_eq-for-I-7}
\|I_t\|_{L^{2}(0,T;(H^{1})')} \leq C_2\,.
\end{equation}
This completes the proof of \eqref{estimate3_diffusion}.

Now 
we introduce the non-empty, convex and weakly compact subspace $W_0$ of $W(0,T)$ by
\begin{align*}
W_0=\Big\{w \in W(0,T): \, \|w\|_{L^2(0,T; H^1)}  \leq C_1\,, & \|w\|_{L^{\infty}(0,T; L^2)} \leq \| I_0 \|_{L^2}\,, \\
&  \| w_t \|_{L^2(0,T; (H^{1})')} \leq C_2\,,  w(0)=I_0\Big\}\,,
\end{align*}
%
and consider the mapping
\begin{align*}
\mathcal{P}:~ & W_0 \goto W_0 \\
& w\mapsto I_w\,.
\end{align*}
We intend to show that the mapping $\mathcal{P} : w \rightarrow I_w $ is weakly continuous from $W_0$ into $W_0$:
Let $w_k$ be a sequence that converges weakly to some $w$ in $W_0$ and let $I_k = I_{w_k}$. We have to show that $\mathcal{P}(w_k):= I_k$ converges weakly to $\mathcal{P}(w): = I_w$. From the classical results of compact inclusion in Sobolev spaces \cite{raadams1975}, we can extract subsequences of $\{w_k\}$ and $\{I_k\}$, still denoted by $\{w_k\}$ and $\{I_k\}$ respectively, such that for some $I\in W_0$ we have, as $k\goto \infty$,
\begin{align}\label{eq:conv-1}
\begin{cases}
 w_{k} \longrightarrow  w\, \hspace{0.2cm} {\rm in}~~L^2(0,T;L^2) \hspace{0.2cm} {\rm and~a.e.~on}~~\Omega_T\,,\\[0.5em]
\displaystyle I_{k} \longrightarrow  I \hspace{0.2cm} \text{weakly}~ \text{in}~~L^{2}(0,T;H^1)\,,\\[1em]
\displaystyle I_{k} \longrightarrow  I \hspace{0.2cm} \text{in}~~L^{2}(0,T; L^2)\,,\\[1em]
\displaystyle {\partial_{x_i} I_k} \longrightarrow  {\partial_{x_i} I}\,\,\,(i=1, 2) \hspace{0.2cm} \text{weakly}~ \text{in}~~L^{2}(0,T; L^2)\,,\\[1em]
\partial_t I_k \longrightarrow \partial_t I \hspace{0.2cm} \text{weakly}~\text{in}~~L^{2}(0,T;(H^{1})')\,.
\end{cases}
\end{align}
Integrating \eqref{eq:P_w_diffusion} from 0 to $T$ yields
\begin{align}\label{eq:conv-first}
  \int_0^T    \int_{\Omega}\partial_t I_k \varphi dx dt
    + \int_0^T \int_{\Omega} \langle \mathcal A_{w_k}\left(x,t \right)\nabla I_k, \nabla \varphi \rangle dx dt = 0,
\end{align}
and by using assumptions \eqref{eq:ass-A2-diffusion}, \eqref{eq:ass-A1-diffusion}, as well as \eqref{eq:conv-1}, we observe that
\begin{align}\label{eq:conv-2}
\bigg|\int_0^T & \int_{\Omega} \langle \mathcal A_{w_k}\left(x,t \right)\nabla I_k, \nabla \varphi \rangle dx dt  
- \int_0^T \int_{\Omega} \langle \mathcal A_{w}\left(x,t \right)\nabla I, \nabla \varphi \rangle dx dt \bigg|\\\notag
 &\leq \int_0^T \int_{\Omega} \| \mathcal A_{w_k}\left(x,t \right) - \mathcal A_{w}\left(x,t \right) \| \langle \nabla I_k, \nabla \varphi \rangle dx dt \\\notag
&+ \int_0^T \int_{\Omega} \| \mathcal A_{w}\left(x,t \right)\| |\nabla (I_k - I)| |\nabla \varphi| dx dt \longrightarrow 0 \quad \text{as}\quad k \to \infty.
\end{align}
Using \eqref{eq:conv-1} and \eqref{eq:conv-2} it follows that 
\begin{align}\label{eq:conv-last}
  \int_0^T    \int_{\Omega}\partial_t I \varphi dx dt
    + \int_0^T \int_{\Omega} \langle \mathcal A_{w}\left(x,t \right)\nabla I, \nabla \varphi \rangle dx dt = 0,
\end{align}
and hence, we can conclude that $I = I_w = \mathcal{P}(w)$.
Since the solution of $\mathcal{P}$ is unique, the whole sequence $I_k = \mathcal{P}(w_k)$ converges weakly in $W_0$ to $I = \mathcal P(w)$,
and thus we have showed that the mapping $\mathcal P$ is weakly continuous from $W_0$ to $W_0$. 
Consequently, thanks to Schauder’s fixed point theorem, there exists $w \in W_0$ such that $w = P(w) = I_w$, which solves the problem \eqref{para-general1}--\eqref{para-general3}.
This proof of existence is complete.\\


\noindent
\emph{Proof of uniqueness}. Let $\bar I$ and $\hat I$ be two solutions of \eqref{para-general1}--\eqref{para-general3}. 
We have, for almost every $t \in (0, T)$,
\begin{align*}
&\bar I_t - \nabla \cdot \left( \bar \alpha(x,t)\, \nabla \bar I\right) = 0, \quad 
 \bar I(x,0) = I_0(x) \quad \text{in}\quad \Omega, \quad\text{and} \quad
 \partial_n \bar I = 0 \quad \text{on} \quad \partial \Omega \times (0,T), \\
&\hat I_t - \nabla \cdot \left( \hat \alpha(x,t)\, \nabla \hat I\right) = 0, \quad 
 \hat I(x,0) = I_0(x) \quad \text{in}\quad \Omega, \quad\text{and} \quad
 \partial_n \hat I = 0 \quad \text{on} \quad \partial \Omega \times (0,T), 
\end{align*}
where $\bar \alpha(x,t) = \mathcal A(x,t,\bar I,\nabla \bar I, \bar I(\cdot),\nabla \bar I(\cdot))$ and
$\hat \alpha(x,t) = \mathcal A(x,t,\hat I,\nabla \hat I, \hat I(\cdot),\nabla \hat I(\cdot))$ are $n \times n$ matrices.
Using the above equations, we obtain
\begin{align}\label{eq:unique-explain}
I_t- \nabla \cdot \left(\bar \alpha(x,t) \nabla I \right) = \nabla \cdot \left(\left(\bar \alpha(x,t) - \hat\alpha(x,t)\right)\nabla \hat I\right),    
\end{align}
where $I=\bar I(x,t) - \hat I(x,t)$. 
Now, multiplying this by $I$, 
integrating over $\Omega$, 
using integration by parts and the boundary condition, we get
\begin{align*}
& \frac{1}{2}\dfrac{d}{dt}\large\int_{\Omega}I^2 dx + \int_{\Omega} \langle\bar  \alpha(x,t) \nabla I, \nabla I\rangle dx = - \int_{\Omega}\langle\left(\bar \alpha(x,t) - \hat\alpha(x,t) \right)\nabla \hat I, \nabla I\rangle dx  
\end{align*}
Using the Cauchy Schwarz inequality and \eqref{eq:ass-A1-diffusion}, we have
\begin{align*}
 & \frac{1}{2}\dfrac{d}{dt} \left( \Vert I \Vert_{L^2}^2 \right) + \theta \Vert \nabla I \Vert_{L^2}^2  \leq  \left(\Vert \bar \alpha(x,t) - \hat\alpha(x,t) \Vert\right)_{L^\infty}\Vert \nabla  \hat I \Vert_{L^2}\Vert \nabla I \Vert_{L^2}.
\end{align*}
Now, using assumption \eqref{eq:ass-A2-diffusion}, the definition of $\bar \alpha$ and $\hat \alpha$ and Young's inequality, we can obtain
\begin{align*}
  \frac{1}{2}\dfrac{d}{dt} \left( \Vert I \Vert_{L^2}^2 \right) + \theta \Vert \nabla I \Vert_{L^2}^2  &\leq \frac{2}{\theta} C \Vert \nabla  \hat I \Vert_{L^2}^2 \Vert I \Vert_{L^2}^2 +     \frac{\theta}{2} \Vert \nabla I \Vert_{L^2}^2,\quad 
\text{or}\\ 
  \dfrac{d}{dt} \left( \Vert I \Vert_{L^2}^2 \right)  &\leq \frac{4}{\theta} C \Vert \nabla  \hat I \Vert_{L^2}^2 \Vert I \Vert_{L^2}^2.
\end{align*}
Since $I(x,0) = \bar I(x,0) - \hat I(x,0) = I_0(x) - I_0(x) = 0$, by using the above relation and Gronwall’s lemma, we can complete the proof of uniqueness. $\hfill\Box$

\begin{cor}\label{corr:diffusion-final}
Problem \eqref{diff_power_constant}  with diffusion function \eqref{diff_edgeg} satisfying \eqref{eq:ass-exponent_p} admits a weak solution in $ W(0,T) $, unique in the same space, if $J \in H^1(\Omega)$ and $0 < \alpha_1 \leq I_0 \leq \beta_1$.
\end{cor}

\noindent
{\bf Proof.}
It suffices to show that Theorem \ref{thm:existence-diffusion} applies to our suggested model \eqref{diff_power_constant}, i.e. to show that the corresponding function in $\mathcal{A}_w$ defined in \eqref{eq:def_A_diffusion} satisfies assumption \eqref{eq:ass-A1-diffusion}.
We have, for any $w \in W_{1}$,
\begin{align}\label{eq:AW-diffusion}
\mathcal{A}_w\left(x,t\right) =
\epsilon
+ \frac{ 2| w_\xi|^\nu}{\left( M^{w}_{\xi}\right)^\nu + |w_\xi|^\nu}  \cdot \dfrac{1}{1+ \left(\frac{\vert \nabla w_\xi \vert }{K} \right)^{p\left(w_\xi, \vert \nabla w_\xi \vert \right)} }\,.
\end{align}
We first observe that
\begin{align}\label{eq:A1-particular_diffusion}
\epsilon \leq \mathcal{A}_w\left(x,t\right) \leq 1,
\end{align}
and hence $A_w$ satisfies \eqref{eq:ass-A1-diffusion} with $\theta = \epsilon$ and $C = 1$.

We now prove that for the suggested model \eqref{diff_power_constant}, with the diffusion coefficients \eqref{diff_constatp} and \eqref{diff_edgeg} satisfies the assumption \eqref{eq:ass-A2-diffusion}. To do so, we first rewrite the left-hand side of \eqref{eq:ass-A2-diffusion} as
\begin{align}\label{eq:est-for-A3}
&\left|\mathcal A(x,t,I_1,\nabla I_1,I_1(\cdot),\nabla I_1(\cdot) ) - \mathcal A(x,t,I_2,\nabla I_2,I_2(\cdot),\nabla I_2(\cdot) )\right|\notag\\
&\leq \frac{ 2|(I_1)_\xi|^\nu}{\left(M^{(I_1)}_{\xi}\right)^\nu + | (I_1)_\xi|^\nu} \left| \frac{1}{1+ \left(\frac{\vert \nabla (I_1)_\xi \vert }{K} \right)^{ p\left((I_1)_\xi, \vert \nabla (I_1)_\xi \vert \right)  } }  - \frac{1}{1+ \left(\frac{\vert \nabla (I_2)_\xi \vert }{K} \right)^{ p\left((I_2)_\xi, \vert \nabla (I_2)_\xi \vert \right) } }\right|\notag\\
&+ \frac{1}{1+ \left(\frac{\vert \nabla (I_2)_\xi \vert }{K} \right)^{p\left((I_2)_\xi, \vert \nabla (I_2)_\xi \vert \right)} } \left| \frac{ 2|(I_1)_\xi|^\nu}{\left(M^{(I_1)}_{\xi}\right)^\nu + |(I_1)_\xi|^\nu}  - \frac{ 2| (I_2)_\xi|^\nu}{\left( M^{(I_2)}_{\xi}\right)^\nu + | (I_2)_\xi|^\nu} \right|\,,
\end{align}
where $I_{1}$ and $I_{2}$ are two weak solutions of \eqref{diff_constatp}. Put $X_i = \frac{\vert \nabla (I_i)_\xi \vert }{K}$ and $p_i = p\left((I_i)_\xi, \vert \nabla (I_i)_\xi \vert \right)$ for $i = 1,2$ and observe that
\begin{align}\label{eq:x1-x2:1}
 \frac{1}{1+ X_1^{p_1}}  - \frac{1}{1+ X_2^{p_2} }
 =\frac{X_1^{p_1} - X_2^{p_2} }{\left( 1+ X_1^{p_1} \right) \left(1+ X_2^{p_2}\right)} 
\leq 
\frac{1}{\left( 1+ X_1^{p_1} \right) \left(1+ X_2^{p_2}\right)} 
\left(\left|X_1^{p_1} - X_2^{p_1}  \right|
+
\left|X_2^{p_1} - X_2^{p_2}  \right|\right)\,.
\end{align}
We have
\begin{align}\label{eq:x1-x2:2}
\left|X_1^{p_1} - X_2^{p_1}  \right|
\leq
\left|X_1 - X_2 \right| \left( 1 + X_1^{p_1} + X_2^{p_1} \right),
\end{align}
in which $\left|X_1 - X_2 \right| \leq C\|\nabla (I_1)_\xi| - |\nabla (I_2)_\xi\| \leq  C|\nabla (I_1)_\xi - \nabla (I_2)_\xi| \leq C |\nabla G_\xi| |I_1 - I_2|$. Next, 
\begin{align}\label{eq:x1-x2:3}
\left|X_2^{p_1} - X_2^{p_2}  \right|
\leq
\log\left(X_2\right) \left(X_2^{p_1} + X_2^{p_2}\right) |p_2 - p_1|
\leq  \left(1 + X_2^{p_1+1} + X_2^{p_2+1}\right) |p_2 - p_1|\,,
\end{align}
where we used $p(\cdot,\cdot) \geq 1$, assumed in \eqref{eq:ass-exponent_p}, in the last inequality.
We now use the differentiability assumption of $p(\cdot,\cdot)$ in \eqref{eq:ass-exponent_p} to conclude (Lipschitz continuity suffices)
\begin{align}\label{eq:p1-p2}
|p_2 - p_1| &= 
|p\left((I_2)_\xi, \vert \nabla (I_2)_\xi \vert \right) - p\left((I_1)_\xi, \vert \nabla (I_1)_\xi \vert \right)|\notag\\
&\leq C \left( |(I_1)_\xi - (I_2)_\xi| + \|\nabla (I_1)_\xi| - |\nabla (I_2)_\xi\|\right)\notag\\
&\leq
C \left( |(I_1)_\xi - (I_2)_\xi|  +  |\nabla (I_1)_\xi - \nabla (I_2)_\xi|\right)\notag\\
&\leq
C \left( |G_\xi * (I_1 - I_2)|  +  |\nabla G_\xi * (I_1 - I_2)|\right)\notag\\
&\leq
C \left( |G_\xi|  +  |\nabla G_\xi| \right)  |I_1 - I_2|.
\end{align}
Using \eqref{eq:x1-x2:1}--\eqref{eq:p1-p2} we conclude
\begin{align}\label{eq:A3-check-1}
\left| \frac{1}{1+ X_1^{p_1}}  - \frac{1}{1+ X_2^{p_2} } \right| 
\leq C \| I_1 - I_2 \|_{L^2},
\end{align}
where $C$ may depend on $\xi$, $\Omega,$ $K$, $I_0$. We can estimate the second term in \eqref{eq:est-for-A3} as follows
\begin{align}\label{eq:A3-check-2}
& \left| \frac{ |  (I_1)_\xi|^\nu  }{\left( M^{(I_1)}_{\xi}\right)^\nu + |  (I_1)_\xi|^\nu}  - \frac{|  (I_2)_\xi|^\nu}{\left( M^{(I_2)}_{\xi}\right)^\nu + |  (I_2)_\xi|^\nu} \right| \notag \\
& \le \left\Vert |(I_1)_\xi|^\nu \left\{ \left( M^{(I_2)}_{\xi}\right)^\nu -\left( M^{(I_1)}_{\xi}\right)^\nu \right\} \right\Vert_{L^{\infty}} + \left\Vert \left( M^{(I_1)}_{\xi}\right)^\nu \left\{ | (I_1)_\xi|^\nu-  | (I_2)_\xi|^\nu  \right\} \right\Vert_{L^{\infty}} \notag\\
&\le  C \left\{ \left\Vert    M^{(I_2)}_{\xi} - M^{(I_1)}_{\xi}  \right\Vert_{L^{\infty}} +   \left\Vert   (I_1)_\xi - (I_2)_\xi  \right\Vert_{L^{\infty}} \right\} \notag\\
&\qquad(\text{where we have used that $0 < \alpha_1 \leq I_1, I_2 \leq \beta_1$ by Lemma \ref{le:bounded_2024}})\notag\\
&\le  C \left\{ \left\Vert  G_\xi \ast I_2(x,t)- G_\xi \ast I_1(x,t)   \right\Vert_{L^{\infty}} +   \left\Vert  I_1 - I_2 \right\Vert_{L^{2}} \right\} \notag\\
&\le  C \left\Vert  I_1 - I_2 \right\Vert_{L^{2}} \,,
\end{align}
where $C$ may depend on $\xi$, the diffusion coefficient $g$, $\alpha_1$ and $\beta_1$. Using the estimates \eqref{eq:A3-check-1} and \eqref{eq:A3-check-2} in \eqref{eq:est-for-A3} shows that we can satisfy assumption \eqref{eq:ass-A2-diffusion}. 
Hence the proof of Corollary \ref{corr:diffusion-final} is completed.
$\hfill \Box$\\

\subsection{Telegraph Diffusion Model}
In this section, we will establish the existence and uniqueness of weak solutions of a class of PDEs, including the proposed model \eqref{maina}--\eqref{mainb} with the diffusion coefficient \eqref{diff_edgeg}. 
Let $\Omega$ denote a bounded domain in $\mathbf{R}^n$ and let $T > 0$.
We consider hyperbolic PDEs on the form
\begin{align}
 &I_{tt} + \lambda I_t - \nabla \cdot\left( \mathcal A\, \nabla I\right) = 0 \quad (x,t) \in \Omega \times (0,T),   \tag{E.1}\label{eq:PDE-general1}\\
 &I(x,0) = I_0(x), \quad I_t(x,0) = 0 \quad x \in \Omega,  \tag{E.2} \label{eq:PDE-general2}\\
 &\langle\mathcal A \nabla I, {\bf n}\rangle = 0  \quad \text{on} \quad \partial \Omega \times (0,T),  \tag{E.3} \label{eq:PDE-general3}
\end{align}
where $\mathcal A = \mathcal A(x,t,I,\nabla I, I(\cdot),\nabla I(\cdot)): \mathbf{R}^n \times \mathbf{R} \times \mathbf{R} \times \mathbf{R}^n \times \mathbf{R} \times \mathbf{R}^n \to \mathbf{R}^{n\times n}$
is a once differentiable symmetric matrix 
satisfying some additional assumptions given below. 
Here, $I(\cdot)$ and $\nabla I(\cdot)$ indicate dependence of $I(x,t)$ and $\nabla I(x,t)$ in a non-local manner.
We will make use of the solution space 
\begin{align*}
\mathcal{W}(0,T)&=\Big\{w\in L^\infty(0,T; H^1)\,, w_t \in L^\infty(0,T; L^2); \,w_{tt} \in L^2(0,T; (H^1)') \Big\}\,,
\end{align*}
%
which is a Hilbert space for the graph norm \cite{jllions1968}, 
and we adopt the following version of weak solutions:
\begin{defi}[Weak solution]\label{defi:weak}
A function $I = I(x,t)$ is called a weak solution of \eqref{eq:PDE-general1}--\eqref{eq:PDE-general3} if
\begin{itemize}
 \item[a)] $I \in \mathcal{W}(0, T) $ and \eqref{eq:PDE-general2} and \eqref{eq:PDE-general3} hold a.e. with derivatives in the sense of distributions.
 \item[b)] For all $\phi \in H^1$ and a.e. $t\in (0,T)$ it holds that
 \begin{align*}
\left \langle I_{tt}, \phi   \right \rangle + {\displaystyle \int_{\Omega}}
\Big(  I_t \phi +  \langle \mathcal A\left(x,t,I,\nabla I, I(\cdot),\nabla I(\cdot)\right)  \nabla I, \nabla \phi\rangle  \Big)\,dx =0.
  \end{align*}
\end{itemize}
\end{defi}
%
\noindent
Next we define, for a constant $M_1 > 0$, the function space
\begin{align}\label{eq:M_1}
 \mathcal{W}_{M_1}= \big\{  v\in \mathcal{W}(0,T):~  \|v\|_{L^\infty(0,T;H^1)} + \|v_t\|_{L^\infty(0,T; L^2)} \le M_1  \|I_0\|_{H^1(\Omega)} \big\},
\end{align}
and for any element $w \in \mathcal{W}_{M_1}$ we consider the function $\mathcal{A}_w : \mathbf{R}^n \times \mathbf{R} \to \mathbf{R}^{n\times n}$ defined as
\begin{align}\label{eq:def_A}
\mathcal{A}_w (x,t) = \mathcal A\left(x,t,w,\nabla w, w(\cdot),\nabla w(\cdot)\right).
\end{align}
To prove existence and uniqueness, we assume that there exist continuous functions $\phi_\mathcal A, \Phi_\mathcal A, \Psi_\mathcal A: \mathbf{R}_+ \to \mathbf{R}_+$,
$\phi_\mathcal A^{-1}, \Phi_\mathcal A, \Psi_\mathcal A$ nondecreasing, such that the following hold, where $z = M_1 \|I_0\|_{H^1(\Omega)}$:
\begin{align}
\phi_\mathcal A(z) |\eta|^2 \leq \langle \mathcal A_w(x,t) \eta, 
 \eta\rangle, \,\, \|\mathcal{A}_w\|_{L^\infty(\Omega)} &\leq \Phi_\mathcal A(z), \,\, \text{whenever}\,\, (x,t) \in \Omega\times (0,T), \; \eta \in \mathbf{R}^n,\; w \in \mathcal{W}. 
 \tag{A.1}\label{eq:ass-A1}\\
  \left\|\partial_t \mathcal A_w(x,t)\right\|_{L^\infty(\Omega)} &\leq  \Psi_\mathcal A(z), \, \,\text{whenever}\,\, (x,t) \in \Omega\times (0,T)\; w \in \mathcal{W}.
  \tag{A.2}\label{eq:ass-A2}
\end{align}
We also assume that there exists a constant $C$ such that, whenever $I_1, I_2 \in \mathcal{W}$, the following holds:
\begin{align}
& \left(\|\mathcal A(x,t,I_1,\nabla I_1,I_1(\cdot),\nabla I_1(\cdot) ) - \mathcal A(x,t,I_2,\nabla I_2,I_2(\cdot),\nabla I_2(\cdot) )\|\right)_{L^\infty(\Omega)} \notag \\
& \hspace{7cm}  \leq C\|I_1(x,t) - I_2(x,t)\|_{L^2(\Omega)} \quad \text{a.e.}\quad  t \in (0,T).\tag{A.3}\label{eq:ass-A3}
\end{align}
Our main result of well-posedness is summarized in the following theorem.
\begin{thm}\label{thm:existence}
Suppose that assumptions \eqref{eq:ass-A1}--\eqref{eq:ass-A3}  hold, that $\Omega \in \mathbf{R}^n$ is bounded and that $I_0 \in H^2(\Omega)$.
Then there exists a weak solution of problem \eqref{eq:PDE-general1}--\eqref{eq:PDE-general3} if at least one of the following statements $(i)-(iii)$ hold true:
\begin{itemize}
  \item[$(i)$] $\|I_0\|_{H^1(\Omega)}$ is small enough. 
  \item[$(ii)$] \eqref{eq:ass-functions-1} holds and $T$ is small enough.
  \item[$(iii)$] \eqref{eq:ass-functions-2} holds.
\end{itemize}
Moreover, if \eqref{eq:ass-A1}--\eqref{eq:ass-A3} hold then any weak solution in $\mathcal{W}$ of problem \eqref{eq:PDE-general1}--\eqref{eq:PDE-general3} is unique.
\end{thm}
To prove Theorem \ref{thm:existence}, we will use the Schauder fixed-point theorem together with properties of solutions of linear PDEs on the form
\begin{equation}\label{eq:linearPDE}
\left.\begin{aligned}
 &I_{tt} + \lambda I_t - \nabla \cdot \left( a(x,t) \nabla I\right) = 0 \quad (x,t) \in \Omega \times (0,T),  \\
 &I(x,0) = I_0(x), \quad  I_t(x,0) = 0 \quad x \in \Omega,\\
 &\langle a(x,t) \nabla I, {\bf n}\rangle = 0 \quad \text{on} \quad \partial \Omega,
\end{aligned}\right\}
\end{equation}
where $a(x,t) : \mathbf{R}^n \times \mathbf{R}  \to \mathbf{R}^{n\times n}$ is an a.e. once differentiable symmetric matrix satisfying \eqref{eq:ass-A1}--\eqref{eq:ass-A2} so that the equations are uniformly hyperbolic. In particular, we will use the following classical result, which can be proved by the Galerkin method, see, e.g., \cite{LCEvans1998}.
\begin{lem}\label{le:exist-linear}
If $I_0 \in H^1(\Omega)$, then the linear problem \eqref{eq:linearPDE} has a unique weak solution $I \in \mathcal{W}$.
In addition, if $I_0 \in H^2(\Omega)$ then $\frac{\partial I}{\partial t} \in L^\infty (0,T;H^1(\Omega))$.
\end{lem}
\noindent
{\bf Proof of Theorem \ref{thm:existence}.}
We begin with the existence part. 
Let $w \in \mathcal{W}_{M_1}$ for some constant $M_1 > 0$ that will be determined later and let $\mathcal{A}_w$ be as in \eqref{eq:def_A}.
Define ($Q_w$) as the problem
\begin{align}\label{eq:P_w}
    \left\langle I_{tt}, \varphi \right\rangle_{H^1(\Omega)' \times H^1(\Omega)}
    + \int_{\Omega} \left(\lambda I_t \varphi + \langle \mathcal A_w\left(x,t \right)\nabla I, \nabla \varphi\rangle \right) dx = 0\,,
\end{align}
whenever $\varphi \in H^1(\Omega)$ and a.e. $t \in (0,T)$, together with $I(x,0) = I_0$, $\frac{\partial I}{\partial t}(x,0) = 0$ and the boundary conditions in \eqref{eq:linearPDE} in the sense of distribution.
For each $w \in \mathcal{W}_{M_1}$, the problem ($Q_w$) defines weak solutions of the linear PDE \eqref{eq:linearPDE}.
Since $w \in \mathcal{W}_{M_1}$ it follows by assumption that \eqref{eq:ass-A1} and \eqref{eq:ass-A2} hold
and we may therefore use Lemma \ref{le:exist-linear} to conclude that ($Q_w$) has a unique weak solution $I_w \in \mathcal{W}$ with
$\frac{\partial I_w}{\partial t} \in  L^\infty(0,T; H^1(\Omega))$.

We will now show that the solution of the linear problem denoted $I = I_w$ for simplicity,
belongs to the space $\mathcal{W}_{M_1}$ defined in \eqref{eq:M_1}.
Since $I_t \in L^\infty(0,T; H^1)$ we can take $\phi = I_t$ in equation \eqref{eq:P_w} to obtain
\begin{align}\label{eq:stet_1}
&\frac{1}{2}\dfrac{d}{dt}\|I_t\|_{L^2 }^2 + \lambda \|I_t\|_{L^2 }^2
+  \large\int_{\Omega} \langle \mathcal A_w(x,t)\nabla I, \nabla \big(I_t\big)\rangle\, dx=0 \,.
\end{align}
Using integration by parts and \eqref{eq:ass-A2} we conclude that
\begin{align}\label{ext_1}
\large\int_{\Omega} \langle \mathcal A_w(x,t) \nabla I, \nabla \big(I_t\big)\rangle \, dx
&\geq  \frac{1}{2}\dfrac{d}{dt}\int_{\Omega} \langle\mathcal A_w(x,t) \nabla I, \nabla I\rangle \,dx- \frac{1}{2}  \Psi_\mathcal A(z) \|\nabla I\|_{L^2 }^2\,,
\end{align}
and thanks to the lower bound in \eqref{eq:ass-A1} we also have
\begin{align}\label{esti:gradient-inters-gw}
\|\nabla I\|_{L^2 }^2 \le \phi_\mathcal A(z)^{-1} \int_{\Omega} \langle\mathcal A_w(x,t) \nabla I, \nabla I\rangle\, dx\,.
\end{align}
Inserting \eqref{ext_1} in \eqref{eq:stet_1} gives
\begin{align*}
&\frac{1}{2}\dfrac{d}{dt}\|I_t\|_{L^2 }^2 + \lambda \|I_t\|_{L^2 }^2
\leq - \frac{1}{2}\dfrac{d}{dt}\int_{\Omega} \langle\mathcal A_w(x,t) \nabla I, \nabla I\rangle \,dx + \frac{1}{2}\Psi_\mathcal A(z) \|\nabla I\|_{L^2 }^2\,,
\end{align*}
and with \eqref{esti:gradient-inters-gw} this becomes
\begin{align*}
&\frac{1}{2}\dfrac{d}{dt}\|I_t\|_{L^2 }^2 + \lambda \|I_t\|_{L^2 }^2
\leq - \frac{1}{2}\dfrac{d}{dt}\int_{\Omega}\langle\mathcal A_w(x,t) \nabla I, \nabla I\rangle \,dx + \frac{1}{2}  \frac{\Psi_\mathcal A(z)}{\phi_\mathcal A(z)} \int_{\Omega}\langle\mathcal A_w(x,t) \nabla I, \nabla I\rangle\, dx\,.
\end{align*}
Thus
\begin{align*}
 \frac{d}{dt} \Big[\| I_t\|_{L^2 }^2  + \int_{\Omega}\langle\mathcal A_w(x,t) \nabla I, \nabla I\rangle\, dx\Big]
\leq  - 2  \lambda \|I_t\|_{L^2 }^2 +   \frac{\Psi_\mathcal A(z)}{\phi_\mathcal A(z)} \int_{\Omega}\langle\mathcal A_w(x,t) \nabla I, \nabla I\rangle\, dx\,.
\end{align*}
In order to use Gronwall's lemma, we observe that the above inequality implies
\begin{align*}
  \frac{d}{dt} \Theta(t)  \leq  \frac{\Psi_\mathcal A(z)}{\phi_\mathcal A(z)} \,\Theta(t)
\end{align*}
in which $\Theta(t) = \| I_t\|_{L^2 }^2  + \int_{\Omega}\langle\mathcal A_w(x,t) \nabla I, \nabla I\rangle\, dx$.
Hence
\begin{align*}
\Theta(t)
  \leq  \Theta(0) e^{ \frac{\Psi_\mathcal A(z)}{\phi_\mathcal A(z)} \,t} \qquad \text{for a.e.} \qquad t\in [0,T].
\end{align*}
The initial conditions imply that $\Theta(0) \leq \left(\|\mathcal A_w(x,0) \|\right)_{L^\infty } \|\nabla I \|_{L^2 }^2 \leq \Phi_\mathcal A(z)\|I_0\|_{H^1 }^2$.
Thus
\begin{align}\label{bound_I_t}
\| I_t\|_{L^2 }^2   \leq  \Phi_\mathcal A(z) \|I_0\|_{H^1 }^2e^{ \frac{\Psi_\mathcal A(z)}{\phi_\mathcal A(z)} \,t} \qquad \text{for a.e.} \qquad t\in [0,T]\,.
\end{align}
Using also the lower bound in \eqref{eq:ass-A1} we conclude
\begin{align}\label{bound_I_t_nabla_I}
\| I_t\|_{L^2 }^2   + \|\nabla I\|_{L^2 }^2 \leq  \frac{\Phi_\mathcal A(z)}{\phi_\mathcal A(z)} \|I_0\|_{H^1 }^2 e^{ \frac{\Psi_\mathcal A(z)}{\phi_\mathcal A(z)} \,t} \qquad \text{for a.e.} \qquad t\in [0,T]\,.
\end{align}
Since $I(x,t)=I_0(x)+ \displaystyle \int_{0}^{t} I_t(x,s)\,ds$ we have, thanks to Young's inequality, H\"older's inequality and \eqref{bound_I_t},
\begin{align}\label{boundIL2}
 \| I\|_{L^2 }^2 = \int_\Omega I(x,t)^2 \, dx &= \int_\Omega \left(I_0(x)\,  +  \int_{0}^{t} I_t(x,s)\,ds\right)^2 \,dx \notag\\
&\leq 2 \int_\Omega I_0(x)^2\,dx  +  2 \int_\Omega \left(\int_{0}^{t} I_t(x,s)\,ds\right)^2 \,dx \notag\\
&\leq 2\left\Vert I_0 \right\Vert^2_{H^1 } +  2  T^2 \| I_t\|_{L^2 }^2 \notag\\
&\leq 2\left\Vert I_0 \right\Vert^2_{H^1 } +  2  T^2  \Phi_\mathcal A(z) \|I_0\|_{H^1 }^2 e^{ \frac{\Psi_\mathcal A(z)}{\phi_\mathcal A(z)} \,T}.
\end{align}
Combining \eqref{bound_I_t_nabla_I} and \eqref{boundIL2} we see that
\begin{align*}
\| I(t)\|_{H^1 }^2  = \| I(t)\|_{L^2 }^2  + \|\nabla I(t)\|_{L^2 }^2
\leq  \|I_0\|_{H^1 }^2 \left( \left(\phi_\mathcal A(z)^{-1} + 2  T^2\right) \Phi_\mathcal A(z) e^{ \frac{\Psi_\mathcal A(z)}{\phi_\mathcal A(z)} \,T}
+ 2 \right)\,
\end{align*}
for a.e. $t \in [0,T]$.
Therefore
$$
\| I\|_{L^{\infty}(0,T;H^{1})} \leq \|I_0\|_{H^1 } \sqrt{ \left(\phi_\mathcal A(z)^{-1} + 2  T^2\right) \Phi_\mathcal A(z) e^{ \frac{\Psi_\mathcal A(z)}{\phi_\mathcal A(z)} \,T}
+ 2 }\,
$$
and so, using \eqref{bound_I_t_nabla_I} once again, we have
\begin{align*}
\| I\|_{L^{\infty}(0,T;H^{1})} + \| I_t\|_{L^{\infty}(0,T;L^2)} &\leq
2\|I_0\|_{H^1 } \sqrt{ \left(\phi_\mathcal A(z)^{-1} + 2  T^2\right) \Phi_\mathcal A(z) e^{ \frac{\Psi_\mathcal A(z)}{\phi_\mathcal A(z)} \,T}
+ 2 }\,.
\end{align*}
If we obtain
\begin{align}\label{eq:find-constants}
2 \sqrt{ \left(\phi_\mathcal A(z)^{-1} + 2  T^2\right) \Phi_\mathcal A(z) e^{ \frac{\Psi_\mathcal A(z)}{\phi_\mathcal A(z)} \,T}
+ 2 }\, \leq M_1,
\end{align}
in which $z = M_1 \|I_0\|_{H^1 }$, then
\begin{align}\label{eq:bound-existence-1}
\| I\|_{L^{\infty}(0,T;H^{1})} + \| I_t\|_{L^{\infty}(0,T;L^2)} &\leq\,
 M_1 \|I_0\|_{H^1 } \,
\end{align}
and hence the solution $I$ of the linear problem belongs to the space $\mathcal{W}_{M_1}$ defined in
\eqref{eq:M_1}.
To get \eqref{eq:find-constants} satisfied we can assume, for example, that $\|I_0\|_{H^1 }$
is small. In particular, if $z$ is held constant when $M_1$ increases, then there exists $M_1$ such that \eqref{eq:find-constants} holds.
Alternatively, we may assume, for example, that
\begin{align}\label{eq:ass-functions-1}
\frac{1}{z^2}\frac{\Phi_\mathcal A(z)}{\phi_\mathcal A(z)} \to 0 \qquad \text{as} \qquad z \to \infty. \tag{A.4}
\end{align}
In such case, we can satisfy \eqref{eq:find-constants} by taking a large enough $M_1$ and then a small enough $T$.
If we assume that there exists a constant $C$, independent of $z$, such that
\begin{align}\label{eq:ass-functions-2}
\frac{\Phi_\mathcal A(z)}{\phi_\mathcal A(z)} \leq C \qquad \text{and} \qquad \frac{\Psi_\mathcal A(z)}{\phi_\mathcal A(z)} \leq C, \tag{A.5}
\end{align}
then we obtain \eqref{eq:find-constants} by increasing $M_1$. 

We will now show that there exists a constant $C$, which may depend on $\mathcal A$, $\Omega$, $I_0$ and $T$, such that the solution of the linear problem \eqref{eq:P_w}, $I = I_w$, satisfies
\begin{align}\label{eq:bound-existence-2}
\| I_{tt}\|_{L^{\infty}(0,T;{(H^1)^\prime})} \leq C.
\end{align}
Taking $\phi \in H^1$ with $ \Vert \phi \Vert_{H^1}\leq 1$ in \eqref{eq:P_w} we have
\begin{align*}
\left \langle I_{tt}, \phi   \right \rangle + {\displaystyle \int_{\Omega}}
\left(  I_t \phi +   \langle \mathcal A_w(x,t)  \nabla I, \nabla \phi \rangle \right)\,dx =0\,.
\end{align*}
We use Cauchy-Schwarz inequality along with \eqref{eq:bound-existence-1} and boundedness of $\mathcal A_w$  to obtain
\begin{align*}
\begin{split}
\big| \left\langle I_{tt}, \phi \right\rangle \big|
\leq & \left(\left\| I_t \right\|_{L^2 } +\left(\left\| \mathcal A_w(x,t)\right\|\right)_{L^{\infty} }\left\Vert \nabla I\right\Vert_{L^2 } \right)\|\phi\|_{H^1 } \\
\leq & \left(2 M_1 + 2 M_1 \Phi_{\mathcal{A}}(z)\right) \|I_0\|_{H^1 } \left\Vert\phi\right\Vert_{H^1 } \leq C \,.
\end{split}
\end{align*}
for a constant $C$ which may depending on $\mathcal A$, $\Omega$, $I_0$ and $T$. Hence, by the definition of norm in $(H^1)^\prime$, we get
$
\left\Vert I_{tt} \right\Vert_{(H^1)^\prime} \leq C 
$
for all $t \in (0,T)$. This implies $\| I_{tt}\|_{L^{\infty}(0,T;{(H^1)^\prime})} \leq C \|I_0\|_{H^1 }$ which is the desired result.
%


With \eqref{eq:bound-existence-1} and \eqref{eq:bound-existence-2} in mind,  we introduce the non-empty, convex, and weakly compact subspace $\mathcal{W}_0$ of $\mathcal{W}(0, T)$ by
\begin{align*}
\mathcal{W}_0=&\Big\{w \in \mathcal{W}(0,T):\, \|w\|_{L^\infty(0,T; H^1)} + \|\partial_t w\|_{L^\infty(0,T; L^2)} \leq C\|I_0\|_{H^1}\,;\\
    &\hspace{7.5cm}\|\partial_{tt} w\|_{L^2\left(0,T; (H^1)^\prime\right)}\leq C
   \Big\}\,
\end{align*}
and consider the mapping
\begin{align*}
\mathcal{P}:~ & \mathcal{W}_0 \goto \mathcal{W}_0 \\
& w\mapsto I_w\,.
\end{align*}
To use Schauder's fixed-point theorem on $\mathcal{P}$, we have to prove that the mapping $\mathcal{P}:w \rightarrow I_w $ is weakly continuous from $\mathcal{W}_0$ into $\mathcal{W}_0$. Let $w_k$ be a sequence that converges weakly to some $w$ in $\mathcal{W}_0$ and let $I_k = I_{w_k}$. We have to show that $\mathcal{P}(w_k):= I_k$ converges weakly to $\mathcal{P}(w): = I_w$. From the classical results of compact inclusion in Sobolev spaces \cite{raadams1975}, we can extract subsequences of $\{w_k\}$ and $\{I_k\}$ still denoted by $\{w_k\}$ and $\{I_k\}$ respectively such that for some $I\in \mathcal{W}_0$, we have, as $k\goto \infty$, 
\begin{align*}
\begin{cases}
w_{k} \longrightarrow  w \hspace{0.2cm} \text{in} \hspace{0.2cm} L^2(0,T;L^2) \hspace{0.2cm} \text{ and a.e. on } \hspace{0.2cm}  \Omega_T,\\[0.5em]
\displaystyle I_{k} \longrightarrow  I\, \hspace{0.2cm} \text{weakly-$\ast$} ~ \text{in}~~L^{\infty}(0,T;H^1)\,,\\[1em]
\displaystyle I_{k} \longrightarrow   I\, \hspace{0.2cm} \text{in}~~L^{2}(0,T; L^2)\,,\\[1em]
\partial_t I_k \longrightarrow   \partial_t I\, \hspace{0.2cm} \text{weakly-$\ast$} ~\text{in}~~L^{\infty}(0,T;L^2)\,,\\[1em]
\partial_{tt} I_k \longrightarrow \partial_{tt}I \hspace{0.2cm} \text{weakly-$\ast$}~ \text{in}~~L^{2}(0,T;(H^1)^\prime)\,,\\[1em]
\end{cases}
\end{align*}
The above convergence, together with an analogue of \eqref{eq:conv-first}--\eqref{eq:conv-last},  allows us to pass to the limit in the problem ($Q_w$) and obtain $I=\mathcal{P}(w)$.  Moreover, since the solution of ($Q_w$) is unique, the whole
sequence $I_k=\mathcal{P}(w_k)$ converges weakly in $\mathcal{W}_0$ to $I=\mathcal{P}(w)$. Hence, $\mathcal{P}$ is weakly continuous. Consequently, thanks to the Schauder fixed point theorem, there exists $w \in \mathcal{W}_0$ such that $w=\mathcal{P}(w)=I_w$.  Thus, the function $I_w$ solves the problem \eqref{eq:PDE-general1}--\eqref{eq:PDE-general3}.
This completes the proof of existence.

\noindent \emph{{Proof of uniqueness}}. Following the method in \cite{LCEvans1998,majee2020gray} we prove the uniqueness of weak solutions of problem \eqref{eq:PDE-general1}--\eqref{eq:PDE-general3}. Let $I_{1}$ and $I_{2}$ be two weak solutions of \eqref{eq:PDE-general1}--\eqref{eq:PDE-general3}.
As in \eqref{eq:unique-explain} we then obtain
\begin{align}
&I_{tt}+ I_t-\text{div} \big(\mathcal{A}_1 \nabla I\big)  = {\rm div}\big( \big(\mathcal{A}_1-\mathcal{A}_2\big) \nabla I_2 \big)\hspace{1.0cm}\text{in}~~\Omega_T\,, \label{uni:maina}\\
&  I(x,0)= 0\,,~ I_t(x,0)=0\, \hspace{4.0cm} {\rm in}~~~\Omega\,, \label{uni:mainb} \\
& \langle\mathcal A_i  \nabla I_i, {\bf n}\rangle = 0  \hspace{5.4cm} {\rm on}~~~\partial \Omega_T\,, \label{uni:mainc}
\end{align}
where $\mathcal{A}_i=\mathcal{A}_i(x,t)= \mathcal A\left(x,t,{I_i},\nabla {I_i}, {I_i}(\cdot),\nabla {I_i}(\cdot)\right)$ for $i=1,2$, and $I=I_1-I_2$. It suffices to show that $ I \equiv 0$. To verify this, fix $ 0 < s < T$ and set, for $i=1,2$, 
\begin{align}\label{relationvi}
v_{i}(\cdot,t)= \begin{cases}
\displaystyle \int_{t}^{s} I_{i}(\cdot, \tau)d\tau, \hspace{0.5cm} 0<t\leq s\,, \\ 
 0, \hspace{2.3cm} s \leq t < T\,.
\end{cases}
\end{align}
Note that, for $t\in (0,s)$,
\begin{align}\label{eq:fact-1}
 \begin{cases}
 \partial_t v_i(x,t)=-I_i(x,t) \quad i=1,2\,, \\
  v_{i}(\cdot,t) \in H^1\,,~~~ \left(\mathcal A_i  \nabla v_i, {\bf n}\right) = 0~~~\text{on}\,\, \partial \Omega\,\,\text{in the sense of distribution}.
 \end{cases}
 \end{align}
Set $v=v_1-v_2$. Then $v(\cdot,s)=0$. Multiplying \eqref{uni:maina} by $v$, integrating over $\Omega \times (0,s)$ and using \eqref{uni:mainb} and \eqref{uni:mainc}, we obtain
\begin{align*}
&\int_{0}^{s}\int_{\Omega} \big( - \partial_t I \partial_t v - I{\partial_t v}+\langle\mathcal{A}_1 \nabla I, \nabla v\rangle \big)dxdt
= -\int_{0}^{s}\int_{\Omega} \langle \left(\mathcal{A}_1-\mathcal{A}_2\right) \nabla I_{2}, \nabla v\rangle \,dx\,dt.
\end{align*}
Now using \eqref{eq:fact-1} and Cauchy-Schwarz inequality in the above equality, we have
\begin{align*}
& \frac{1}{2}\int_{0}^{s}\int_{\Omega} {\partial_t} (I^2)\, dxdt+\int_{0}^{s}\int_{\Omega} I^2 \,dxdt-\int_{0}^{s}\int_{\Omega}\langle\mathcal{A}_1{\partial_t}{\nabla v}, \nabla v\rangle\,dx\,dt \\
&\hspace{5cm} \le \int_0^s \left(\|(\mathcal{A}_1-\mathcal{A}_2)(t)\|\right)_{L^\infty} \|\nabla I_2(t)\|_{L^2}\|\nabla v(t)\|_{L^2}\,dt.
\end{align*}
Using the fact that
\begin{align*}
 \langle\mathcal{A}_1\partial_t \nabla v, \nabla v\rangle = \frac{1}{2} \partial_t\langle \mathcal{A}_1 \nabla v, \nabla v \rangle - \frac{1}{2} \langle \partial_t \mathcal{A}_1 \nabla v, \nabla v \rangle, \quad
 \nabla v(x,s)=0\,,
\end{align*}
and \eqref{uni:mainb}, we obtain
\begin{align}\label{unique9}
&\frac{1}{2}\|I(s)\|_{L^2}^2+\int_{0}^{s}\|I(t)\|_{L^2}^2\,dt  + \frac{1}{2}\int_{\Omega} \langle\mathcal{A}_1(x,0) \nabla v(x,0), \nabla v(x,0)\rangle\, dx \nonumber  \\
&\leq \frac{1}{2} \Big|\int_{0}^{s}\int_{\Omega} |\nabla v|^2 \left(\|\partial_t \mathcal{A}_1\|\right)_{L^{\infty}}\, dx\,dt\Big| + \int_0^s \left(\|(\mathcal{A}_1-\mathcal{A}_2)(t)\|\right)_{L^\infty} \|\nabla I_2(t)\|_{L^2}\|\nabla v(t)\|_{L^2}\,dt\,.
\end{align}
Now, assumptions \eqref{eq:ass-A1}, \eqref{eq:ass-A2} and 
\eqref{eq:ass-A3} imply the existence of a constant $C$ such that
%
\begin{align*}
\frac{1}{2}\|I(s)\|_{L^2}^2+\int_{0}^{s}\| I(t)\|_{L^2}^2\,dt  + C^{-1} \|\nabla v(0)\|_{L^2}^2 
&\le   C\, \int_0^s \big( \|\nabla v(t)\|_{L^2}^2 + \|I(t)\|_{L^2}^2\big)\,dt\,,
\end{align*}
%
and since $ \|v(0)\|_{L^2}^2 \le T \int_0^s \| I(t)\|_{L^2}^2\,dt $, we also have 
\begin{align}\label{unique15_1}
\frac{1}{2}\|I(s)\|_{L^2}^2+\int_{0}^{s}\| I(t)\|_{L^2}^2\,dt  + C^{-1} \|v(0)\|_{H^1}^2 
\le C\,\int_0^s \big( \|v(t)\|_{H^1}^2 + \|I(t)\|_{L^2}^{2} \big)\,dt\,.
\end{align}
Set
\begin{align*}
w_{i}(\cdot,t)=& \int_{0}^{t} I_{i}(\cdot,\tau)d\tau\, ; \quad w(\cdot,t)=(w_1-w_2)(\cdot,t)\,, \hspace{0.5cm} 0<t\leq T,
\end{align*}
and observe that then 
\begin{align*}
 v(x,0)= w(x,s) \quad {\rm and}~~
 v(x,t)= w(x,s)-w(x,t)~~{\rm for}~~ 0<t\le s\,.
\end{align*}
Hence, \eqref{unique15_1} reduces to
\begin{align}\label{unique16}
&\frac{1}{2}\|I(s)\|_{L^2}^2+\int_{0}^{s}\| I(t)\|_{L^2}^2\,dt  + C \|w(s)\|_{H^1}^2
 \le \tilde{C} s\,\|w(s)\|_{H^1}^2 + C\,\int_0^s \Big( \|w(t)\|_{H^1}^2 + \|I(t)\|_{L^2}^{2} \Big)\,dt\,.
\end{align}
Choose $T_1$ sufficiently small such that  $C^{-1}-\tilde{C} T_1 >0$. 
Then, for $0<s\leq T_1,$ we have, from \eqref{unique16},
\begin{align}
\| I(s)\|_{L^2}^2 + \|w(s)\|_{H^1}^2 \le C \int_0^s\Big( \|w(t)\|_{H^1}^2 + \|I(t)\|_{L^2}^{2}\Big)\,dt\,. \label{unique17}
\end{align}
Consequently, an application of Gronwall's lemma implies $ I \equiv 0$ on $[0,T_1]$.
Finally, we utilize a similar logic on the intervals $(T_1, 2T_1]$, $(2T_1,3T_1],\ldots$ step by step, and eventually deduce that $I_{1} = I_{2}$ on $(0,T)$. This finishes the proof of the uniqueness part, and hence also the proof of Theorem \ref{thm:existence}.
$\hfill\Box$
\begin{cor}\label{corr:hyp-final}
Problem \eqref{maina}--\eqref{mainb} with diffusion function \eqref{diff_edgeg} satisfying \eqref{eq:ass-exponent_p} and with $0 < \alpha_1 \leq I_0 \leq \beta_1$ admits a weak solution in $\mathcal{W}$ if $\|J\|_{H^1(\Omega)}$ is small enough, 
or if $T$ is small enough. 
The solution is unique among functions in $\mathcal{W}$.
\end{cor}

To prove Corollary \ref{corr:hyp-final}, we need the following simple Lemma:

\begin{lem}\label{le:bounded_2024}
Let I be a weak solution of the telegraph diffusion model \eqref{maina}--\eqref{mainb}, or the diffusion model \eqref{diff_power_constant}, with diffusion function given by \eqref{diff_edgeg}.  
Then $\alpha_1 \leq I(x,t) \leq \beta_1$ for a.e. $(x,t) \in  \Omega \times [0,T]$, where
$\alpha_1 = \inf_{x \in \Omega} I_0(x)$ and $\beta_1 = \sup_{x \in \Omega} I_0(x)$.
\end{lem}

\noindent
{\bf Proof.}
We proceed as in \cite[Lemma 3.3]{majee2020gray}.
Integrating the equation \eqref{maina} w.r.t. time variable and using \eqref{mainc}, we get that
\begin{align}\label{eq:bounded_2024}
I_t + \gamma (I - I_0) - \int_{0}^t \text{div} \left(g (I_\xi, |\nabla I_\xi |)    \nabla I\right) ds = 0 \quad \text{for all}\quad (x,t) \in \Omega_T.
\end{align}
Note that $(I - \beta_1)_+ \in H^1$, in which $(\theta)_+ = \max\{0, \theta\}$. 
Multiplying  \eqref{eq:bounded_2024} with $(I - \beta_1)_+$ and integrating over $\Omega$ we obtain
\begin{align}\label{eq:2024_tjohej1}
\frac{1}{2} \frac{d}{dt} \int_{\Omega} \left((I -\beta_1)_+\right)^2 dx +
\gamma \int_{\Omega} (I -\beta_1)_+ (I - I_0) dx + 
\int_{0}^t \int_{I \geq \beta_1} g (I_\xi, |\nabla I_\xi |)    |\nabla I|^2 dx ds = 0.
\end{align}
Since $g \geq 0$ and $(I-J)(I-\beta_1)_+ \geq 0$ it must hold that
$\frac{1}{2} \frac{d}{dt} \int_{\Omega} \left((I -\beta_1)_+\right)^2 dx \leq 0$.
Therefore, since $J \leq \beta_1$ we obtain
$\int_{\Omega}  |(I - \beta_1)_+|^2 dx \leq 0$, for a.e. $t \in [0, T]$, 
so that $I(x,t) \leq \beta_1$ for a.e. $(x,t) \in 
\Omega_T$.

Similarly, multiplying \eqref{eq:bounded_2024} with $(I - \alpha_1)_- \in H^1$, where $(\theta)_- = \min\{0, \theta\}$, and integrating over $\Omega$ we have
\begin{align}\label{eq:2024_tjohej2}
\frac{1}{2} \frac{d}{dt} \int_{\Omega} \left((I -\alpha_1)_-\right)^2 dx +
\gamma \int_{\Omega} (I -\alpha_1)_- (I - I_0) dx + 
\int_{0}^t \int_{I \leq \alpha_1} g (I_\xi, |\nabla I_\xi |)    |\nabla I|^2 dx ds = 0.
\end{align}
Since $g \geq 0$ and $(I-J)(I-\alpha_1)_- \geq 0$ it must hold that
$\frac{1}{2} \frac{d}{dt} \int_{\Omega} \left((I -\alpha_1)_-\right)^2 dx \leq 0$, and since $\alpha_1 \leq I_0$, we also conclude that  
$0 < \alpha_1 \leq I(x,t)$ for a.e.$(x,t) \in \Omega_T$. 
This completes the proof of Lemma \ref{le:bounded_2024} in the case of the telegraph diffusion model.

In the case of the diffusion model \eqref{diff_power_constant}
we also obtain the result from \eqref{eq:2024_tjohej1} and \eqref{eq:2024_tjohej2} but now without the first term and with $\gamma = 1$.
In particular, it follows that 
$$
\int_{\Omega} (I -\beta_1)_+ (I - I_0) dx =0 = \int_{\Omega} (I -\alpha_1)_- (I - I_0) dx
$$
from which we conclude the desired result.   
The proof of Lemma \ref{le:bounded_2024} is complete. $\hfill\Box$

\bigskip

\noindent
{\bf Proof of Corollary \ref{corr:hyp-final}.}
It suffices to show that Theorem \ref{thm:existence} applies to our suggested model \eqref{maina}--\eqref{mainb}, i.e. to show that
the function in $\mathcal{A}_w$ defined in \eqref{eq:def_A} satisfies assumptions \eqref{eq:ass-A1}--\eqref{eq:ass-A3}.
We have, for any $w \in \mathcal{W}_{M_1}$,
\begin{align}\label{eq:AW}
\mathcal{A}_w\left(x,t\right) =
\epsilon
+ \frac{ 2|  w_\xi|^\nu}{\left(M^{w}_{\xi}\right)^\nu + | w_\xi|^\nu}  \cdot \dfrac{1}{1+ \left(\frac{\vert \nabla w_\xi \vert }{K} \right)^{p\left(w_\xi, \vert \nabla w_\xi \vert \right)} }\,.
\end{align}
We first observe that
\begin{align}\label{eq:A1-particular}
\epsilon \leq \mathcal{A}_w\left(x,t\right) \leq 1
\end{align}
and hence $A_w$ satisfies \eqref{eq:ass-A1} with $\phi_\mathcal A = \epsilon$ and $\Phi_\mathcal A = 1$.

To show that $\mathcal A_w$ satisfies \eqref{eq:ass-A2} we first note that
\begin{align*}
\left|\frac{d}{d t} \mathcal A_w(x,t)\right|
\leq \left| \frac{d}{d t} \left( \dfrac{ 2|  w_\xi|^\nu}{\left( M^{w}_{\xi}\right)^\nu + |  w_\xi|^\nu}\right) \right|
+   \left| \frac{d}{d t} \left( \frac{1}{1 + \left( \frac{|\nabla w_\xi|}{K} \right)^{p(w_\xi, |\nabla w_\xi|)} } \right) \right|
= A + B.
\end{align*}
Next, observe that by properties of convolution, we have
$$
|w_\xi| + |G_\xi * w_t| + |\nabla w_\xi| + |\nabla G_\xi * w_t| \leq C_\xi M_1 ||I_0||_{H^1} \qquad \text{in $\Omega \times (0,T]$}
$$
where we denote $w_{\xi,t} = G_\xi * w_t$ and $\nabla w_{\xi,t} = \nabla G_\xi * w_t$.
Since $w_\xi, M_\xi^{w} > \alpha_1$ we have
\begin{align*}
A &=2 \frac{\nu w_\xi^{\nu-1} (M_\xi^{w})^{\nu-1}(M_\xi^{w}  w'_\xi -  w_\xi(M_\xi^{w})' }{(w_\xi^\nu + (M_\xi^{w})^\nu)^2}
\end{align*}
and thus
\begin{align*}
A &\leq 2
\frac{\nu w_\xi^{\nu-1} (M_\xi^{w})^{\nu-1} 2 \left( M_1 C_\xi \|I_0\|_{H^1} \right)^2 }{(\alpha_1^\nu + \alpha_1^\nu)^2} \\
&\leq 4 \nu \left( \alpha_1^{\nu-1} + \left( M_1 C_\xi \|I_0\|_{H^1} \right)^{\nu-1}\right)^2 \frac{\left( M_1 C_\xi \|I_0\|_{H^1} \right)^2}{4 \alpha_1^{2\nu}} \\
&\leq  \nu \left( 2 \alpha_1^{2(\nu-1)} + 2\left( M_1 C_\xi \|I_0\|_{H^1} \right)^{2(\nu-1)}\right) \frac{\left( M_1 C_\xi \|I_0\|_{H^1} \right)^2}{\alpha_1^{2\nu}} \\
&\leq  2\nu \left( \frac{\left( M_1 C_\xi \|I_0\|_{H^1} \right)^2}{\alpha_1} + \frac{\left( M_1 C_\xi \|I_0\|_{H^1} \right)^{2\nu}}{\alpha_1^{2\nu}}\right) \\
&\leq  C \left( M_1 \|I_0\|_{H^1} \right)^{2(\nu + 1)}
\end{align*}
where $C$ depends only on $\nu, \xi$ and $\alpha_1$.
To bound the second term, first note that
\begin{align*}
B &= \left| \frac{d}{dt} \left( \frac{1}{1 + \left( \frac{|\nabla w_\xi|}{K} \right)^{p(w_\xi, |\nabla w_\xi|)} } \right) \right|
\leq \left| \frac{d}{dt} \left( \frac{|\nabla w_\xi|}{K}\right)^{p(w_\xi, |\nabla w_\xi|)} \right|\\
= &  \left(  \frac{|\nabla w_\xi|}{K}\right)^{p(w_\xi, |\nabla w_\xi|)}  \left| \frac{d}{dt} p(w_\xi, |\nabla w_\xi|) \log\left(\frac{|\nabla w_\xi|}{K}\right)  +  p(w_\xi, |\nabla w_\xi|)  \frac{d}{dt} \log\left(\frac{|\nabla w_\xi|}{K}\right) \right|.
\end{align*}
Using that $x \log x \to 0$ as $x \to 0$ with $p \geq 1$ gives
\begin{align*}
B&\leq \left(C_\xi M_1 K^{-1} \|I_0\|_{H^1}\right)^{p^+}
L \left(|G_\xi * w_t| + |\nabla G_\xi * w_t|\right) C_\xi M_1 K^{-1} \|I_0\|_{H^1}\\
&+ \left(  \frac{|\nabla w_\xi|}{K}\right)^{p(w_\xi, |\nabla w_\xi|)} p^+ \frac{K |\nabla G_\xi * w_t|}{|\nabla w_\xi|} \\
&\leq
2 L\left(C_\xi M_1 K^{-1} \|I_0\|_{H^1}\right)^{p^+ + 2}
+ \left(  \frac{|\nabla w_\xi|}{K}\right)^{p(w_\xi, |\nabla w_\xi|) - 1} p^+ C_\xi M_1 K^{-1} \|I_0\|_{H^1}\\
&\leq
2 L\left(C_\xi M_1 K^{-1} \|I_0\|_{H^1}\right)^{p^+ + 2} +
p^+ \left( C_\xi M_1 K^{-1} \|I_0\|_{H^1} \right)^{p^+}\\
&\leq
C \left( M_1  \|I_0\|_{H^1}\right)^{p^+ + 2}
\end{align*}
where $C$ is independent of $M_1$ and $I_0$.
Thus, by summing $A + B$ we realize, by the above displays, that
\begin{align}\label{eq:A2-particular}
\left|\frac{d}{d t} \mathcal A_w(x,t)\right| \leq C \left( M_1  \|I_0\|_{H^1}\right)^{C}
\end{align}
where $C$ depends only on $\alpha_1, \xi, \nu$ and $K$.
This proves that $\mathcal{A}_w$ satisfies \eqref{eq:ass-A2} with $\Psi(z) = C z^C$.

To show that $\mathcal{A}$ satisfies \eqref{eq:ass-A3}, we may proceed as in the proof of Corollary \ref{corr:diffusion-final}. 
Thus, to complete the proof of Corollary \ref{corr:hyp-final}, it only remains to verify that at least one of the three statements $(i)$, $(ii)$ or $(iii)$ in Theorem \ref{thm:existence} holds.
We have
$$
\frac{\Phi_\mathcal A}{\phi_\mathcal A} \leq \frac{1}{\epsilon} \qquad \text{and so} \qquad \frac{1}{z^2} \frac{\Phi_\mathcal A}{\phi_\mathcal A} \to 0 \qquad \text{as} \qquad z \to \infty, 
$$
and hence \eqref{eq:ass-functions-1} holds.
On the other hand,
$$
\frac{\Psi_\mathcal A}{\phi_\mathcal A} \leq \frac{C z^C}{\epsilon} \to \infty \qquad \text{as}\qquad z \to \infty,
$$
and so \eqref{eq:ass-functions-2} does not hold. 
In conclusion, existence follows if $\|I_0\|_{H^1 }$ is small enough,
or if $T$ is small enough, depending only on $\mathcal A$ and $\|I_0\|_{H^1 }$.
The proof is complete.
$\hfill\Box$\\
%


\section{Numerical Method}
\label{sec:numerical}
To solve the system \eqref{maina}--\eqref{mainc} with \eqref{diff_edgeg} numerically, we construct a weighted-$\theta$ finite difference scheme \cite{jovanovic2013analysis}. We replace the derivative terms in \eqref{maina}--\eqref{mainc} using the following finite difference approximations:
\begin{align*}
& \dfrac{\partial I_{i,j}^n }{\partial t} \approx  \displaystyle\frac{I_{i,j}^{n+1}-I_{i,j}^{n-1}}{2 \tau}\,, \quad 
\dfrac{\partial^2 I_{i,j}^n }{\partial t^2} \approx \displaystyle\frac{I_{i,j}^{n+1}-2I_{i,j}^n+I_{i,j}^{n-1}}{\tau ^2}\,,\\
& \nabla_x I_{i,j}^n \approx \displaystyle\frac{I_{i+1,j}^n-I_{i-1,j}^n}{2h}\,, \quad 
\nabla_y I_{i,j}^n \approx \displaystyle\frac{I_{i,j+1}^n-I_{i,j-1}^n}{2h}\,, \\
& \Delta_x I_{i,j}^n \approx  \frac{{I_{i+1,j}^n-2I_{i,j}^n+I_{i-1,j}^n}}{h^2}\,,  \quad 
 \Delta_y I_{i,j}^n \approx \frac{{I_{i,j+1}^n-2I_{i,j}^n+I_{i,j-1}^n}}{h^2}\,, \\
& |\nabla I_{i,j}^n| \approx \sqrt{(\nabla_x I_{i,j}^n)^2 + (\nabla_y I_{i,j}^n)^2}\,.
\end{align*}
In the above, $\tau$  and $h$ denote the time step size and the spatial step size, respectively. $I^n_{i,j}=I(t_n,x_i,y_j),$ where $x_i=ih~(i=0,1,2,...,M-1)$, 
$y_j=jh~(j=0,1,2,...,N-1)$, $t_n=n\tau~(n=0,1,2,\ldots),$ where $n$ is the number of iterations and $M \times N$ is the image dimension. Then using the weighted scheme \cite{jovanovic2013analysis}, the discrete form of the equation \eqref{maina} could be written as 
\begin{align}\label{disc:I}
&(1+ 0.5\gamma \tau) I^{n+1}_{i,j}- \tau^2 \theta_1 \left[\nabla(g \nabla I) \right]^{n+1}_{i,j}
= 2I^{n}_{i,j} + \tau^2 (1-\theta_1-\theta_2) \left[\nabla(g \nabla I) \right]^n_{i,j}\notag  \\
&\hspace{8cm}  + \tau^2 \theta_2  \left[\nabla(g \nabla I) \right]^{n-1}_{i,j} +(0.5 \gamma \tau -1)I^{n-1}_{i,j}\,,
\end{align}
where $\theta_1$ and $\theta_2$ are non negative `weights'. The superscript `$n$' denotes the value at the $n^{th}$ time level $t_n$ and 
\begin{align}\label{rhs_discretization_short}
\left[\nabla(g \nabla I) \right]_{i,j}&=\dfrac{0.5}{\tilde{h}^2}\Big[ ( g_{i,j} + g_{i+1,j}  )  I_{i+1, j}+ ( g_{i,j} + g_{i-1,j}  ) I_{i-1, j} 
    - \left( g_{i+1,j}+2g_{i,j} + g_{i-1,j} \right) I_{i, j}   \Big] \notag  \\
&         + \dfrac{0.5}{\tilde{h}^2}\Big[ (g_{i,j} + g_{i,j+1}) I_{i, j+1} + (g_{i,j} + g_{i,j-1}) I_{i, j-1} 
     - \left( g_{i,j+1}+2g_{i,j} + g_{i,j-1} \right) I_{i, j}   \Big],
\end{align}
where
\begin{align*}
g_{i,j}^n=\epsilon
+b(s_{i,j}) \cdot \dfrac{1}{1+ \left(\frac{|\nabla G_{\xi} \ast I_{i,j} |}{K} \right)^{p\left(G_\xi \ast I_{i,j}, \vert \nabla G_\xi  \ast I_{i,j} \vert \right)} }.
\end{align*}
The initial and boundary conditions are given as follows:
\begin{align*}
& I_{i,j}^0 =I_0(x_i,y_j)\,,I_{i,j}^1=I_{i,j}^0,\,0 \leq i \leq M-1\,, 0 \leq j \leq N-1\,,\\
& I_{0,j}^n=I_{1,j}^n,~I_{M-1,j}^n=I_{M-2,j}^n\,,~I_{i,0}^n=I_{i,1}^n\,,~I_{i,N-1}^n=I_{i,N-2}^n\,.
\end{align*}
We solve the system \eqref{disc:I} using the Gauss-Seidel iterative method \cite{hoffman2018numerical}. The numerical stability condition for the model \eqref{maina}--\eqref{mainc} is $\tau \leq h / \sqrt{\text{max}\, g (x, t) }$  according to the Courant–Friedrichs–Lewy stability criterion \cite{zauderer2011partial}, where $h$ denotes the length of the spatial intervals. We choose a uniform time step size $\tau = 0.25$, spatial step size $h=1$, and $\xi=1$ for all our computations. Moreover, we require a stopping criterion to stop the numerical simulation process. We use two distinct stopping rules: when the clear image is available and another for the real-life images. When the clean data is available (i.e., in the case of artificially noised images), we stop the simulation after getting the best possible peak signal-to-noise ratio (PSNR) \cite{gonzalez2002digital}  value between the clean image ($I$) and the restored image ($I^k$) calculated by the formula
\begin{align}
\text{PSNR} = 10\, \text{log}_{10} \left(\frac{\text{max}(I)^2}{\frac{1}{\text{MN}} \sum\limits_{i=1}^\text{M} \sum\limits_{j=1}^\text{N} (I(i,j)-I^{k}(i,j))^2  }\right),
\end{align}
where $I$ denotes the ground truth image of size $\text{M}\times \text{N}$ and max($I$) is the maximum possible pixel value of $I,$ and $I^k$ denotes the restored image at the $k^{th}$ iteration. For the case of real images, when the clean image is not available, we stop the simulation process after satisfying the following relation
\begin{equation}\label{eq:4stopping}
\frac{{\Vert I^{k+1}-I^k\Vert^2_2}}{{\Vert I^k \Vert^2_2}}\leq \varepsilon,
\end{equation} 
where $\varepsilon > 0$ is a fixed threshold. In \eqref{eq:4stopping}, $I^k$ and $I^{k+1}$ denote the restored images at the $k^{th}$ and ${(k+1)}^{th}$ iteration, respectively. For our simulations, we have used $\varepsilon = 10^{-4}$.

\section{Computational Results}
\label{computation_discussion}


This section discussed the image despeckling ability of the presented model in terms of visual quality and quantitative measures. We discuss the results in three different subsections; \ref{sec:artificial}, \ref{sec:real}, and \ref{sec:benchmark}. The first subsection \ref{sec:artificial} describes the despeckling results using three standard test images degraded by the multiplicative speckle noise with different noise levels. We artificially generate speckle noise with ``looks" $ L = \{ 1,3,5,10 \} $ using MATLAB built-in gamma random noise generator function ``gamrnd" as $\eta=\text{gamrnd}(L,1/L,\text{M, N})$ and then multiply with the clean images,  where  $\text{M}\times \text{N}$ is the image dimensions. In subsection \ref{sec:real}, we describe the results for two real SAR images. Finally, in subsection \ref{sec:benchmark}, we discuss the despeckling results for standard benchmark images introduced in \cite{di2013benchmarking}. We compute the results using the proposed model and compare the results with the results of a nonlocal speckle removal approach. All the numerical tests are performed under Windows $10$ and MATLAB version $R2021b$ running on a desktop computer with an Intel Core $i5$ dual-core CPU at $2.11$ GHz with $16$ GB of RAM.

\subsection{Results for the Artificially Noised Images}
\label{sec:artificial}
To check the image despeckling ability of the present model, we compute three performance measures, i.e., the PSNR, Mean Structural Similarity Index (MSSIM), and despeckling gain measure (DG) \cite{argenti2013tutorial,di2013benchmarking}. Higher PSNR, MSSIM, and DG values suggest that the despeckled result is closer to the corresponding clean image. Also, we report two bias indicators; mean of the ratio image (MoR) and variance of the ratio image (VoR) (ratio image is defined as the point-by-point ratio between the degraded and the despeckled image) \cite{argenti2013tutorial,di2013benchmarking}. Since the ratio image contains only speckles, it has a unit mean and variance equal to $1/L$. Considering MoR $\cong$ 1, VoR indicates insight about under/over smoothing phenomena. A VoR $<$ 1 shows under smoothing; that is, part of the speckle remains in the filtered image, whereas VoR $>$ 1 means over smoothing; that is, the filter also eliminates some details of the underlying image.

\begin{figure}
       \centering
       \begin{subfigure}[b]{0.24\textwidth}           
                \includegraphics[scale=0.36]{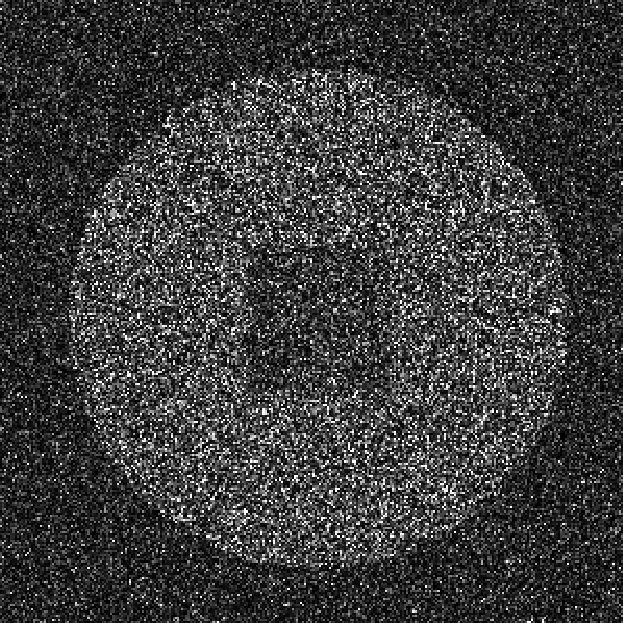}               
                \caption{Noisy $L=1$}
                \label{circle_noisy_l1}
       \end{subfigure}
       \begin{subfigure}[b]{0.24\textwidth}           
                \includegraphics[scale=0.36]{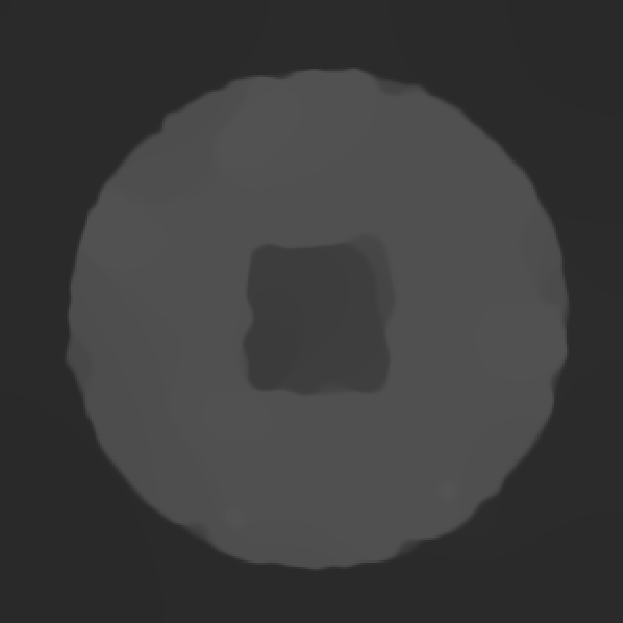}               
                \caption{$p=p_0-p_1$}
                \label{circlel1_avggray_tdmp}
       \end{subfigure}
            \begin{subfigure}[b]{0.24\textwidth}           
                \includegraphics[scale=0.36]{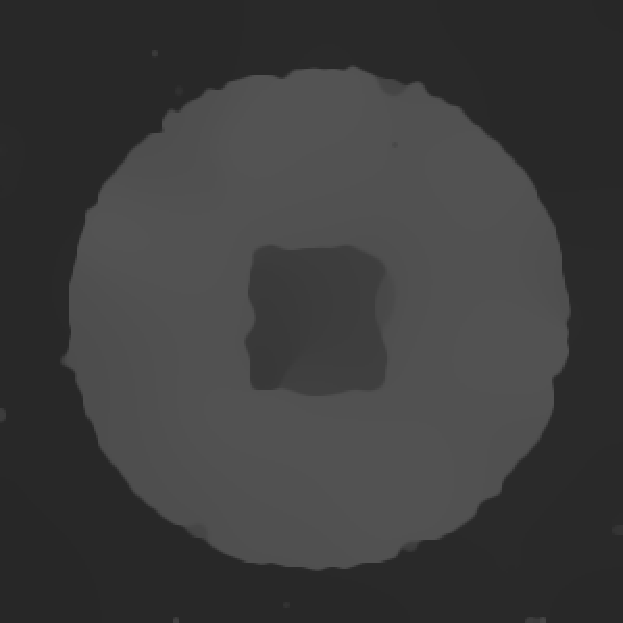}               
                \caption{$p=p_0-p_2$}
                \label{circlel1_zzdb_tdmp}
       \end{subfigure}
       \begin{subfigure}[b]{0.24\textwidth}           
                \includegraphics[scale=0.36]{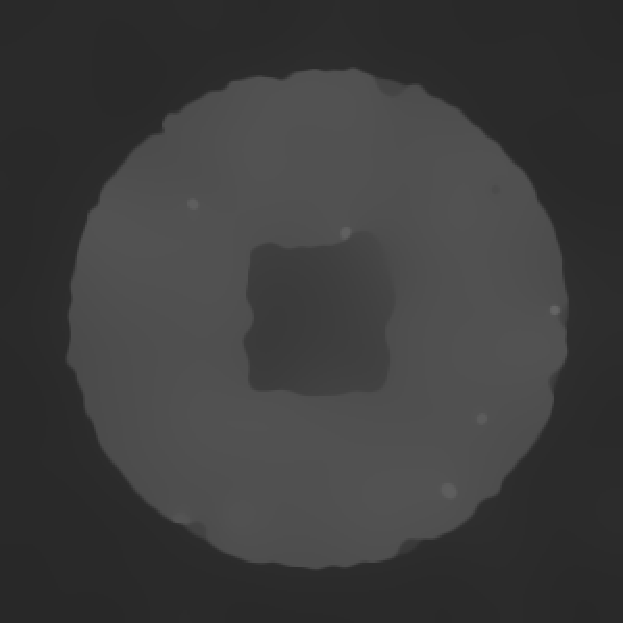}               
                \caption{$p=p_0-p_3$}
                \label{circlel1_grad_tdmp}
       \end{subfigure}
       
      \begin{subfigure}[b]{0.24\textwidth}           
                \includegraphics[scale=0.36]{circle_speckle_look_10}               
                \caption{Noisy $L=10$}
                \label{circle_noisy_l10}
       \end{subfigure}
       \begin{subfigure}[b]{0.24\textwidth}           
                \includegraphics[scale=0.36]{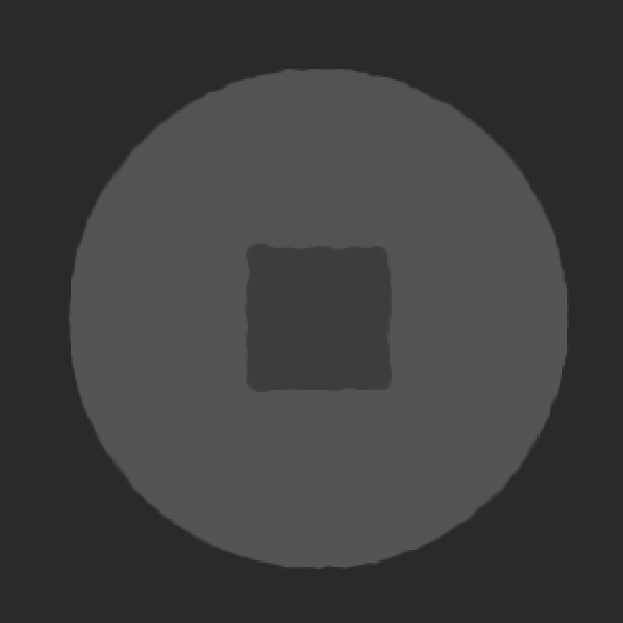}               
                \caption{$p=p_0-p_1$}
                \label{circlel10_avggray_tdmp}
       \end{subfigure}
            \begin{subfigure}[b]{0.24\textwidth}           
                \includegraphics[scale=0.36]{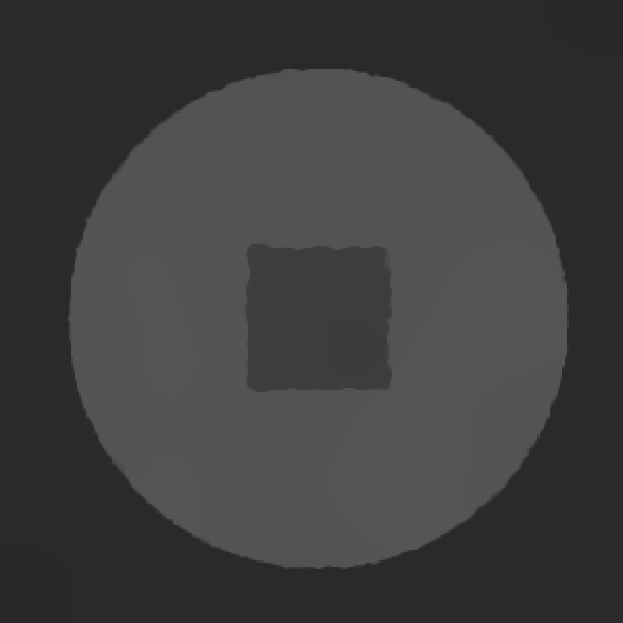}               
                \caption{$p=p_0-p_2$}
                \label{circlel10_zzdb_tdmp}
       \end{subfigure}
       \begin{subfigure}[b]{0.24\textwidth}           
                \includegraphics[scale=0.36]{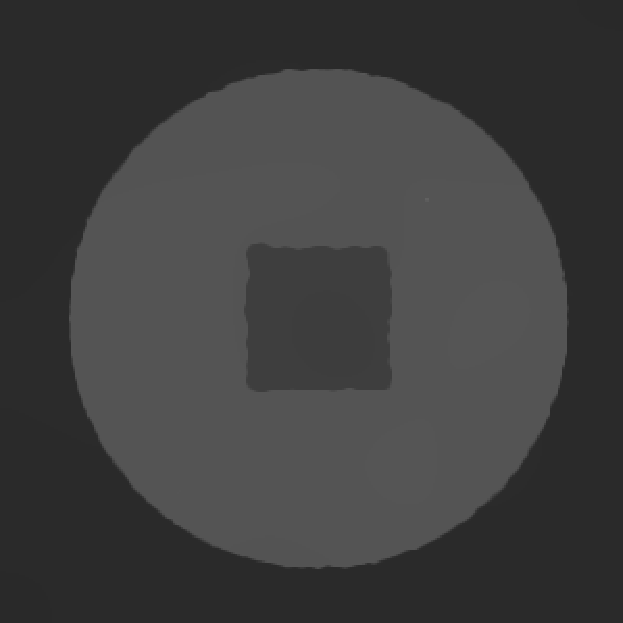}               
                \caption{$p=p_0-p_3$}
                \label{circlel10_grad_tdmp}
       \end{subfigure}
       
  \hspace*{-1cm}    \begin{subfigure}[b]{0.25\textwidth}           
                \includegraphics[scale=0.335]{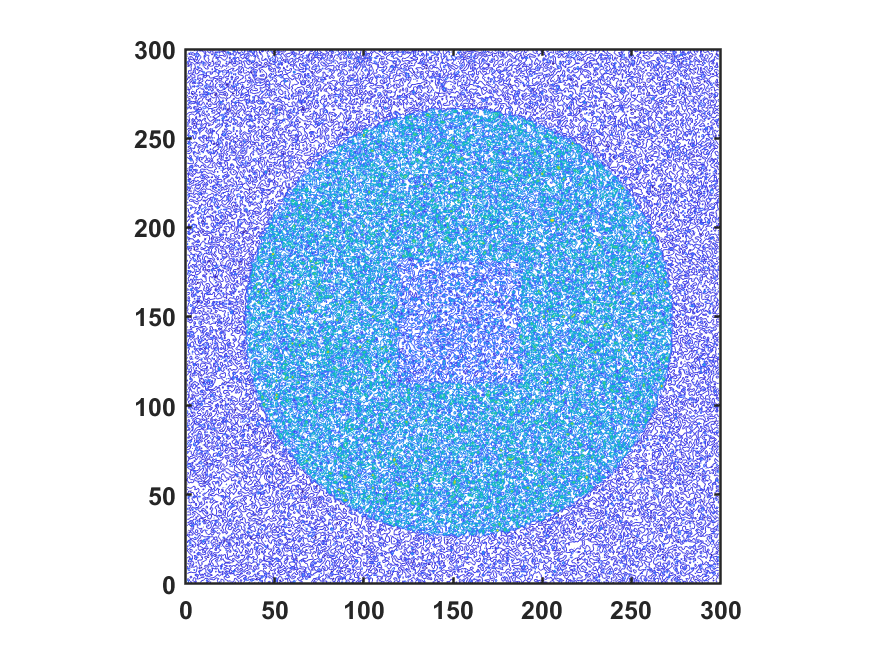}               
                \caption{Noisy $L=10$}
                \label{circle_noisy_l13d}
       \end{subfigure}
       \begin{subfigure}[b]{0.25\textwidth}           
                \includegraphics[scale=0.335]{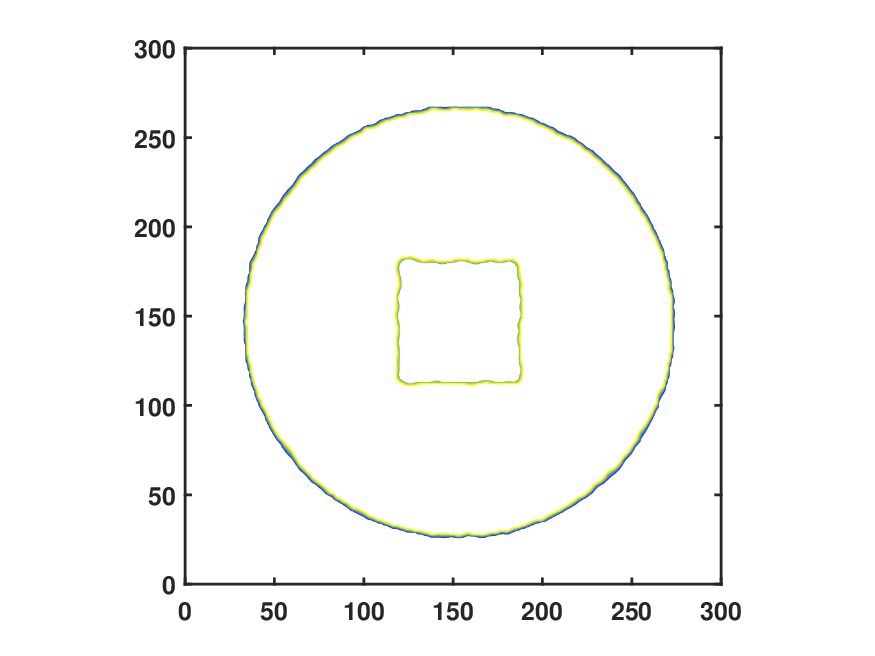}               
              \caption{$p=p_0-p_1$}
                \label{circlel1_avggray_tdmp3d}
       \end{subfigure}
            \begin{subfigure}[b]{0.25\textwidth}           
                \includegraphics[scale=0.335]{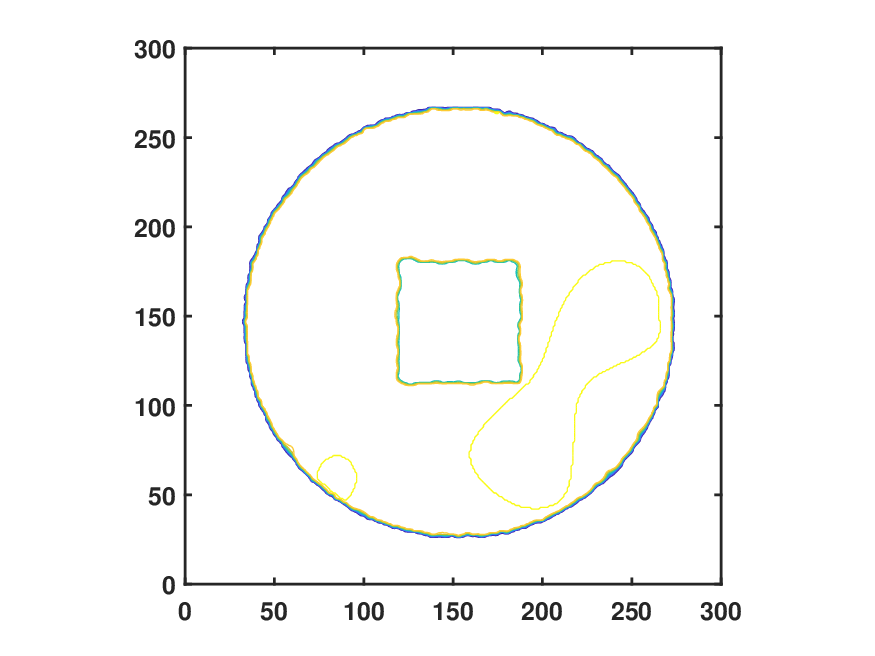}               
                \caption{$p=p_0-p_2$}
                \label{circlel1_zzdb_tdmp3d}
       \end{subfigure}
       \begin{subfigure}[b]{0.25\textwidth}           
                \includegraphics[scale=0.335]{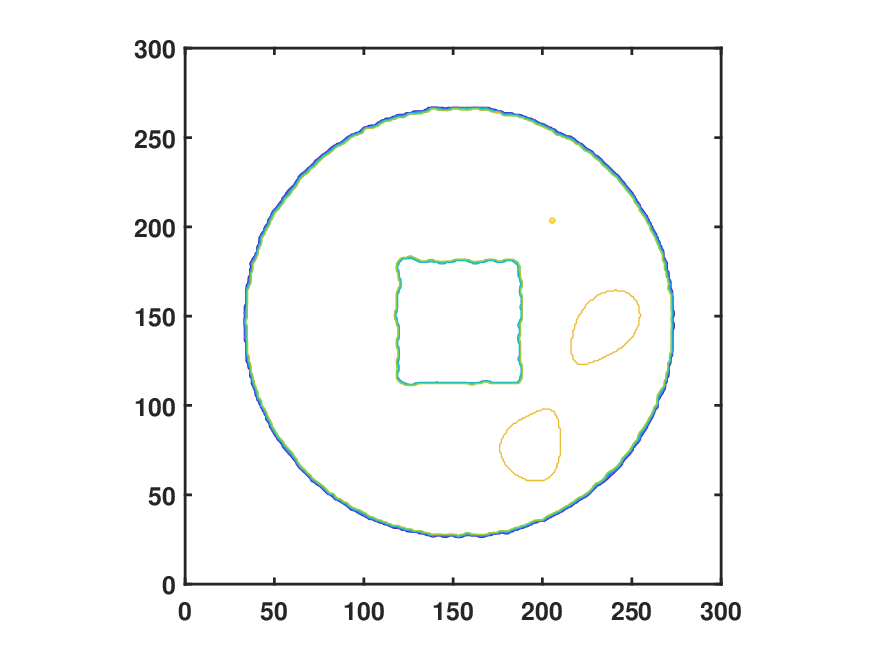}               
                \caption{$p=p_0-p_3$}
                \label{circlel1_grad_tdmp3d}
       \end{subfigure}
       
\caption{\footnotesize Restored Images using \eqref{maina}--\eqref{mainc}.}\label{circle_1_10_restored_tdmp}
\end{figure}

\begin{figure}
       \centering
       \begin{subfigure}[b]{0.24\textwidth}           
                \includegraphics[scale=0.36]{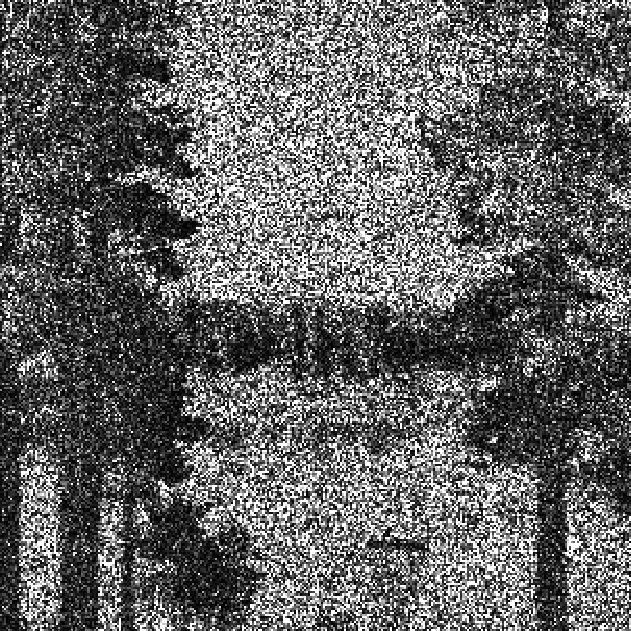}               
                \caption{Noisy $L=1$}
                \label{lake1_noisy_l1}
       \end{subfigure}
       \begin{subfigure}[b]{0.24\textwidth}           
                \includegraphics[scale=0.36]{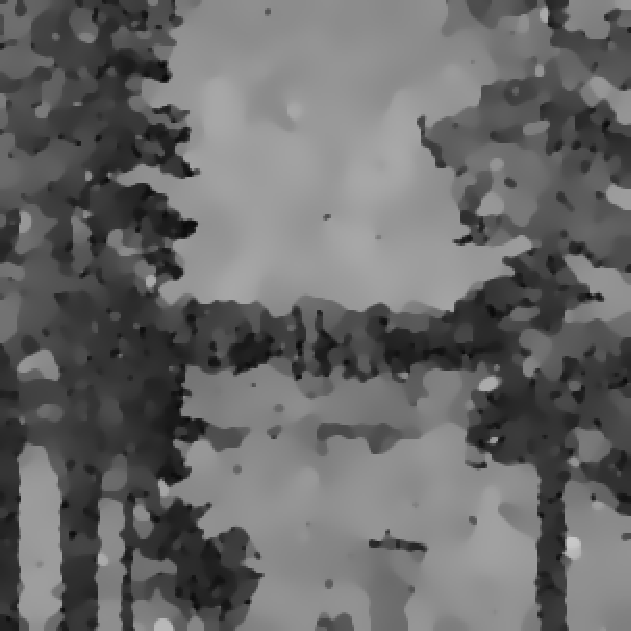}               
                \caption{$p=p_0-p_1$}
                \label{lake1l1_avggray_tdmp}
       \end{subfigure}
            \begin{subfigure}[b]{0.24\textwidth}           
                \includegraphics[scale=0.36]{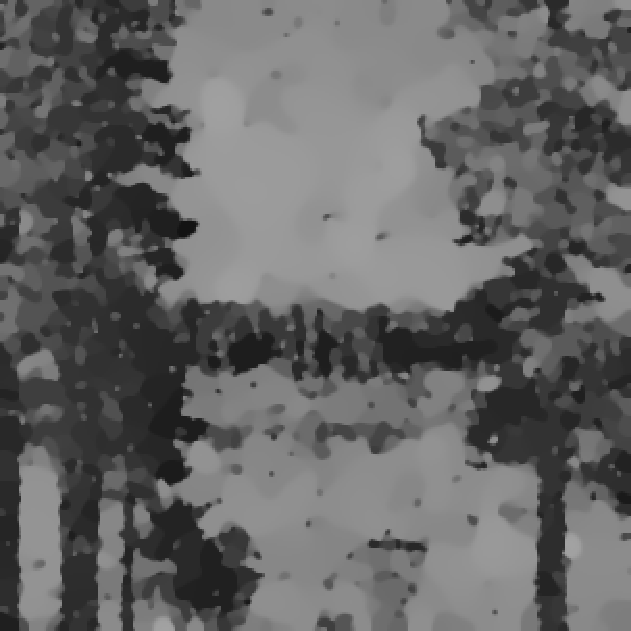}               
                \caption{$p=p_0-p_2$}
                \label{lake11_zzdb_tdmp}
       \end{subfigure}
       \begin{subfigure}[b]{0.24\textwidth}           
                \includegraphics[scale=0.36]{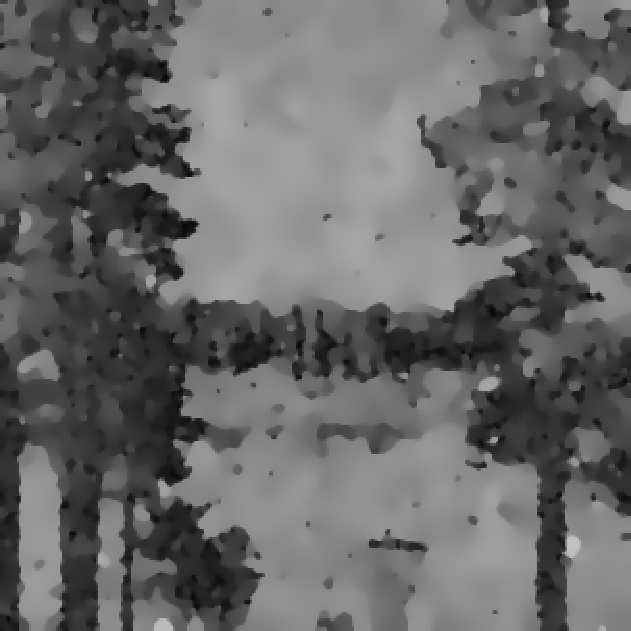}               
                \caption{$p=p_0-p_3$}
                \label{lake11_grad_tdmp}
       \end{subfigure}
       
      \begin{subfigure}[b]{0.24\textwidth}           
                \includegraphics[scale=0.36]{lake1_speckle_look_10}               
                \caption{Noisy $L=10$}
                \label{lake1_noisy_l10}
       \end{subfigure}
       \begin{subfigure}[b]{0.24\textwidth}           
                \includegraphics[scale=0.36]{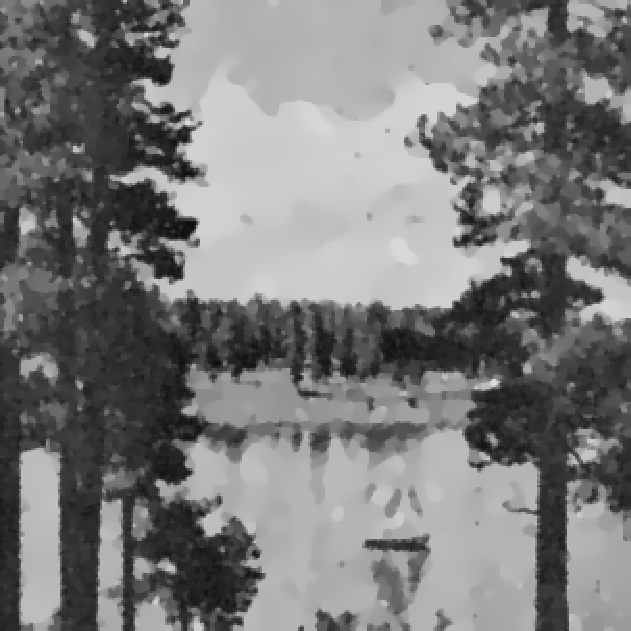}               
                \caption{$p=p_0-p_1$}
                \label{lake1l10_avggray_tdmp}
       \end{subfigure}
            \begin{subfigure}[b]{0.24\textwidth}           
                \includegraphics[scale=0.36]{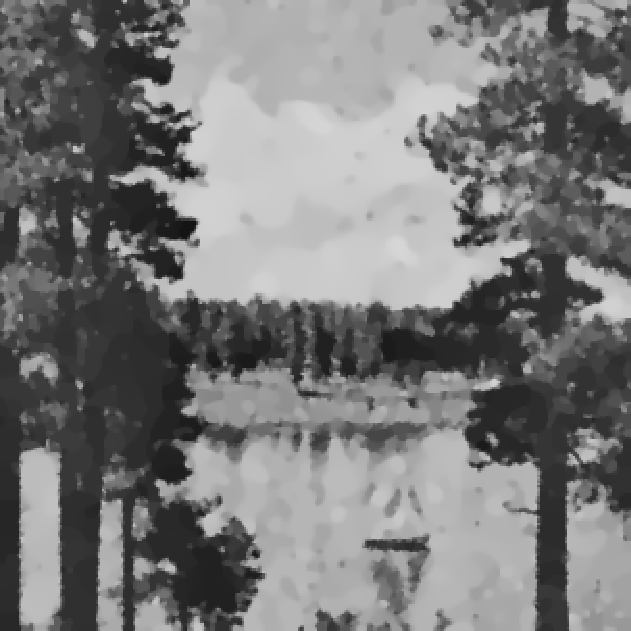}               
                \caption{$p=p_0-p_2$}
                \label{lake1l1_zzdb_tdmp}
       \end{subfigure}
       \begin{subfigure}[b]{0.24\textwidth}           
                \includegraphics[scale=0.36]{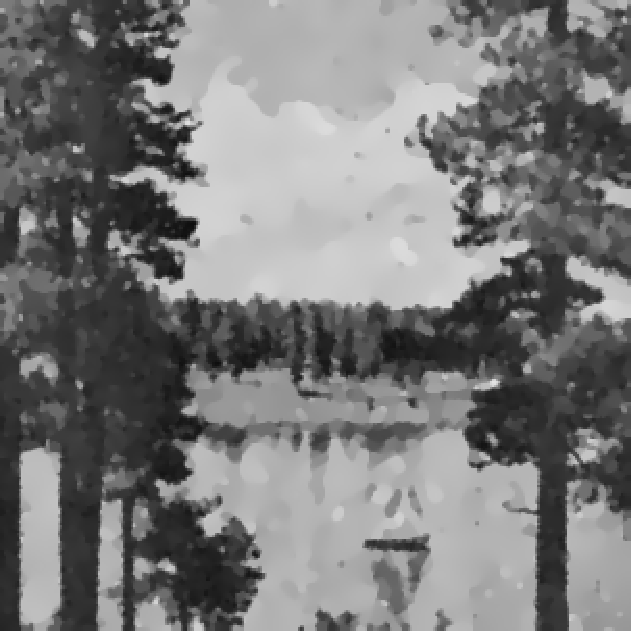}               
                \caption{$p=p_0-p_3$}
                \label{lake1l1_grad_tdmp}
       \end{subfigure}
       
\caption{\footnotesize Restored Images using \eqref{maina}--\eqref{mainc}.}\label{lake_1_10_restored_tdmp}
\end{figure}

\begin{figure}
       \centering
       \begin{subfigure}[b]{0.24\textwidth}           
                \includegraphics[scale=0.21]{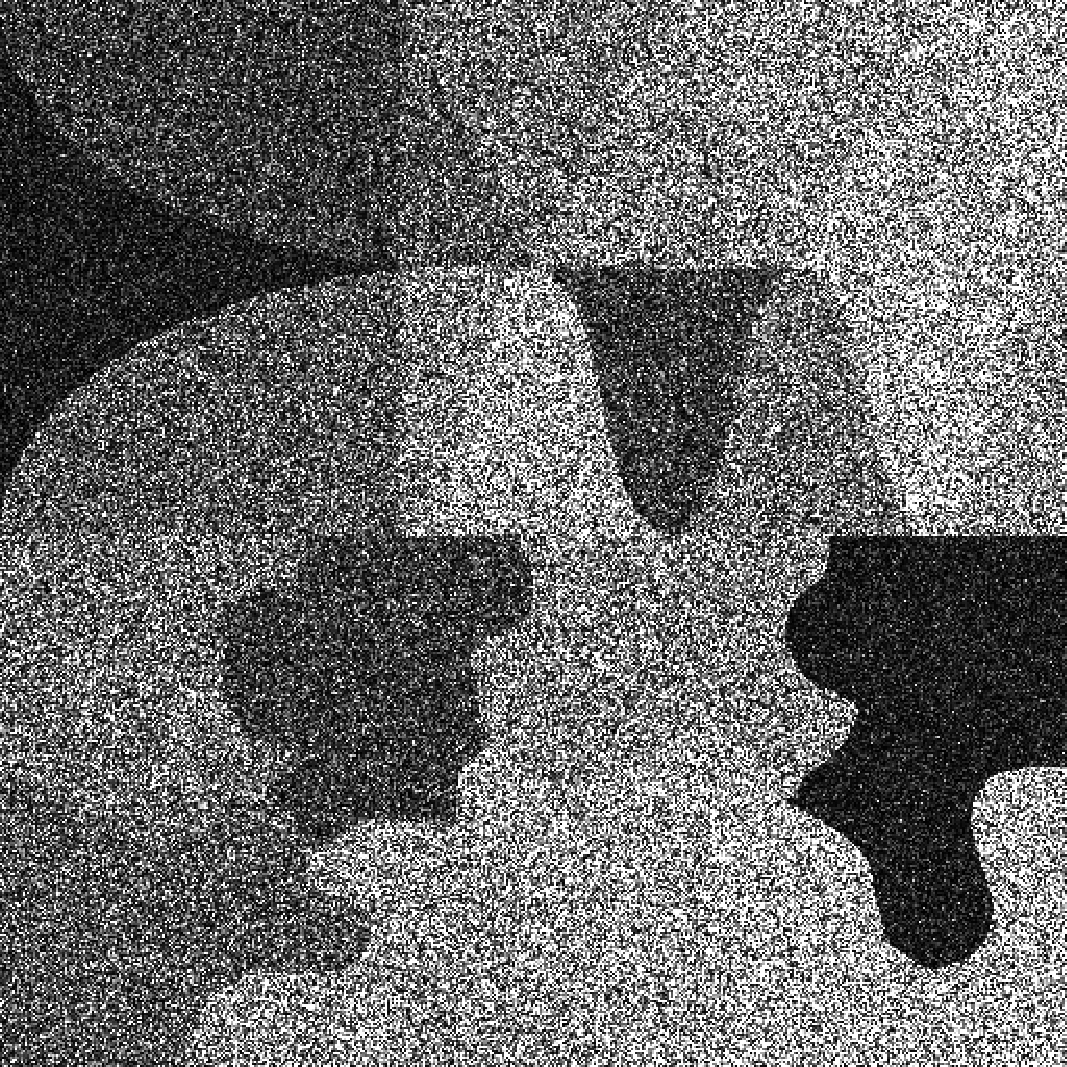}               
                \caption{Noisy $L=1$}
                \label{texturel1_noisy_l1}
       \end{subfigure}
       \begin{subfigure}[b]{0.24\textwidth}           
                \includegraphics[scale=0.21]{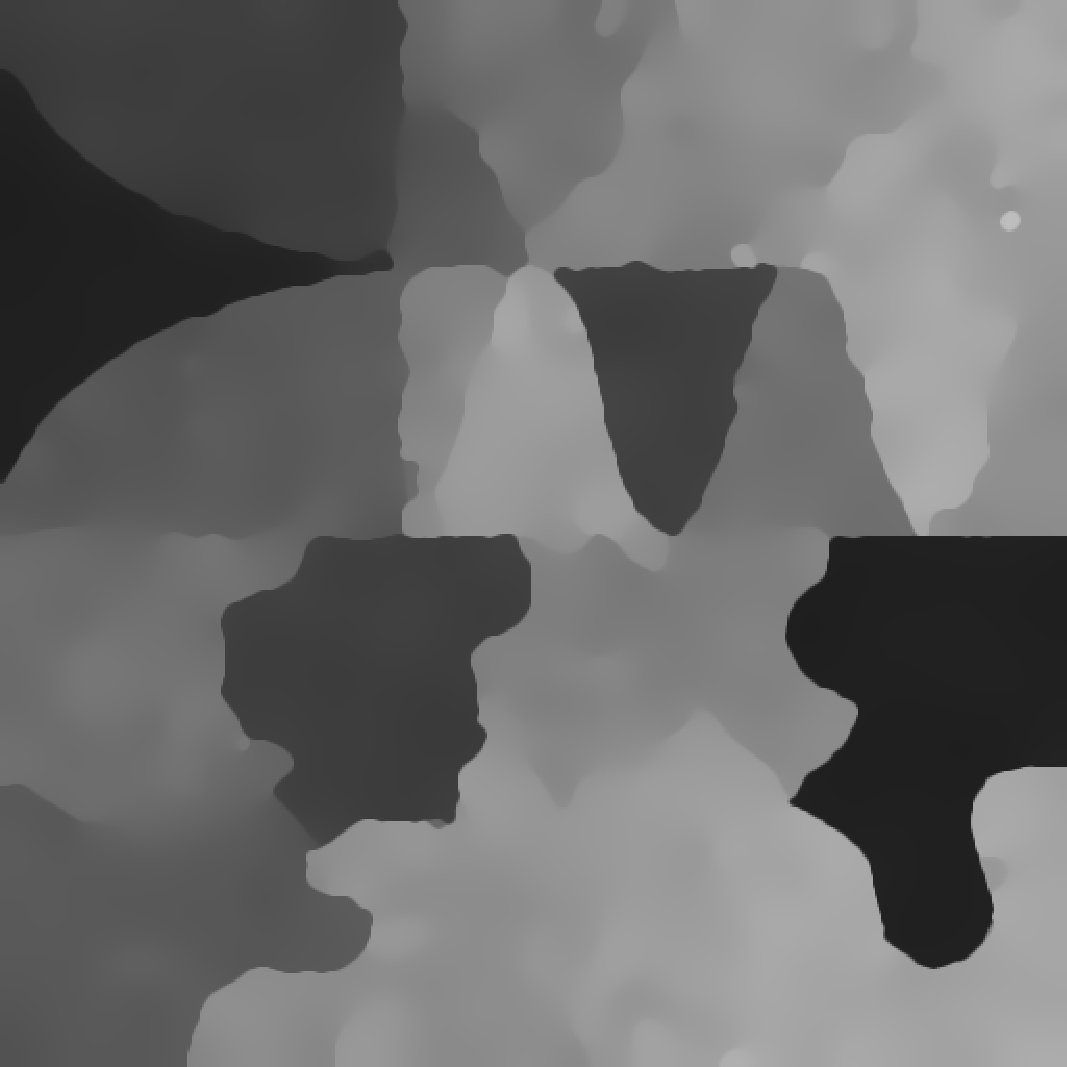}               
                \caption{$p=p_0-p_1$}
                \label{texturel1_avggray_tdmp}
       \end{subfigure}
            \begin{subfigure}[b]{0.24\textwidth}           
                \includegraphics[scale=0.21]{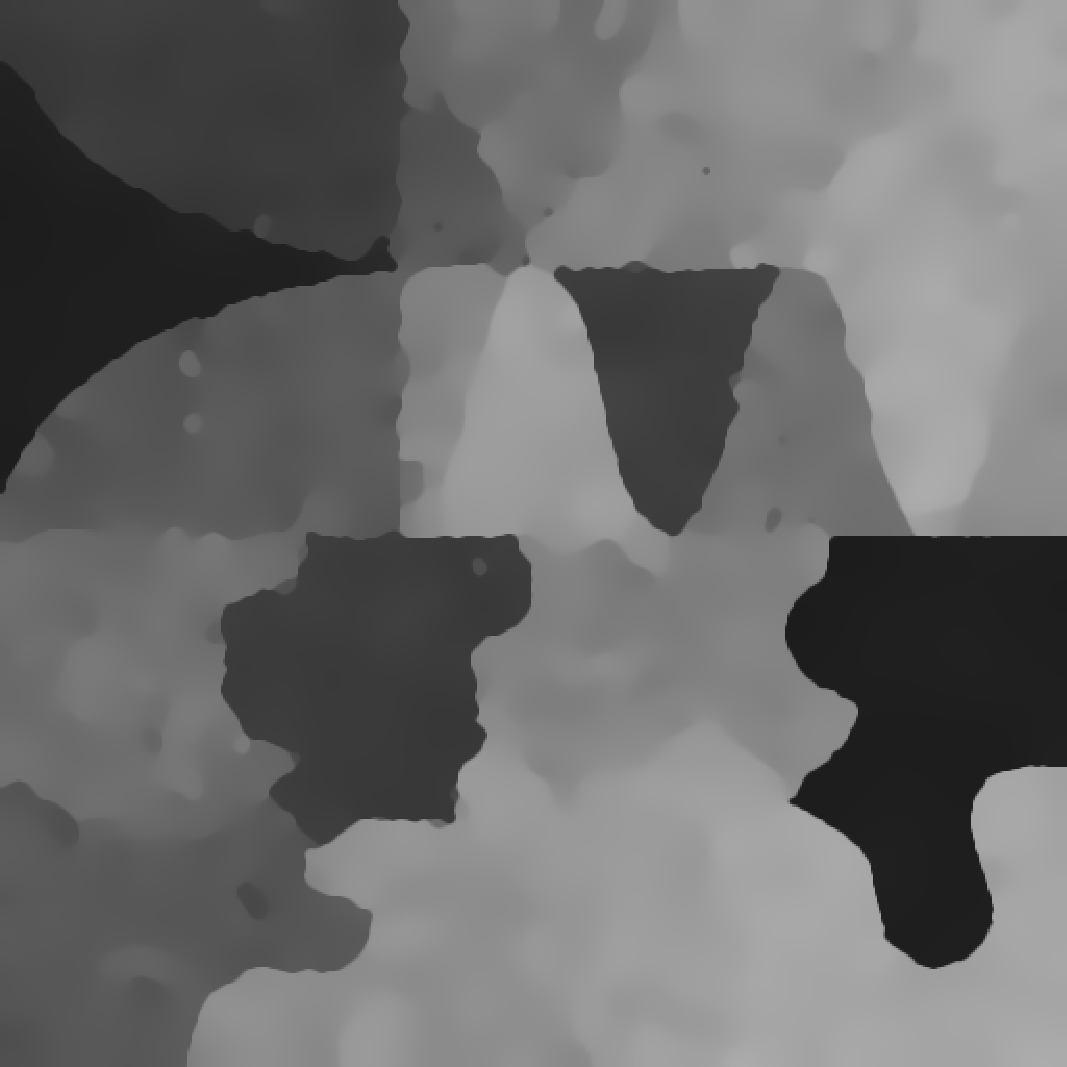}               
                \caption{$p=p_0-p_2$}
                \label{texturel1_zzdb_tdmp}
       \end{subfigure}
       \begin{subfigure}[b]{0.24\textwidth}           
                \includegraphics[scale=0.21]{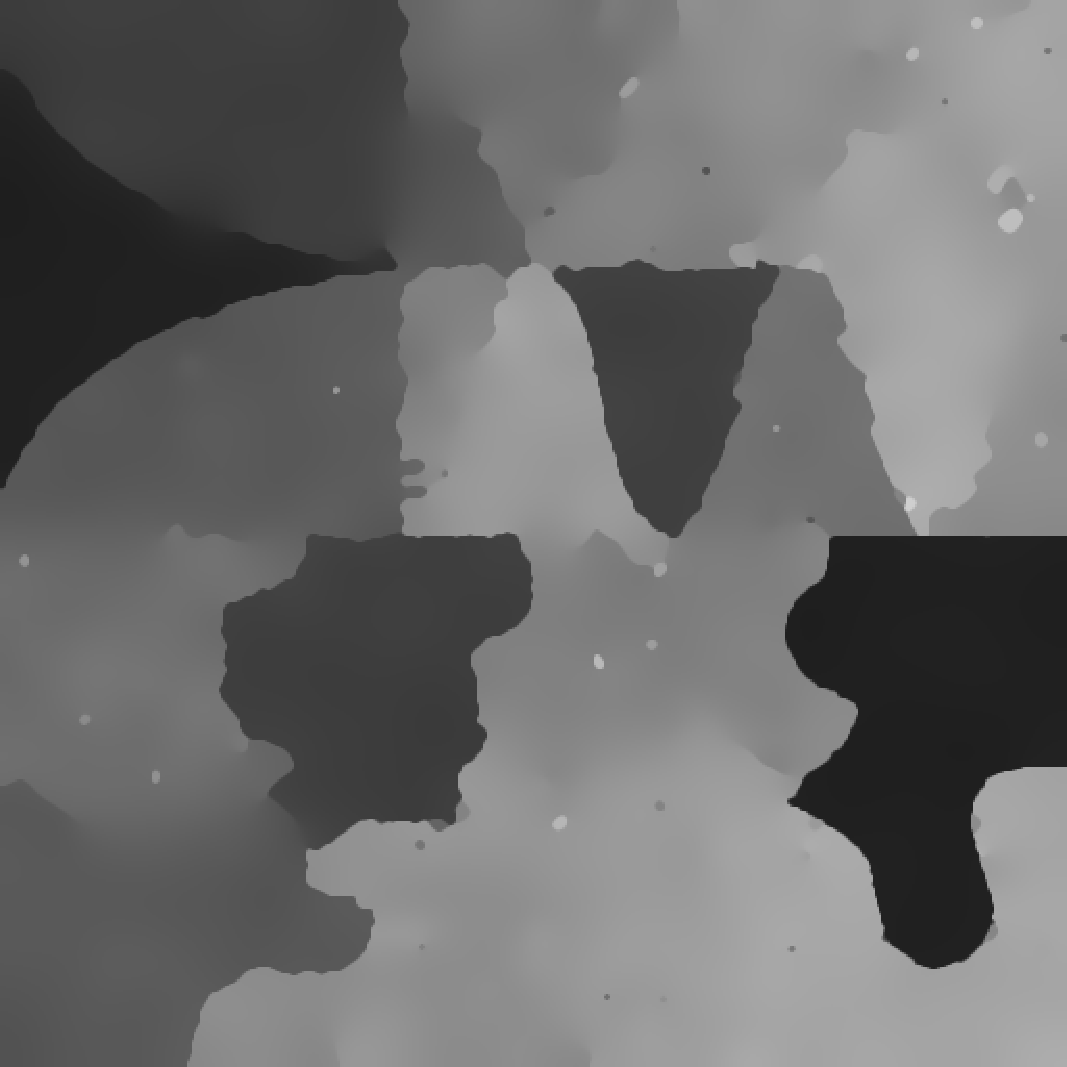}               
                \caption{$p=p_0-p_3$}
                \label{texturel1_grad_tdmp}
       \end{subfigure}
       
      \begin{subfigure}[b]{0.24\textwidth}           
                \includegraphics[scale=0.21]{texture_speckle_look_10}               
                \caption{Noisy $L=10$}
                \label{texturel1_noisy_l10}
       \end{subfigure}
       \begin{subfigure}[b]{0.24\textwidth}           
                \includegraphics[scale=0.21]{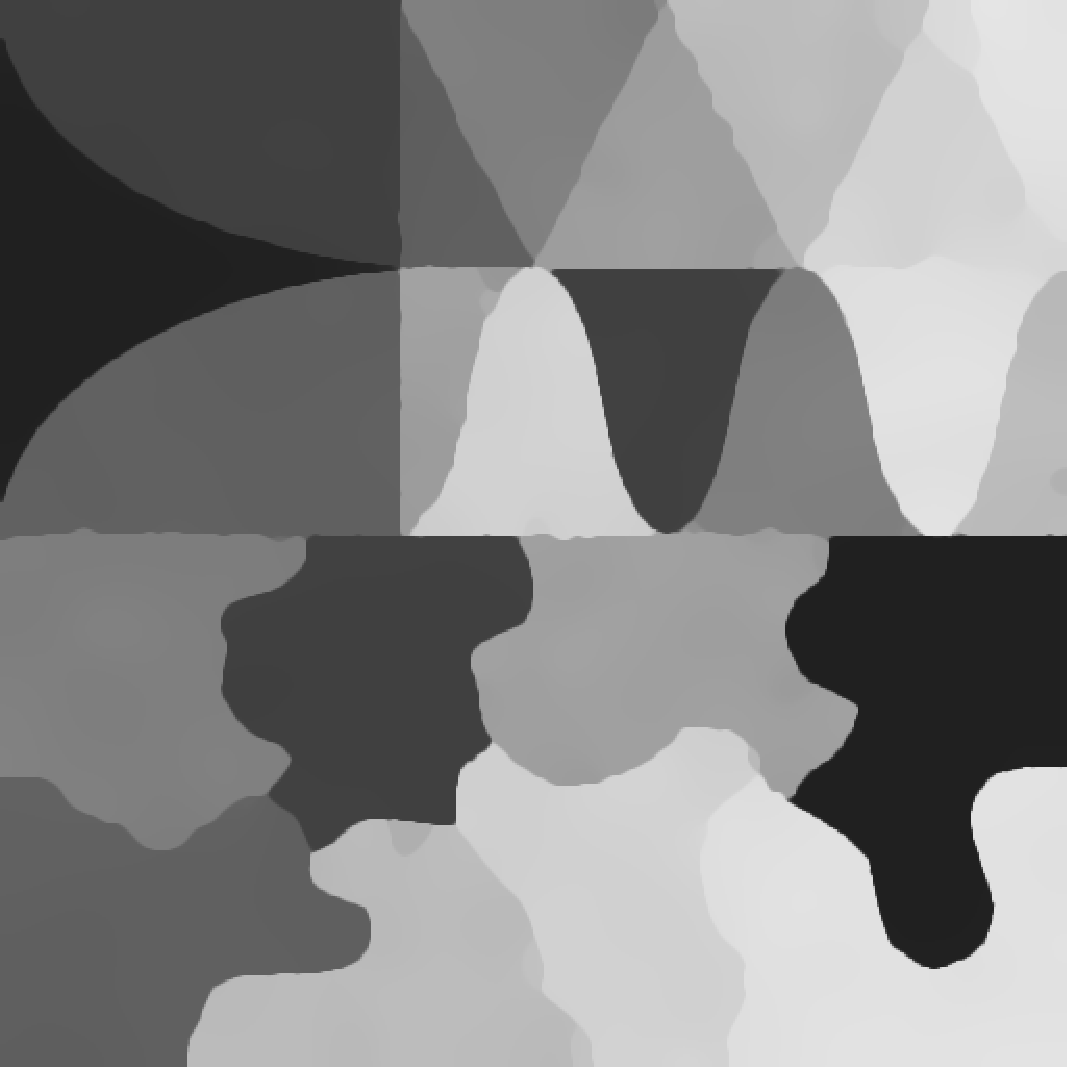}               
                \caption{$p=p_0-p_1$}
                \label{texturel10_avggray_tdmp}
       \end{subfigure}
            \begin{subfigure}[b]{0.24\textwidth}           
                \includegraphics[scale=0.21]{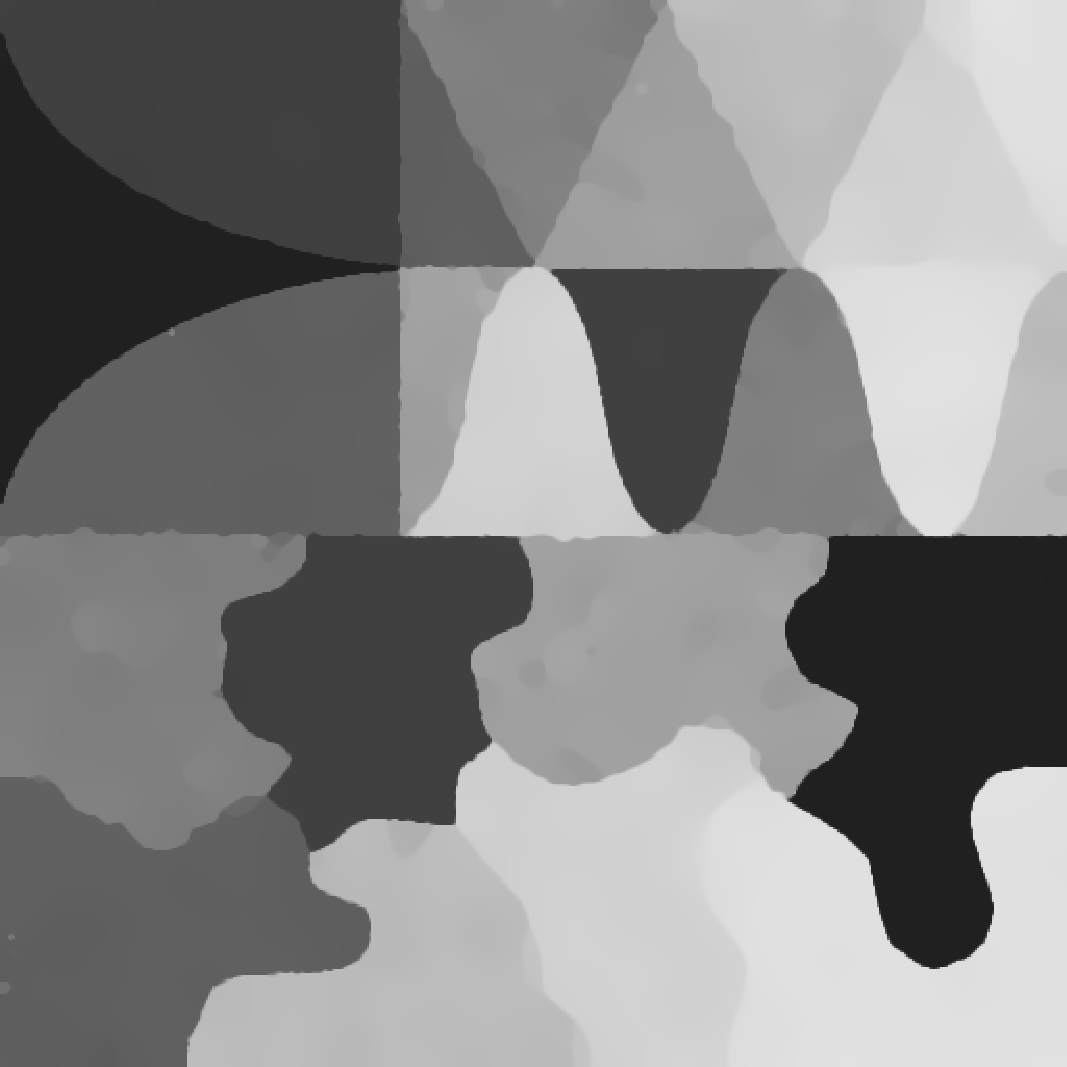}               
                \caption{$p=p_0-p_2$}
                \label{texturel10_zzdb_tdmp}
       \end{subfigure}
       \begin{subfigure}[b]{0.24\textwidth}           
                \includegraphics[scale=0.21]{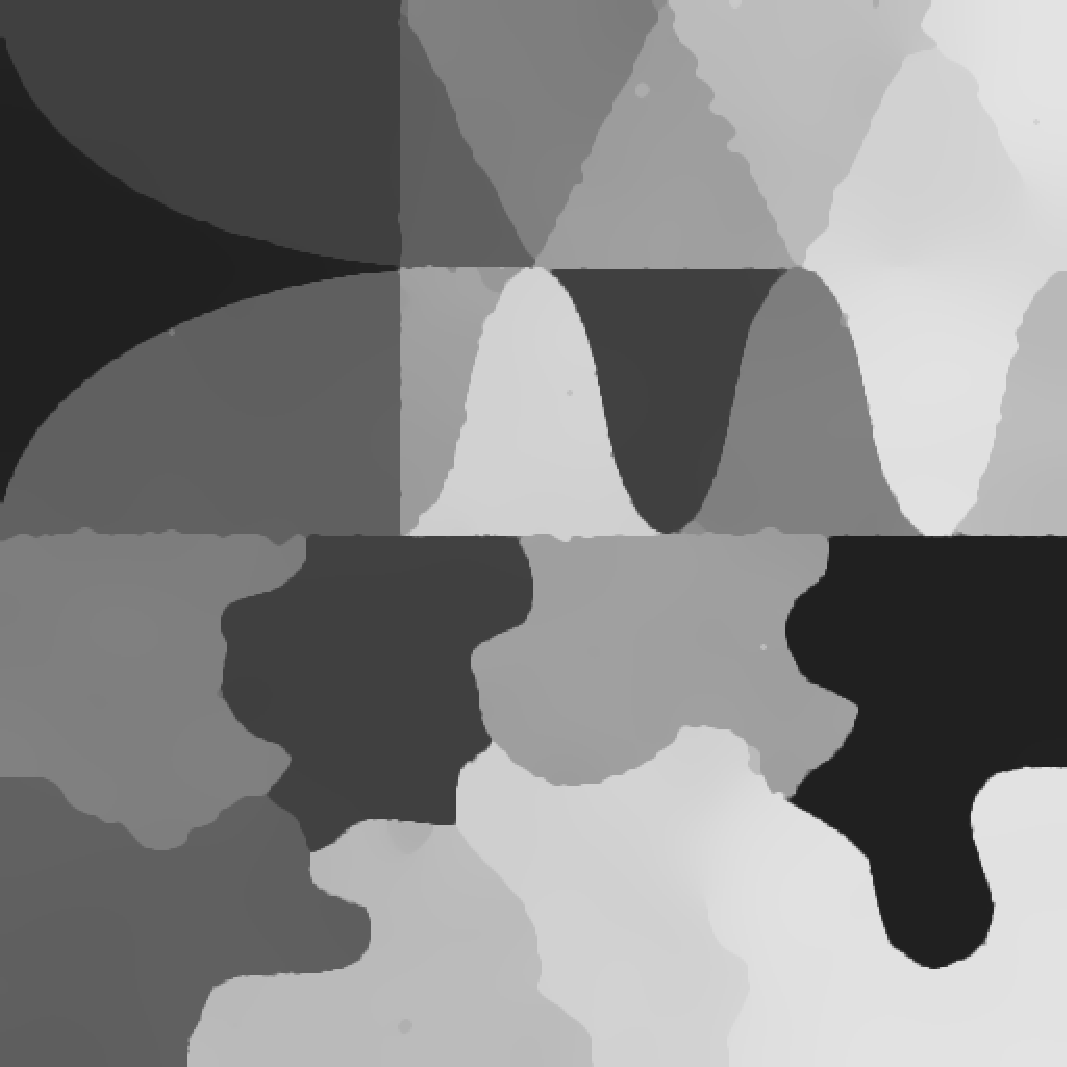}               
                \caption{$p=p_0-p_3$}
                \label{texturel10_grad_tdmp}
       \end{subfigure}
       
\caption{\footnotesize Restored Images using \eqref{maina}--\eqref{mainc}.}\label{texture_1_10_restored_tdmp}
\end{figure}

Figure \ref{circle_1_10_restored_tdmp} represents the despeckled results of a circle image initially degraded by speckle noise with $L=\{1, 10 \}$. The first row describes the image degraded by speckle noise with $L=1$ and the restored images by \eqref{maina}--\eqref{mainb} using the three exponents, $p=p_0-p_1$, $p=p_0-p_2$, and $p=p_0-p_3$, recall \eqref{p_1_relation}--\eqref{p_3_relation}. From the results, one can conclude that the first exponent outperforms the other two exponents. The third exponent destroys the edges when removing the image noise. The second row shows the speckled image degraded by $L=10$ and the restored images. Despeckled images indicate that all the exponents can restore the image quite accurately. However, from the 2D contour plots for the corresponding images, we can conclude that the first exponent can remove image noise better than the second and third exponents. In Figure \ref{lake_1_10_restored_tdmp}, we show the despeckling results for a lake image. From the results in the second row, it can be clearly shown that the first and third exponents can enhance the edges better than the second exponent. Figure \ref{texture_1_10_restored_tdmp} represents the despeckling results for a mosaic image initially degraded by speckle noise with $L=\{1, 10 \}$. Restored results indicate that all the exponents can adequately remove the speckle noise; however, the first exponent preserves the edges better than the other two. The last exponent introduced some undesired effects into the restored image for the high noise level. The computed values of corresponding quantitative measures are presented in Table \ref{tab:psnr_ssim_mor_vor_dg}, and the best results are highlighted in bold shapes.

\begin{figure}
       \centering
       
        \begin{subfigure}[b]{0.19\textwidth}           
                \includegraphics[scale=0.17]{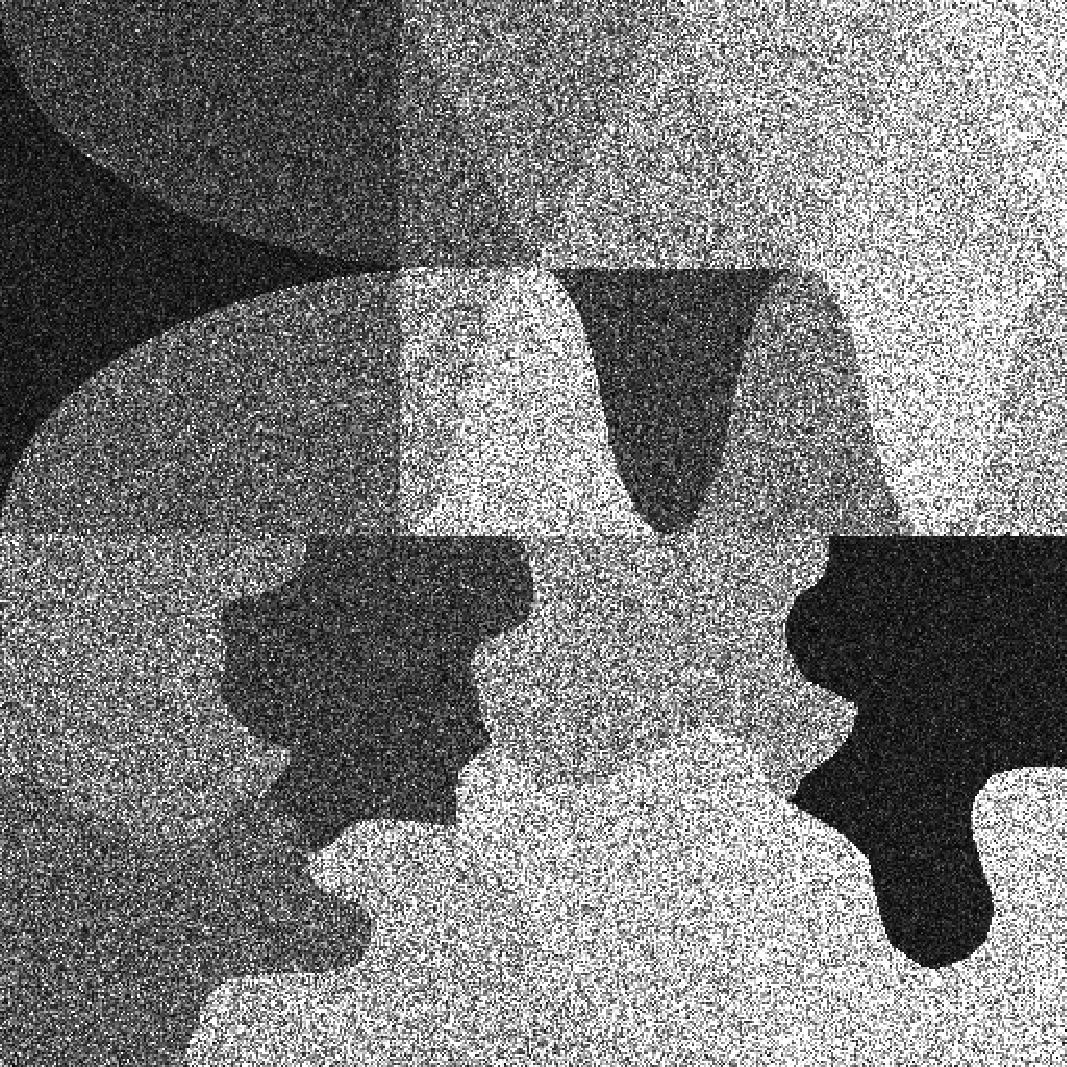}               
                \caption{Noisy}
                \label{texture_3}
       \end{subfigure}
       \begin{subfigure}[b]{0.19\textwidth}           
                \includegraphics[scale=0.17]{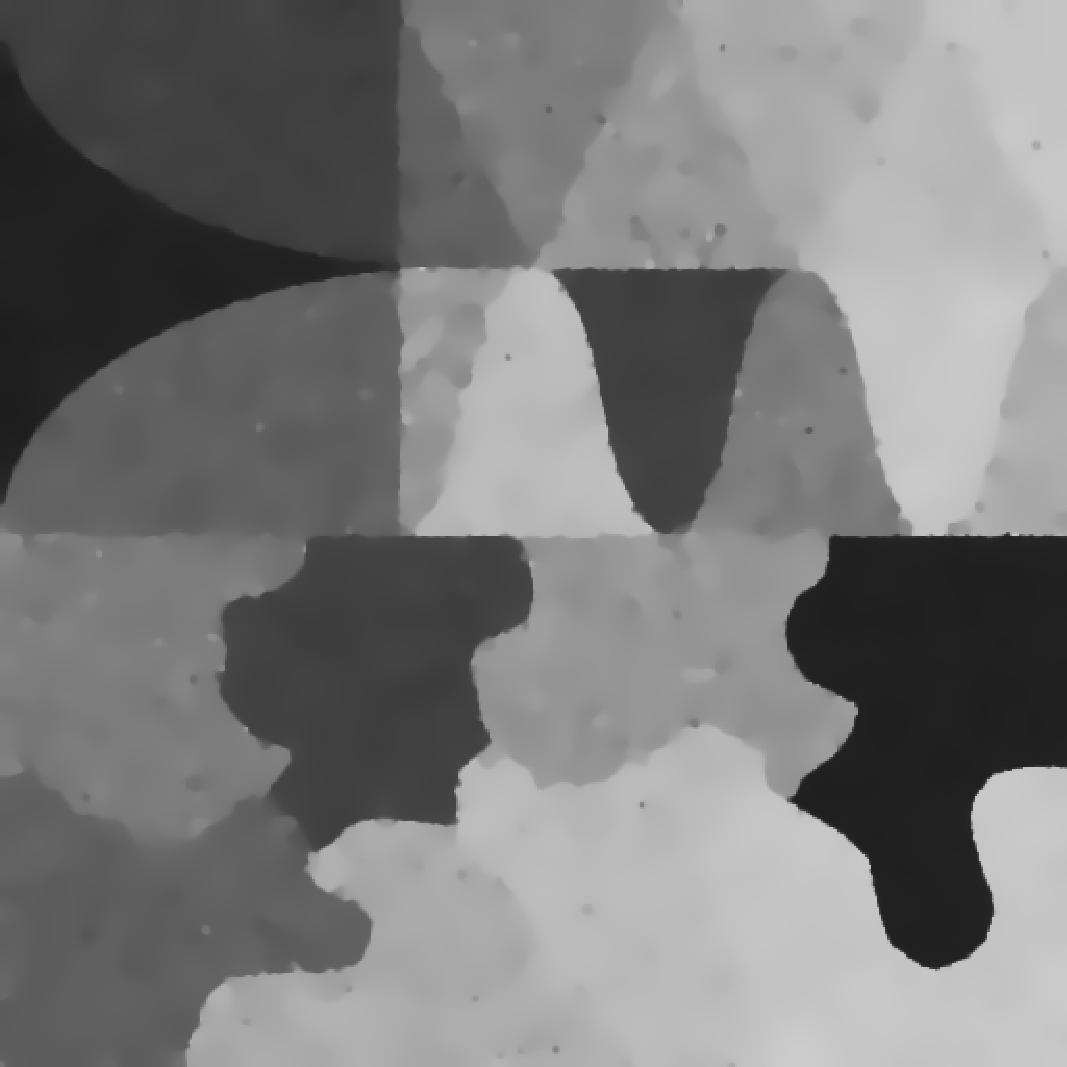}               
                \caption{DCE}
                \label{texture_3_diffconst}
       \end{subfigure}
       \begin{subfigure}[b]{0.19\textwidth}           
                \includegraphics[scale=0.17]{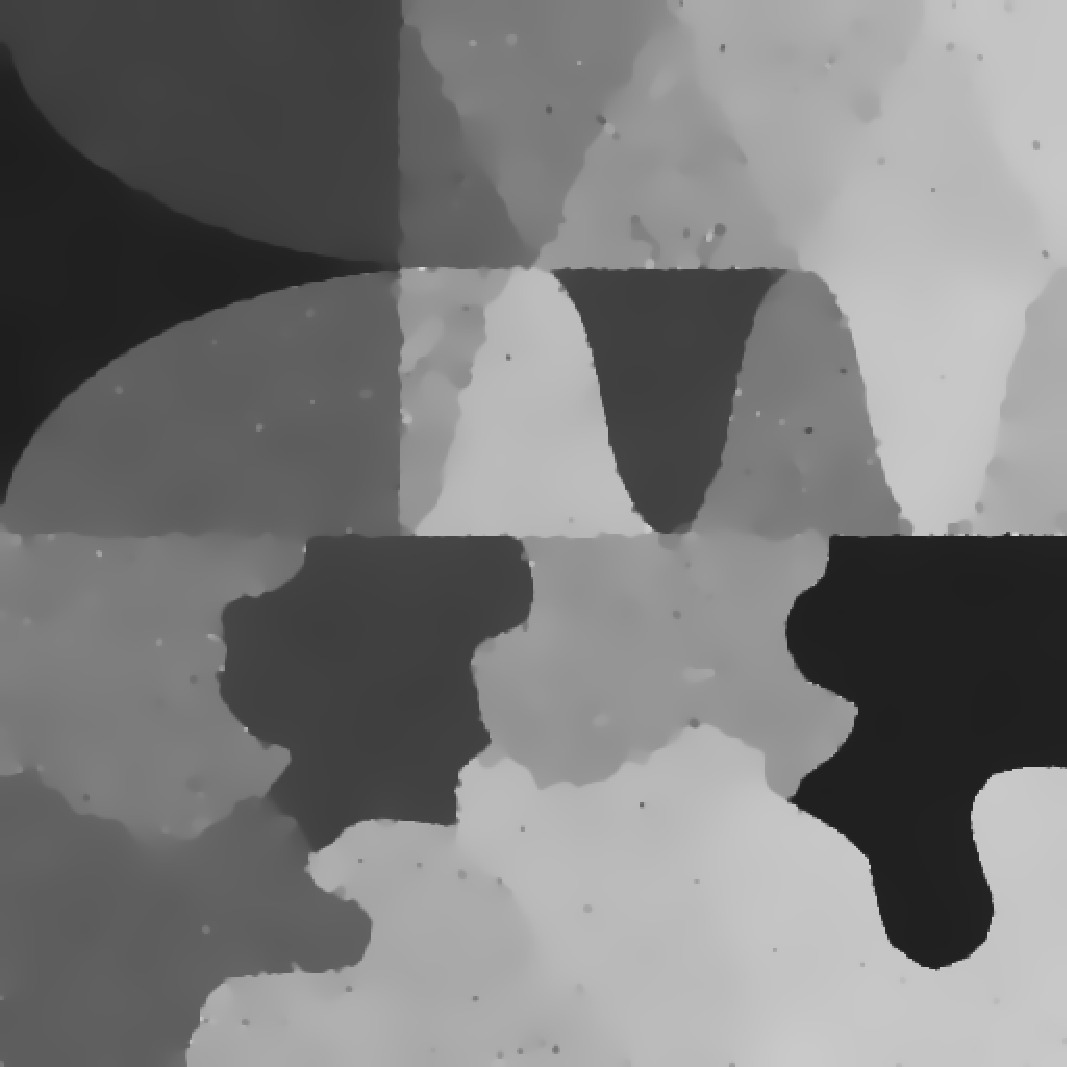}               
                \caption{DVE}
                \label{texture_3_diffgray}
       \end{subfigure}
            \begin{subfigure}[b]{0.19\textwidth}           
                \includegraphics[scale=0.17]{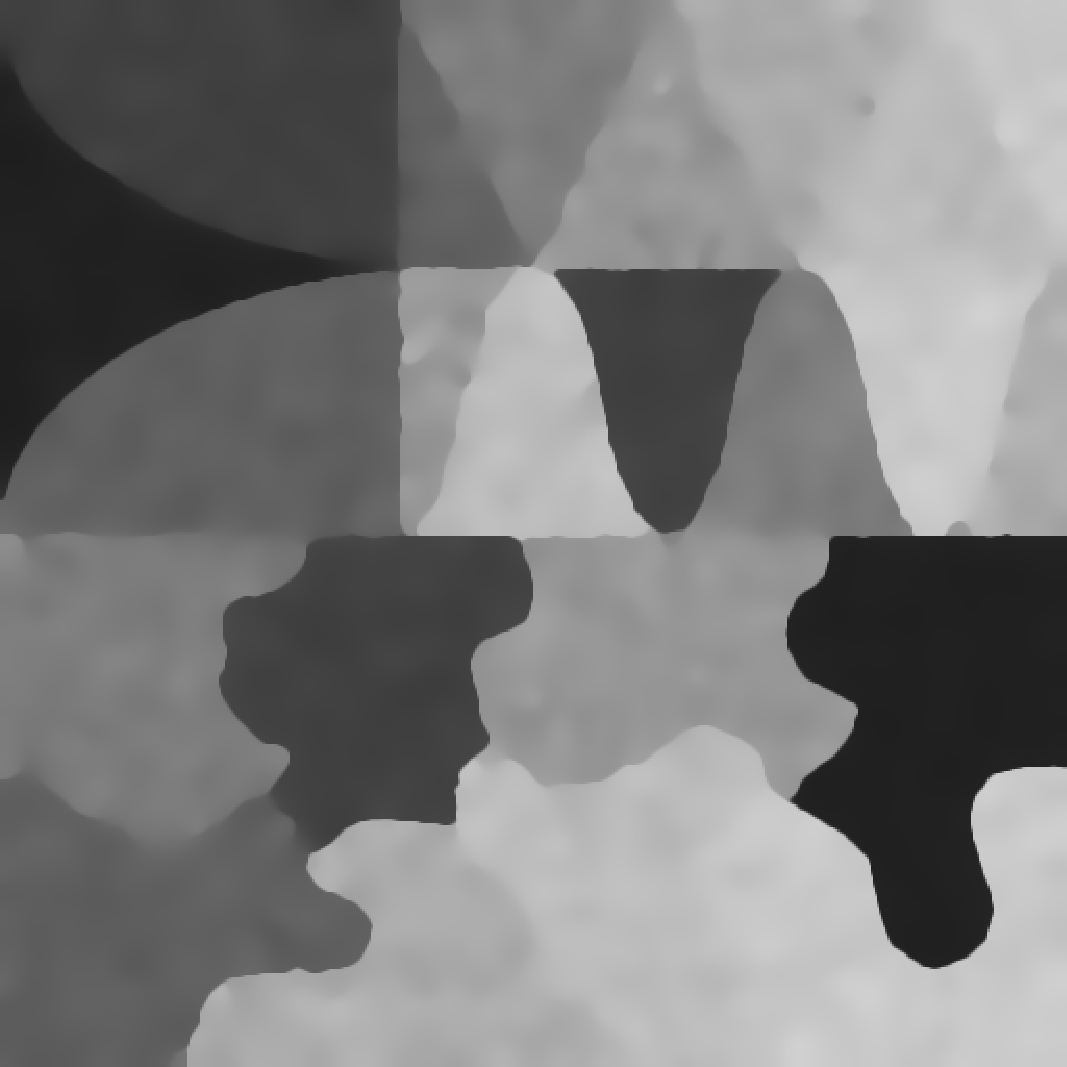}               
                \caption{TCE}
                \label{texture_3_tdmpcons}
       \end{subfigure}
       \begin{subfigure}[b]{0.19\textwidth}           
                \includegraphics[scale=0.17]{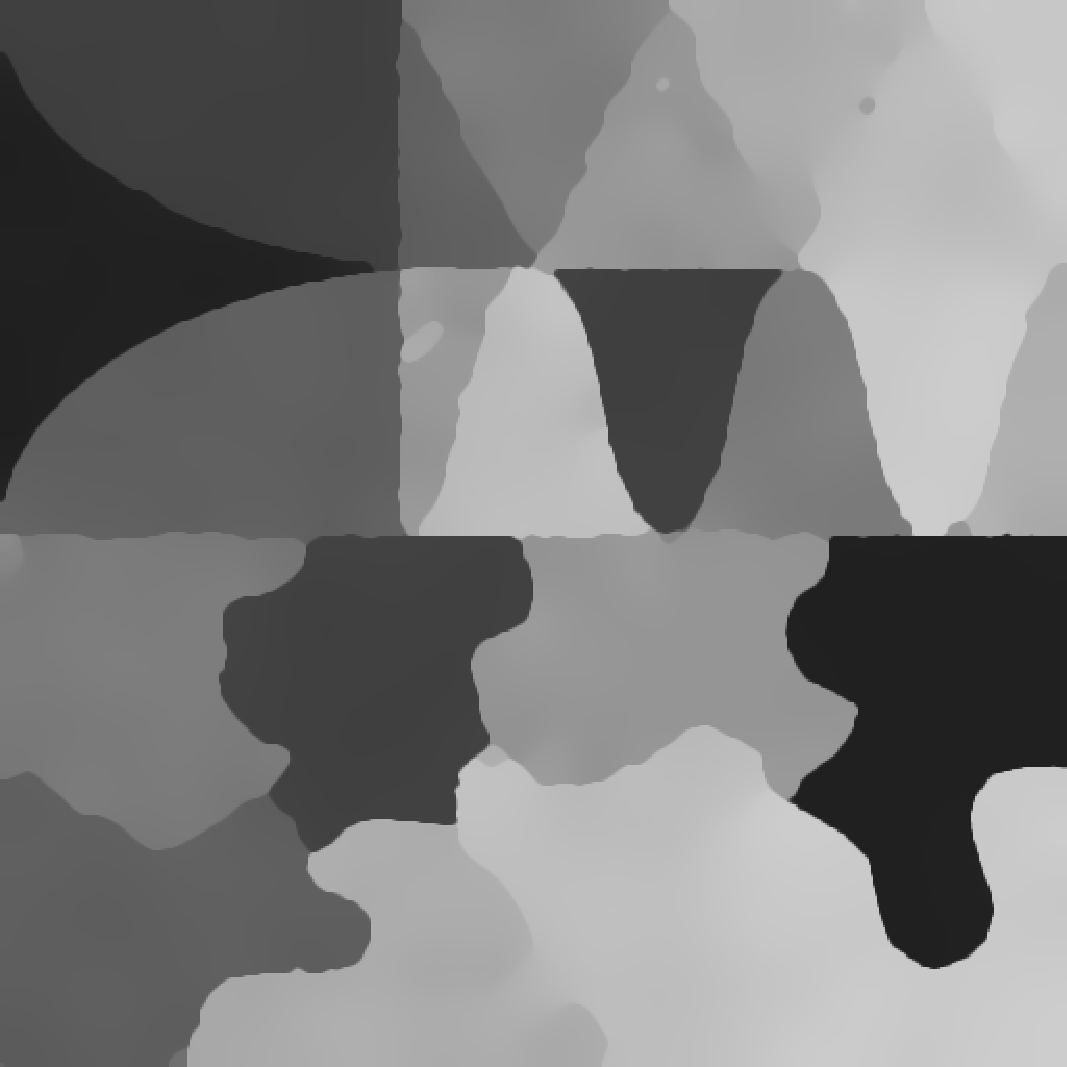}               
                \caption{TVE}
                \label{texture_3_tdmavggray}
       \end{subfigure}
       
       \begin{subfigure}[b]{0.19\textwidth}           
                \includegraphics[scale=0.17]{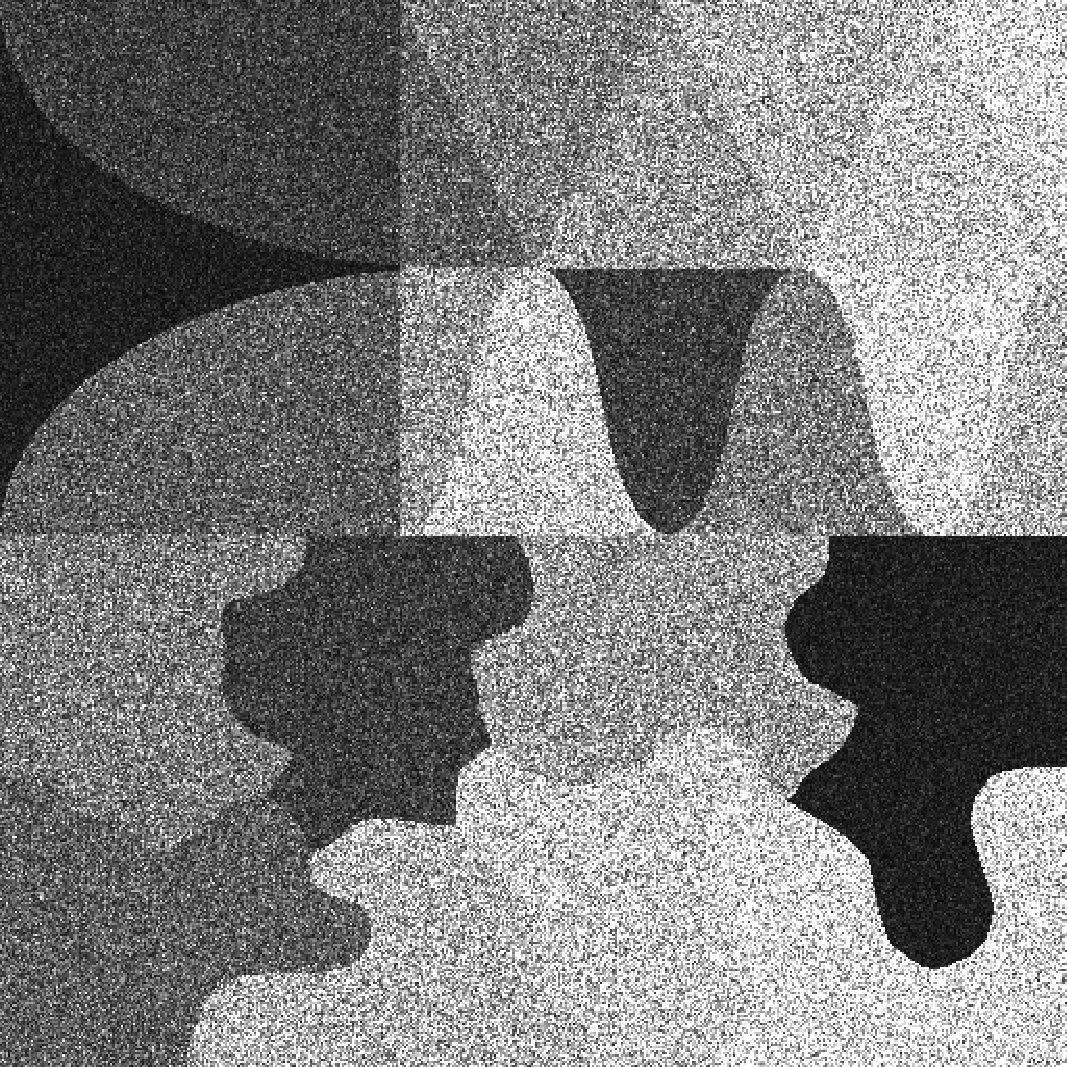}               
                \caption{Noisy}
                \label{texture_5}
       \end{subfigure}
       \begin{subfigure}[b]{0.19\textwidth}           
                \includegraphics[scale=0.17]{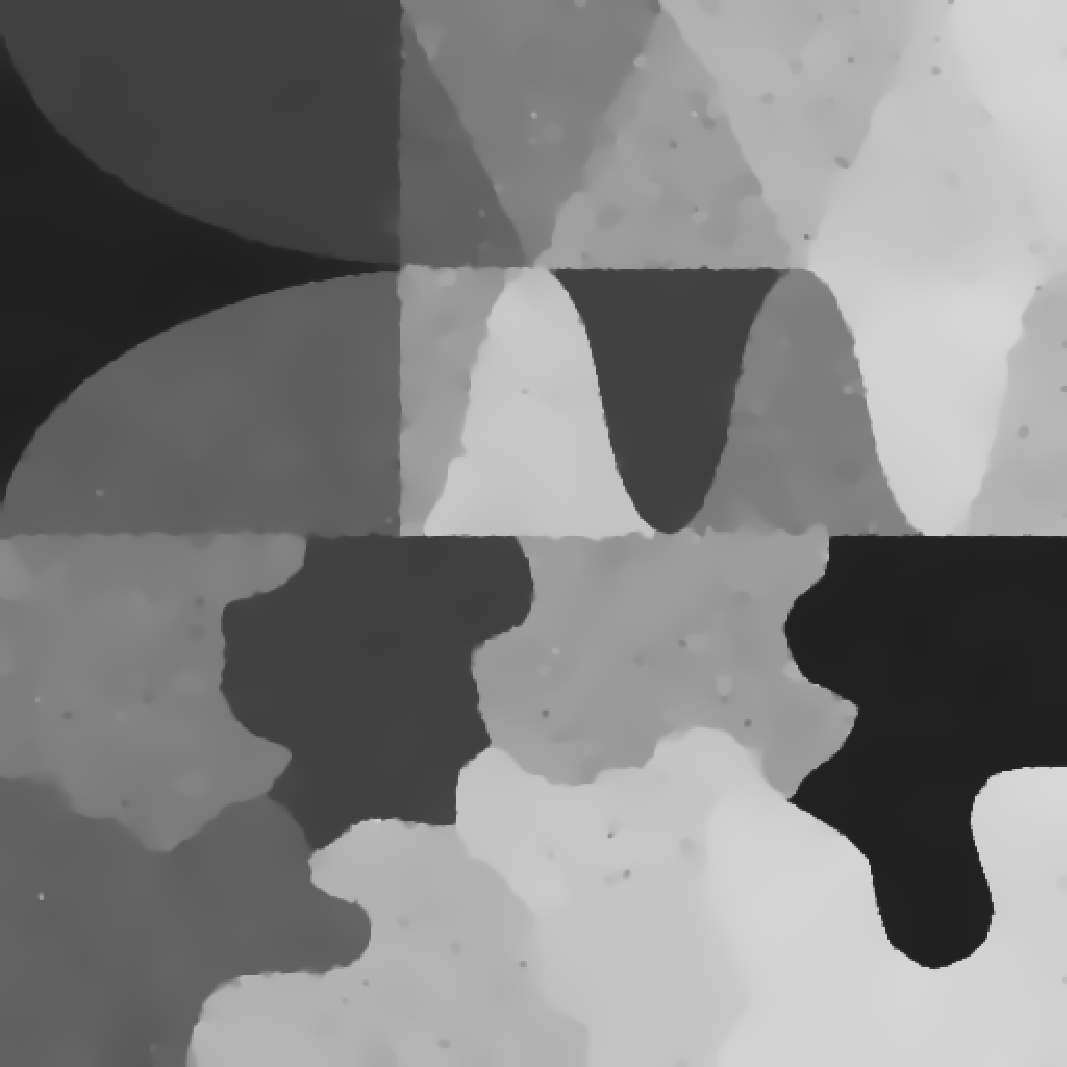}               
                \caption{DCE}
                \label{texture_5_diffconst}
       \end{subfigure}
       \begin{subfigure}[b]{0.19\textwidth}           
                \includegraphics[scale=0.17]{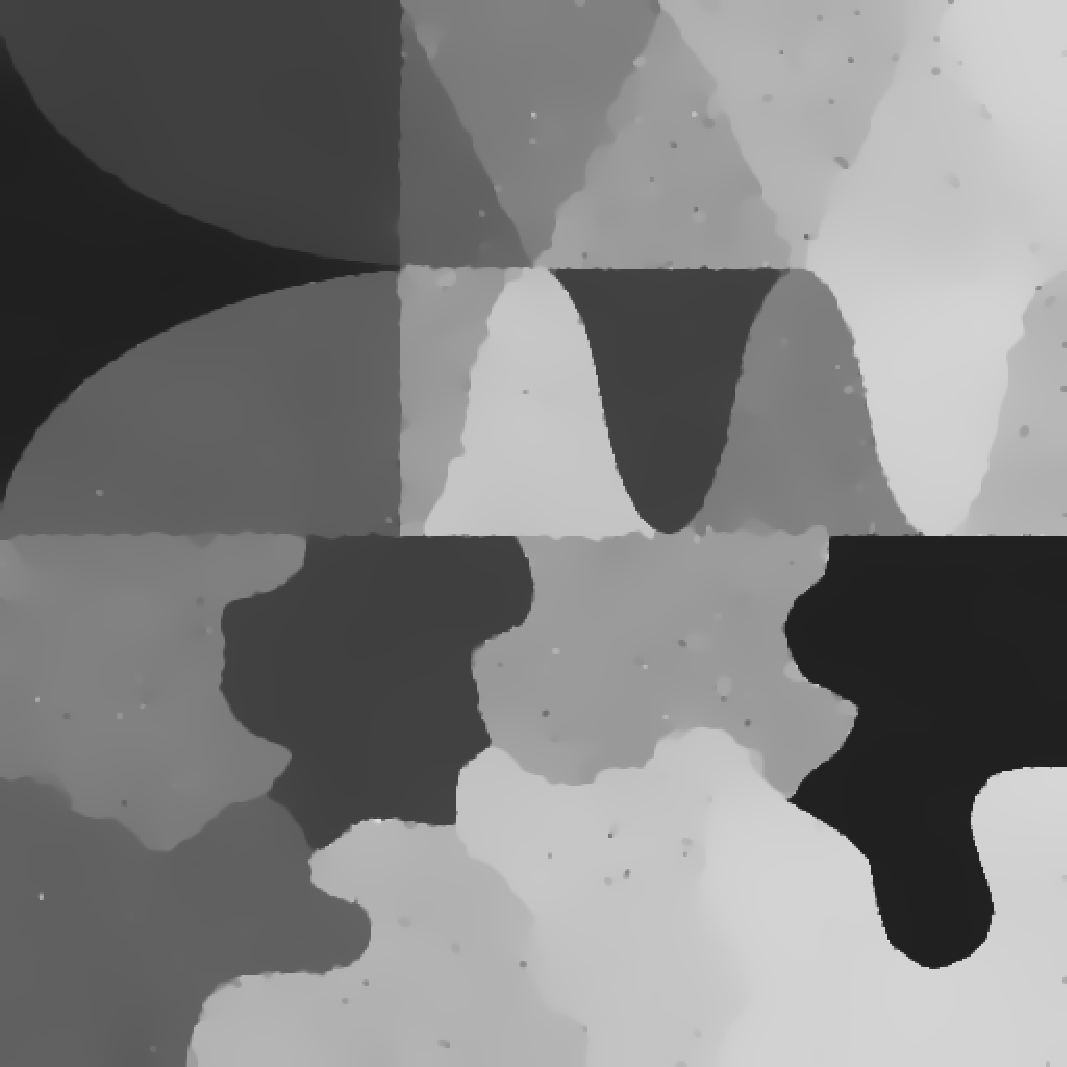}               
                \caption{DVE}
                \label{texture_5_diffgray}
       \end{subfigure}
            \begin{subfigure}[b]{0.19\textwidth}           
                \includegraphics[scale=0.17]{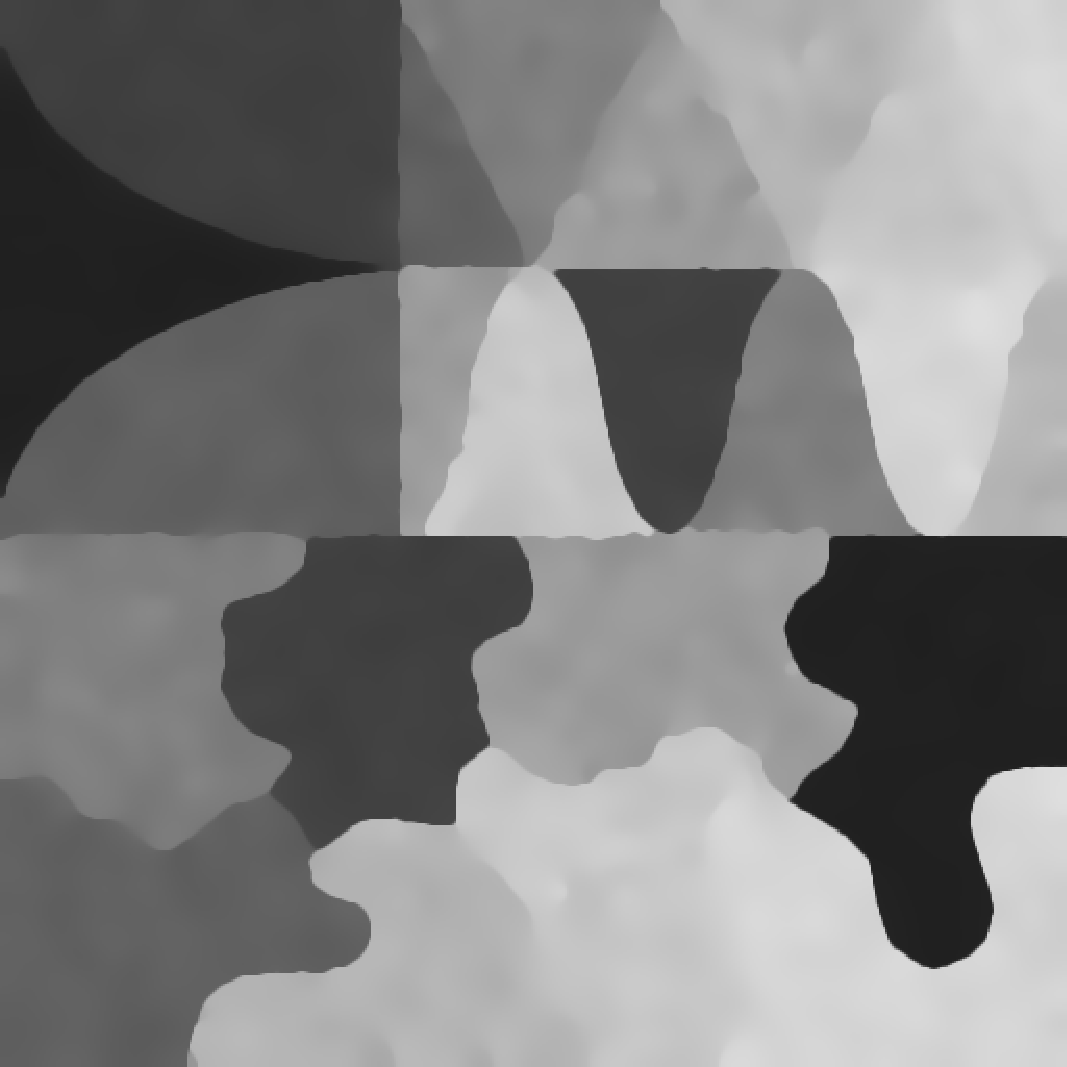}               
                \caption{TCE}
                \label{texture_5_tdmpcons}
       \end{subfigure}
       \begin{subfigure}[b]{0.19\textwidth}           
                \includegraphics[scale=0.17]{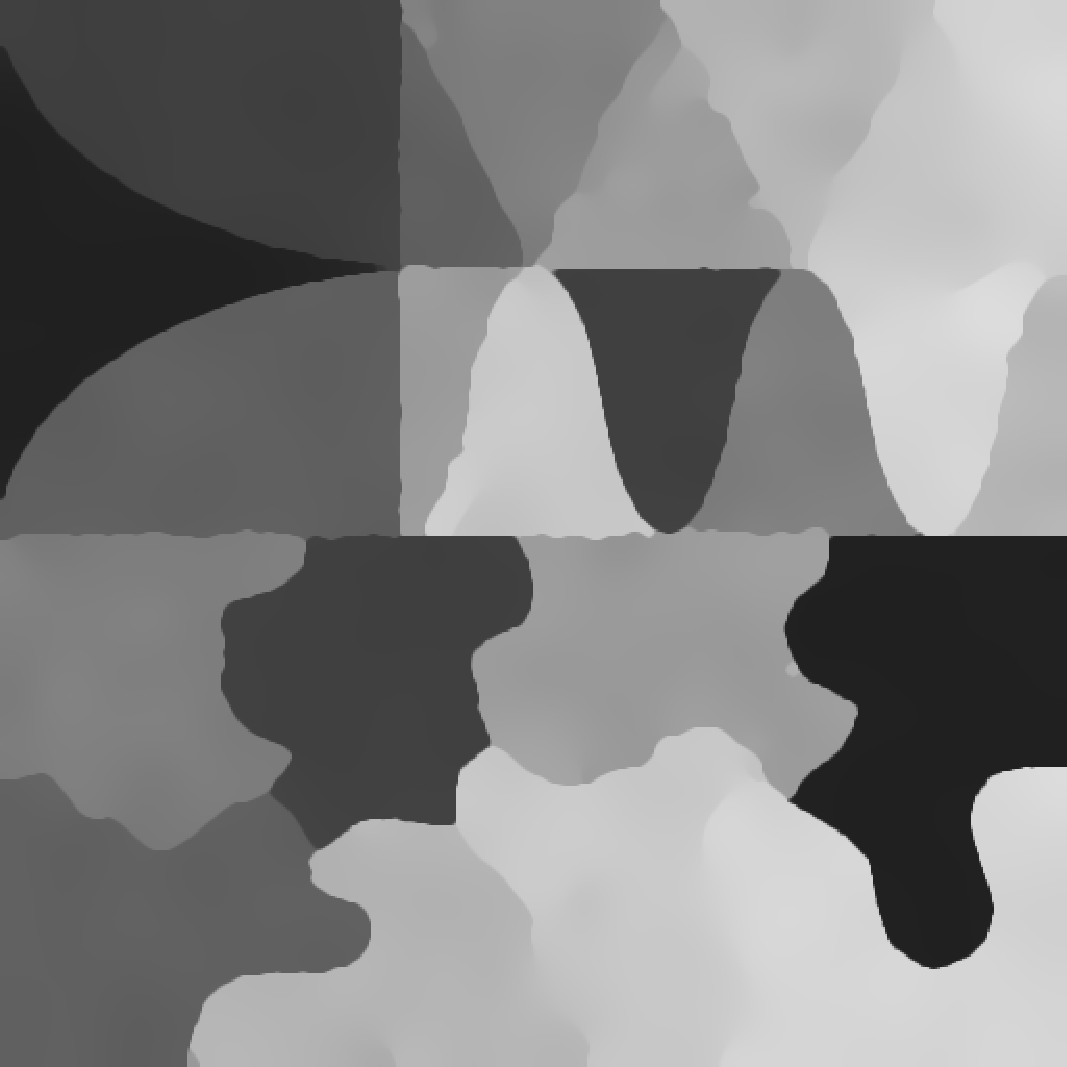}               
                \caption{TVE}
                \label{texture_5_tdmavggray}
       \end{subfigure}
       
\caption{\footnotesize \textbf{First Row:} (a) $L=3$, (b) PSNR=19.90, DG=7.70, (c) PSNR=19.92, DG=7.73, (d) PSNR=20.36, DG=8.17, (e) PSNR=20.61, DG=8.42. \textbf{Second Row:} (f) $L=5$, (g) PSNR=22.35, DG=8.49, (h) PSNR=22.39, DG=8.51, (i) PSNR=22.87, DG=9.00, (j) PSNR=23.14, DG=9.26.}\label{texture_3_5_restored_all_pde}
\end{figure}

\begin{figure}
           \begin{subfigure}[b]{0.9\textwidth}           
                \includegraphics[scale=1.0]{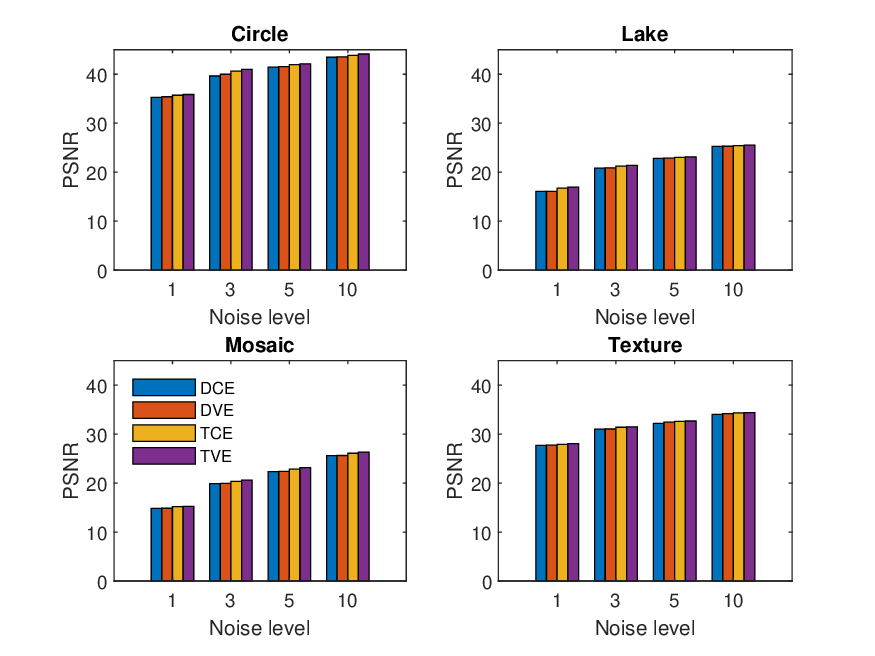}               
                \caption{Color legend shown for the mosaic image also works for the other images.}
                \label{bar_graph_psnr}
       \end{subfigure}

         \begin{subfigure}[b]{0.9\textwidth}           
                \includegraphics[scale=1.0]{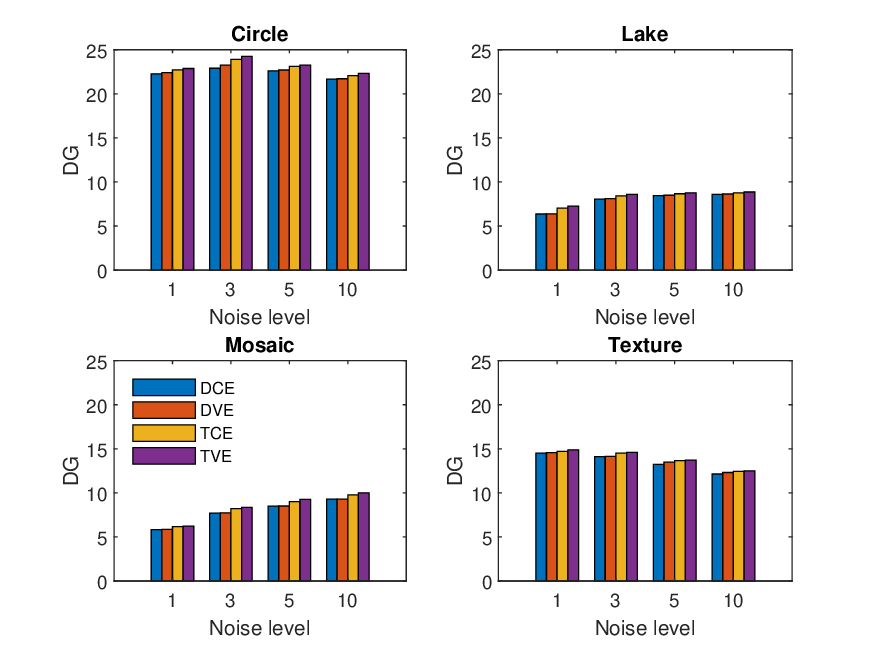}               
                \caption{Color legend shown for the mosaic image also works for the other images.}
                \label{bar_graph_dg}
       \end{subfigure}
 \caption{\footnotesize Comparison of PSNR and DG values for different images in different noise levels $L=\{1, 3, 5, 10\}$.}\label{bar_graph:psnr_dg}
\end{figure}

In Figure \ref{texture_3_5_restored_all_pde}, we compare the despeckling results computed by the approaches: diffusion model with the constant exponent (DCE), diffusion model with the variable exponent (DVE), telegraph model with constant exponent (TCE), and the telegraph model with the variable exponent (TVE). We show the results for the Mosaic image with two different noise levels, $L=\{3, 5\}$. Also, in the caption of the figure, we mention the quantitative measures for the restored images. From the results, one can easily conclude the performances of the diffusion and telegraph models with constant and variable exponents. Also, in Figure \ref{bar_graph:psnr_dg}, we compare the PSNR and despeckling gain values for the restored images computed by the approaches: DCE, DVE, TCE, and TVE. From the results, one can conclude that TVE always wins over the other techniques.

\begin{table*}
\caption{\footnotesize Comparison of quantitative results.}
\label{tab:psnr_ssim_mor_vor_dg}
\centering
\scalebox{0.75}{
\begin{tabular}[t]{llcccc}
\toprule
\multirow{2}[8]{*}{Image} & \multirow{2}[8]{*}{$L$}  & \multicolumn{1}{c}{$p=p_0-p_1$} & \multicolumn{1}{c}{ $p=p_0-p_2$ } & \multicolumn{1}{c}{$p=p_0-p_3$}  &  

\multicolumn{1}{c}{Clean}  \\
\cmidrule(r){3-3}
\cmidrule(r){4-4}
\cmidrule(r){5-5}
\cmidrule(r){6-6}  
       
&   & \multicolumn{1}{c}{PSNR/MSSIM/MoR/VoR/DG} & \multicolumn{1}{c}{PSNR/MSSIM/MoR/VoR/DG} & \multicolumn{1}{c}{PSNR/MSSIM/MoR/VoR/DG}  & \multicolumn{1}{c}{MoR/VoR}\\
      \midrule
Circle   & 1   &  35.86/\textbf{0.9722}/0.998/0.870/22.87   & \textbf{35.93}/0.9667/1.008/0.890/\textbf{22.94}  & 35.44/0.9691/0.996/0.865/22.44 & 0.977/0.839 \\
         & 3   &  \textbf{40.99}/\textbf{0.9826}/1.001/0.330/\textbf{24.25}   & 40.39/0.9800/1.002/0.330/23.65  & 40.23/0.9810/0.997/0.327/23.49 & 0.999/0.328\\
         & 5   &  \textbf{42.11}/\textbf{0.9842}/1.000/0.200/\textbf{23.25}   & 41.56/0.9830/1.001/0.201/22.70  & 41.50/0.9830/0.998/0.199/22.64 & 1.000/0.200\\
         & 10  &  \textbf{44.13}/\textbf{0.9879}/1.000/0.100/\textbf{22.31}   & 43.46/0.9870/1.000/0.100/21.64  & 43.62/0.9872/0.999/0.100/21.80 & 1.001/0.100\\
         &     &                                            &                                         &               &            \\
Lake     & 1   &  \textbf{16.96}/0.5674/0.934/0.542/\textbf{7.26}    & 16.41/0.5816/0.973/0.582/6.71   & 16.71/0.5678/0.947/0.549/7.01  & 0.861/0.555\\
         & 3   &  \textbf{21.39}/0.6945/0.965/0.220/\textbf{8.60}    & 21.08/0.7085/0.982/0.228/8.30   & 21.26/0.6984/0.970/0.222/8.48  & 0.946/0.240\\
         & 5   &  \textbf{23.14}/0.7510/0.977/0.140/\textbf{8.76}    & 23.01/0.7523/0.984/0.145/8.64   & 23.13/0.7468/0.976/0.143/\textbf{8.76}  & 0.966/0.158\\
         & 10  &  \textbf{25.53}/0.8067/0.982/0.071/\textbf{8.86}    & 25.40/0.8026/0.988/0.077/8.74   & 25.41/0.8010/0.983/0.0767/8.74 & 0.984/0.085\\
         &     &                                            &                                          &               &            \\
Mosaic   & 1   & 15.25/\textbf{0.8822}/0.989/0.609/6.23   & \textbf{15.29}/0.8797/1.004/0.637/\textbf{6.27}   & 15.20/0.8718/0.992/0.612/6.18  & 0.839/0.510 \\
         & 3   & \textbf{20.61}/\textbf{0.9478}/0.992/0.241/\textbf{8.42}     & 20.44/0.9426/0.999/0.244/8.25   & 20.50/0.9440/0.993/0.241/8.32  & 0.931/0.231\\
         & 5   & \textbf{23.14}/\textbf{0.9615}/0.993/0.154/\textbf{9.26}     & 22.87/0.9573/0.999/0.155/8.98   & 23.03/0.9598/0.994/0.154/\textbf{9.15}  & 0.954/0.151\\
         & 10  & \textbf{26.33}/\textbf{0.9782}/0.997/0.082/\textbf{10.00}    & 26.13/0.9741/0.999/0.083/9.80   & 26.26/0.9746/0.997/0.082/9.94  & 0.973/0.082\\
         &     &                                            &                                         &               &        \\  
Texture  & 1   & \textbf{28.07}/0.8404/0.967/0.827/\textbf{14.89}   & 27.89/\textbf{0.8409}/1.003/0.889/14.71  & 27.84/0.8348/0.980/0.840/14.67 & 0.981/0.861\\
         & 3   & \textbf{31.50}/\textbf{0.9009}/0.985/0.3201/\textbf{14.60}   & 31.32/0.8977/0.996/0.328/14.42  & 31.34/0.8990/0.988/0.322/14.44 & 0.998/0.392\\
         & 5   & \textbf{32.68}/\textbf{0.9211}/0.989/0.197/\textbf{13.72}    & 32.45/0.9178/0.996/0.200/13.49  & 32.58/0.9205/0.990/0.197/13.62 & 0.999/0.200\\
         & 10  & \textbf{34.39}/\textbf{0.9393}/0.992/0.100/\textbf{12.49}    & 34.20/0.9374/0.996/0.102/12.31  & 34.28/0.9391/0.993/0.102/12.38 & 1.001/0.100\\
         &     &                                   &                                      &                           &         \\     
\midrule          
\end{tabular}
}
\end{table*}

\subsection{Results for the Real Images}
\label{sec:real}
\begin{figure}
       \centering
        \begin{subfigure}[b]{0.3\textwidth}           
                \includegraphics[scale=0.55]{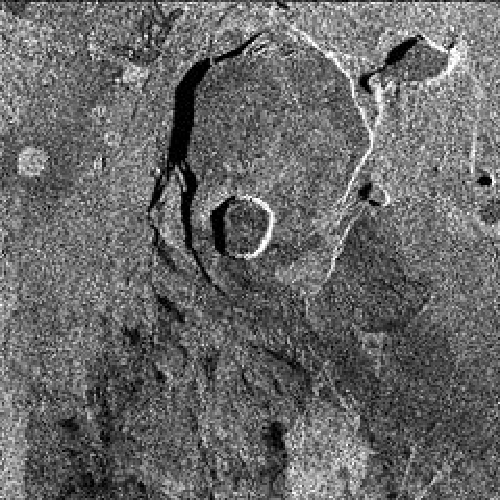}               
               \caption{Radar Image}
                \label{Image1_real}
       \end{subfigure}
        \begin{subfigure}[b]{0.3\textwidth}           
                \includegraphics[scale=0.55]{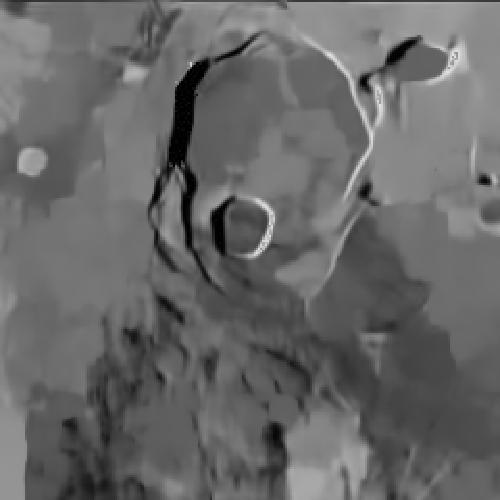}               
                \caption{SAR-BM3D}
                \label{Image1_sarbm3d}
       \end{subfigure}
       
        \begin{subfigure}[b]{0.3\textwidth}           
                \includegraphics[scale=0.55]{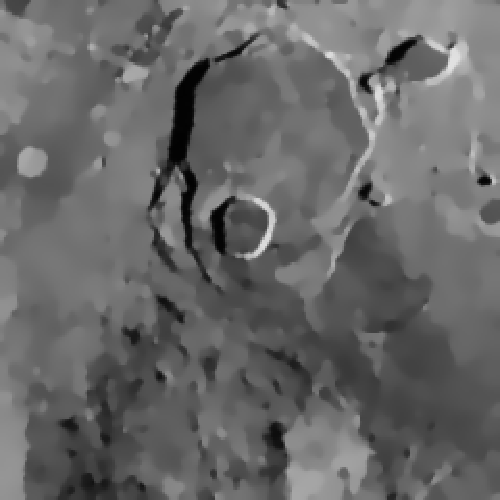}               
                \caption{$p=p_0-p_1$}
                \label{Image1_avggray1}
       \end{subfigure}
       \begin{subfigure}[b]{0.3\textwidth}           
                \includegraphics[scale=0.55]{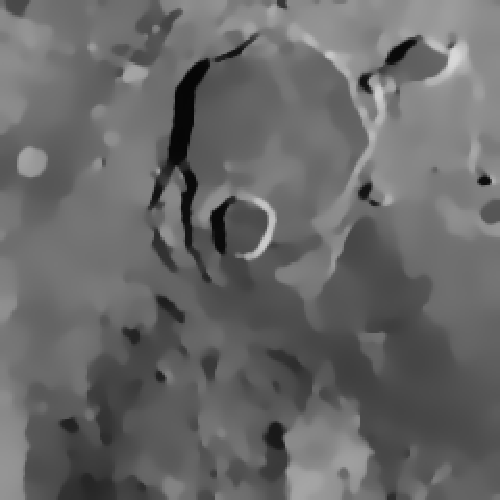}               
                \caption{$p=p_0-p_2$}
                \label{Image1_zzdb1}
       \end{subfigure} 
        \begin{subfigure}[b]{0.3\textwidth}           
                \includegraphics[scale=0.55]{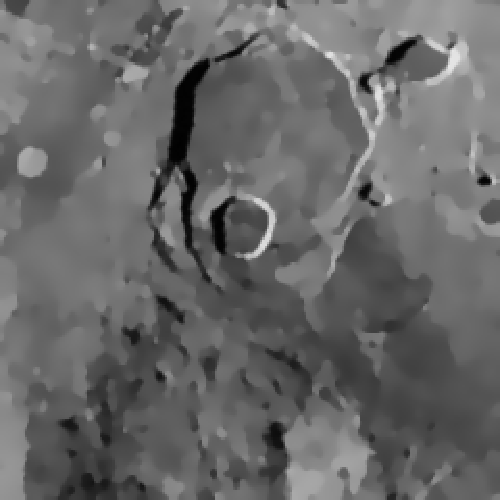}               
                \caption{$p=p_0-p_3$}
                \label{Image1_grad1}
       \end{subfigure}

\caption{\footnotesize (a) 240 $\times$ 240 Radar Image. (b) SI=0.2811. (c) SI=0.2788. (d) SI=0.2782. (e) SI=0.2982. }\label{Image1_tdmp}
\end{figure}

This section illustrates the image filtering ability of the proposed model in real SAR images. We have no prior information about the original noise-free image for the real-life images. Therefore, to compare the quantitative results, we calculate the value of the Speckle Index (SI) \cite{dewaele1990comparison}. A smaller value of SI indicates efficient speckle elimination. Figure \ref{Image1_tdmp} shows the despeckling results of the SAR image computed through \eqref{maina}--\eqref{mainc} and SAR-BM3D \cite{parrilli2011nonlocal}. From the restored images \ref{Image1_sarbm3d}--\ref{Image1_grad1}, one can conclude that all the methods can remove noise adequately though the first exponent preserves the texture and edges better than the others.

\subsection{Benchmarking Framework for SAR Despeckling}
\label{sec:benchmark}
\begin{figure}
       \centering
        \begin{subfigure}[b]{0.28\textwidth}           
                \includegraphics[scale=0.49]{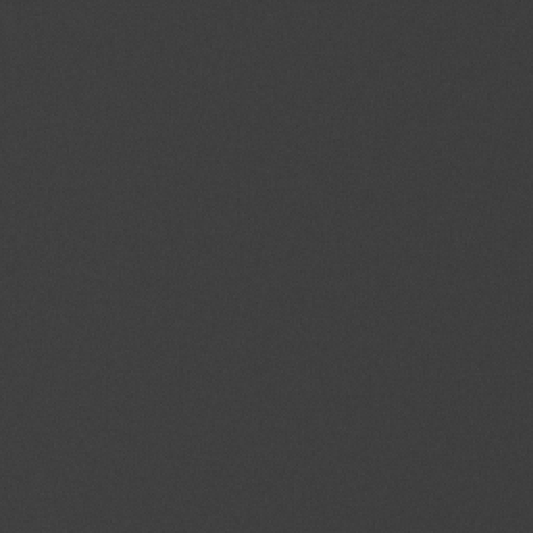}               
               \caption{}
                \label{homogeneous_clean1}
       \end{subfigure}
         \begin{subfigure}[b]{0.28\textwidth}           
                \includegraphics[scale=0.49]{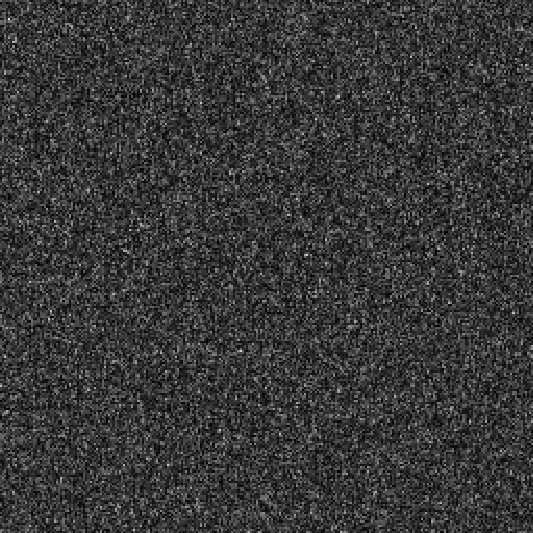}               
                \caption{}
                \label{homogeneous_noisy1}
       \end{subfigure}
       \begin{subfigure}[b]{0.28\textwidth}           
                \includegraphics[scale=0.49]{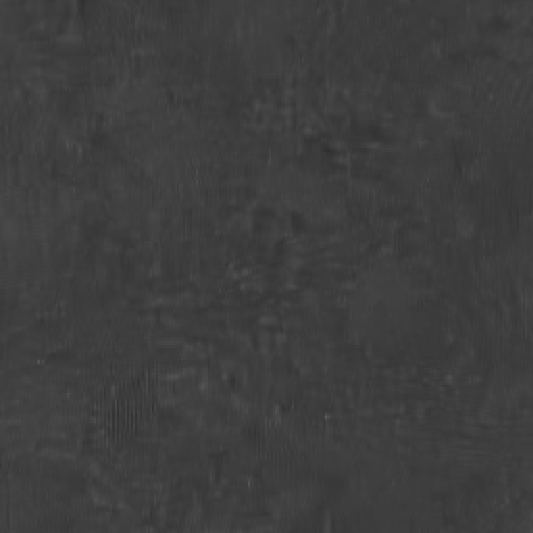}               
                \caption{}
                \label{homogeneous_sarbm3d1}
       \end{subfigure}
       
             \begin{subfigure}[b]{0.28\textwidth}           
                \includegraphics[scale=0.49]{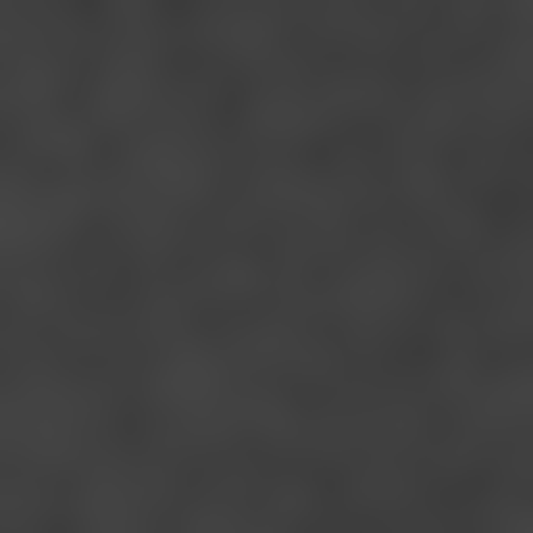}               
                \caption{}
                \label{homogeneous_avggraytdm1}
       \end{subfigure}
           \begin{subfigure}[b]{0.28\textwidth}           
                \includegraphics[scale=0.49]{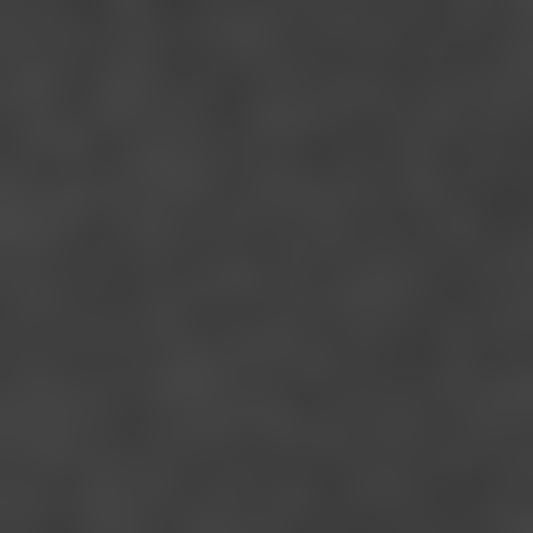}               
                \caption{}
                \label{homogeneous_zzdbtdm1}
       \end{subfigure}
            \begin{subfigure}[b]{0.28\textwidth}           
                \includegraphics[scale=0.49]{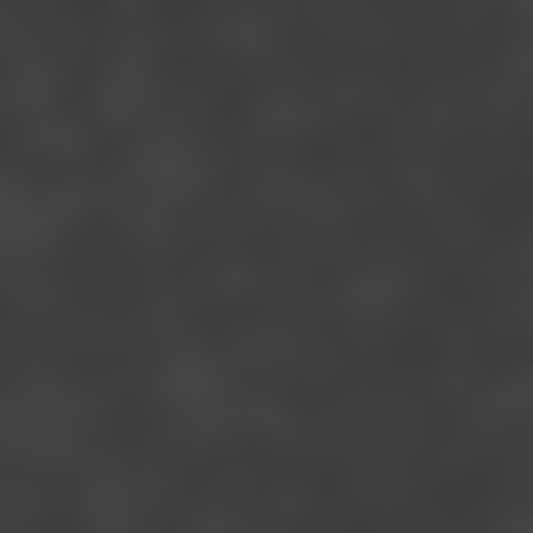}               
                \caption{}
                \label{homogeneous_gradtdm1}
       \end{subfigure}
       
\caption{\footnotesize 256 $\times$ 256 Homogeneous image. (a) 512-look. (b) One-look. (c) SAR-BM3D. (d) $Q=Q_0-p_1$. (e) $Q=Q_0-p_2$. (f) $Q=Q_0-p_3$.}\label{homogeneous_tdmpnew}
\end{figure}

\begin{figure}
       \centering
        \begin{subfigure}[b]{0.28\textwidth}           
                \includegraphics[scale=0.245]{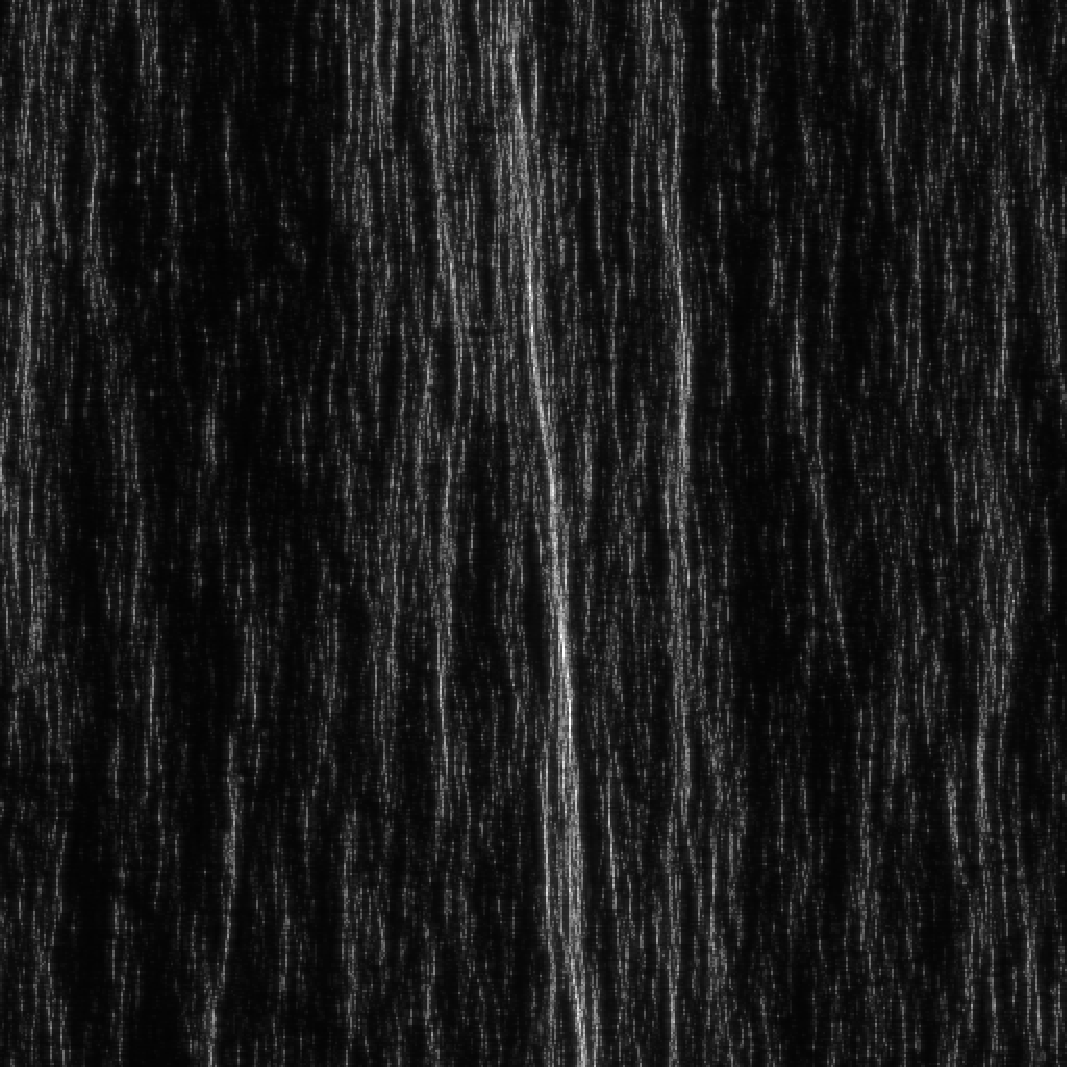}               
               \caption{}
                \label{dem_clean1}
       \end{subfigure}
         \begin{subfigure}[b]{0.28\textwidth}           
                \includegraphics[scale=0.245]{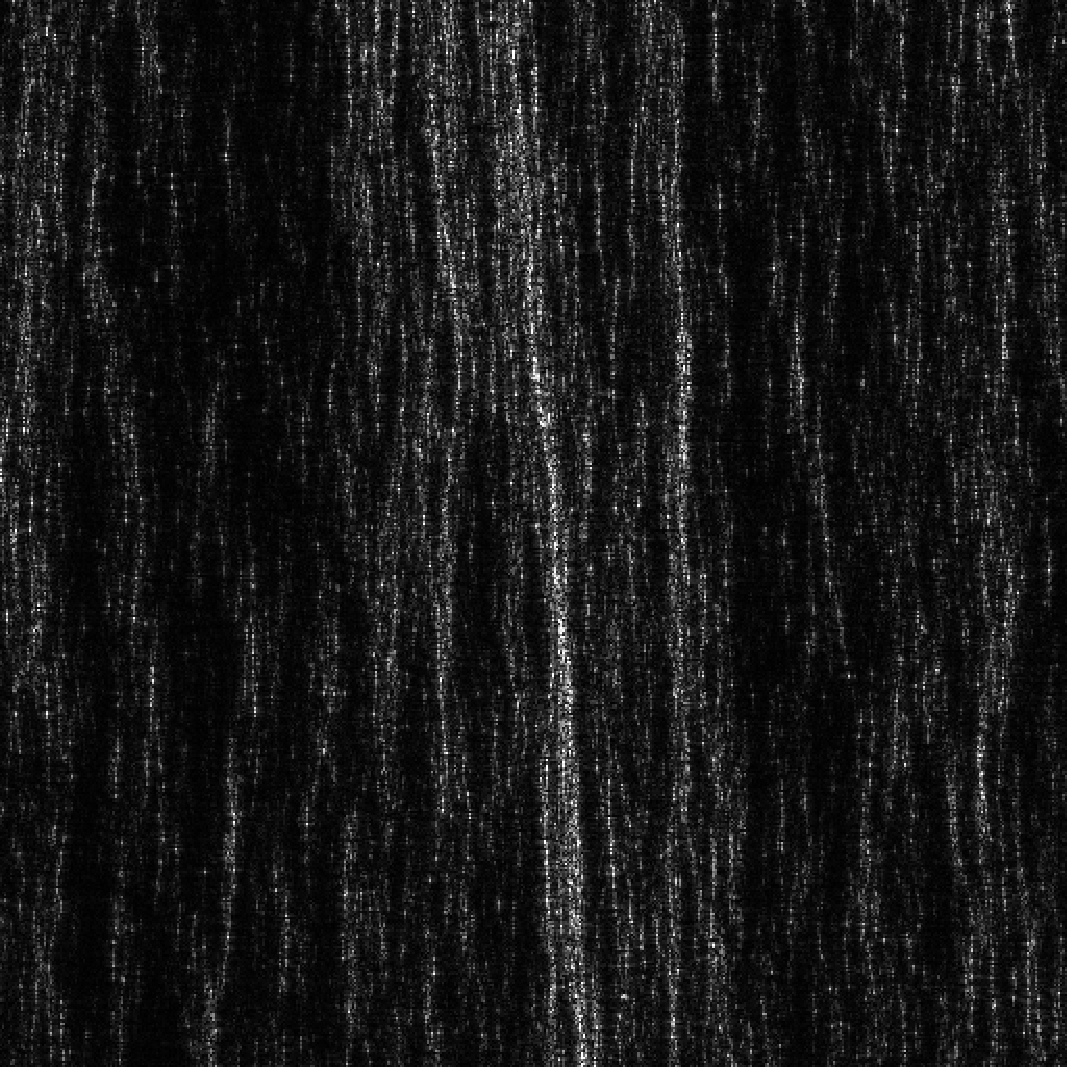}               
                \caption{}
                \label{dem_noisy1}
       \end{subfigure}
       \begin{subfigure}[b]{0.28\textwidth}           
                \includegraphics[scale=0.245]{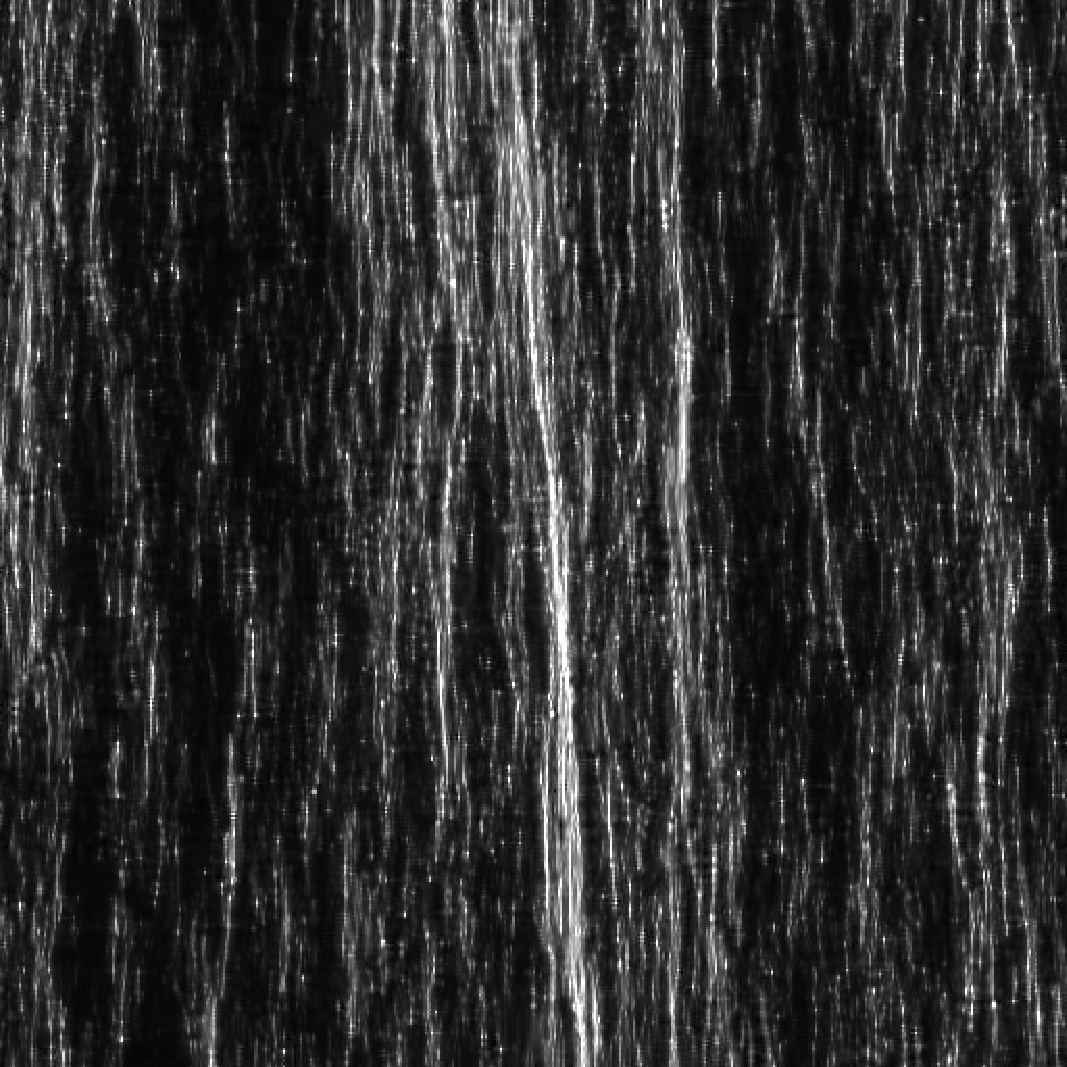}               
                \caption{}
                \label{dem_sarbm3d1}
       \end{subfigure}
       
             \begin{subfigure}[b]{0.28\textwidth}           
                \includegraphics[scale=0.245]{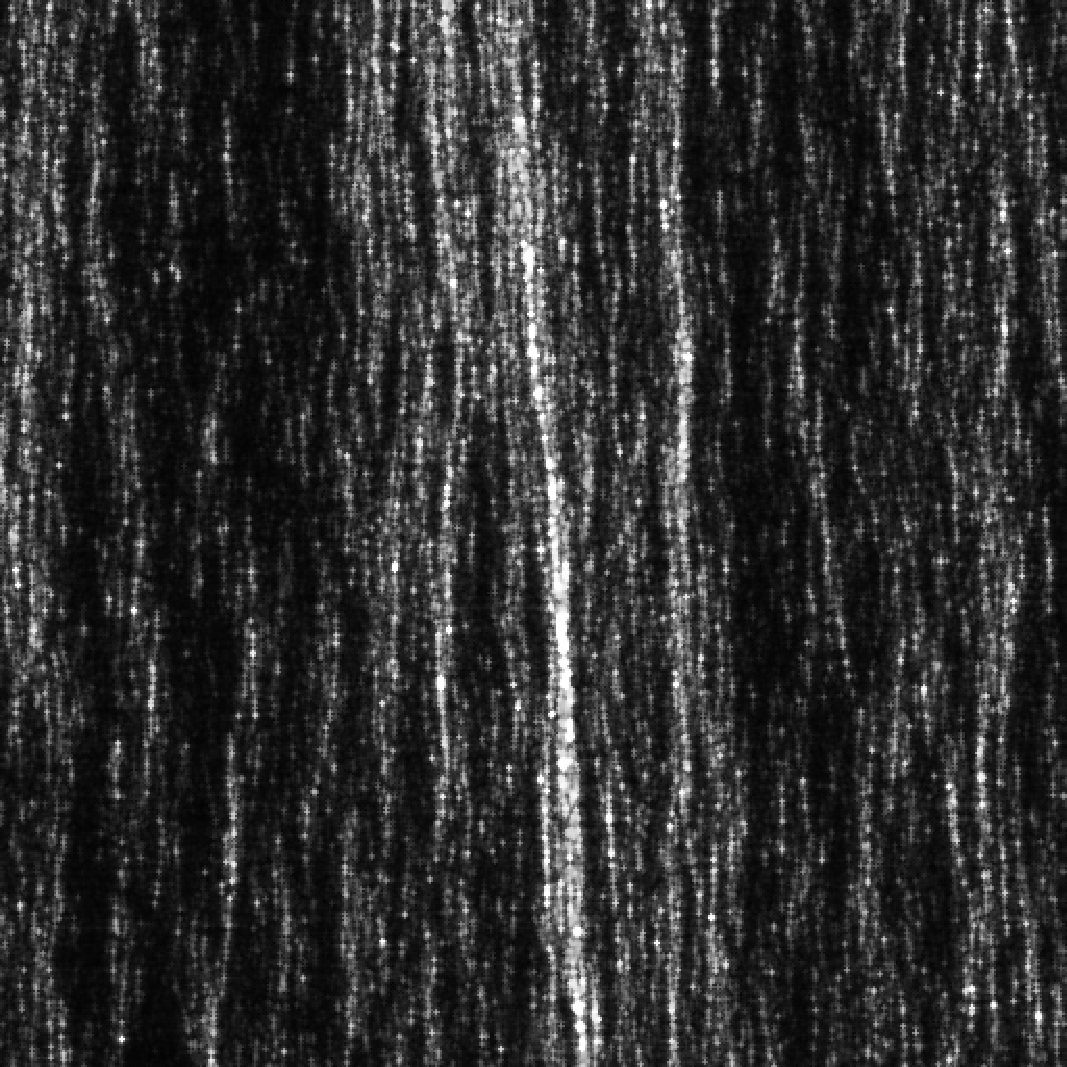}               
                \caption{}
                \label{dem_avggraytdm1}
       \end{subfigure}
           \begin{subfigure}[b]{0.28\textwidth}           
                \includegraphics[scale=0.245]{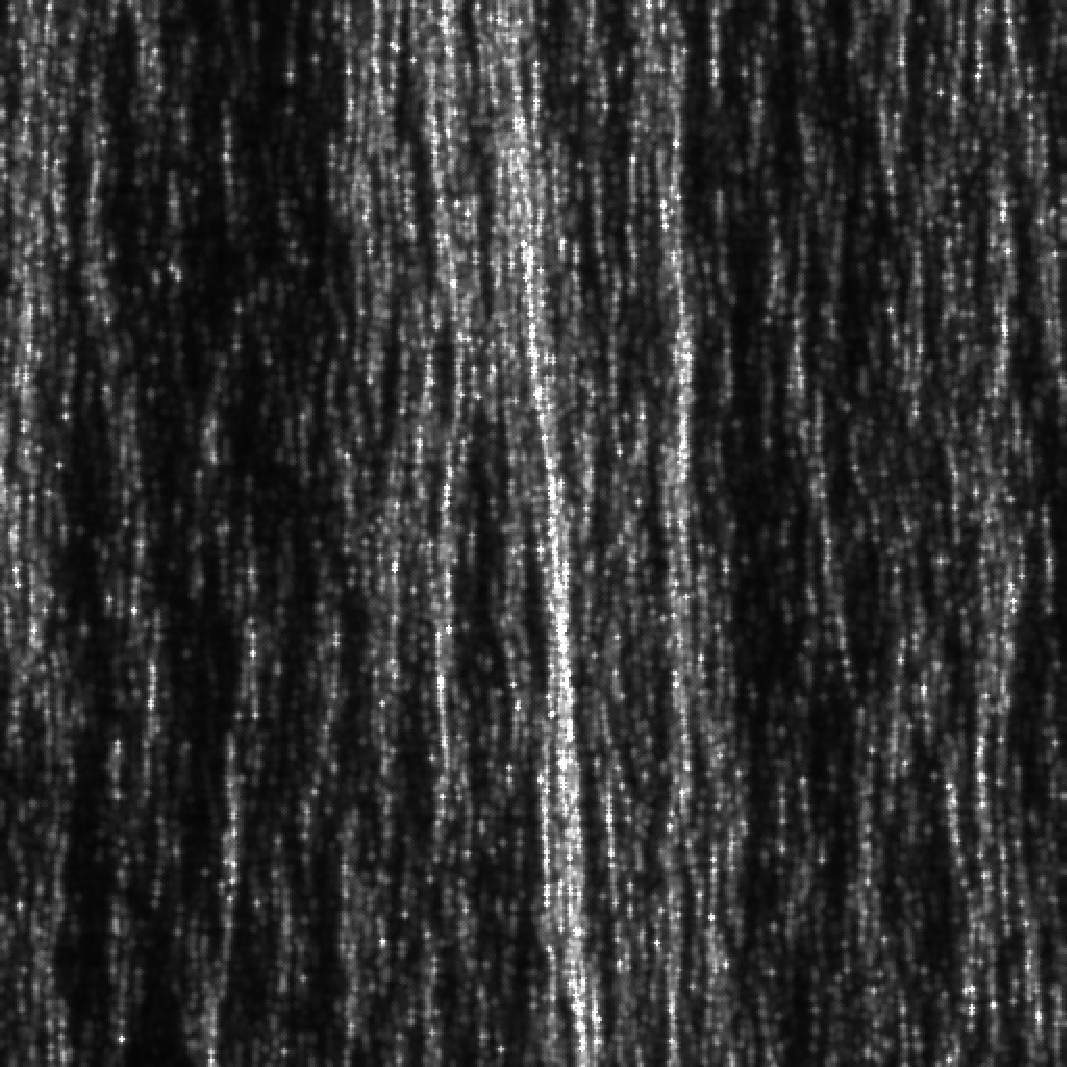}               
                \caption{}
                \label{dem_zzdbtdm1}
       \end{subfigure}
            \begin{subfigure}[b]{0.28\textwidth}           
                \includegraphics[scale=0.245]{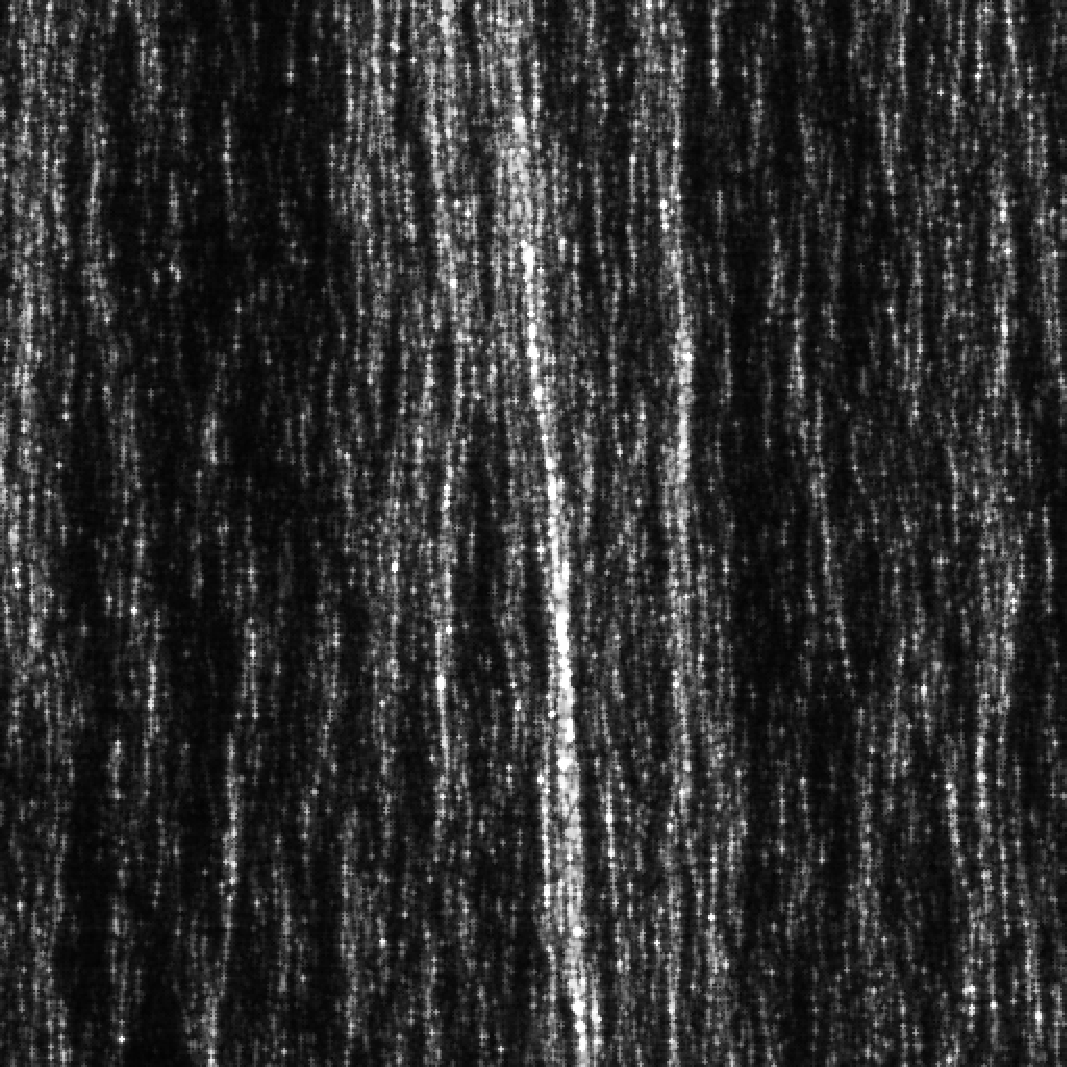}               
                \caption{}
                \label{dem_gradtdm1}
       \end{subfigure}
\caption{\footnotesize 512 $\times$ 512 DEM image. (a) 512-look. (b) One-look. (c) SAR-BM3D. (d) $Q=Q_0-p_1$. (e) $Q=Q_0-p_2$. (f) $Q=Q_0-p_3$.}\label{dem_tdmpnew}
\end{figure}

\begin{figure}
       \centering
        \begin{subfigure}[b]{0.28\textwidth}           
                \includegraphics[scale=0.245]{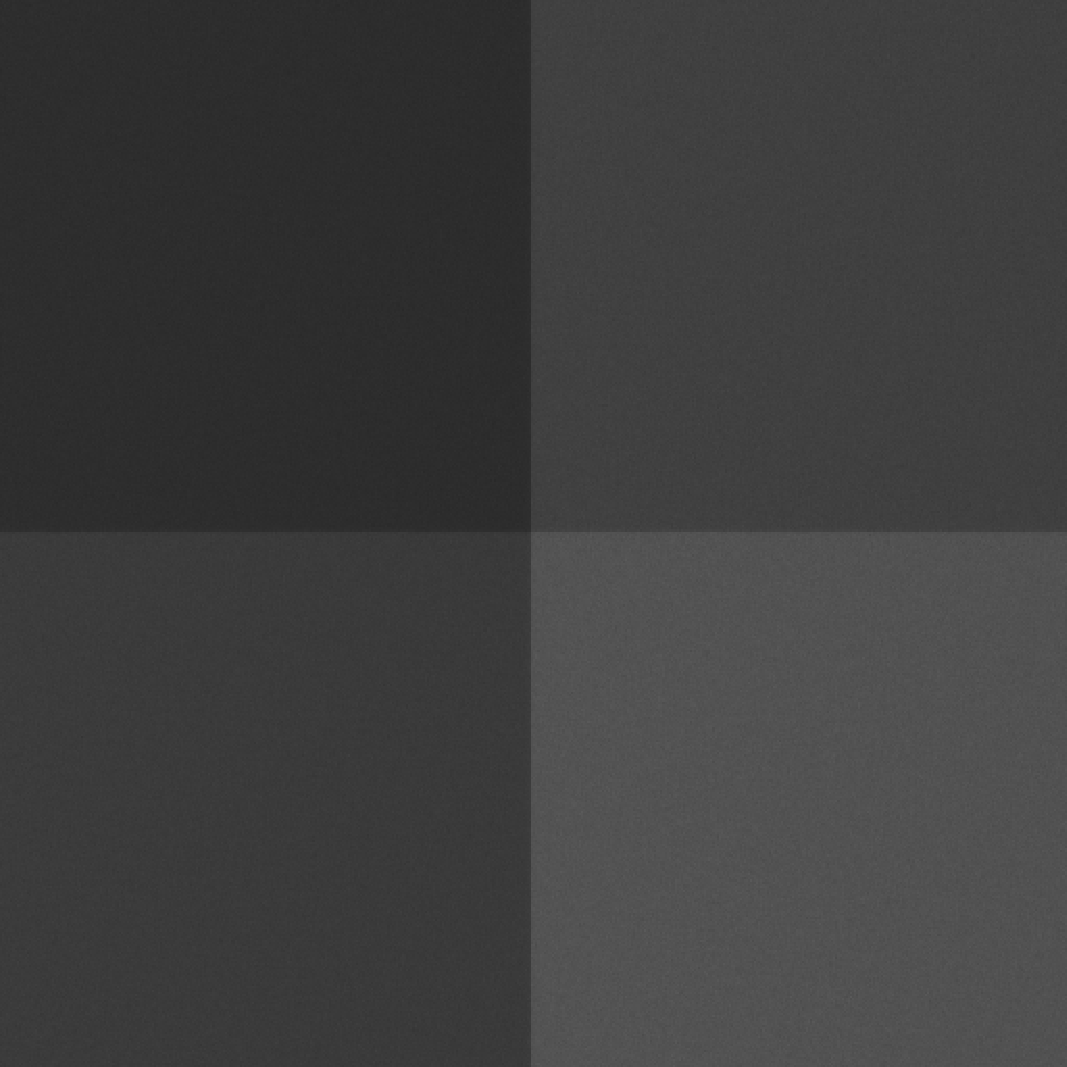}               
               \caption{FOM=0.993}
                \label{squares_clean1}
       \end{subfigure}
         \begin{subfigure}[b]{0.28\textwidth}           
                \includegraphics[scale=0.245]{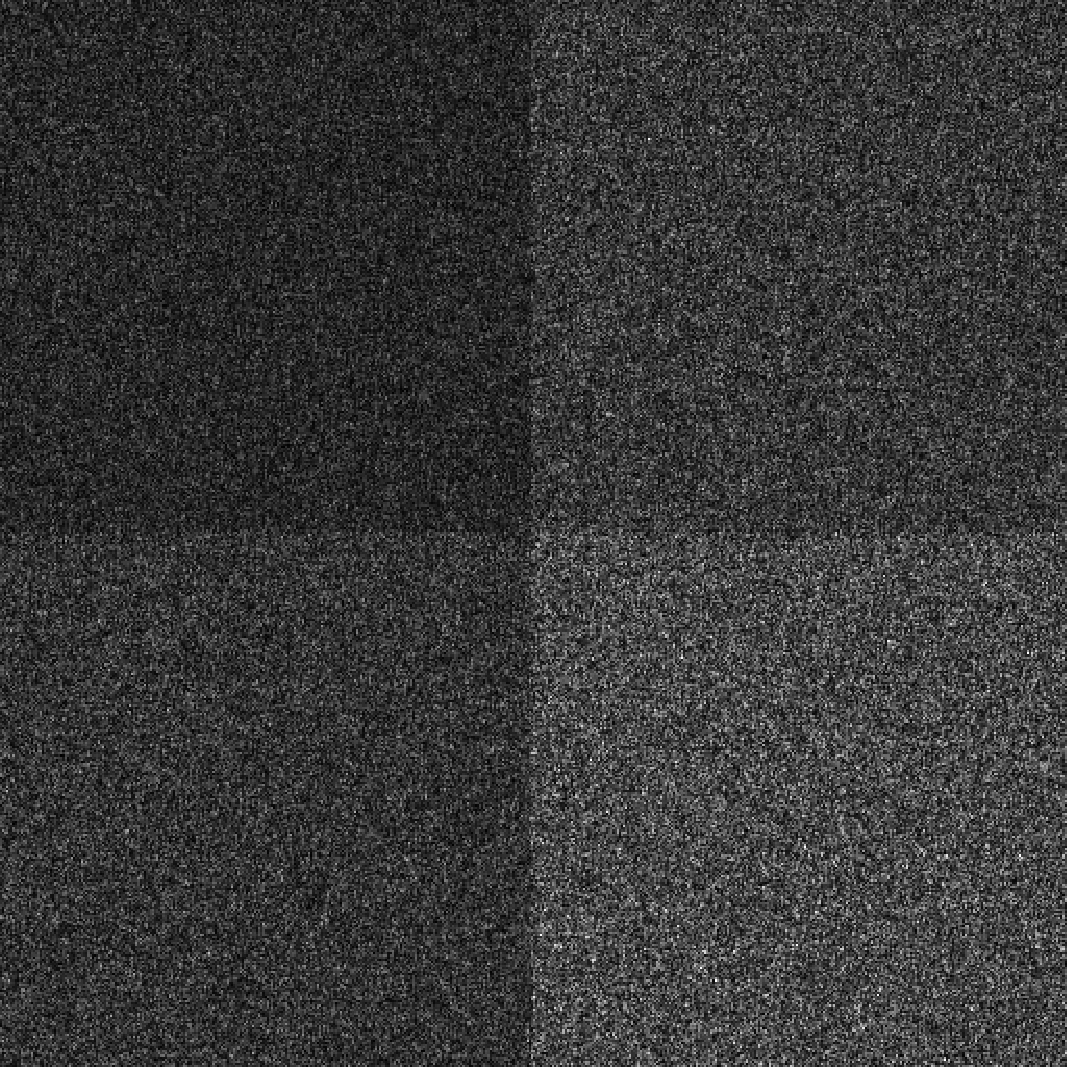}               
                \caption{FOM=0.792}
                \label{squares_noisy1}
       \end{subfigure}
       \begin{subfigure}[b]{0.28\textwidth}           
                \includegraphics[scale=0.245]{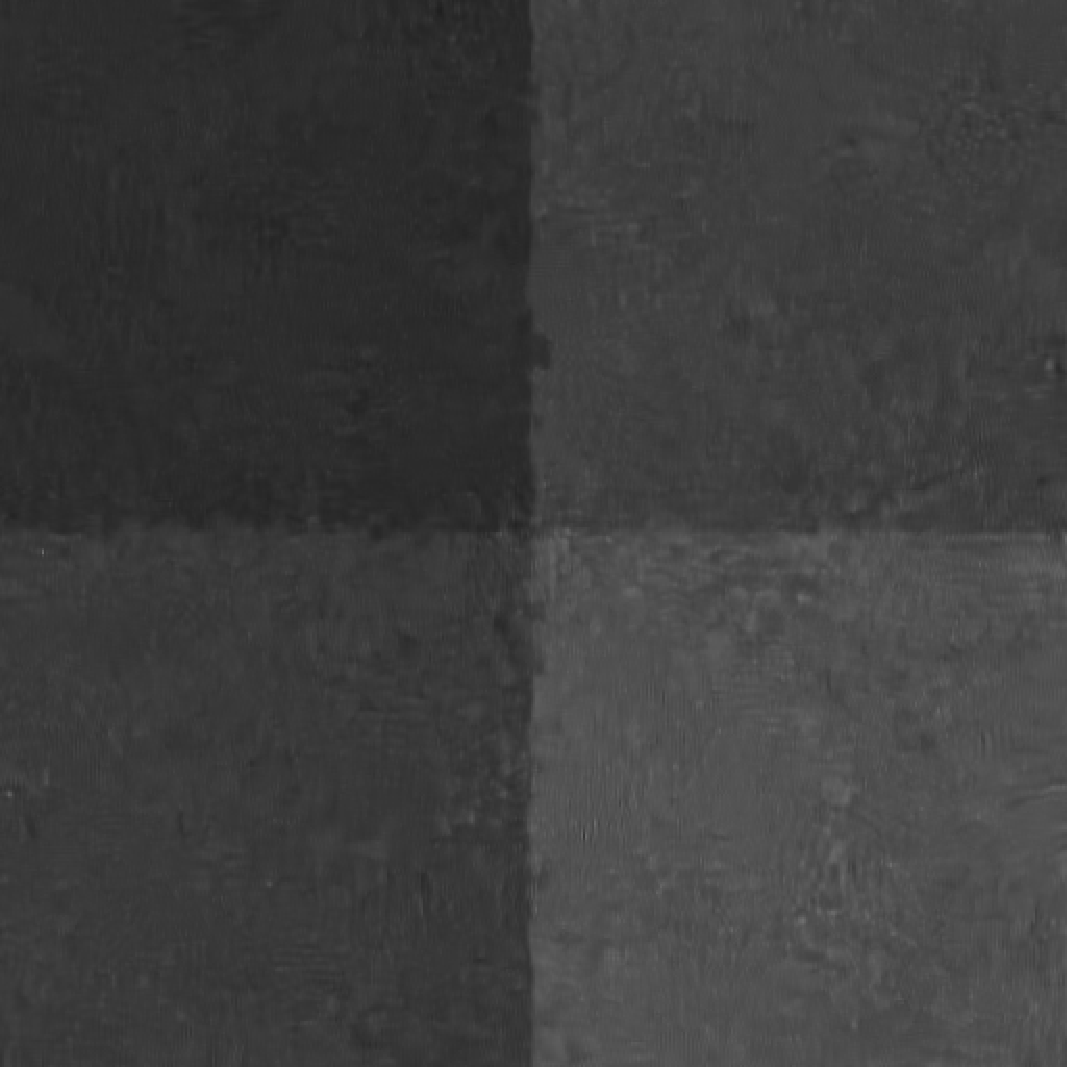}               
                \caption{FOM=0.847}
                \label{squares_sarbm3d1}
       \end{subfigure}
       
        \begin{subfigure}[b]{0.28\textwidth}           
                \includegraphics[scale=0.245]{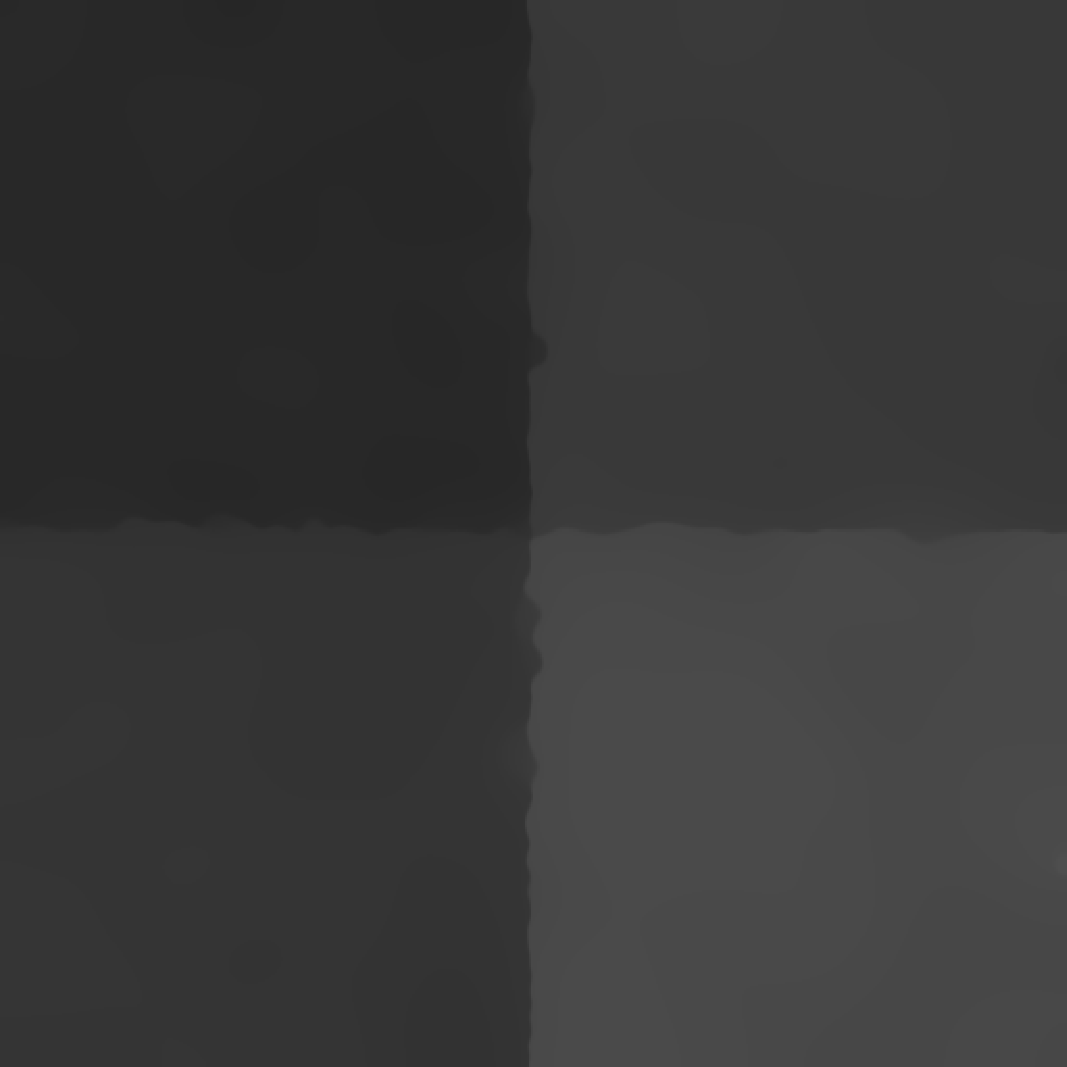}               
                \caption{FOM=0.864}
                \label{squares_avggray1}
       \end{subfigure}
        \begin{subfigure}[b]{0.28\textwidth}           
                \includegraphics[scale=0.245]{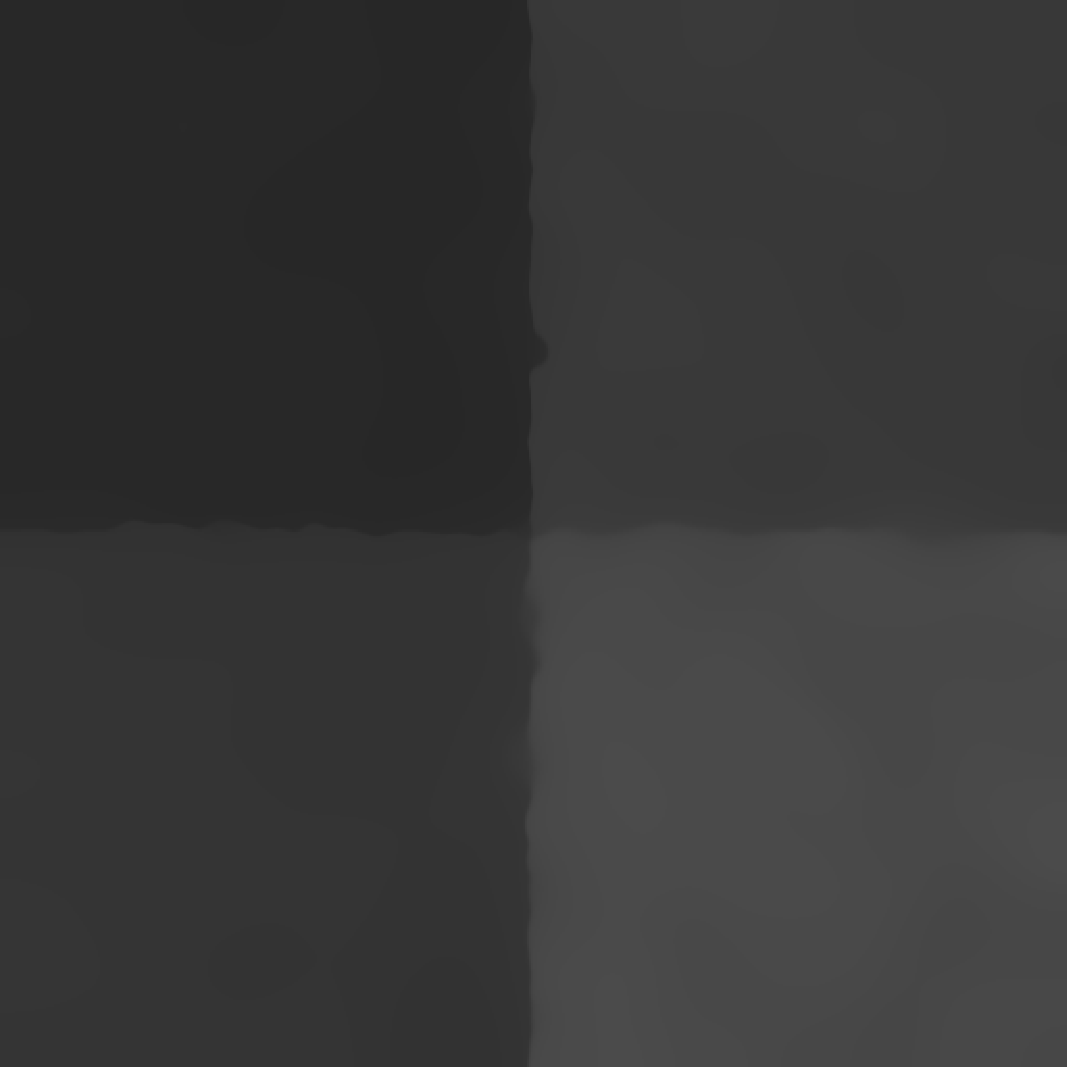}               
                \caption{FOM=0.811}
                \label{squares_zzdbtdm1}
       \end{subfigure}
        \begin{subfigure}[b]{0.28\textwidth}           
                \includegraphics[scale=0.245]{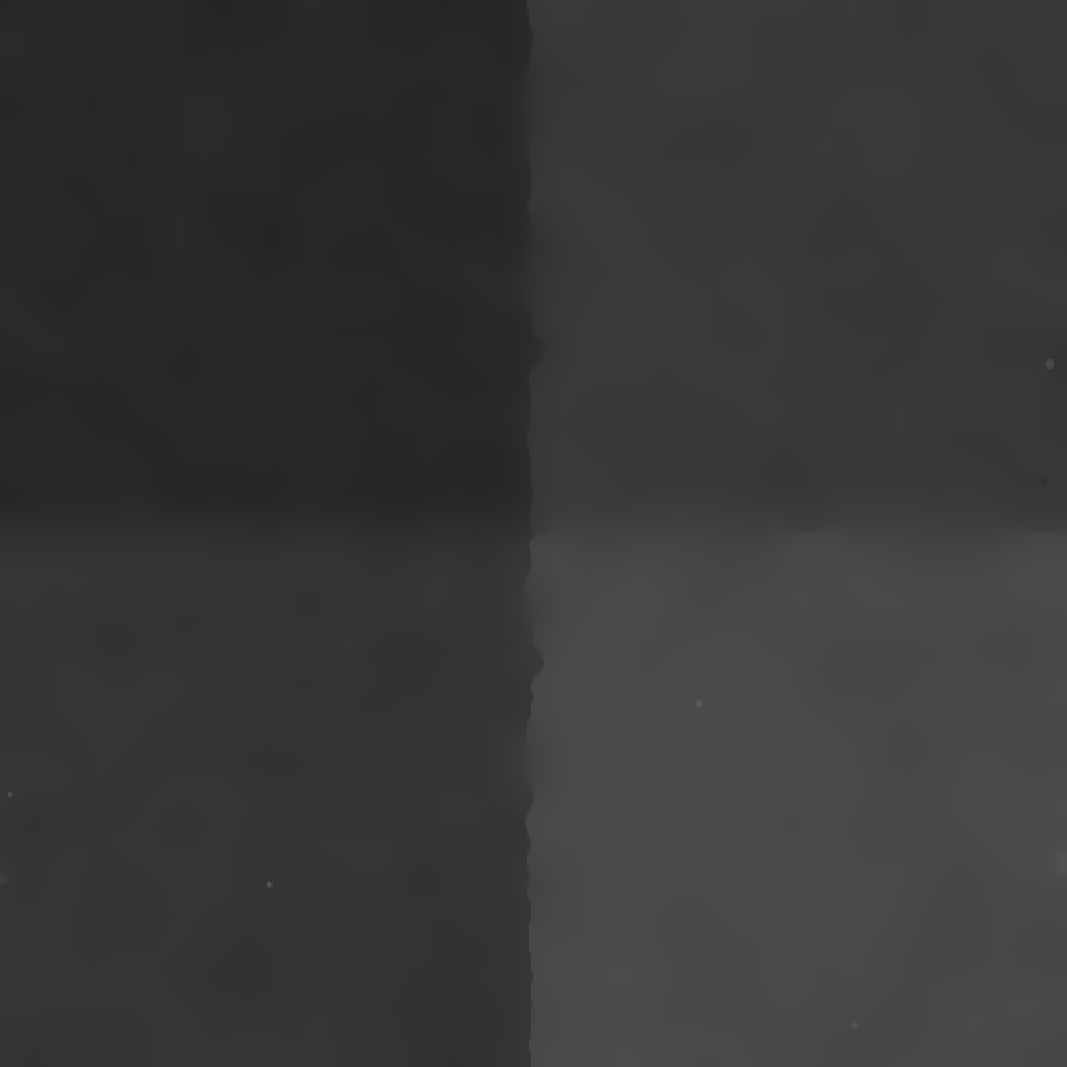}               
                \caption{FOM=0.836}
                \label{squares_gradtdm1}
       \end{subfigure}
     
\caption{\footnotesize 512 $\times$ 512 Square image. (a) 512-look. (b) One-look. (c) SAR-BM3D. (d) $Q=Q_0-p_1$. (e) $Q=Q_0-p_2$. (f) $Q=Q_0-p_3$.}\label{squares_tdmpnew}
\end{figure}

\begin{figure}
       \centering
        \begin{subfigure}[b]{0.28\textwidth}           
                \includegraphics[scale=0.49]{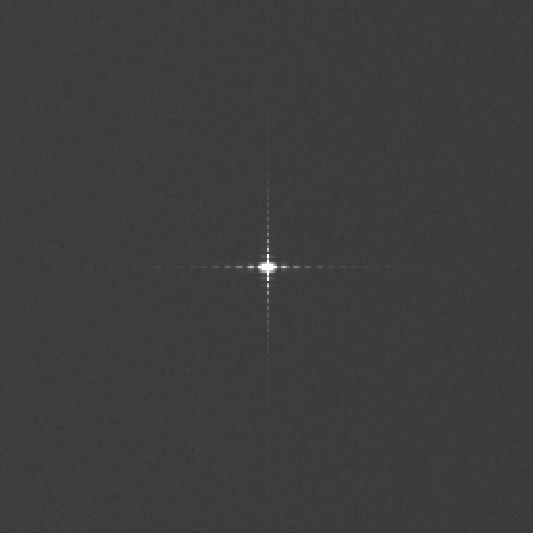}               
               \caption{$C_{\text{NN}}=7.75$, $C_{\text{BG}}=36.56$}
                \label{corner_clean1}
       \end{subfigure}
         \begin{subfigure}[b]{0.28\textwidth}           
                \includegraphics[scale=0.49]{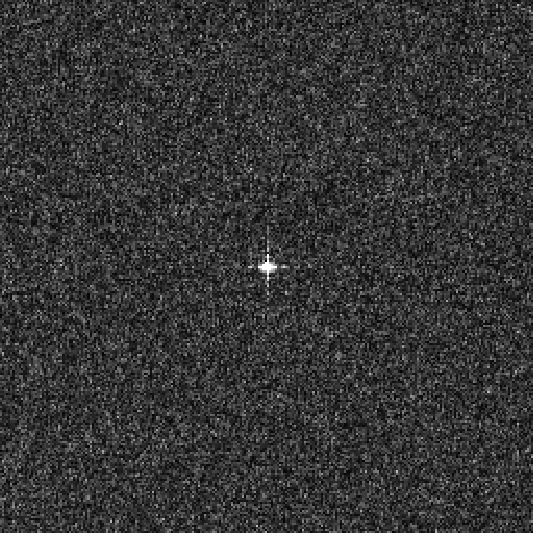}               
                \caption{$C_{\text{NN}}=7.77$, $C_{\text{BG}}=36.50$}
                \label{corner_noisy1}
       \end{subfigure}
        \begin{subfigure}[b]{0.28\textwidth}           
                \includegraphics[scale=0.49]{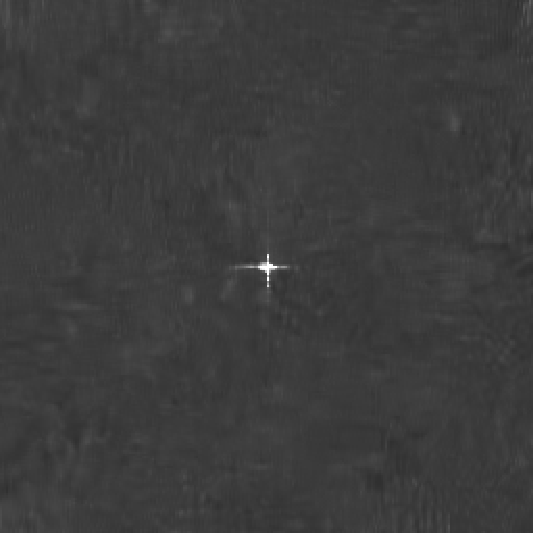}               
               \caption{$C_{\text{NN}}=7.39$, $C_{\text{BG}}=35.45$}
                \label{corner_sarbm3d1}
       \end{subfigure}

     \begin{subfigure}[b]{0.28\textwidth}           
                \includegraphics[scale=0.49]{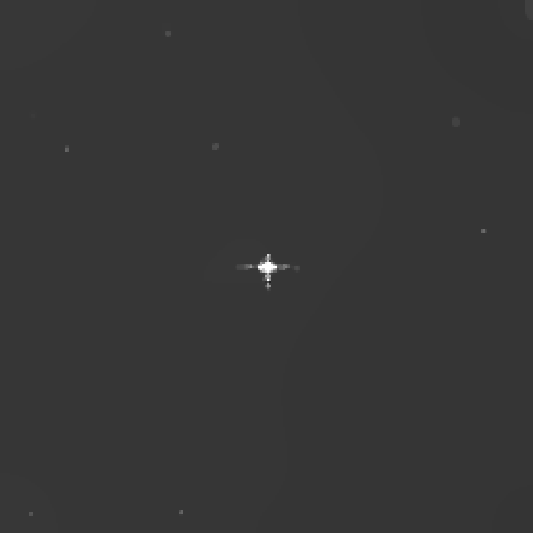}               
                \caption{$C_{\text{NN}}=7.42$, $C_{\text{BG}}=36.59$}
                \label{corner_avggray}
       \end{subfigure}
         \begin{subfigure}[b]{0.28\textwidth}           
                \includegraphics[scale=0.49]{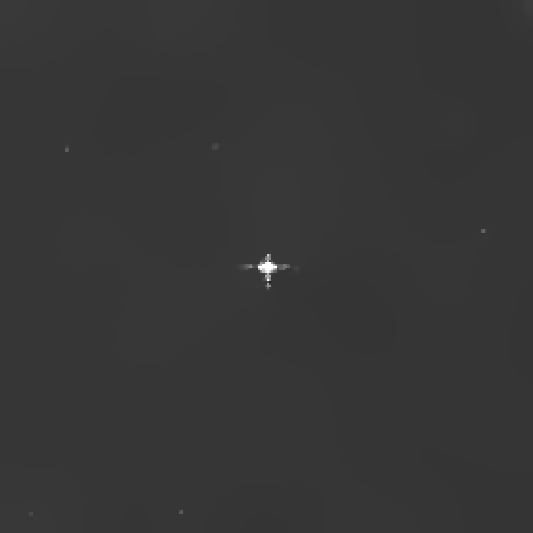}               
                \caption{$C_{\text{NN}}=7.40$, $C_{\text{BG}}=36.54$}
                \label{corner_zzdbtdm1}
       \end{subfigure}
         \begin{subfigure}[b]{0.28\textwidth}           
                \includegraphics[scale=0.49]{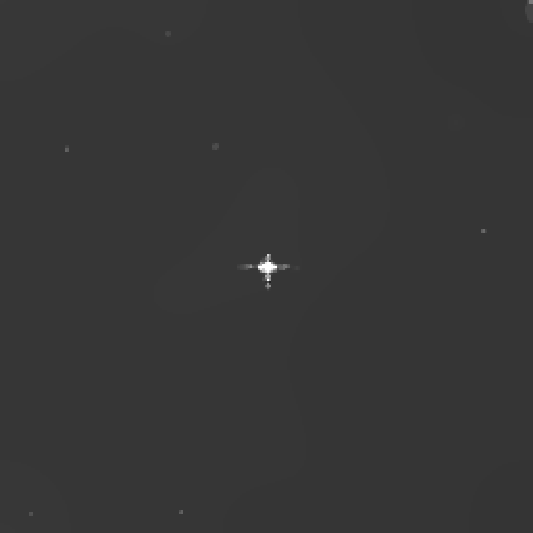}               
                \caption{$C_{\text{NN}}=7.42$, $C_{\text{BG}}=36.58$}
                \label{corner_gradtdm1}
       \end{subfigure}
     
\caption{\footnotesize 256 $\times$ 256 Corner image. (a) 512-look. (b) One-look. (c) SAR-BM3D. (d) $Q=Q_0-p_1$. (e) $Q=Q_0-p_2$. (f) $Q=Q_0-p_3$.}\label{corner_tdmpnew}
\end{figure}

\begin{figure}
       \centering
        \begin{subfigure}[b]{0.32\textwidth}           
                \includegraphics[scale=0.38]{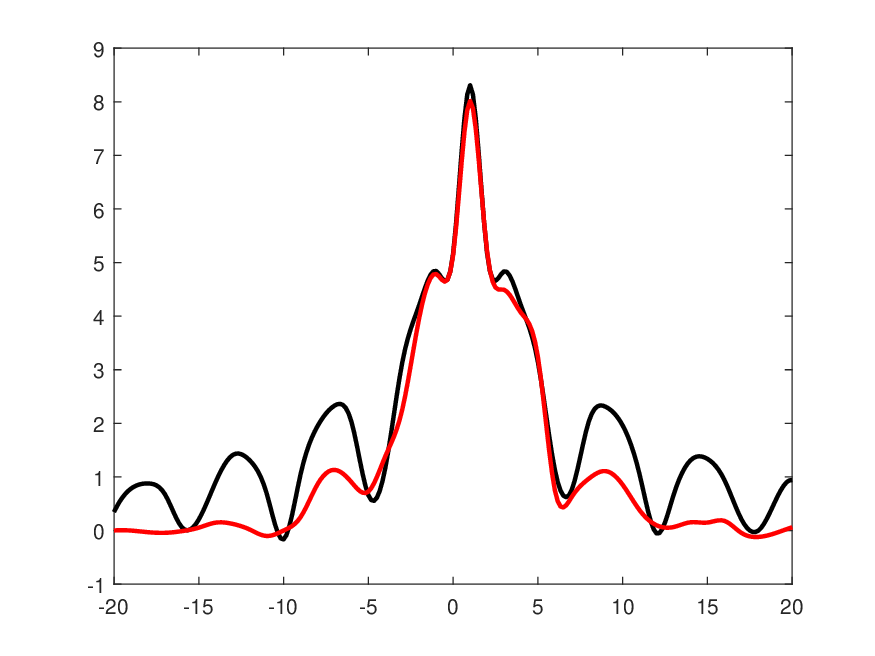}               
               \caption{}
                \label{corner_prof_sarbm3d}
       \end{subfigure}
         \begin{subfigure}[b]{0.32\textwidth}           
                \includegraphics[scale=0.38]{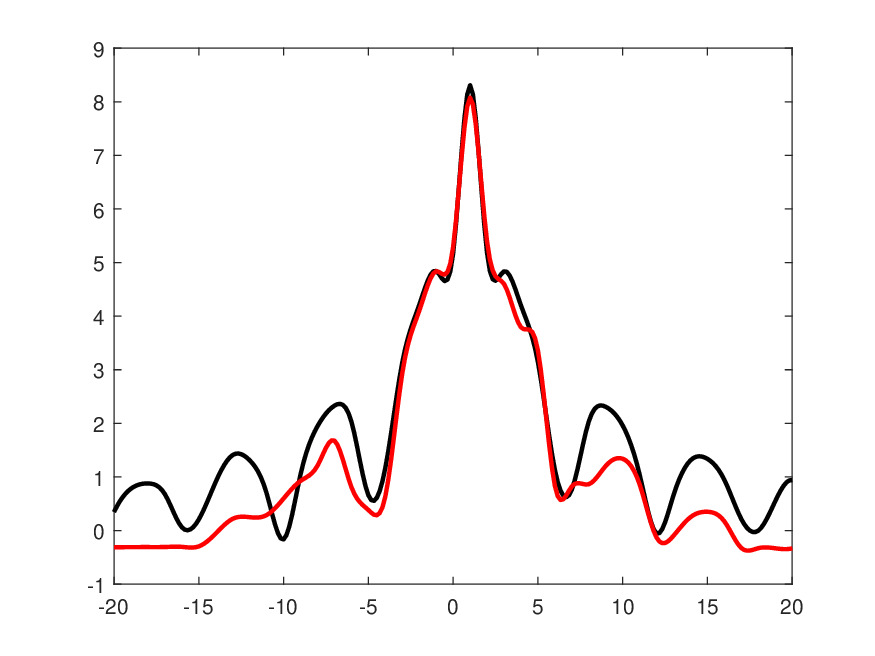}               
                \caption{}
                \label{corner_prof_tdm_avggray}
       \end{subfigure}
       
          \begin{subfigure}[b]{0.32\textwidth}           
                \includegraphics[scale=0.38]{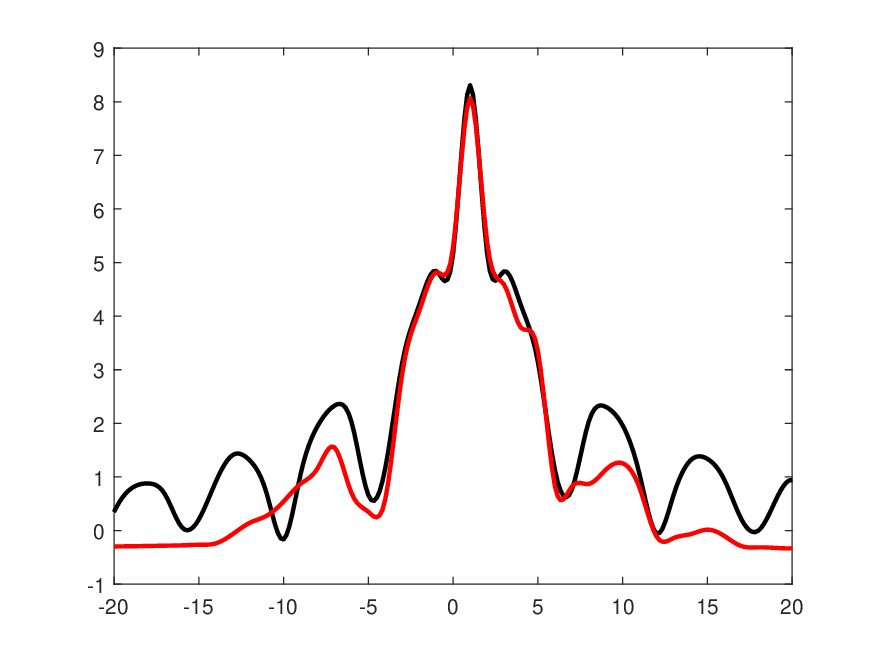}               
                \caption{}
                \label{corner_prof_tdm_zzdb}
       \end{subfigure}
           \begin{subfigure}[b]{0.32\textwidth}           
                \includegraphics[scale=0.38]{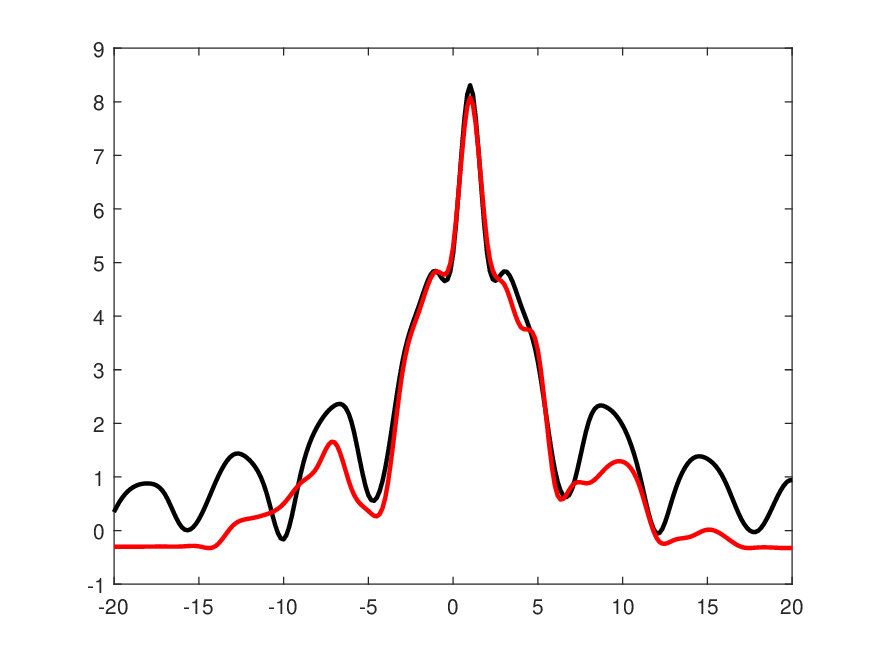}               
                \caption{}
                \label{corner_prof_tdm_grad}
       \end{subfigure}
       
 \caption{\footnotesize Corner reflector: restored range profiles (red) compared with the clean one (black). (a) SAR-BM3D. (b) $Q=Q_0-p_1$. (c) $Q=Q_0-p_2$. (d) $Q=Q_0-p_3$.}\label{corner_profile_tdmpnew}
\end{figure}

\begin{figure}
       \centering
        \begin{subfigure}[b]{0.28\textwidth}           
                \includegraphics[scale=0.49]{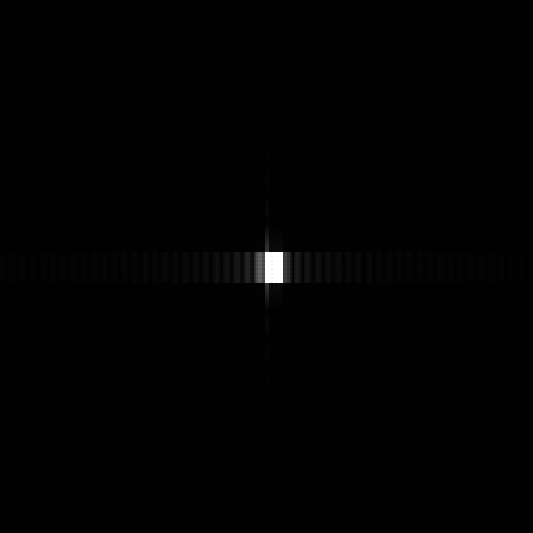}               
               \caption{$C_{\text{DR}}=65.90$}
                \label{building_clean1}
       \end{subfigure}
         \begin{subfigure}[b]{0.28\textwidth}           
                \includegraphics[scale=0.49]{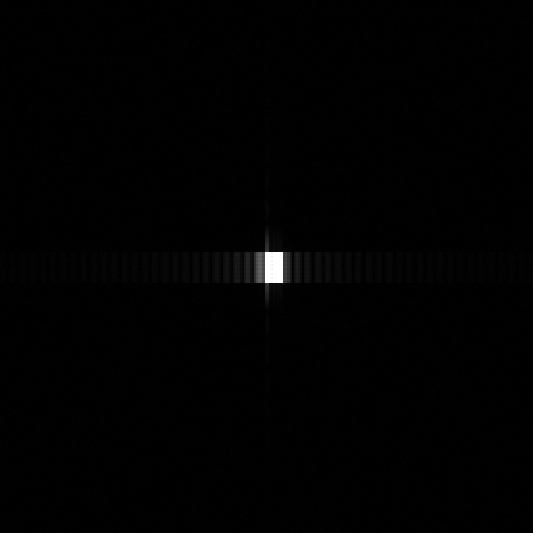}               
                \caption{$C_{\text{DR}}=65.90$, BS=0.09}
                \label{building_noisy1}
       \end{subfigure}
        \begin{subfigure}[b]{0.28\textwidth}           
                \includegraphics[scale=0.49]{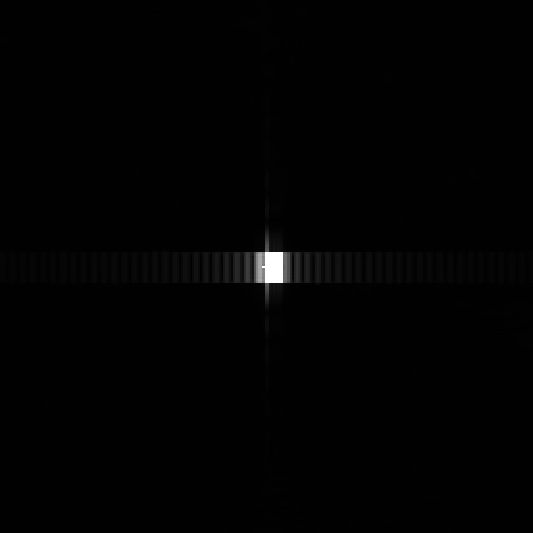}               
                \caption{$C_{\text{DR}}=65.91$, BS=1.46}
                \label{building_sarbm3d1}
       \end{subfigure}
       
       \begin{subfigure}[b]{0.28\textwidth}           
                \includegraphics[scale=0.49]{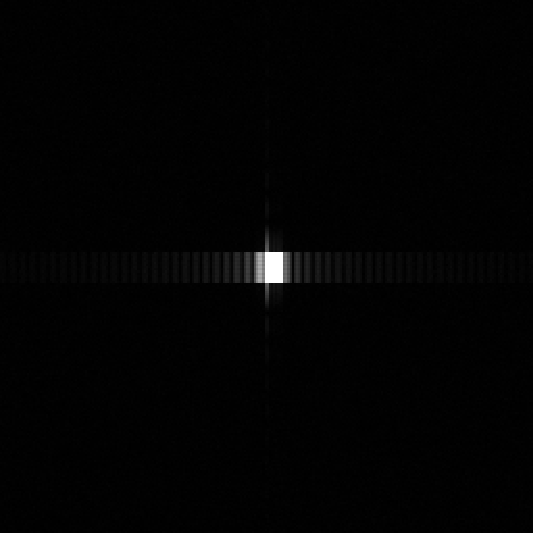}               
                \caption{$C_{\text{DR}}=65.89$, BS=0.1091}
                \label{building_avggray1}
       \end{subfigure}
           \begin{subfigure}[b]{0.28\textwidth}           
                \includegraphics[scale=0.49]{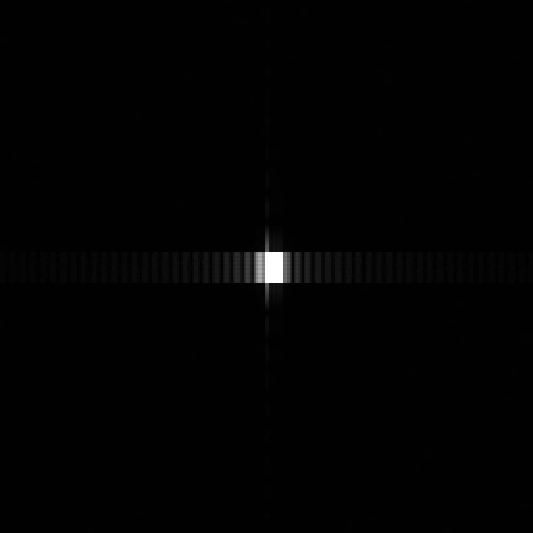}               
                \caption{$C_{\text{DR}}=66.92$, BS=0.2218}
                \label{building_zzdb1}
       \end{subfigure}
           \begin{subfigure}[b]{0.28\textwidth}           
                \includegraphics[scale=0.49]{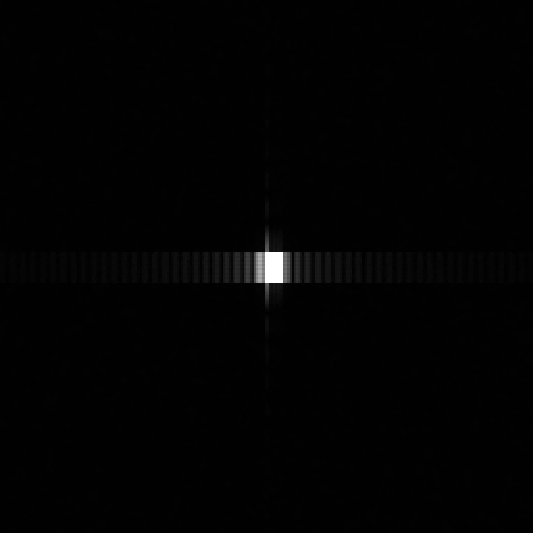}               
                \caption{$C_{\text{DR}}=66.81$, BS=0.2267}
                \label{building_gradtdm1}
       \end{subfigure}
     
\caption{\footnotesize 256 $\times$ 256 Building image. (a) 512-look. (b) One-look. (c) SAR-BM3D. (d) $Q=Q_0-p_1$. (e) $Q=Q_0-p_2$. (f) $Q=Q_0-p_3$.}\label{building_tdmpnew}
\end{figure}

\begin{figure}
       \centering
        \begin{subfigure}[b]{0.32\textwidth}           
                \includegraphics[scale=0.38]{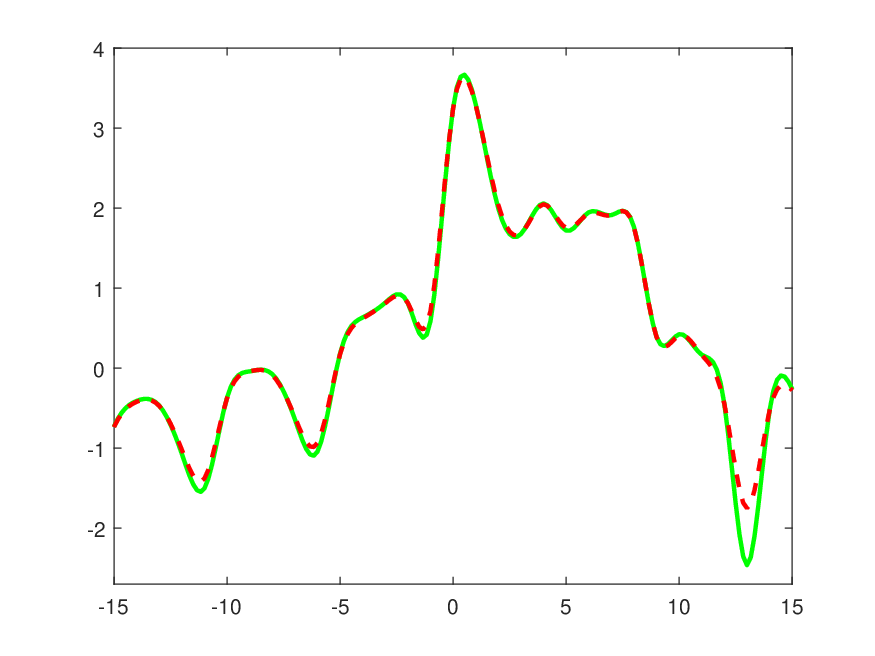}               
               \caption{}
                \label{building_prof_sarbm3d}
       \end{subfigure}
         \begin{subfigure}[b]{0.32\textwidth}           
                \includegraphics[scale=0.38]{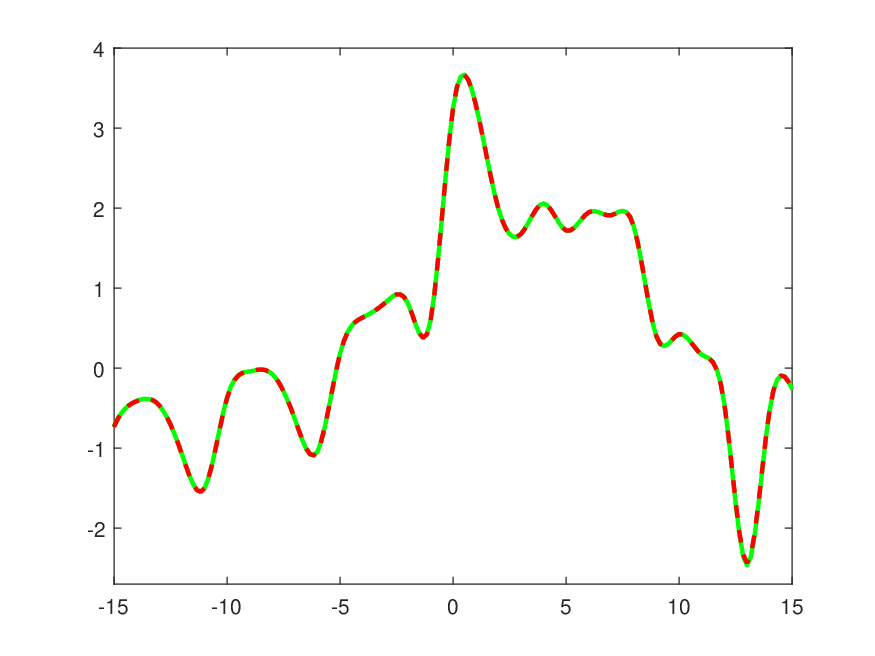}               
                \caption{}
                \label{building_prof_avggray1}
       \end{subfigure}
       
               \begin{subfigure}[b]{0.32\textwidth}           
                \includegraphics[scale=0.38]{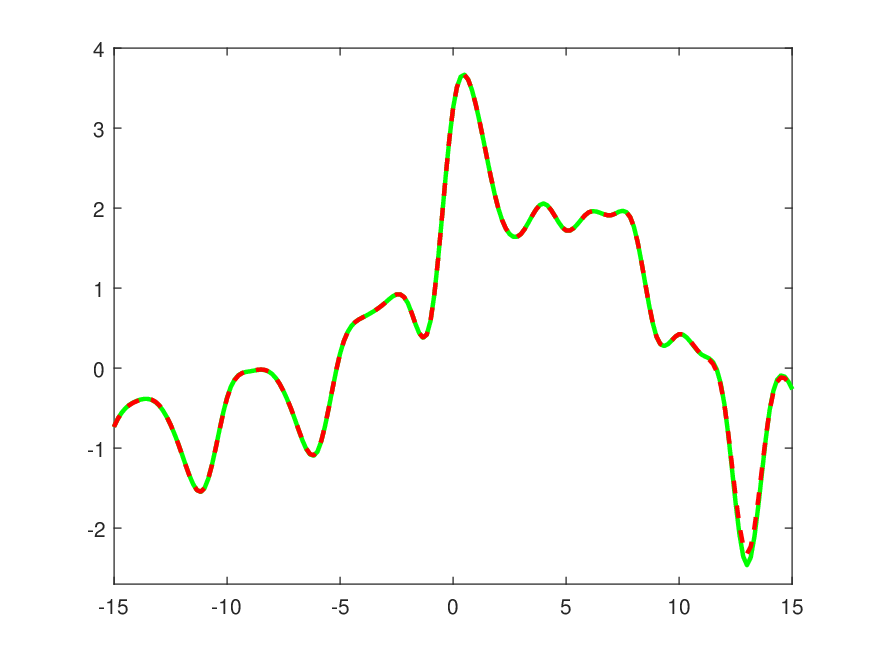}               
               \caption{}
                \label{building_prof_zzdb1}
       \end{subfigure}
         \begin{subfigure}[b]{0.32\textwidth}           
                \includegraphics[scale=0.38]{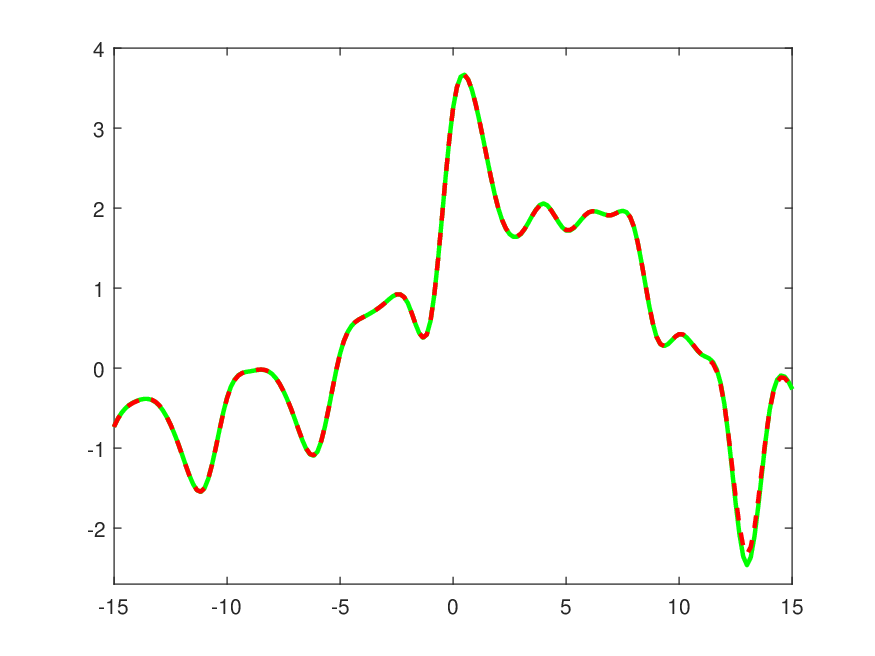}               
                \caption{}
                \label{building_prof_grad1}
       \end{subfigure}
     
\caption{\footnotesize Restored building profiles (red dashed) compared with the clean one (green). (a) SAR-BM3D. (b) $Q=Q_0-p_1$. (c) $Q=Q_0-p_2$. (d) $Q=Q_0-p_3$.}\label{building_profile_tdmpnew}
\end{figure}

This section discusses the results of a benchmark framework to show the despeckling ability of the present technique. Also, we compare the results of the present model with the results of SAR-BM3D \cite{parrilli2011nonlocal}. We compute the despeckled results for the five standard images (Homogeneous: Figure \ref{homogeneous_noisy1}, DEM: Figure \ref{dem_noisy1}, Squares: Figure \ref{squares_noisy1}, Corner: Figure \ref{corner_noisy1}, Building: Figure \ref{building_noisy1}) using the suggested model and compare the outcomes with the results of SAR-BM3D. We evaluate different quality measurement parameters \cite{di2013benchmarking} for the restored images. Each parameter is computed w.r.t. a 512-look reference image of the same scene. To gain a concrete idea of these images and quantitative measures, we refer the reader to \cite{di2013benchmarking}. Figure \ref{homogeneous_tdmpnew} illustrates the restored results for the Homogeneous image. It is simplest and most important as it provides knowledge on the speckle removal efficiency of an approach. Table \ref{measures_homogeneous} reports bias indicators mean value of the image (MoI), MoR, VoR, and performance measures ENL, ENL$^{*}$, and DG for the despeckled images. From the quantitative measures, one can observe that none of the methods introduces significant bias on MoI and MoR. Comparatively, a smaller value for VoR indicates imperfect speckle filtering for the case of SAR-BM3D. Also, in terms of ENL, ENL$^{*}$, and DG, the first two exponents work better than the third exponent and SAR-BM3D. Also, from the qualitative results, one can notice that \eqref{maina}--\eqref{mainc} can remove speckle noise better than the SAR-BM3D. Figure \ref{dem_tdmpnew} shows the despeckling results for the DEF image; it is a highly textured image. Observing the filtered images, one can conclude that despeckling results are not significantly biased. Also, the quantitative results MoI, MoR, VoR, coefficient of variation ($C_x$), and DG are presented in table \ref{measures_dem}. In terms of MoI and MoR, all the approaches can produce satisfactory despeckling results. A higher VoR value indicates the second exponent can remove the speckle better than the others. In terms of $C_x$ and DG, the results for the proposed method are comparable to SAR-BM3D. In Figure \ref{squares_tdmpnew}, we present the filtered images for the Squares image. We calculated the edge-preserving parameter, the figure of merit (FOM), for each image and mentioned the values in the captions of images. From the despeckled images, one can conclude that the first exponent can remove the speckle efficiently and improve the edges better than the other exponents and SAR-BM3D. Also, the values of quantitative measure FOM support the qualitative results. Figure \ref{corner_tdmpnew} shows the restored results for the Corner image. Also, the called contrast measures $C_{\text{NN}}$ (for the eight nearest pixels) and $C_{\text{BG}}$ (for the background) are computed and mentioned in the captions of the images. For the case of the proposed model, the Contrast values are closer to the contrast values of the reference image. Figure \ref{corner_profile_tdmpnew} describes a range profile of the corner response, and one can notice that the restored profile using the first exponent is more intimate to the profile for the clean image as compared to the other results. Figure \ref{building_tdmpnew} illustrates the filtering outcomes for the Building image. We computed two parameters: contrast measure $C_{\text{DR}}$ and building smearing (BS) figure and mentioned them in the caption of the images. All the strategies can preserve the double reflection line quite accurately. Also, Figure \ref{building_profile_tdmpnew} describes average range profiles (calculated averaging over the range lines interested by the presence of the building) on a logarithmic scale. The first exponent correctly specifies the whole building region and preserves it integrally, with BS close to zero.

\begin{table}
\caption{\footnotesize MEASURES FOR HOMOGENEOUS} 
\centering 
\scalebox{0.75}{
\begin{tabular}{c c c c c c c} 
\hline\hline 
 & MoI & MoR & VoR & ENL & ENL$^{*}$ & DG \\ [0.5ex] 
\hline 
Clean        & 1.000 & 0.998 & 0.987 & 436.97 & 510.36 & $-$ \\ 
Noisy        & 0.998 &  $-$  & $-$   & 1.00   & 1.02   &  0 \\
SAR-BM3D     & 0.978 &  0.979 & 0.814 & 102.44 & 111.91 & 19.40\\
$Q_0-p_1$    & 1.000 &  0.989 & 0.960 & 365.85 & 507.20 & 23.65\\
$Q_0-p_2$    & 0.989 &  0.999 & 0.979 & 369.08 & 508.20 & 23.60\\
$Q_0-p_3$    & 0.990 &  0.999 & 0.977 & 348.94 & 475.61 & 23.46\\ [1ex] 
\hline 
\end{tabular}
}
\label{measures_homogeneous} 
\end{table}

\begin{table}
\caption{\footnotesize MEASURES FOR DEM} 
\centering 
\scalebox{0.85}{
\begin{tabular}{c c c c c c c} 
\hline\hline 
 & MoI & MoR & VoR & $C_x$  & DG \\ [0.5ex] 
\hline 
Clean      & 1.000 &  1.001 & 0.999 & 2.40  & $-$ \\ 
Noisy      & 1.003 &  $-$   & $-$   & 3.54  &  0 \\
SAR-BM3D   & 0.968 &  0.833 & 0.415 & 2.43  & 5.32\\
$Q_0-p_1$  & 0.988 &  0.849 & 0.433 & 2.20  & 4.02\\
$Q_0-p_2$  & 0.988 &  0.820 & 0.557 & 2.24  & 3.72\\
$Q_0-p_3$  & 0.988 &  0.850 & 0.429 & 2.26  & 4.00\\ [1ex] 
\hline 
\end{tabular}
}
\label{measures_dem} 
\end{table}

\section{Conclusion}
\label{sec:Conclusion}
This work discusses a class of diffusion and telegraph diffusion-based approaches for image despeckling. Instead of a constant exponent, we use three distinct variable exponents for the edge indicator function, which depend on either the image's gray level or its gradient. Though the discussed diffusion model behaves as some existing PDE approaches in some specific conditions, our main interest is in the variable exponent-based telegraph diffusion equation. To the best of our knowledge, the variable exponent-based telegraph diffusion model has never been used for speckle noise elimination. Also, we rigorously study the existence and uniqueness of a weak solution to the discussed model using the Schauder fixed point theorem. Then we run numerical experiments on both artificial noise and real SAR images. Further, to check the SAR despeckling ability of the proposed model, we use a standard benchmarking framework. We compare the despeckled results of the present approach with the results of a nonlocal technique. The study confirms that the proposed hybrid model is capable of performing image despeckling efficiently. As a result, the proposed method is a significant alternative in the domain of image despeckling.




\end{document}